%% file: elliptic.tex
\documentclass[11pt]{article}

\usepackage[margin=2.5cm]{geometry}
\usepackage[utf8]{inputenc}
\usepackage{graphicx}
\usepackage{amssymb}
\usepackage{amsmath}
\usepackage{mathtools}
\usepackage[usenames,dvipsnames,svgnames,table]{xcolor}
\usepackage{enumerate}
\usepackage{enumitem}
\setlist[enumerate, 1]{label=(\roman*)}
\setlist[itemize]{leftmargin=1.5em}
\setlist[description]{leftmargin=1em}
\setlist[itemize, 1]{label=$\blacktriangleright$}
\setlist[itemize, 2]{label=$\bullet$}
\usepackage[normalem]{ulem} % for striking out text with \sout

\usepackage{tikz}
\usetikzlibrary{positioning}

\usepackage{units}
\usepackage{bbm}

\usepackage{ifthen}

\usepackage{graphicx}
\usepackage{float}
\usepackage{pgfplots} 

\usepackage[hidelinks]{hyperref}  % create links
\usepackage[format=plain,labelsep=period,justification=justified,font=small,labelfont=bf, margin=1em]{caption}

\usepackage{amsthm}

\theoremstyle{plain}
\newtheorem{theorem}{Theorem}[section]

\newtheorem{lemma}[theorem]{Lemma}
\newtheorem{proposition}[theorem]{Proposition}

\newtheorem{definition}[theorem]{Definition}

\theoremstyle{definition}

\theoremstyle{remark}

\newcommand{\Z}{\mathbb{Z}}
\newcommand{\N}{\mathbb{N}}
\newcommand{\R}{\mathbb{R}}

\newcommand{\setbar}{\ | \ }
\newcommand{\set}[2]{\left\{#1 \setbar #2 \right\}}

\newcommand{\halfDelta}{\Delta^{\nicefrac{1}{2}}}

\renewcommand{\eqref}[1]{(\refeq{#1})}

\newcommand{\bela}[1]{\begin{equation}\label{#1}}
\newcommand{\ela}{\end{equation}}
\newcommand{\bear}[1]{\begin{array}{#1}}
\newcommand{\ear}{\end{array}}
\newcommand{\bx}{\mbox{\boldmath $x$}}

\newcommand{\bn}{\boldsymbol{n}}

\newcommand{\bs}{\mbox{\boldmath $s$}}
\newcommand{\be}{\boldsymbol{e}}
\newcommand{\bc}{\mbox{\boldmath $c$}}

\newcommand{\q}{\mbox{\boldmath $x$}}
\newcommand{\bsigma}{\boldsymbol{\sigma}}
\newcommand{\bdelta}{\boldsymbol{\delta}}
\newcommand{\sqr}[1]{{\left(#1\right)}_{1/2}}
\newcommand{\as}{\\[.6em]}

\newcommand{\AS}{\\[1.2em]}

\newcommand{\dis}{\displaystyle}

\newcommand{\ignore}[1]{}
\newcommand{\jac}[1]{\operatorname{#1}}

\title{On a discretization of confocal quadrics. II.\linebreak A geometric approach to general parametrizations}
\author{Alexander I. Bobenko$^1$, Wolfgang K. Schief$^2$,\\ Yuri B. Suris$^1$, Jan Techter$^1$ \bigskip\\  
$^1$Institut f\"ur Mathematik, TU Berlin, \\ Str.\@ des 17.\@ Juni 136, 10623 Berlin, Germany\bigskip\\
$^2$School of Mathematics and Statistics,\\ The University of New South Wales, Sydney, NSW 2052, Australia}
\date{\today}

\begin{document}

\maketitle

\begin{abstract}
We propose a discretization of classical confocal coordinates. It is based on a novel characterization thereof as factorizable orthogonal coordinate systems. Our geometric discretization leads to
factorizable discrete nets with a novel discrete analog of the orthogonality property. A discrete confocal coordinate system may be constructed geometrically via polarity with respect to a sequence of classical confocal quadrics. Various sequences correspond to various discrete parametrizations. 
The coordinate functions of discrete confocal quadrics are computed explicitly.
The theory is illustrated with a variety of examples in two and three dimensions. These include confocal coordinate systems parametrized in terms of Jacobi elliptic functions. Connections with incircular (IC) nets and a generalized Euler-Poisson-Darboux system are established.
\end{abstract}

\newpage
\tableofcontents
\newpage

\section{Introduction}

Confocal quadrics have played a prominent role in classical mathematics due to their beautiful geometric properties and numerous relations and applications to various branches of mathematics. Optical properties of quadrics and their confocal families were already discovered by the ancient Greeks and continued to fascinate mathematicians for many centuries, culminating in the famous Ivory and Chasles theorems from 19th century given a modern interpretation by Arnold~\cite[Appendix 15]{Arnold}. 
Geodesic flows on quadrics and billiards in quadrics are classical examples of integrable systems \cite{jacobi, Moser, Veselov, Fuchs-Tabachnikov}.
Gravitational properties of ellipsoids were studied in detail by Newton, Ivory and others, see \cite[Part 8]{Fuchs-Tabachnikov}, and are based to a large extent on the geometric properties of confocal quadrics. Quadrics in general and confocal systems of quadrics in particular constitute popular objects in geometry. Poncelet and Ivory theorems play a central role there \cite{Dragovic-Radnovic, Izmestiev-Tabachnikov}.
In differential geometry quadrics provide non-trivial examples of isothermic surfaces which form one of the most interesting classes of ``integrable'' surfaces, that is, surfaces which are governed by integrable differential equations and possess a rich theory of transformations with remarkable permutability properties \cite{bobenko-suris}. Importantly, confocal quadrics also lie at the heart of confocal coordinates systems which give rise to separation of variables in the Laplace operator. As such, they support a rich theory of special functions, including Lam\'e functions and their generalizations \cite{Bateman-Erdelyi, Whittaker-Watson}. 

%{\bf orthogonal coordinate sytems: smooth and discrete}
 
In general, coordinate systems are instances of {\em smooth nets}, that is, maps $\mathbb R^M\supset U\to\mathbb R^N$.
In this paper we present a novel characterization  of confocal coordinate system: a coordinate system $\mathbb R^N \supset  U\ni\bs \mapsto \q(\bs)\in \mathbb R^N$ is confocal if and only if it is orthogonal and the coordinates $x_i$ factorize as functions of the parameters $s_j$, that is,
$$
 x_i(\bs)=f_1^{i}(s_1)f_2^{i}(s_2)\cdots f_N^{i}(s_N),\quad i = 1,\ldots,N
$$
(see Section~\ref{s.characterization}). 

Orthogonal coordinate systems constitute a classical topic in differential geometry. They were extensively treated in the fundamental monograph by Darboux \cite{darboux}.
From the viewpoint of the theory of
integrable systems they were investigated in 
\cite{Zakharov_1998}. Algebro-geometric orthogonal
coordinate systems were constructed in \cite{Krichever_1997}. Although it is natural to expect that confocal coordinate systems belong to this class, it remains an open problem to include them in Krichever's construction (see \cite{Mironov-Taimanov}).  

Discretizing coordinate systems consists of finding suitable approximating {\em discrete nets}, that is, maps $\mathbb Z^M\supset U\to \mathbb R^N$. Various discretizations of orthogonal coordinate systems have been proposed. The most investigated variant is the class of circular nets  \cite{Bobenko_1999, Cieslinski_Doliwa_Santini_1997, KS1998}, where all elementary quadrilaterals are inscribed in circles. This class inherits the property of orthogonal coordinate systems to be invariant under M\"obius transformations. A special case of Darboux-Egorov metrics was discretized in \cite{Akhmetshin_Volvovski_Krichever_1999} as circular nets whose quadrilaterals have two opposite right angles. Another discretization of orthogonal nets is given by conical nets \cite{Liu_Pottmann_Wallner_Yang_Wang_2006}, which are characterized by the property that any four neighboring planar quadrilaterals are tangent to a sphere. This class is preserved by Laguerre transformations. Circular and conical nets may be unified in the context of Lie geometry as principal contact element nets \cite{bobenko-suris, Pottmann_Wallner_2008}.

Recently, in \cite{BSST16}, we have proposed an integrability-preserving discretization of systems of confocal quadrics or, equivalently, systems of confocal coordinates in $\R^N$. The discretization is based on a discrete version \cite{konopelchenko-schief} of the Euler-Poisson-Darboux system, which is known to encode algebraically classical systems of confocal quadrics. This algebraic approach resulted in particular in a new geometric orthogonality condition, which is of central importance in the current paper. In this new sense, a discrete coordinate system  ${(\frac{1}{2}\mathbb Z)}^N\supset U\ni\bn \mapsto \q(\bn)\in \mathbb R^N$ is discrete orthogonal if any edge of any of the $\Z^N$ sublattices is orthogonal to the dual facet of the dual $\Z^N$ sublattice (see Section~\ref{s.orthogonality}). 

In Section~\ref{s.discrete-confocal} we define {\em discrete confocal coordinate systems} as orthogonal (in the new sense) nets
${(\frac{1}{2}\mathbb Z)}^N\supset U\ni\bn \mapsto \q(\bn)\in \mathbb R^N$ such that the coordinates $x_i$ factorize as functions of the lattice arguments $n_j$, that is,
$$
 x_i(\bn)=f_1^{i}(n_1)f_2^{i}(n_2)\cdots f_N^{i}(n_N),\quad i = 1,\ldots,N.
$$
We provide an explicit description of discrete confocal coordinate systems in Theorems~\ref{Th discr confocal}, \ref{th resolution}, \ref{th funct eqs}.

In Section~\ref{s.geometric} we show that discrete confocal coordinate systems admit a geometric characterization in terms of polarity with respect to quadrics of a classical confocal family. The connection with the particular case discussed in \cite{BSST16} is set down in Section~\ref{s.article1}. 

Sections~\ref{section_2d}--\ref{section_3d} contain an extensive collection of examples of discrete confocal coordinates in the cases $N=2$ and $N=3$. We begin by presenting in Section \ref{s: 2d trig/hyp} the discrete analogue of the classical parametrization of systems of confocal conic sections in terms of trigonometric and hyperbolic functions. Then, in Sections \ref{sect elliptic IC}--\ref{sect hyperbolic IC} we record a novel parametrization of confocal coordinate systems in $\R^2$, both continuous and discrete, in terms of Jacobi elliptic functions. The discrete confocal coordinate systems of these families are intimately related to incircular (IC) nets studied in \cite{Akopyan-Bobenko}. In Sections \ref{s: akopyan}--\ref{s: 2 concentric}, several geometrically remarkable parametrizations are discretized, related to 3-webs comprized by conic sections, circles and lines. In Section \ref{section_3d}, we show that the classical confocal coordinate systems in $\R^3$ parametrized in terms of Jacobi elliptic functions admit a natural discrete analogue.  Finally, in the Appendix, we present a generalized discrete Euler-Poisson-Darboux systems which algebraically encodes discrete confocal coordinate systems.

\paragraph{Acknowledgement}

This research was supported by the DFG Collaborative Research Center TRR 109 ``Discretization in Geometry and Dynamics''. W.S. was also supported by the Australian Reserach Council (DP1401000851).

%%%%%%%%%%%%%%%%%%%%%%%%%%%%%%%%%%%

\section{Classical confocal coordinates}

For given $a_1 > a_2 > \cdots > a_N > 0$,
we consider the one-parameter family of confocal quadrics in $\R^N$ given by
\begin{equation}   \label{eq:confocal-family}
 Q(\lambda)=\left\{ \q=(x_1,\ldots,x_N)\in\mathbb R^N:\ \sum_{k=1}^N \frac{x_k^2}{a_k + \lambda} = 1\right\}, \quad \lambda \in \R.
\end{equation}
Note that the quadrics of this family are centered at the origin and have the principal axes aligned along the coordinate directions.
For a given point $\q=(x_1, \ldots, x_N) \in \R^N$ with $x_1x_2\ldots x_N\neq 0$, equation $\sum_{k=1}^N x_k^2/(a_k + \lambda) = 1$ is, after clearing the denominators, a polynomial equation of degree $N$ in $\lambda$, with $N$ real roots $u_1,\ldots, u_N$ lying in the intervals
\begin{equation}\label{eq: intervals}
  -a_1 < u_1 < -a_2 < u_2 < \cdots < -a_N < u_N,
\end{equation}
so that
\begin{equation} \label{eq: equation factorized}
  \sum_{k=1}^N \frac{x_k^2}{\lambda+a_k}  - 1= - \frac{\prod_{m=1}^N (\lambda - u_m)}{\prod_{m=1}^N (\lambda+a_m) }.
\end{equation}
These $N$ roots correspond to the $N$ confocal quadrics of the family \eqref{eq:confocal-family} that intersect at the point $\q=(x_1, \ldots, x_N)$:
\begin{equation}   \label{eq:confocal-quadrics}
  \sum_{k=1}^N \frac{x_k^2}{a_k + u_i} = 1, \quad i=1,\ldots,N \quad\Leftrightarrow\quad \q\in \bigcap_{i=1}^N Q(u_i).
\end{equation}
The $N$ quadrics $Q(u_i)$ are all of different signatures. Evaluating the residue of the right-hand side of \eqref{eq: equation factorized} at $\lambda=-a_k$, one can easily express $x_k^2$ through $u_1, \ldots, u_N$:
\begin{equation} \label{eq:elliptic-coordinates-squares}
  x_k^2 = \frac{\prod_{i=1}^N (u_i+a_k)}{\prod_{i \neq k} (a_k - a_i)}, \quad k = 1, \ldots, N.
 \end{equation}
Thus, for each point $(x_1, \ldots, x_N) \in \R^N$ with $x_1x_2\ldots x_N\neq 0$, there is exactly one solution $(u_1, \ldots, u_N) \in \mathcal{U}$ of \eqref{eq:elliptic-coordinates-squares},
where
\begin{equation}
  \mathcal{U} = \set{ (u_1,\ldots,u_N) \in \R^N }{ -a_1 < u_1 < -a_2 < u_2 < \ldots < -a_N < u_N }.
  \label{eq:domain}
\end{equation}
On the other hand, for each $(u_1, \ldots, u_N) \in \mathcal{U}$ there are exactly $2^N$ solutions $(x_1, \ldots, x_N) \in \R^N$,
which are mirror symmetric with respect to the coordinate hyperplanes. In what follows, when we refer to a solution of \eqref{eq:elliptic-coordinates-squares}, we always mean the solution with values in
\begin{equation*}
  \R_+^N = \set{ (x_1, \ldots, x_N) \in \R^N }{ x_1 > 0, \ldots, x_N > 0 }.
\end{equation*}
Thus, we are dealing with a parametrization of the first hyperoctant of $\R^N$,
$\q : \mathcal{U}\ni (u_1, \ldots, u_N) \mapsto (x_1,\ldots,x_N)\in\R_+^N$, given by
\begin{equation} \label{eq:elliptic-coordinates}
  x_k = \frac{\prod_{i=1}^{k-1} \sqrt{-(u_i+a_k)}\prod_{i=k}^N\sqrt{u_i+a_k}}{\prod_{i = 1}^{k-1} \sqrt{a_i - a_k}\prod_{i=k+1}^N\sqrt{a_k-a_i}}, 
  \quad k = 1, \ldots, N,
 \end{equation}
such that the coordinate hyperplanes $u_i = {\rm const}$ are mapped to (parts of) the respective quadrics given by \eqref{eq:confocal-quadrics}.
The coordinates $(u_1,\ldots,u_N)$ are called \emph{confocal coordinates} (or {\em elliptic coordinates}, following Jacobi \cite[Vorlesung 26]{jacobi}).

For various applications, it is often useful to re-parametrize the coordinate lines according to $u_i=u_i(s_i)$, $i=1,\ldots,N$. One of the reasons of the usefulness of this procedure is the possibility to uniformize the square roots in the above formulas, that is, to present them as single-valued functions of the new coordinates $s_i$. A classical example in the dimension $N=2$, where
$$
x_1=\frac{\sqrt{u_1+a_1}\sqrt{u_2+a_1}}{\sqrt{a_1-a_2}}, \quad x_2=\frac{\sqrt{-(u_1+a_2)}\sqrt{u_2+a_2}}{\sqrt{a_1-a_2}},
$$
is to set
$$
u_1=-a_1\sin^2 s_1-a_2\cos^2 s_1, \quad u_2=a_1\sinh^2 s_2-a_2\cosh^2 s_2,
$$
so that
$$
u_1+a_1=(a_1-a_2)\cos^2 s_1, \quad -(u_1+a_2)=(a_1-a_2)\sin^2 s_1,
$$
and
$$
u_2+a_1=(a_1-a_2)\cosh^2 s_2, \quad u_2+a_2=(a_1-a_2)\sinh^2 s_2.
$$
Accordingly, one obtains a version of elliptic coordinates in the plane free from branch points and naturally periodic with respect to $s_1$:
$$
x_1=\sqrt{a_1-a_2}\cos s_1\cosh s_2, \quad x_2=\sqrt{a_1-a_2}\sin s_1\sinh s_2.
$$
Such a re-parametrization, being a relatively trivial operation for classical coordinate systems, does not have a simple counterpart in the discrete context. Actually, the lack of the notion of a re-parametrization is one of the main and fundamental differences between discrete differential geometry and discrete analysis, on the one hand, and their classical analogs, on the other hand. It is one of the principal goals of this paper to present a natural geometric construction of a general parametrization  for discrete confocal coordinate systems. 

%%%%%%%%%%%%%%%%%%%%%%%%%%%%%%%

\section{Characterization of confocal coordinate systems}
\label{s.characterization}

Our main subject in this paper are coordinate systems, i.e., maps $\q:\mathbb R^N\supset U \to\mathbb R^N$ on open sets $U$ such that $\det(\partial x_i/\partial s_j)_{i,j=1}^N\neq 0$. We now demonstrate that the two properties, factorization and orthogonality, are sufficient  to characterize confocal coordinates. For the sake of simplicity, we restrict ourselves to coordinate systems satisfying the additional condition $\partial x_i/\partial s_j\neq 0$, which excludes degenerate cases like cylindric, spherical coordinates etc.
  
\begin{theorem} \label{theorem confocal charact N}
If a coordinate system $\q:\mathbb R^N\supset U\to\mathbb R^N$ satisfies two conditions:
\begin{enumerate}
\item[{\rm i)}] $\q(\bs)$ {\em factorizes}, in the sense that
\begin{equation}\label{factorized N}
\renewcommand{\arraystretch}{1.5}
\left\{ \begin{array}{l} x_1(\bs)=f_1^1(s_1)f_2^1(s_2)\cdots f_N^1(s_N), \\ 
x_2(\bs)=f_1^2(s_1)f_2^2(s_2)\cdots f_N^2(s_N), \\ \ldots \\
x_N(\bs)=f_1^N(s_1)f_2^N(s_2)\cdots f_N^N(s_N),
\end{array} \right.
\end{equation}
with all $f_i^k(s_i)\neq 0$ and $\big(f_i^k\big)'(s_i)\neq 0$;
\item[{\rm ii)}] $\q$ is {\em orthogonal}, that is, 
\begin{equation}\label{ortho N}
\langle \partial_i \q , \partial _j \q\rangle=0 \quad for \quad i\neq j,
\end{equation}
\end{enumerate}
then all coordinate hypersurfaces are confocal quadrics.
\end{theorem}
\begin{proof}
One easily computes that the orthogonality condition \eqref{ortho N} for a factorized net \eqref{factorized N} is equivalent to
$$
\sum_{k=1}^N \big(f_i^k(s_i)\big)'f_i^k(s_i)\big(f_j^k(s_j)\big)'f_j^k(s_j)\prod_{\ell\neq i,j}\big(f_\ell^k(s_\ell)\big)^2=0,
$$
or
\begin{equation} \label{F'F'F N}
\sum_{k=1}^N \big(F_i^k(s_i)\big)'\big(F_j^k(s_j)\big)'\prod_{\ell\neq i,j}F_\ell^k(s_\ell)=0,
\end{equation}
where $F_i^k(s_i)=\big(f_i^k(s_i)\big)^2$. 
\begin{lemma} \label{lemma 1}
Equation \eqref{F'F'F N} is equivalent to
\begin{equation}\label{FFF N}
\sum_{k=1}^NF_1^k(s_1)F_2^k(s_2) \cdots F_N^k(s_N)=A_1(s_1)+A_2(s_2)+ \ldots +A_N(s_N),
\end{equation}
with some functions $A_i(s_i)$. 
\end{lemma}
\begin{proof}
Equation \eqref{F'F'F N} reads: $\partial^2 F/\partial s_i\partial s_j=0$, where the function $F$ is the left-hand side of \eqref{FFF N}. Induction with respect to $N$ shows that this is equivalent to $F$ being a sum of functions of single arguments.
\end{proof}

\begin{lemma}\label{lemma 2}
Assume that all $f_i^k\neq 0$ and $\big(f_i^k\big)'\neq 0$. Then, for each $i=1,\ldots,N$ there exists a function $F_i(s_i)$ such that 
\begin{equation}\label{GH thru F N}
F_i^k(s_i)=\alpha_i^kF_i(s_i)+\beta_i^k,  \quad k=1,\ldots, N,
\end{equation}
for some constants $\alpha_i^k\neq 0$ and $\beta_i^k$.
\end{lemma}
\begin{proof} Note that the assumption of lemma is equivalent to $\big(F_i^k\big)'\neq 0$. We will prove that for each $i=1,\ldots,N$ we have
\begin{equation}\label{F'G' N}
\begin{pmatrix}\big(F_i^1\big)'(s_i) \\ \ldots \\ \big(F_i^N\big)'(s_i)\end{pmatrix}\in \mathbb R\begin{pmatrix}\alpha_i^1 \\ \ldots \\ \alpha_i^N\end{pmatrix}.
\end{equation}
For any fixed $i$, equation \eqref{F'F'F N} can be formulated as the following $N-1$ orthogonality conditions:
\begin{equation}\label{N-1 vectors}
\renewcommand{\arraystretch}{1.5}
\begin{pmatrix} \big(F_i^1\big)' \\ \ldots \\ \big(F_i^N\big)' \end{pmatrix} \perp
\begin{pmatrix} \big(F_j^1\big)'\prod_{\ell\neq i,j} F_\ell^1  \\ \ldots \\  \big(F_j^N\big)'\prod_{\ell\neq i,j} F_\ell^N 
\end{pmatrix}, \quad j\neq i.
\end{equation}
Multiplying the $N-1$ vectors on the right-hand side from the left by the non-degenerate matrix
$$
{\rm diag}\Big(\prod_{\ell\neq i} F_\ell^1, \ldots,  \prod_{\ell\neq i} F_\ell^N \Big)^{-1},
$$
we obtain vectors
$$
\renewcommand{\arraystretch}{1.5}
\begin{pmatrix} \big(F_j^1\big)'/F_j^1  \\ \ldots \\  \big(F_j^N\big)'/ F_j^N 
\end{pmatrix}=
\frac{1}{2}\begin{pmatrix} \big(f_j^1\big)'/f_j^1  \\ \ldots \\  \big(f_j^N\big)'/ f_j^N 
\end{pmatrix}, \quad j\neq i,
$$
Multiplying the latter vectors from the left by the non-degenerate matrix
$$
{\rm diag}\Big(\prod_{\ell=1}^N f_\ell^1, \ldots,  \prod_{\ell=1}^N f_\ell^N \Big),
$$
we obtain vectors
$$
\renewcommand{\arraystretch}{1.5}
\frac{1}{2}\begin{pmatrix} \big(f_j^1\big)'\prod_{\ell\neq j} f_\ell^1  \\ \ldots \\  \big(f_j^N\big)'\prod_{\ell\neq j} f_\ell^N 
\end{pmatrix}=\frac{1}{2}\partial_j\q, \quad j\neq i.
$$
The latter vectors are linearly independent and span an $(N-1)$-dimensional subspace of $\mathbb R^N$.  Thus, the vector on the left-hand side of \eqref{N-1 vectors} lies in the orthogonal complement of an $(N-1)$-dimensional subspace which is manifestly independent of $s_i$. This orthogonal complement is a one-dimensional space which does not depend on $s_i$. This proves \eqref{F'G' N}.
\end{proof}

%(In the case $N=3$ we have set $F_i(u_i)=F_i^1(u_i)$, so that $\alpha_i^1=1$, $\beta_i^1=0$, and then $\alpha_i^2=\alpha_i$, $\beta_i^2=\beta_i$, and $\alpha_i^3=\gamma_i$, $\beta_i^3=\delta_i$). 
Substituting \eqref{GH thru F N} into the left-hand side of equation \eqref{FFF  N}, we arrive at an expression which may be represented as the polynomial 
$$
\sum_{k=1}^N(\alpha_1^kz_1+\beta_1^k)\cdots(\alpha_N^kz_N+\beta_N^k)
$$
of degree $N$ in $N$ formal variables $z_1, \ldots,z_N$, evaluated at $z_i=F_i(s_i)$. It is easy to deduce that the result is a sum of functions of single variables, as in \eqref{FFF N}, if and only if in the above polynomial all monomials of degree $\ge 2$ vanish, leaving us with
\begin{equation}\label{polynom z}
\sum_{k=1}^N\big(\alpha_1^kz_1+\beta_1^k\big)\cdots\big(\alpha_N^kz_N+\beta_N^k\big)=\sum_{i=1}^N \rho_i z_i+c.
\end{equation}
We can identify the coefficients of the monomials of degree $\le 1$:
\begin{equation}\label{rho N}
\rho_i=\sum_{k=1}^N \alpha_i^k\prod_{\ell\neq i} \beta_{\ell}^k, \qquad c=\sum_{k=1}^N \prod_{\ell=1}^N\beta_\ell^k,
\end{equation}
while the vanishing of the coefficients of all monomials of degree $\ge 2$ can be expressed as a certain set of equations for the coefficients $\alpha_i^k$, $\beta_i^k$. As a result,  \eqref{FFF N} adopts the concrete form
\begin{equation}
\sum_{k=1}^N\big(\alpha_1^kF_1(s_1)+\beta_1^k\big)\cdots\big(\alpha_N^kF_N(s_N)+\beta_N^k\big)=\sum_{i=1}^N \rho_iF_i(s_i)+c.
\end{equation}
It turns out that the above mentioned equations for the coefficients $\alpha_i^k$, $\beta_i^k$ imply certain identities involving functions of $N-1$ variables.

\begin{lemma} \label{lemma coord surfaces}
The following formulas hold true for all $i=1,2,\ldots,N$:
\begin{equation}\label{for coord surfaces N}
\sum_{k=1}^N \alpha_i^k \prod_{\ell\neq i} F_\ell^k(s_\ell)=\rho_i.
\end{equation}
\end{lemma}
\begin{proof}
Differentiate equation
$$
\sum_{k=1}^N \prod_{\ell=1}^N F_\ell^k(s_\ell)=\sum_{i=1}^N \rho_i F_i(s_i)+c
$$
with respect to $s_i$:
$$
\sum_{k=1}^N \big(F_i^k(s_i)\big)' \prod_{\ell\neq i} F_\ell^k(s_\ell)=\rho_i\big(F_i(s_i)\big)'.
$$
Taking into account equation \eqref{GH thru F N} and dividing by $F_i'\neq 0$, we arrive at \eqref{for coord surfaces N}.
\end{proof}
Equations  \eqref{for coord surfaces N} describe coordinate hypersurfaces $s_i={\rm const}$. Indeed, observe that, according to \eqref{factorized N} and to $F_i^k=\big(f_i^k\big)^2$,  these equations can be expressed as quadrics:
$$
\sum_{k=1}^N \frac{\alpha_i^k}{F_i^k(s_i)}x_k^2=\rho_i,
$$
or, equivalently, due to \eqref{GH thru F N},
\begin{equation}\label{coord surfaces N}
\sum_{k=1}^N \frac{x_k^2}{\rho_i F_i(s_i)+\rho_i\beta_i^k/\alpha_i^k}=1.
\end{equation}
In order to show that these quadrics for all $i=1,\ldots, N$ and for any values of $s_i$  belong to a confocal family
$$
\sum_{k=1}^N \frac{x_k^2}{\lambda+a_k}=1,
$$
it remains to show that, for any two indices $k\neq m$ from $1,\ldots,N$, the expressions 
$$
\rho_i\left(\frac{\beta_i^k}{\alpha_i^k}-\frac{\beta_i^m}{\alpha_i^m}\right),
$$
which should be equal to $a_k-a_m$,  do not depend on $i$. Upon setting in identity \eqref{polynom z}
$$
z_1=-\frac{\beta_1^{(k_1)}}{\alpha_1^{(k_1)}}, \quad \ldots, \quad z_N=-\frac{\beta_N^{(k_N)}}{\alpha_N^{(k_N)}}
$$
for an arbitrary permutation $(k_1,\ldots,k_N)$ of $(1,\ldots,N)$, we arrive at
$$
\rho_1\frac{\beta_1^{(k_1)}}{\alpha_1^{(k_1)}}+\ldots+\rho_N\frac{\beta_N^{(k_N)}}{\alpha_N^{(k_N)}}=c.
$$
Subtracting two such equations for two permutations, differing at only two positions $i$, $j$, where they take values $k$, $m$ and $m$, $k$, respectively, we arrive at
$$
\rho_{i}\frac{\beta_{i}^{(k)}}{\alpha_{i}^{(k)}}+\rho_{j}\frac{\beta_{j}^{(m)}}{\alpha_{j}^{(m)}}=
\rho_{i}\frac{\beta_{i}^{(m)}}{\alpha_{i}^{(m)}}+\rho_{j}\frac{\beta_{j}^{(k)}}{\alpha_{j}^{(k)}},
$$ 
or
$$
\rho_{i}\left(\frac{\beta_{i}^{(k)}}{\alpha_{i}^{(k)}}-\frac{\beta_{i}^{(m)}}{\alpha_{i}^{(m)}}\right)=
\rho_{j}\left(\frac{\beta_{j}^{(k)}}{\alpha_{j}^{(k)}}-\frac{\beta_{j}^{(m)}}{\alpha_{j}^{(m)}}\right).
$$ 
This is the desired result, since $i$ and $j$ are arbitrary.
\end{proof}

We have demonstrated that the equations of the coordinate hypersurfaces of a factorized orthogonal coordinate system \eqref{factorized N} can be put as
\begin{equation}\label{eq: confocal coords hypersurfaces}
\sum_{k=1}^N \frac{x_k^2}{u_i+a_k}=1, \quad i=1,\ldots,N,
\end{equation}
where the parameters $a_k$ are given by
\begin{equation}\label{eq: ak}
a_k=\rho_i\frac{\beta_i^k}{\alpha_i^k}+c_i, \quad k=1,\ldots,N,
\end{equation}
with suitable constants $c_i$ (which ensure that the right-hand side of \eqref{eq: ak} does not depend on $i$), while the quantities 
\begin{equation}
u_i=u_i(s_i)=\rho_iF_i(s_i)-c_i,\quad i=1,\ldots,N
\end{equation} 
can be considered as the confocal coordinates of the points $\q(\bs)$. We remark that, since confocal quadrics of the same signature do not intersect, the $N$ confocal coordinates should belong to $N$ disjoint intervals \eqref{eq: intervals} (possibly, upon a re-numbering). 

Conversely, we know that for any confocal coordinate system, equation \eqref{eq: confocal coords hypersurfaces} is equivalent to
\begin{equation} \label{eq: confocal coords squares}
  x_k^2 = \frac{\prod_{i=1}^N (u_i+a_k)}{\prod_{i \neq k} (a_k - a_i)}, \quad k = 1, \ldots, N,
 \end{equation}
 and positivity of these expressions is equivalent to \eqref{eq: intervals}. Thus, formulas for $x_k$ contain square roots $\sqrt{\pm(u_i+a_k)}$ (see \eqref{eq:elliptic-coordinates}). Suppose that these square roots are uniformized by the re-parametrization
 \begin{equation}  \label{eq: ff=u+a smooth}
  \renewcommand{\arraystretch}{1.4}
  \big(f_i^k(s_i)\big)^2=\left\{\begin{array}{ll} u_i+a_k, & k\le i,\\ 
                                                     -\big(u_i+a_k\big), & k>i.
                                                   \end{array}\right.                            
\end{equation}
The latter equations are consistent, if for any $1\le i\le N$ the squares of the functions $f_i^k(s_i)$, $1\le k\le N$, satisfy a system of $N-1$ linear equations:
\begin{equation}\label{eq: squares smooth}
  \renewcommand{\arraystretch}{1.6}
\left\{
\begin{array}{l}
\big(f_i^1(s_i)\big)^2-\big(f_i^k(s_i)\big)^2= a_1-a_k, \quad k\le i, \\
\big(f_i^1(s_i)\big)^2+\big(f_i^k(s_i)\big)^2= a_1-a_k, \quad k> i.
\end{array} \right.
\end{equation}
Under such a re-parametrization, formulas for confocal coordinates can be written as
\begin{equation}\label{eq: x=fff smooth}
x_k(\bs)= \frac{\prod_{j=1}^N f_j^k(s_j)}{\prod_{i = 1}^{k-1} \sqrt{a_i - a_k}\prod_{i=k+1}^N\sqrt{a_k-a_i}}, \quad k = 1, \ldots, N
\end{equation}
and, hence, the coordinate system factorizes. Note that \eqref{eq: x=fff smooth} is equivalent to \eqref{factorized N} modulo a scaling of the functions $f_i^k$.

%%%%%%%%%%%%%%%%%%%%%%%%%%%%%%%%%%%%%%%%

\section{Discrete orthogonality}
\label{s.orthogonality}

\begin{figure}[H]
  \centering
  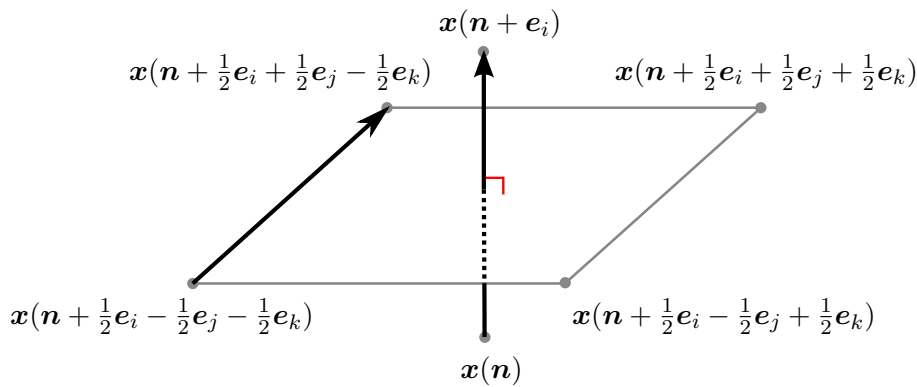
  \caption{Discrete orthogonality in dimension $N=3$.}
  \label{fig:orthogonality}
\end{figure}

We will use the discrete version of the characteristic properties from Theorem \ref{theorem confocal charact N} to define discrete confocal coordinate systems. These will be special nets defined on the square lattice of stepsize $1/2$,
\begin{equation}\label{eq: dense net}
\q: (\tfrac{1}{2}\mathbb Z)^N\supset \mathcal U \to\mathbb R^N.
\end{equation}
The suitable notion of orthogonality is a novel one, introduced in \cite{BSST16}. We denote by $\be_i$ the unit vector of the coordinate direction $i$.
\begin{definition}\label{def ortho}
A net \eqref{eq: dense net} is called {\em orthogonal} if for each edge $[\bn,\bn+\be_i]$, all $2^{N-1}$ vertices of the dual facet,
$$
\q(\bn+\tfrac{1}{2}\bsigma) \;\;{ for\;\; all} \;\;\bsigma=(\sigma_1,\ldots,\sigma_N)\in\{\pm 1\}^N \;\;{ with} \;\; \sigma_i=1,
$$ 
lie in a hyperplane orthogonal to the line $(\q(\bn),\q(\bn+\be_i))$ (see Fig. \ref{fig:orthogonality}).
\end{definition}
\bigskip

Note that the original definition from \cite{BSST16} referred to pairs of nets defined on two dual lattices $\mathbb Z^N$ and $(\mathbb Z+\tfrac{1}{2})^N$. The lattice $(\tfrac{1}{2}\mathbb Z)^N$ contains $2^{N-1}$ pairs of dual sublattices of this type, namely
$$
\mathbb Z^N+\tfrac{1}{2}\bdelta\quad {\rm and} \quad \mathbb Z^N+\tfrac{1}{2}\bar{\bdelta},
$$
for any $\bdelta=(\delta_1,\ldots,\delta_N)\in\{0,1\}^N$ and $\bar{\bdelta}=(1-\delta_1,\ldots,1-\delta_N)\in\{0,1\}^N$.

\begin{proposition}
All elementary quadrilaterals
\begin{equation}\label{elem quad}
(\q(\bn),\q(\bn+\be_j),\q(\bn+\be_j+\be_k),\q(\bn+\be_k))
\end{equation}
of a generic orthogonal net are planar.
\end{proposition}
\begin{proof} An elementary quadrilateral \eqref{elem quad} can be considered as the intersection of $N-2$ facets dual to the edges 
$$
\big[\q(\bn+\tfrac{1}{2}\sum_{\ell\neq i}\be_\ell-\tfrac{1}{2}\be_i),\q(\bn+\tfrac{1}{2}\sum_{\ell\neq i}\be_\ell+\tfrac{1}{2}\be_i)\big],  \quad i\neq j,k
$$ 
of the dual sublattice. Each of these $N-2$ facets lies in a hyperplane. The intersection of $N-2$ hyperplanes in $\mathbb R^N$ is generically a two-dimensional plane. 
\end{proof}

Clearly, the definition of orthogonality can be equivalently formulated as follows: the two lines containing any pair of dual edges are orthogonal:
\begin{equation}\label{eq: dual edges}
\big(\q(\bn),\q(\bn+\be_i)\big)\perp \big(\q(\bn+\tfrac{1}{2}\bsigma), \q(\bn+\tfrac{1}{2}\bsigma+\be_j)\big),
\end{equation}
where $\bsigma\in\{\pm 1\}^N$ is any $N$-tuple of signs with $\sigma_i=1$ and $\sigma_j=-1$ (and orthogonality is understood in the sense of orthogonality of the direction vectors). 
From this it is easy to see that pairs of dual sublattices actually play symmetric roles in the definition of orthogonality.

\section{Discrete confocal coordinate systems}
\label{s.discrete-confocal}

For discrete nets $\q:\mathbb Z^N\supset \mathcal U\to\mathbb R^N$, at any point $\bn\in\mathcal U$ and for any coordinate direction $j=1,\ldots,N$, there exist two natural discrete tangent vectors,  $\Delta_j\q(\bn)=\q(\bn+\be_j)-\q(\bn)$ and $\bar\Delta_j\q(\bn)=\q(\bn)-\q(\bn-\be_j)$. We call such a net a {\em discrete coordinate system} if at any $\bn\in \mathcal U$, the $N$ discrete tangent vectors (arbitrarily chosen among $\Delta_j\q(\bn)$ and $\bar\Delta_j\q(\bn)$ for any $j$) are linearly independent. 

A net $\q:\left(\tfrac{1}{2}\mathbb Z\right)^N\to\mathbb R^N$ defined on the lattice of a half stepsize can be considered as consisting of $2^N$ subnets defined on sublattices $\mathbb Z^N+\tfrac{1}{2}\bdelta$ for $\bdelta\in\{0,1\}^N$, and we call it a discrete coordinate system if all $2^N$ subnets satisfy the above condition.

\begin{definition}
A discrete coordinate system $\q:\left(\tfrac{1}{2}\mathbb Z\right)^N\supset \mathcal U\to\mathbb R^N$ is called {\em a discrete confocal coordinate system} if it satisfies two conditions:
\begin{enumerate}
\item[{\rm i)}] $\q(\bn)$ {\em factorizes}, in the sense that for any $\bn\in\mathcal U$
\begin{equation}\label{d factorized}
\renewcommand{\arraystretch}{1.5}
\left\{ \begin{array}{l} x_1(\bn)=f_1^1(n_1)f_2^1(n_2)\cdots f_N^1(n_N), \\ x_2(\bn)=f_1^2(n_1)f_2^2(n_2)\cdots f_N^2(n_N), \\ 
\ldots \\ x_N(\bn)=f_1^N(n_1)f_2^N(n_2)\cdots f_N^N(n_N),
\end{array} \right.
\end{equation}
with $f_i^k(n_i)\neq 0$ and $\bar\Delta f_i^k(n_i)=f_i^k(n_i)-f_i^k(n_i-1)\neq 0$;
\item[{\rm ii)}] $\q$ is {\em orthogonal} in the sense of Definition \ref{def ortho}.
\end{enumerate}
\end{definition}
\begin{theorem}\label{Th discr confocal}
For a discrete confocal coordinate system, there exist $N$ real numbers $a_k$, $1\le k\le N$, and $N$ sequences $u_i:\tfrac{1}{2}\mathbb Z+\tfrac{1}{4}\to\mathbb R$ such that the following equations are satisfied for any $\bn\in\mathcal U$ and for any $\bsigma\in\{\pm 1\}^N$:
\begin{equation}\label{discr confocal gen}
\sum_{k=1}^N \frac{x_k(\bn)x_k(\bn+\frac{1}{2}\bsigma)}{a_k+u_i}=1,\quad u_i=u_i(n_i+\tfrac{1}{4}\sigma_i), \quad i=1,\ldots,N.
\end{equation}
Equivalently, 
\begin{equation}\label{discr confocal gen x thru u}
x_k(\bn)x_k(\bn+\tfrac{1}{2}\bsigma)= \frac{\prod_{j=1}^N (u_j+a_k)}{\prod_{j \neq k} (a_k - a_j)}, \quad u_j=u_j(n_j+\tfrac{1}{4}\sigma_j), \quad k = 1, \ldots, N.
\end{equation}
\end{theorem}
\begin{proof}
Orthogonality condition \eqref{eq: dual edges} written in full reads:
\begin{eqnarray}\label{eq: d ortho long}
\lefteqn{\sum_{k=1}^N(f_i^k(n_i+1)-f_i^k(n_i))f_i^k(n_i+\tfrac{1}{2})\ \cdot}\nonumber\\
&&(f_j^k(n_j+\tfrac{1}{2})-f_j^k(n_j-\tfrac{1}{2}))f_j^k(n_j)\cdot\prod_{\ell\neq i,j} f_\ell^k(n_\ell)f_\ell^k(n_\ell+\tfrac{1}{2}\sigma_\ell)=0.\nonumber\\
\end{eqnarray}
We introduce the quantities
\begin{equation}
F_i^k(n_i+\tfrac{1}{4})=f_i^k(n_i)f_i^k(n_i+\tfrac{1}{2}),
\end{equation}
assigned to the points of the lattice $\tfrac{1}{2}\mathbb Z +\tfrac{1}{4}$, and the difference operator
\begin{equation}
\Delta^{\nicefrac{1}{2}} F(n)=F(n+\tfrac{1}{4})-F(n-\tfrac{1}{4}).
\end{equation}
With this notation, relation \eqref{eq: d ortho long} takes the form
\begin{equation}\label{eq: d ortho short}
\sum_{k=1}^N \Delta^{\nicefrac{1}{2}}F_i^k(n_i+\tfrac{1}{2})\cdot
\Delta^{\nicefrac{1}{2}}F_j^k(n_j)\cdot \prod_{\ell\neq i,j} F_\ell^k(n_\ell+\tfrac{1}{4}\sigma_\ell)=0.
\end{equation}
Since it is supposed that this relation holds true for all $\bn\in(\tfrac{1}{2}\mathbb Z)^N$, we write it, omitting all arguments due to their arbitrariness, as 
\begin{equation} \label{d F'F'F}
\sum_{k=1}^N \Delta^{\nicefrac{1}{2}} F_i^k\cdot  \Delta^{\nicefrac{1}{2}} F_j^k \cdot \prod_{\ell\neq i,j}F_\ell^k=0.
\end{equation}
Now one sees immediately that the following analogues of Lemmas \ref{lemma 1}, \ref{lemma 2} hold true in the discrete context {\em mutatis mutandis}.
\begin{lemma} \label{d lemma 1}
Equation \eqref{d F'F'F} is equivalent to
\begin{equation}\label{d FFF}
\sum_{k=1}^NF_1^k(n_1+\tfrac{1}{4}) \cdot\ldots\cdot F_N^k(n_N+\tfrac{1}{4})=A_1(n_1)+ \ldots +A_N(n_N),
\end{equation}
with some functions $A_i(n_i)$. 
\end{lemma}

\begin{lemma}\label{d lemma 2}
Assume that all $f_i^k\neq 0$ and $\bar\Delta f_i^k\neq 0$. Then, for each $i=1,\ldots,N$ there exists a function $F_i(n_i+\tfrac{1}{4})$ such that 
\begin{equation}\label{d GH thru F N}
F_i^k(n_i+\tfrac{1}{4})=\alpha_i^kF_i(n_i+\tfrac{1}{4})+\beta_i^k,  \quad k=1,\ldots, N
\end{equation}
for some constants $\alpha_i^k\neq 0$ and $\beta_i^k$.
\end{lemma}
 \begin{proof} Note that the assumption of lemma is equivalent to  $\Delta^{\nicefrac{1}{2}}F_i^k\neq 0$. The statement of lemma is equivalent to
\begin{equation}\label{d F'G' N}
\begin{pmatrix}\Delta^{\nicefrac{1}{2}} F_i^1(n_i) \\ \ldots \\ \Delta^{\nicefrac{1}{2}}F_i^N(n_i)\end{pmatrix}\in \mathbb R\begin{pmatrix}\alpha_i^1 \\ \ldots \\ \alpha_i^N\end{pmatrix}.
\end{equation}
To prove this, take equation \eqref{d F'F'F} with all $\sigma_\ell=1$ and observe that, for any fixed $i$,  it can be formulated as the following $N-1$ orthogonality conditions:
\begin{equation}\label{d N-1 vectors}
\renewcommand{\arraystretch}{1.5}
\begin{pmatrix} \Delta^{\nicefrac{1}{2}} F_i^1 \\ \ldots \\ \Delta^{\nicefrac{1}{2}} F_i^N\end{pmatrix} \perp
\begin{pmatrix} \Delta^{\nicefrac{1}{2}}  F_j^1\cdot\prod_{\ell\neq i,j} F_\ell^1  \\ \ldots \\  \Delta^{\nicefrac{1}{2}} F_j^N\cdot \prod_{\ell\neq i,j} F_\ell^N 
\end{pmatrix}, \quad j\neq i.
\end{equation}
Multiplying the $N-1$ vectors on the right-hand side from the left by the non-degenerate matrix
$$
{\rm diag}\Big(\prod_{\ell\neq i} F_\ell^1, \ldots,  \prod_{\ell\neq i} F_\ell^N \Big)^{-1},
$$
we obtain vectors
$$
\renewcommand{\arraystretch}{1.5}
\begin{pmatrix} \big(\Delta^{\nicefrac{1}{2}}F_j^1 \big)/F_j^1  \\ \ldots \\  \big(\Delta^{\nicefrac{1}{2}}F_j^N \big)/ F_j^N 
\end{pmatrix}, \quad j\neq i.
$$
We have:
$$
 \frac{\Delta^{\nicefrac{1}{2}}F_j^1}{F_j^1}=\frac{(f_j^k(n_j+\tfrac{1}{2})-f_j^k(n_j-\tfrac{1}{2}))f_j^k(n_j)}{f_j^k(n_j+\tfrac{1}{2})f_j^k(n_j)}=
\frac{\bar\Delta_jf_j^k(n_j+\tfrac{1}{2})}{f_j^k(n_j+\tfrac{1}{2})}.
$$
Multiplying the latter vectors from the left by the non-degenerate matrix
$$
{\rm diag}\Big(\prod_{\ell=1}^N f_\ell^1(n_\ell+\tfrac{1}{2}), \ldots,  \prod_{\ell=1}^N f_\ell^N(n_\ell+\tfrac{1}{2}) \Big),
$$
we obtain vectors
$$
\bar\Delta_j\q(\bn+\tfrac{1}{2}\bsigma), \quad j\neq i.
$$
The latter vectors are linearly independent and span an $(N-1)$-dimensional subspace of $\mathbb R^N$.  Thus, the vector on the left-hand side of \eqref{d N-1 vectors} lies in the orthogonal complement of an $(N-1)$-dimensional subspace which is manifestly independent of $n_i$. This orthogonal complement is a one-dimensional space which does not depend on $n_i$. This proves \eqref{d F'G' N}. 
\end{proof}

As a result, a discrete analogue of Lemma \ref{lemma coord surfaces} holds true:
\begin{equation}\label{d for coord surfaces N}
\sum_{k=1}^N \alpha_i^k \prod_{\ell\neq i} F_\ell^k(n_\ell+\tfrac{1}{4}\sigma_\ell)=\rho_i.
\end{equation} 
Now observe that, according to $F_i^k(n_i+\tfrac{1}{4})=f_i^k(n_i)f_i^k(n_i+\tfrac{1}{2})$, to \eqref{d factorized}, and to \eqref{d GH thru F N}, equations  \eqref{d for coord surfaces N} can be expressed as follows:
\begin{equation}\label{d coord surfaces N}
\sum_{k=1}^N \frac{x_k(\bn)x_k(\bn+\tfrac{1}{2}\bsigma)}{\rho_iF_i(n_i+\tfrac{1}{4}\sigma_i)+\rho_i\beta_i^k/\alpha_i^k} = 1.
\end{equation}
The same arguments as after equation \eqref{coord surfaces N} show that the expressions 
$$
\rho_i\left(\frac{\beta_i^k}{\alpha_i^k}-\frac{\beta_i^m}{\alpha_i^m}\right),
$$
which should be equal to $a_k-a_m$,  do not depend on $i$. This finishes the proof, by setting
$$
a_k=\rho_i\frac{\beta_i^k}{\alpha_i^k}+c_i,
$$
with suitable constants $c_i$, and 
$$
u_i(n_i+\tfrac{1}{4})=\rho_i F_i(n_i+\tfrac{1}{4})-c_i.
$$
\end{proof}

Upon a re-numbering, we can assume that $a_1>a_2>\cdots>a_N>0$. Formula \eqref{discr confocal gen x thru u} shows that, as long as the points $\q(\bn)$ and $\q(\bn+\tfrac{1}{2}\bsigma)$ stay in one hyperoctant like $\mathbb R_+^N$, the quantities $u_i=u_i(n_i+\tfrac{1}{2}\sigma_i)$ lie in the intervals \eqref{eq: intervals}. If the points $\q(\bn)$ and $\q(\bn+\tfrac{1}{2}\bsigma)$ lie on different sides of a coordinate hyperplane $x_i=0$ of $\q$, the corresponding  quantity $u_i$ is outside the corresponding interval.

It is convenient to re-scale the functions $f_j^k(n_j)$ in \eqref{d factorized} by certain constant factors so that it takes the form 
\begin{equation}\label{eq: x=fff}
x_k(\bn)= \frac{\prod_{j=1}^N f_j^k(n_j)}{\prod_{i = 1}^{k-1} \sqrt{a_i - a_k}\prod_{i=k+1}^N\sqrt{a_k-a_i}}, \quad k = 1, \ldots, N.
\end{equation}
Thus, relations \eqref{discr confocal gen x thru u} give rise to the following theorem.

\begin{theorem}\label{th resolution}
For given sequences $u_i:\tfrac{1}{2}\mathbb Z+\tfrac{1}{4}\to\mathbb R$, $1\le i\le N$, consider functions $f_i^k(n_i)$ as solutions of the respective difference equations
\begin{equation}  \label{eq: ff=u+a}
  \renewcommand{\arraystretch}{1.4}
  f_i^k(n_i)f_i^k(n_i+\tfrac{1}{2})=\left\{\begin{array}{ll} u_i(n_i+\tfrac{1}{4})+a_k, & k\le i,\\ 
                                                     -\big(u_i(n_i+\tfrac{1}{4})+a_k\big), & k>i.
                                                   \end{array}\right.                            
\end{equation}
Given a sequence $u_i$, equations \eqref{eq: ff=u+a} define the functions $f_i^k$, $k=1,\ldots,N$ uniquely by prescribing their values at one point. Then, $\q$ defined by \eqref{eq: x=fff} constitutes a discrete confocal coordinate system. The right-hand sides of equations \eqref{eq: ff=u+a} are positive as long as $\q(\bn)$ and $\q(\bn+\tfrac{1}{2}(\be_1+\ldots+\be_N))$ stay in one hyperoctant of $\mathbb R^N$. 
\end{theorem}

Formulas  \eqref{eq: ff=u+a}  may be regarded as a discrete parametrization of the variables $u_i$.
Another interpretation of equations \eqref{eq: ff=u+a} (and the corresponding {\em modus operandi}) is as follows. 
\begin{theorem}\label{th funct eqs}
For $i=1,\ldots,N$, consider a system of $N-1$ functional equations for functions $f_i^k(n_i)$, $1\le k\le N$:
\begin{equation}\label{eq: funct eqs for discrete squares}
  \renewcommand{\arraystretch}{1.6}
\left\{
\begin{array}{l}
f_i^1(n_i)f_i^1(n_i+\tfrac{1}{2})-f_i^k(n_i)f_i^k(n_i+\tfrac{1}{2})= a_1-a_k, \quad k\le i, \\
f_i^1(n_i)f_i^1(n_i+\tfrac{1}{2})+f_i^k(n_i)f_i^k(n_i+\tfrac{1}{2})= a_1-a_k, \quad k> i.
\end{array} \right.
\end{equation}
For any solution of these $N$ systems, the function $\q$ defined by \eqref{eq: x=fff} constitutes a discrete confocal coordinate system. The corresponding values $u_i(n_i+\tfrac{1}{4})$ are determined by equations \eqref{eq: ff=u+a}.
\end{theorem}
\begin{proof}
We arrive at equations \eqref{eq: funct eqs for discrete squares} by eliminating $u_i$ between equations  \eqref{eq: ff=u+a}. Note that equations \eqref{eq: funct eqs for discrete squares} do not depend on $u_i=u_i(n_i+\tfrac{1}{4})$. The latter can be determined {\em a posteriori} from any of equations \eqref{eq: ff=u+a}.
\end{proof}

 It is well known that functional equations \eqref{eq: funct eqs for discrete squares} admit solutions in terms of trigonometric/hyperbolic functions if $N=2$, and in terms of elliptic functions if $N=3$. We discuss these solutions in Sections \ref{section_2d}, \ref{section_3d}, respectively.

\section{Geometric interpretation}
\label{s.geometric}

The main formula from Theorem \ref{Th discr confocal}, 
\begin{equation}\label{discr confocal gen again}
\sum_{k=1}^N \frac{x_k(\bn)x_k(\bn+\frac{1}{2}\bsigma)}{a_k+u_i}=1,\quad u_i=u_i(n_i+\tfrac{1}{4}\sigma_i), \quad i=1,\ldots,N,
\end{equation}
admits a remarkable geometric interpretation. Recall that the {\em polarity} with respect to a non-degenerate quadric is a projective transformation between the points $\q\in \mathbb P^N$ and the hyperplanes $\Pi\in (\mathbb P^N)^*$. In non-homogeneous coordinates, if the quadric $Q$ is given by a quadratic form $Q(\q)=0$, then the hyperplane $\Pi$ polar to a point $\q=(x_1,\ldots,x_N)\in \mathbb P^N$ with respect to $Q$ consists of all points $\boldsymbol y=(y_1,\ldots,y_N)\in\mathbb P^N$ satisfying $\bar Q(\q, \boldsymbol y)=0$, where $\bar Q$ is the symmetric bilinear form corresponding to the quadratic form $Q$. We write $\q=P_Q(\Pi)$ and $\Pi=P_Q(\q)$.
Thus, formula \eqref{discr confocal gen again} is equivalent to saying that
\begin{itemize}
\item[] {\em the point  $\q(\bn+\frac{1}{2}\bsigma)$ lies in the intersection of the polar hyperplanes of $\q(\bn)$ with respect to the confocal quadrics $Q(u_i)$,  $i=1,\ldots,N$:}
\begin{equation}\label{eq: geom constr}
\q(\bn+\tfrac{1}{2}\bsigma)=\bigcap_{i=1}^N P_{Q(u_i)}(\q(\bn)), \quad u_i=u_i(n_i+\tfrac{1}{4}\sigma_i).
\end{equation}
\end{itemize}
Of course, the roles of $\q(\bn)$ and $\q(\bn+\tfrac{1}{2}\bsigma)$ in this formula are completely symmetric.

This interpretation can be used to give a {\em geometric construction} of a discrete confocal coordinate system $\q: (\frac{1}{2}\mathbb Z)^N\supset\mathcal U\to\mathbb R^N$, or, better, of its restriction to two dual sublattices like $\mathbb Z^N$ and $(\mathbb Z+\frac{1}{2})^N$. Suppose that for each $i=1,\ldots,N$ a sequence of quadrics of the confocal family \eqref{eq:confocal-family} is chosen, with the parameters 
\begin{equation}\label{uk}
u_i:\big(\tfrac{1}{2}\mathbb Z+\tfrac{1}{4}\big)\cap \mathcal I_i\to \mathbb R,
\end{equation}
indexed by a discrete variable $n_i+\tfrac{1}{4}\in \mathcal I_i$, where $n_i\in\tfrac{1}{2}\mathbb Z$. It is convenient to think of $u_i(n_i+\tfrac{1}{4})$ as being assigned 
to the interval $[n_i,n_i+\frac{1}{2}]$, for which  $n_i+\frac{1}{4}$ is the midpoint. We denote by $\mathcal V$, $\mathcal V^*$ the parts of the respective lattices $\mathbb Z^N$, $(\mathbb Z+\tfrac{1}{2})^N$ lying in the region $\prod_{i=1}^N \mathcal I_i$. We construct a discrete net $\q:\mathcal V\cup\mathcal V^*\to\mathbb R^N$ recurrently, starting with an arbitrary point $\q(\bn_0)$, as long as the components of $\q(\bn)$ are non-vanishing.
\medskip 

\begin{figure}[H]
  \centering
  \includegraphics[width=0.33\textwidth]{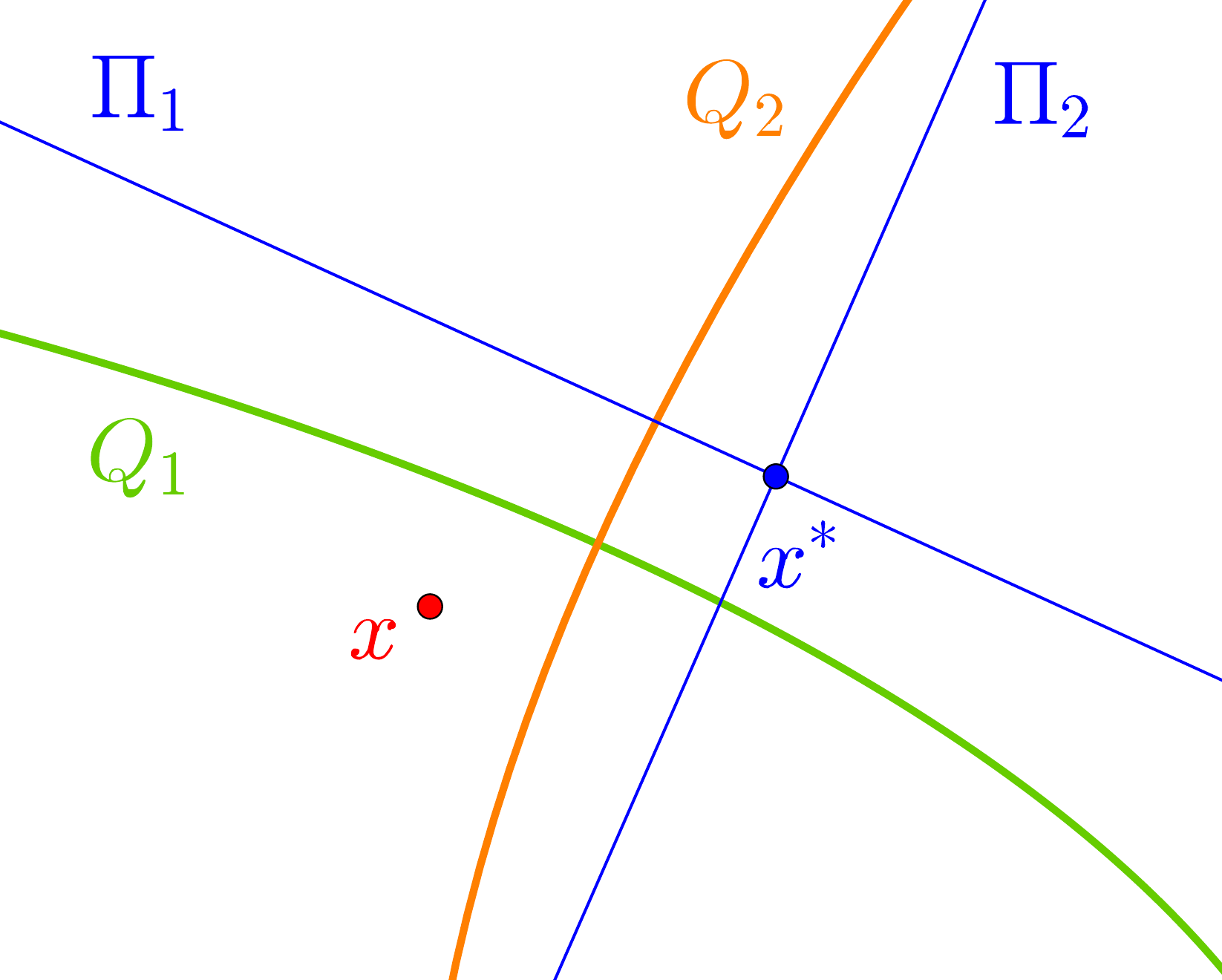}
    \caption{Geometric construction of $\bx^*$ in the case $N=2$ as the intersection of the polar lines $\Pi_1$ and $\Pi_2$ 
    of $\bx$ with respect to the confocal conics $Q_1$ and $Q_2$.    
  }
  \label{Fig geom construction}
\end{figure}

{\bf Construction} (cf. Figure \ref{Fig geom construction}). {\em Let $\bn$ and $\bn^*$ be two neighboring points in the two dual sublattices, in the sense that
$$
\bn^*=\bn+\tfrac{1}{2}\bsigma, \quad \bsigma=(\sigma_1,\ldots,\sigma_N), \quad \sigma_i=\pm 1.
$$
Suppose that $\q(\bn)=\q$ is already known. Then $\q(\bn^*)=\q^*$ is constructed as the intersection point of the $N$ polar hyperplanes}
\begin{equation}
\q^*=C_{\bn,\frac{1}{2}\bsigma}(\q):
=\bigcap_{i=1}^N P_{Q(u_i)}(\q), \quad u_i=u_i(n_i+\tfrac{1}{4}\sigma_i).
\end{equation}

\medskip

In order to show that this construction is well defined, the following statement is required.

\begin{proposition} \label{Prop well defined}
The following diagram is commutative for any $\bsigma,\tilde{\bsigma}\in\{\pm 1\}^N$:
\begin{center}
\begin{tikzpicture}
  % Tell it where the nodes are
  \node (A) {$\quad\mathbb R^N\quad$};
  \node (B) [below=of A] {$\quad\mathbb R^N\quad$};
  \node (C) [right=of A] {$\quad\mathbb R^N\quad$};
  \node (D) [right=of B] {$\quad\mathbb R^N\quad$};
  % Tell it what arrows to draw
  \draw[-stealth] (A)-- node[left] {\small $C_{\bn,\frac{1}{2}\bsigma}$} (B);
  \draw[-stealth] (B)-- node [below] {\small $C_{\bn+\frac{1}{2}\bsigma,\frac{1}{2}\tilde\bsigma}$} (D);
  \draw[-stealth] (A)-- node [above] {\small $C_{\bn,\frac{1}{2}\tilde\bsigma}$} (C);
  \draw[-stealth] (C)-- node [right] {\small $C_{\bn+\frac{1}{2}\tilde\bsigma, \frac{1}{2}\bsigma}$} (D);
\end{tikzpicture}
\end{center}
Thus, applying the above construction along a path depends only on the initial and the end points of the path and not on the path itself. 
\end{proposition} 
\begin{proof} 
Denote the right-hand side of \eqref{discr confocal gen x thru u} by
\begin{equation}
  x_k(\bn)x_k(\bn+\tfrac{1}{2}\bsigma) = B_k(\bn,\bsigma):=\dfrac{\prod_{j=1}^N (u_j(n_j+\tfrac{1}{4}\sigma_j)+a_k)}{\prod_{j \neq k} (a_k - a_j)}.
\ela
Then the commutativity of the diagram is equivalent to
\begin{equation}\label{BB=BB}
  \frac{B_k(\bn+\frac{1}{2}\bsigma,\tilde{\bsigma})}{B_k(\bn,\bsigma)}=\frac{B_k(\bn+\tfrac{1}{2}\tilde{\bsigma},\bsigma)}{B_k(\bn,\tilde{\bsigma})},
\end{equation}
or
\begin{equation}\label{eq for comm}
  \frac{\prod_{j=1}^N (u_j(n_j+\tfrac{1}{2}\sigma_j+\tfrac{1}{4}\tilde\sigma_j)+a_k)}{\prod_{j=1}^N (u_j(n_j+\tfrac{1}{4}\sigma_j)+a_k)}=
   \frac{\prod_{j=1}^N (u_j(n_j+\tfrac{1}{2}\tilde\sigma_j+\tfrac{1}{4}\sigma_j)+a_k)}{\prod_{j=1}^N (u_j(n_j+\tfrac{1}{4}\tilde\sigma_j)+a_k)}.
\end{equation}
For each $j$, we have either $\tilde\sigma_j=\sigma_j$, or $\tilde\sigma_j=-\sigma_j$. In the first case, the corresponding factors in the numerators of both sides of \eqref{eq for comm} are equal, as well as the corresponding factors in the denominators. In the second case, the corresponding factors in the numerator and in the denominator on the left-hand side are equal, and the same holds true for the corresponding factors in the numerator and in the denominator on the right-hand side. This proves \eqref{BB=BB}.
\end{proof}

\begin{figure}[H]
  \centering
  \includegraphics[width=0.34\textwidth]{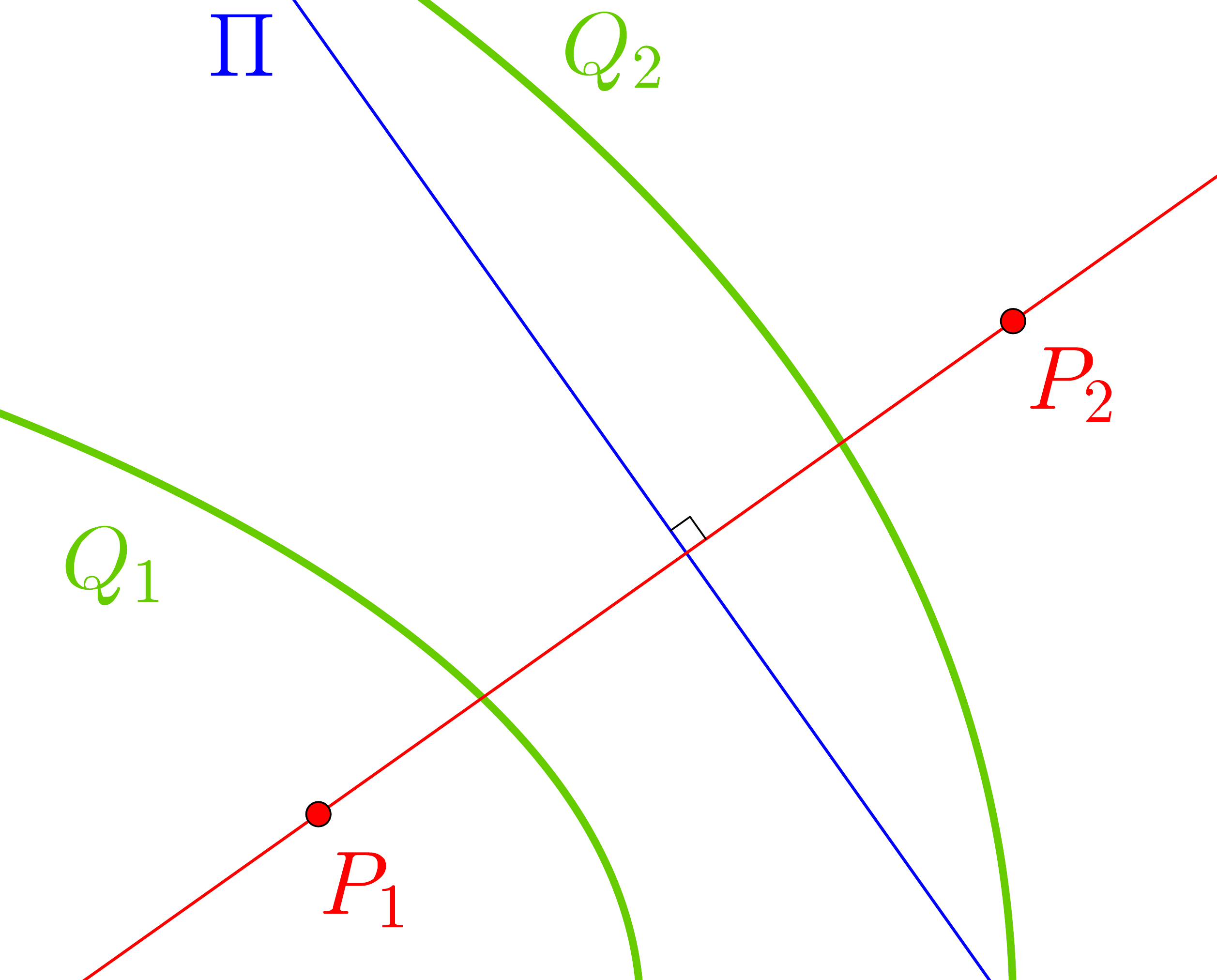}
  \caption{Orthogonality in the case $N=2$: If $P_2$ is related to $P_1$ via polarity in two confocal conics, 
    that is, $\Pi = P_{Q_1}(P_1)$ and $P_2 = P_{Q_2}(\Pi)$, then the line through $P_1$ and $P_2$
    is orthogonal to $\Pi$.
  }
  \label{Fig ortho confocal}
\end{figure}

The discrete orthogonality property is now a consequence of the following lemma (cf. Figure \ref{Fig ortho confocal}).

\begin{lemma} \label{Lemma ortho}
Let $\Pi$ be a hyperplane. Then the poles of $\Pi$ with respect to all quadrics of the confocal family \eqref{eq:confocal-family} lie on a line $\ell$. This line $\ell$ is orthogonal to $\Pi$.
\end{lemma}

\begin{proof}
Let the equation of the hyperplane $\Pi$ be $\sum_{k=1}^N c_kx_k=1$, where $\bc=(c_1,\ldots,c_N)$ is a normal vector for $\Pi$. Take two quadrics of the confocal family, $Q_1=Q(u)$ and $Q_2=Q(v)$. Set
$$
P_1=P_{Q_1}(\Pi)=(y_1, \ldots,y_N)\quad {\rm and} \quad P_2=P_{Q_2}(\Pi)=(z_1, \ldots,z_N). 
$$ 
Then we get the following two forms of the equation of the hyperplane $\Pi$:
\begin{equation}
\sum_{i=1}^N \frac{x_iy_i}{a_i+u}=1 \quad {\rm and}\quad \sum_{i=1}^N \frac{x_iz_i}{a_i+v}=1.
\end{equation}
Thus,
$$
c_i=\frac{y_i}{a_i+u}=\frac{z_i}{a_i+v},
$$
and, hence,
$$
y_i-z_i=c_i(a_i+u)-c_i(a_i+v)=c_i(u-v),
$$
so that the vector $P_2-P_1$ is proportional to $\bc$ and therefore is orthogonal to $\Pi$. Thus, denoting by $\ell$ the line which passes through $P_1$ orthogonally to $\Pi$, we see that the pole $P_2$ of the hyperplane $\Pi$ with respect to the quadric $Q_2$ lies on $\ell$. It remains to note that $Q_2$ is an arbitrary quadric of the confocal family of $Q_1$.
\end{proof}

\begin{theorem}
The nets $\q(\mathcal V)$ and $\q(\mathcal V^*)$ are orthogonal in the sense of Definition \ref{def ortho}.
\end{theorem}
\begin{proof}
From \eqref{eq: geom constr} it follows that, for a fixed index $k$, all points $\q(\bn+\tfrac{1}{2}\bsigma)$ with $\sigma_k=1$ lie in the hyperplane $\Pi=P_{Q_1}(\q(\bn))$,
where 
$
Q_1=Q(u_k(n_k+\tfrac{1}{4})).
$
These points are exactly the vertices of the facet of $\q(\mathcal V^*)$ dual to the edge $[\q(\bn),\q(\bn+\be_k)]$ of $\q(\mathcal V)$. Now, since
$$
\q(\bn)=P_{Q_1}(\Pi), \quad \q(\bn+\be_k)=P_{Q_2}(\Pi),
$$
where
$
Q_2=Q(u_k(n_k+\tfrac{3}{4})), 
$
it follows from Lemma \ref{Lemma ortho} that the line $(\q(\bn),\q(\bn+\be_k))$ is orthogonal to the hyperplane $\Pi$. 
\end{proof}

%%%%%%%%%%%%%%%%%%%%%%%%%%%%%%%%%%%%%%%%%%%%%%%%%%%%%%%%%%%

\section{Discrete confocal coordinates in terms of gamma functions}
\label{sect gamma}
\label{s.article1}

There exists an important particular case when the difference equations \eqref{eq: ff=u+a} admit an explicit solution, namely by the choice 
\begin{equation*}
  u_i(n_i + \tfrac{1}{4}) = n_i + \epsilon_i, \quad i = 1,\ldots,N,  
\end{equation*}
where $\epsilon_i \in \R$ are some fixed shifts.
This can be considered to correspond to the smooth case \eqref{eq:elliptic-coordinates}
where we take the quantities $u_i$ as coordinates without further re-parametrization. With this choice, equations \eqref{eq: ff=u+a} turn into
\begin{equation}\label{Q26}
  \renewcommand{\arraystretch}{1.4}
  f_i^k(n_i)f_i^k(n_i+\tfrac{1}{2})=\left\{
    \begin{array}{ll} 
      n_i + a_k + \epsilon_i, & k\le i,\\
      -\big(n_i +a_k + \epsilon_i\big), & k>i.
    \end{array}
  \right.                            
\end{equation}
These equations can be solved in terms of the ``discrete square root'' function defined as
\bela{Q27}
  \sqr{u} = \frac{\Gamma(u+\frac{1}{2})}{\Gamma(u)},
\ela
which satisfies the identities
\bela{Q28}
  \textstyle\sqr{u}\sqr{u+\frac{1}{2}} = u,\quad \sqr{-u}\sqr{-u-\frac{1}{2}} = -u-\frac{1}{2}.
\ela
We can write solutions of \eqref{Q26} as
\begin{equation}
\renewcommand{\arraystretch}{1.4}
  f_i^k(n_i) =
  \left\{ \begin{array}{ll}    
      \sqr{n_i + a_k + \epsilon_i} \quad & \text{for} \quad i \geq k, \\ 
      \sqr{-n_i - a_k - \epsilon_i + \frac{1}{2}}  \quad & \text{for} \quad i < k.
    \end{array}\right.
\end{equation}

One can impose boundary conditions
\begin{alignat}{2}
  \label{eq:discrete-boundary-conditions1}
  &x_k |_{n_k=-\alpha_k} = 0  &&\quad \text{for} \quad k= 1, \ldots, N,\\
  \label{eq:discrete-boundary-conditions2}
  & x_k |_{n_{k-1} = -\alpha_k} = 0 &&\quad \text{for} \quad k= 2, \ldots, N,
\end{alignat}
on the integer lattice $\Z^N$ for certain integers $\alpha_1> \cdots >\alpha_N$, which imitate the corresponding property of the continuous confocal coordinates.
These boundary conditions are satisfied provided that
$$
a_k-\alpha_k+\epsilon_k=0, \quad a_k-\alpha_k+\epsilon_{k-1}=\frac{1}{2},
$$
for which the shifts $\epsilon_k$ should satisfy $\epsilon_{k-1} -\epsilon_k = \frac{1}{2}$. Choosing $\epsilon_k = - \frac{k}{2}$ and $a_k =\alpha_k+\frac{k}{2}$,
we finally arrive at the solutions 
\begin{equation}\label{eq: solution gamma}
\renewcommand{\arraystretch}{1.4}
  f_i^k(n_i) =
  \left\{ \begin{array}{ll}    
      \sqr{n_i + \alpha_k + \frac{k-i}{2}} \quad & \text{for} \quad i \geq k, \\ 
      \sqr{-n_i - \alpha_k - \frac{k-i}{2} + \frac{1}{2}}  \quad & \text{for} \quad i < k.
    \end{array}\right.
\end{equation}
These are the functions introduced and studied in \cite{BSST16}, as solutions of the discrete Euler-Darboux-Poisson equations (cf. Appendix). 

%%%%%%%%%%%%%%%%%%%%%%%%%%%%%%%%%%%%%%%%%%%%

\section{The case \boldmath $N=2$}\label{section_2d}

\subsection{Classical confocal coordinate system}

\begin{figure}[H]
  \begin{center}
    \scalebox{1.2}{\input{trigonometric-symmetric-cont1_image.pgf}}
  \end{center}
  \caption{
    Two-dimensional classical confocal coordinate system \eqref{P7}
    in terms of trigonometric functions with $a=2$, $b=1$.}
    \label{Fig 2d elliptic coord}
\end{figure}
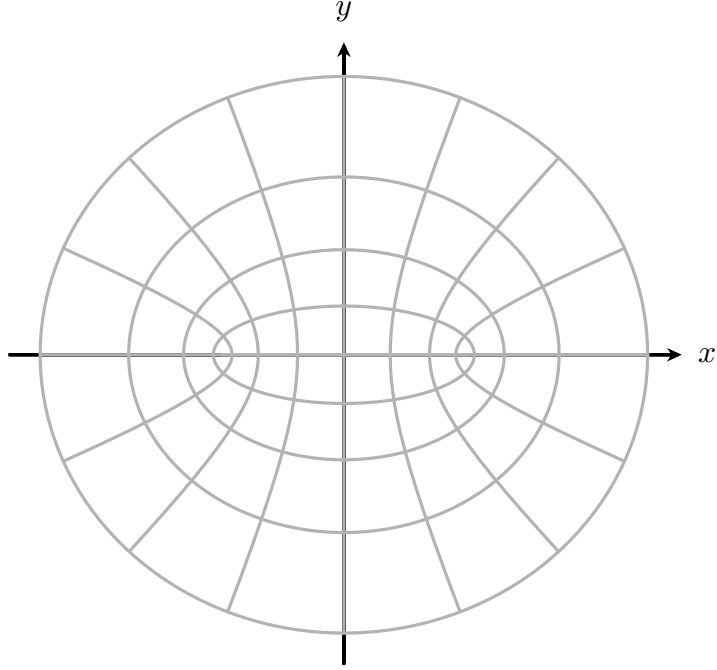

We have seen that, given the family of confocal conics
\begin{equation}\label{eq: 2d confocal}
\frac{x^2}{a+\lambda} +\frac{y^2}{b+\lambda} =1,
\end{equation}
the defining equations
\bela{P1}
  \frac{x^2}{u+a} + \frac{y^2}{u+b} = 1,\quad  \frac{x^2}{v+a} + \frac{y^2}{v+b} = 1
\ela
of confocal coordinates $\big\{(u,v): -a< u < -b < v\big\}$ on the plane give rise to the expressions
\bela{P2}
  x^2 = \frac{(u+a)(v+a)}{a-b},\quad y^2 = \frac{(u+b)(v+b)}{b-a}.
\ela
For an arbitrary re-parametrization of the coordinate lines, $u=u(s_1)$, $v=v(s_2)$, we obtain
\bela{P3}
  x = \frac{f_1(s_1)f_2(s_2)}{\sqrt{a-b}},\quad y = \frac{g_1(s_1)g_2(s_2)}{\sqrt{a-b}},
\ela
where
\begin{equation}\label{P4}
\renewcommand{\arraystretch}{1.4}
\left\{\begin{array}{l}
  (f_1(s_1))^2=u + a , \\ (g_1(s_1))^2 =-(u+ b),
  \end{array}\right. \quad 
\left\{\begin{array}{l}  
  (f_2(s_2))^2 = v + a, \\
  (g_2(s_2))^2 = v+ b.
 \end{array}\right. 
 \end{equation}
Elimination of $u$ and $v$ leads to 
\begin{align}
 ( f_1(s_1))^2 + (g_1(s_1))^2 & =  a-b, \label{P5a} \\ 
  (f_2(s_2))^2 - (g_2(s_2))^2 & =  a-b. \label{P5b}
\end{align}
The probably most obvious parametrization of solutions of these functional equations is by means of trigonometric/hyperbolic functions:
\begin{eqnarray}
  f_1(s_1) =\sqrt{a-b}\,\cos s_1,& & g_1(s_1) = \sqrt{a-b}\ \sin s_1,\label{P6a}\\
  f_2(s_2)= \sqrt{a-b}\,\cosh s_2\,,& & g_2(s_2)  = \sqrt{a-b}\,\sinh s_2. \label{P6b}
\end{eqnarray}
Accordingly, we obtain the representation
\bela{P7}
  \begin{pmatrix} x \\ y \end{pmatrix} = \sqrt{a-b} \begin{pmatrix} \cos s_1\,\cosh s_2\\  \sin s_1\,\sinh s_2 \end{pmatrix}
\ela
of the confocal system of coordinates on the plane with the relation between $(u,v)$ and $(s_1,s_2)$ given by \eqref{P4}. This coordinate system is depicted in Figure \ref{Fig 2d elliptic coord}.

\subsection{Discrete confocal coordinate systems}
\label{s: 2d trig/hyp}

For any discrete set of confocal quadrics \eqref{eq: 2d confocal}, indexed by $u(n_1+\tfrac{1}{4})\in \mathbb R$ and $v(n_2+\tfrac{1}{4})\in \mathbb R$ with $n_1,n_2\in\frac{1}{2}\mathbb Z$, we have introduced the discrete confocal quadrics defined by the equations of polarity relating nearest neighbors $\bx(\bn)$ and $\bx(\bn+\tfrac{1}{2}\bsigma)$:
\begin{eqnarray}
  \frac{x(\bn)x(\bn+\tfrac{1}{2}\bsigma)}{u(n_1+\tfrac{1}{4}\sigma_1)+a} 
  + \frac{y(\bn)y(\bn+\tfrac{1}{2}\bsigma)}{u(n_1+\tfrac{1}{4}\sigma_1)+b} & = &  1, \nonumber\\ 
  \frac{x(\bn)x(\bn+\tfrac{1}{2}\bsigma)}{v(n_2+\tfrac{1}{4}\sigma_2)+a} 
  + \frac{y(\bn)y(\bn+\tfrac{1}{2}\bsigma)}{v(n_2+\tfrac{1}{4}\sigma_2)+b} & = & 1. \label{P8}
\end{eqnarray}
This is equivalent to
\begin{align}\label{P9}
 x(\bn)x(\bn+\tfrac{1}{2}\bsigma) & =  \frac{\big(u(n_1+\tfrac{1}{4}\sigma_1)+a\big)\big(v(n_2+\tfrac{1}{4}\sigma_2)+a\big)}{a-b}, \nonumber\\
 y(\bn)y(\bn+\tfrac{1}{2}\bsigma) & =  \frac{\big(u(n_1+\tfrac{1}{4}\sigma_1)+b\big)\big(v(n_2+\tfrac{1}{4}\sigma_2)+b\big)}{a-b}.
\end{align}
According to Theorem \ref{th resolution}, we can resolve this as follows:
\bela{P10}
  x(\bn) = \frac{f_1(n_1)f_2(n_2)}{\sqrt{a-b}},\quad y(\bn) = \frac{g_1(n_1)g_2(n_2)}{\sqrt{a-b}},
\ela
where
\begin{equation}\label{P11a}
\renewcommand{\arraystretch}{1.4}
\left\{\begin{array}{l}
  f_1(n_1)f_1(n_1+\tfrac{1}{2})=u(n_1+\tfrac{1}{4})+a, \\
  g_1(n_1)g_1(n_1+\tfrac{1}{2})=-\big(u(n_1+\tfrac{1}{4})+ b\big),
\end{array} \right.   
\end{equation}
\begin{equation}\label{P11b}
\renewcommand{\arraystretch}{1.4}
 \left\{\begin{array}{l} 
 f_2(n_2)f_2(n_2+\tfrac{1}{2})=v(n_2+\tfrac{1}{4})+ a, \\
 g_2(n_2)g_2(n_2+\tfrac{1}{2}) =v(n_2+\tfrac{1}{4})+ b.
\end{array} \right.
\end{equation}

\paragraph{Parametrization in terms of gamma functions.} A solution of equations \eqref{P11a}, \eqref{P11b} found in \cite{BSST16} is given by $a=\alpha+\tfrac{1}{2}$, $b=\beta+1$,
\begin{equation*}
\renewcommand{\arraystretch}{1.4}
\left\{\begin{array}{l}
  f_1(n_1)=(n_1+\alpha)_{1/2}, \\
  g_1(n_1)=(-n_1-\beta)_{1/2},
\end{array} \right.   \quad 
 \left\{\begin{array}{l} 
 f_2(n_2)=(n_2+\alpha-\tfrac{1}{2})_{1/2}, \\
 g_2(n_2)=(n_2+\beta)_{1/2},
\end{array} \right.
\end{equation*}
so that 
\begin{equation*}
\renewcommand{\arraystretch}{1.4}
 \begin{pmatrix} x(\bn)\\ y(\bn) \end{pmatrix}=\frac{1}{\sqrt{\alpha-\beta-\tfrac{1}{2}}}
   \begin{pmatrix}   (n_1+\alpha)_{1/2}\ (n_2+\alpha-\tfrac{1}{2})_{1/2} \\  (-n_1-\beta)_{1/2}\ (n_2+\beta)_{1/2}  \end{pmatrix}.
\end{equation*}

\paragraph{Parametrization in terms of trigonometric/hyperbolic functions.}
We obtain functional equations satisfied by the functions $f_i$, $g_i$ by eliminating $u(n_1+\tfrac{1}{4})$ and $v(n_2+\tfrac{1}{4})$ from equations \eqref{P11a}, \eqref{P11b}:
\begin{align}
 & f_1(n_1)f_1(n_1+\tfrac{1}{2}) + g_1(n_1)g_1(n_1+\tfrac{1}{2}) =  a-b,    \label{P12a}\\
 & f_2(n_2)f_2(n_2+\tfrac{1}{2}) - g_2(n_2)g_2(n_2+\tfrac{1}{2})  =a-b.     \label{P12b}
\end{align}
By virtue of the addition theorems for trigonometric and hyperbolic functions, one easily finds solutions to these functional equations which approximate functions \eqref{P6a}, \eqref{P6b}:
\begin{equation}\label{P13a}
  f_1(n_1) =\dis\sqrt{\frac{a-b}{\cos\frac{\delta_1}{2}}}\,\cos (\delta_1 n_1+c_1),\quad
  g_1(n_1) =\dis\sqrt{\frac{a-b}{\cos\frac{\delta_1}{2}}}\,\sin (\delta_1 n_1+c_1),
\end{equation}
and
\begin{equation}\label{P13b}
  f_2(n_2) =\dis \sqrt{\frac{a-b}{\cosh\frac{\delta_2}{2}}}\,\cosh (\delta_2n_2+c_2), \quad
  g_2(n_2) =\dis \sqrt{\frac{a-b}{\cosh\frac{\delta_2}{2}}}\,\sinh (\delta_2n_2+c_2).
\end{equation}
Thus,
\bela{P14}
  \begin{pmatrix} x(\bn)\\ y(\bn) \end{pmatrix} = \sqrt{\frac{a-b}{\cos\frac{\delta_1}{2}\,\cosh\frac{\delta_2}{2}}}\begin{pmatrix} \cos (\delta_1 n_1+c_1)\,\cosh (\delta_2n_2+c_2)\\    \sin (\delta_1 n_1+c_1)\,\sinh (\delta_2n_2+c_2) \end{pmatrix}.
\ela
The discrete coordinate curves $n_2={\rm const}$ are to be interpreted as discrete ellipses. In order that they be closed curves,  it is necessary to choose the lattice parameter $\delta_1$ according to
\bela{P15}
  \delta_1 = \frac{2\pi}{m},\quad m\in\N.
\ela
One obtains a picture which is symmetric with respect to the coordinate axes if $c_1=c_2=0$.
The parameters $u(n_1+\tfrac{1}{4})$ and $v(n_2+\tfrac{1}{4})$ of the associated lattice of continuous confocal quadrics \eqref{P1} are obtained from \eqref{P11a}, \eqref{P11b} and \eqref{P13a}, \eqref{P13b}.

Figures \ref{Fig discr 2d trig full}--\ref{Fig polarity conics} display a discrete confocal coordinate system for $a=2$, $b=1$, $m=2$ and $\delta_2=\delta_1$. In the continuous case encoded in the parametrisation \eqref{P7}, the foci on the $x$-axis correspond to $(s_1,s_2)=(0,0)$ and $(s_1,s_2)=(\pi,0)$. Their discrete analogs in the sublattice $\mathbb Z^2$ in Figure \ref{iihh} (top) correspond to $(n_1,n_2)=(0,0)$ resp. $(n_1,n_2)=(4,0)$. The valence of these points is 2, as opposed to the regular points of valence 4. In the sublattice $(\mathbb Z+\frac{1}{2})^2$, the analogs of the foci are the ``focal edges'' connecting pairs of neighboring points of valence 3. For instance, the analog of the right focus is the edge  $[\q(\frac{1}{2},\frac{1}{2}),\q(\frac{1}{2},-\frac{1}{2})]$. In the sublattices $\mathbb Z\times (\mathbb Z+\frac{1}{2})$ and $(\mathbb Z+\frac{1}{2})\times \mathbb Z$, the analogs of the foci are the double points like $\q(0,\frac{1}{2})=\q(0,-\frac{1}{2})$ and $\q(\frac{1}{2},0)=\q(-\frac{1}{2},0)$, both having  valence 3 (see  Figure \ref{iihh}, bottom). Figure \ref{Fig polarity conics} shows the confocal conics participating in the polarity relations of a discrete confocal coordinate system \eqref{P14}.

\begin{figure}[H]
  \begin{center}
    \scalebox{1.2}{\input{trigonometric-symmetric1_image.pgf}}
  \end{center}
  \caption{
    Two-dimensional discrete confocal coordinate system \eqref{P14} on $\left(\tfrac{1}{2}\Z\right)^2$ in terms of trigonometric/hyperbolic  functions
     with $a=2$, $b=1$, $m=8$, $\delta_1 = \delta_2 = \frac{2\pi}{m}$, $c_1 = c_2 = 0$.
  }
  \label{Fig discr 2d trig full}
\end{figure}
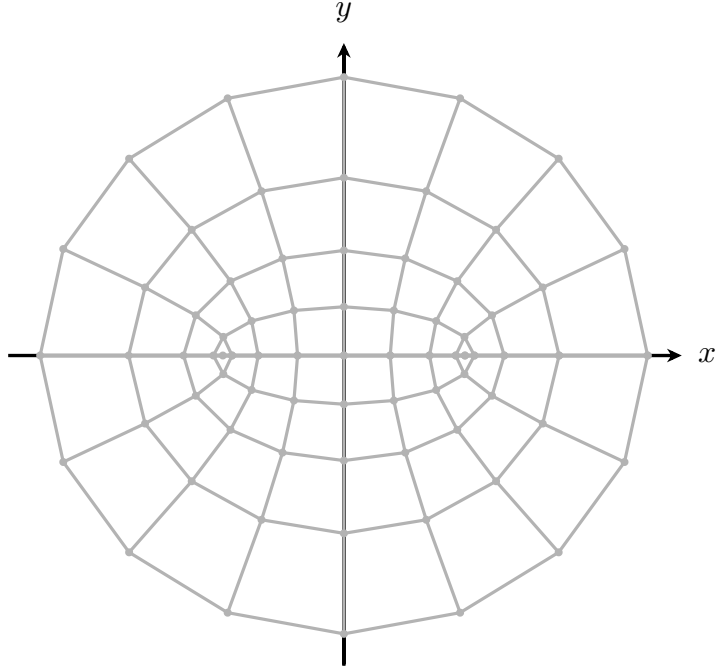

\begin{figure}[H]
  \begin{center}
    \scalebox{1.1}{\input{trigonometric-symmetric5_image.pgf}}
    \scalebox{1.1}{\input{trigonometric-symmetric7_image.pgf}}
  \end{center}
  \caption{
    Pairs of dual orthogonal sublattices.
    Gray points show discrete confocal coordinates \eqref{P14} on $\left(\tfrac{1}{2}\Z\right)^2$
    with $a=2$, $b=1$, $m=8$, $\delta_1 = \delta_2 = \frac{2\pi}{m}$, $c_1 = c_2 = 0$.
    (top) Sublattice on $\Z^2$ in blue and on $\left(\Z + \tfrac{1}{2}\right)^2$ in red.
    (bottom) Sublattice on $\Z \times (\Z+\tfrac{1}{2})$ in blue and on $(\Z+\tfrac{1}{2}) \times \Z$ in pink.
  }
  \label{iihh}
\end{figure}
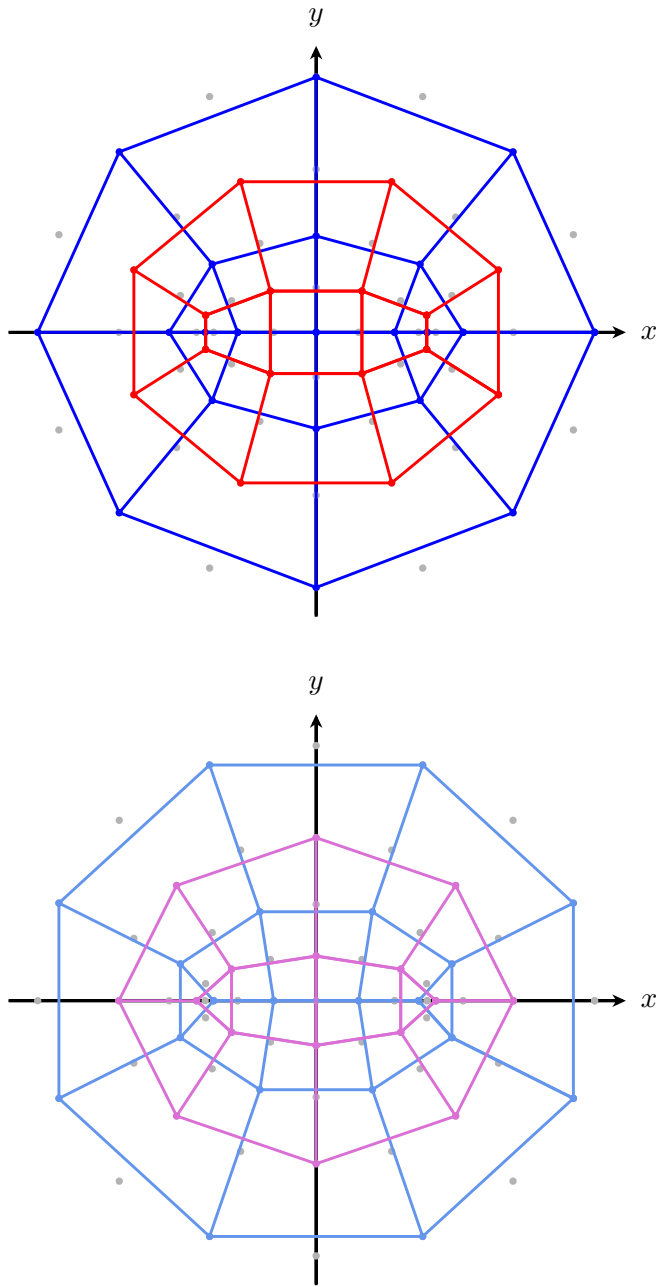

\begin{figure}[H]
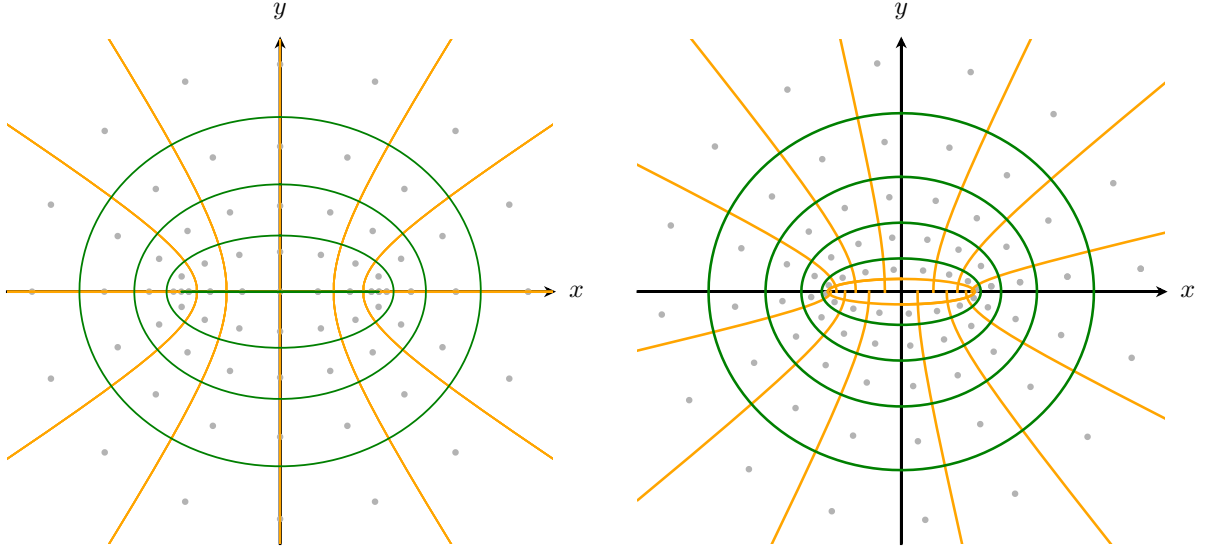

  \begin{center}
    \scalebox{0.98}{\input{trigonometric-symmetric3_image.pgf}}
    \scalebox{0.98}{\input{trigonometric-unsymmetric3_image.pgf}}
  \end{center}
  \caption{
    Polarity relation for discrete confocal conics.
    Gray points show discrete confocal coordinates \eqref{P14} on $\left(\tfrac{1}{2}\Z\right)^2$
    with $a=2$, $b=1$, $m=8$, $\delta_1 = \delta_2 = \frac{2\pi}{m}$.
    The corresponding classical confocal conics which give rise to the polarity relation
    between gray points are shown in orange
    (for the values $u\left(n_1 + \frac{1}{4}\right)$)
    and green (for the values $v\left(n_2 + \frac{1}{4}\right)$).
    (top) Symmetric case with $c_1 = c_2 = 0$.
    All orange conics are hyperbolas and all green conics are ellipses.
    Note that near the coordinate axes those conics become degenerate
    and the polarity relation is not injective anymore.
    (bottom) Asymmetric case with $c_1 = 0.1$, $c_2 = 0.3$.
    Moving along the $n_2$-direction, the polarity across the $y$-axis
    is established by a conic with value $u\left(n_1 + \frac{1}{4}\right) < -a$, which is purely imaginary,
    while the polarity across the $x$-axis is established by a conic
    with value $u\left(n_1 + \frac{1}{4}\right) > -b$, which is an ellipse. 
  }
  \label{Fig polarity conics}
\end{figure}

% \begin{figure}[H]
%   \begin{center}
%     \input{pics/trigonometric-symmetric2_image.pgf}
%   \end{center}
%   \caption{$a=2$, $b=1$, $m=8$, $\delta_1 = \delta_2 = \frac{2\pi}{m}$}
%   \label{symmetric}
% \end{figure}

%% \begin{figure}[H]
%%   \begin{center}
%%     \scalebox{1.1}{\input{pics/trigonometric-unsymmetric5_image.pgf}}
%%     \scalebox{1.1}{\input{pics/trigonometric-unsymmetric7_image.pgf}}
%%   \end{center}
%%   \caption{}
%%   \label{ihih}
%% \end{figure}

\subsection{Parametrization by elliptic functions}

The trigonometric/hyperbolic parametrization \eqref{P6a}, \eqref{P6b} is not the only explicit solution to the functional equations \eqref{P5a}, \eqref{P5b}. One can find further ones in terms of elliptic functions. 
For instance, equation \eqref{P5a} admits the solution
\begin{equation}\label{ell param 1}
f_1(s_1)=\sqrt{a-b}\ \jac{cn}(s_1,k_1), \quad g_1(s_1)=\sqrt{a-b}\ \jac{sn}(s_1,k_1)
\end{equation}
with an arbitrary modulus $k_1$ (with \eqref{P6a} being the limiting case $k_1\to 0$), or the solution
\begin{equation}\label{ell param 2}
f_1(s_1)=\sqrt{a-b}\ \jac{dn}(s_1,k_1), \quad g_1(s_1)=\sqrt{a-b}\ k_1\jac{sn}(s_1,k_1).
\end{equation}
Similarly, equation \eqref{P5b} admits the solution
\begin{equation}\label{ell param 2 1}
f_2(s_2)=\sqrt{a-b}\ \frac{1}{\jac{dn}(s_2,k_2)}, \quad g_2(s_2)=\sqrt{a-b}\ \frac{k_2\jac{sn}(s_2,k_2)}{\jac{dn}(s_2,k_2)}
\end{equation}
with an arbitrary modulus $k_2$ (with the limiting case \eqref{P6b} as $k_2\to 1$). Further examples of solutions of  \eqref{P5b} are:
\begin{equation}\label{ell param 2 2}
f_2(s_2)  =\sqrt{a-b}\  \dfrac{1}{\jac{sn}(s_2,k_2)}, \quad g_2(s_2) = \sqrt{a-b}\ \dfrac{\jac{cn}(s_2,k_2)}{\jac{sn}(s_2,k_2)},
\end{equation}
or
\begin{equation}\label{ell param 2 3}
f_2(s_2)=\sqrt{a-b}\ \frac{1}{k_2'}\jac{dn}(s_2,k_2), \quad g_2(s_2)=\sqrt{a-b}\ \frac{k_2}{k_2'}\jac{cn}(s_2,k_2),
\end{equation}
where $k_2'=\sqrt{1-k_2^2}$. All such solutions can be seen as based on relations between squares of theta functions, and are connected by simple transformations in the complex domain, but they have rather different properties in the real domain. For instance, in \eqref{ell param 2 1} one of the participating functions is odd and another is even, while in \eqref{ell param 2 2} both functions are odd and in \eqref{ell param 2 3} both functions are odd. On the other hand, in \eqref{ell param 2 1} and in \eqref{ell param 2 2} both participating functions have no singularities on the real axis, while in \eqref{ell param 2 3} both have simple poles at $s_2=2\mathsf{K}(k_2)$. Thus, the corresponding parametrizations of the confocal coordinates cover different regions of the plane $\mathbb R^2$ and have, in principle, different geometric features. 

It turns out that any solution of the quadratic relations \eqref{P5a}, \eqref{P5b} admits a corresponding solution of the bilinear relations \eqref{P12a}, \eqref{P12b}, the latter approximating the former in the continuum limit. These solutions can be derived with the help of the addition formulas for the theta functions (or for the Jacobi elliptic functions). As an example, we mention the addition formulas
\begin{align}
 \jac{cn}(s,k)\jac{cn}(s+\eta,k)+\jac{sn}(s,k)\jac{sn}(s+\eta,k)\jac{dn}(\eta,k) & =  \jac{cn}(\eta,k), \qquad \label{Q13a} \\
 \jac{dn}(s,k)\jac{dn}(s+\eta,k) +k^2\jac{sn}(s,k)\jac{sn}(s+\eta,k)\jac{cn}(\eta,k)  & =   \jac{dn}(\eta,k),  \qquad \label{Q13b}
 \end{align}
 which constitute bilinear analogs of the identities
 \begin{align*}
 \jac{cn}^2(s,k)+\jac{sn}^2(s,k) & =  1,\\
\jac{dn}^2(s,k)+k^2\jac{sn}^2(s,k) & =  1.
\end{align*}
As a consequence, we find the following two solutions of the functional equation \eqref{P12a}:
$$
f_1(n_1) = \alpha\jac{cn}(\delta n_1 + c_1,k_1),\quad  g_1(n_1)  = \beta\jac{sn}(\delta n_1 + c_1,k_1),
 $$
 where
$$
\alpha = \textstyle\sqrt{\dfrac{a-b}{\jac{cn}(\frac{\delta}{2},k_1)}},
 \quad \beta = \textstyle\sqrt{(a-b)\dfrac{\jac{dn}(\frac{\delta}{2},k_1)}{\jac{cn}(\frac{\delta}{2},k_1)}},
$$ 
and 
$$
f_1(n_1) = \alpha\jac{dn}(\delta n_1 + c_1,k_1),\quad  g_1(n_1)  = \beta\jac{sn}(\delta n_1 + c_1,k_1),
 $$
 where
$$
\alpha = \textstyle\sqrt{\dfrac{a-b}{\jac{dn}(\frac{\delta}{2},k)}},
 \quad \beta = k_1\textstyle\sqrt{(a-b)\dfrac{\jac{cn}(\frac{\delta}{2},k)}{\jac{dn}(\frac{\delta}{2},k)}}.
$$ 
They approximate solutions \eqref{ell param 1}, resp. \eqref{ell param 2} in the continuum limit $\delta\to 0$.

In the following two sections, we will consider in detail two parametrizations of the continuous and discrete confocal coordinate systems of this kind with very remarkable geometric properties.

\subsection{Confocal coordinates outside of an ellipse, \\ diagonally related to a straight line coordinate system}
\label{sect elliptic IC}

\subsubsection{Continuous case}

Consider a coordinate system \eqref{P3} with
\bela{eq: 2d ell par}
 \bear{r@{}lr@{}l}
 f_1(s_1) &{} = \alpha_1\jac{sn}(s_1,k),\quad & g_1(s_1) &{} = \beta_1\jac{cn}(s_1,k),\as
f_2(s_2) &{} = \alpha_2\dfrac{\jac{dn}(s_2,k)}{\jac{cn}(s_2,k)}, \quad &g_2(s_2) &{} = \beta_2\dfrac{1}{\jac{cn}(s_2,k)},
 \ear
\ela
where $s_1\in[0,2\mathsf{K}(k)]$ and $s_2\in[0,\mathsf{K}(k))$, and the amplitudes $\alpha_1$, $\beta_1$, $\alpha_2$ and $\beta_2$ are chosen as follows:
\bela{eq: 2d ell par ampl}
\alpha_1=\beta_1=\sqrt{a-b}, \quad \alpha_2=\frac{1}{k}\sqrt{a-b}, \quad \beta_2=\frac{k'}{k}\sqrt{a-b},
\ela
where $k'=\sqrt{1-k^2}$. Observe that  the modulus $k$ in both pairs $(f_1,g_1)$ and $(f_2,g_2)$ is chosen to be the same. The remarkable geometric property mentioned above is this (cf. Figure \ref{Fig cont system outside ellipse}):
\begin{proposition}
In the coordinate system \eqref{P3}, \eqref{eq: 2d ell par}, the points $(x,y)$ with $s_1+s_2=\xi={\rm const}$ lie on straight lines. The same is true for points $(x,y)$ with $s_2-s_1=\eta={\rm const}$. Moreover, all these lines are tangent to the ellipse 
$$
\mathcal E_0: \quad \frac{x^2}{a_0} +\frac{y^2}{b_0} =1,
$$
where
$$
a_0=\frac{1}{k^2}(a-b), \quad b_0=\frac{(k')^2}{k^2}(a-b).
$$
This ellipse belongs to the confocal family \eqref{eq: 2d confocal}.
\end{proposition}
\begin{proof}
Due to the fact that the functions $f_2,g_2$ are even with respect to $s_2$, it is enough to demonstrate the second statement. We set $s_2=s_1+\eta$ and use addition theorems for elliptic functions to derive:
\begin{align*}
x(s_1,s_1+\eta) &  =  \frac{\sqrt{a-b}}{k}\ \frac{\jac{sn}(s_1)\jac{dn}(s_1)\jac{dn}(\eta)-k^2\jac{sn}^2(s_1)\jac{cn}(s_1)\jac{sn}(\eta)\jac{cn}(\eta)}
                             {\jac{cn}(s_1)\jac{cn}(\eta)-\jac{sn}(s_1)\jac{dn}(s_1)\jac{sn}(\eta)\jac{dn}(\eta)}, \\
y(s_1,s_1+\eta) & =  \frac{k'\sqrt{a-b}}{k}\ \frac{\jac{cn}(s_1)-k^2\jac{sn}^2(s_1)\jac{cn}(s_1)\jac{sn}^2(\eta)}
                             {\jac{cn}(s_1)\jac{cn}(\eta)-\jac{sn}(s_1)\jac{dn}(s_1)\jac{sn}(\eta)\jac{dn}(\eta)}.
\end{align*}
For these points, equation $Ax+By=C$ is satisfied with
$$
\frac{\sqrt{a-b}}{k}\ A=-C\jac{sn}(\eta), \quad \frac{k'\sqrt{a-b}}{k}\ B=C\jac{cn}(\eta).
$$
Obviously, for any $\eta$ the coefficients $A,B,C$ satisfy
$$
a_0\left(\frac{A}{C}\right)^2+b_0\left(\frac{B}{C}\right)^2=1
$$
and the quantities $a_0$, $b_0$ obey $a_0-b_0=a-b$.
\end{proof}

\begin{figure}[H]
  \begin{center}
    \scalebox{1.2}{\input{ic-symmetric-cont2_image.pgf}}
  \end{center}
  \caption{
    Two-dimensional continuous confocal coordinate system \eqref{P3}, \eqref{eq: 2d ell par}
    with $a=2$, $b=1$, $k=0.9$.
    Points with $s_1 + s_2 = \text{const}$ as well as points with $s_1 - s_2 = \text{const}$
    lie on straight lines which are tangent to an ellipse $\mathcal E_0$.
    The parametrization is only defined outside $\mathcal E_0$.
  }
  \label{Fig cont system outside ellipse}
\end{figure}

%\note[JT]{$\nu$'s for Koenigs here?}

\subsubsection{Discrete case and ``elliptic'' IC-nets}

A solution of the functional equations \eqref{P12a} and \eqref{P12b} which approximates \eqref{eq: 2d ell par} in the continuum limit $\delta\to 0$, is given by
\bela{S13}
 \bear{r@{}lr@{}l}
 f_1(n_1) &{} = \hat\alpha_1\jac{sn}(\delta n_1 + c_1,k),\quad & g_1(n_1) &{} = \hat\beta_1\jac{cn}(\delta n_1 + c_1,k),\as
 f_2(n_2) &{} = \hat\alpha_2\dfrac{\jac{dn}(\delta n_2 + c_2,k)}{\jac{cn}(\delta n_2 + c_2,k)},\quad & g_2(n_2) &{} = \hat\beta_2\dfrac{1}{\jac{cn}(\delta n_2 + c_2,k)}.
 \ear
\ela
Using addition theorems for elliptic functions, we easily see that this is a solution if  
\begin{equation}\label{S14}
\hat\alpha_1 = \textstyle\sqrt{(a-b)\dfrac{\jac{dn}(\frac{\delta}{2},k)}{\jac{cn}(\frac{\delta}{2},k)}}, \quad \hat\beta_1 = \textstyle\sqrt{\dfrac{a-b}{\jac{cn}(\frac{\delta}{2},k)}},
\end{equation}
and
\begin{equation}\label{S14a}
 \hat\alpha_2 =\displaystyle \frac{1}{k}\hat\alpha_1,\quad 
 \hat\beta_2 = \displaystyle\frac{k'}{k}\hat\beta_1.
\end{equation}
Here, the constants $\delta,k,c_1,c_2$ are arbitrary except that $0<k^2<1$ and $\jac{cn}(\tfrac{\delta}{2},k)>0$.  However, for reasons of symmetry and closure, one should choose $c_1=c_2=0$ and
\bela{S15}
  \delta = \frac{\mathsf{K}(k)}{m},\quad m\in \tfrac{1}{2}\N,
\ela
so that the parameters $n_i$ may be restricted to $n_1\in[ -2m,2m]$ and $n_2\in[0,m-\frac{1}{2}]$. 
The same computation as in the previous section allows us to show (cf. Figure \ref{Fig elliptic IC}):
\begin{proposition}
The points $(x,y)$ with $n_1+n_2=\xi={\rm const}$ lie on straight lines. The same is true for points $(x,y)$ with $n_1-n_2=\eta={\rm const}$. Moreover, all these lines are tangent to the ellipse 
\bela{S15a}
\hat{\mathcal E}_0: \quad \frac{x^2}{\hat a_0} +\frac{y^2}{\hat b_0} =1,
\ela
where
$$
\hat a_0=(a-b)\frac{1}{k^2}\dfrac{\jac{dn}^2(\frac{\delta}{2},k)}{\jac{cn}^2(\frac{\delta}{2},k)}, \quad 
\hat b_0=(a-b)\frac{(k')^2}{k^2}\dfrac{1}{\jac{cn}^2(\frac{\delta}{2},k)}.
$$
This ellipse belongs to the confocal family \eqref{eq: 2d confocal}, since $\hat a_0-\hat b_0=a-b$.
\end{proposition}

\begin{figure}[H]
  \begin{center}
    \scalebox{1.2}{\input{ic-symmetric2_image.pgf}}
  \end{center}
  \caption{
    Two-dimensional discrete confocal coordinate system \eqref{P10}, \eqref{S13}
    on $\left(\tfrac{1}{2}\Z\right)^2$ with $a=2$, $b=1$, $k=0.9$, $m=3$, $c_1 = c_2 = 0$.
    Points with $n_1 + n_2 = \text{const}$ as well as points with $n_1 - n_2 = \text{const}$
    lie on straight lines which are tangent to an ellipse.
    The parametrization is only defined outside this ellipse.
  }
  \label{Fig elliptic IC}
\end{figure}

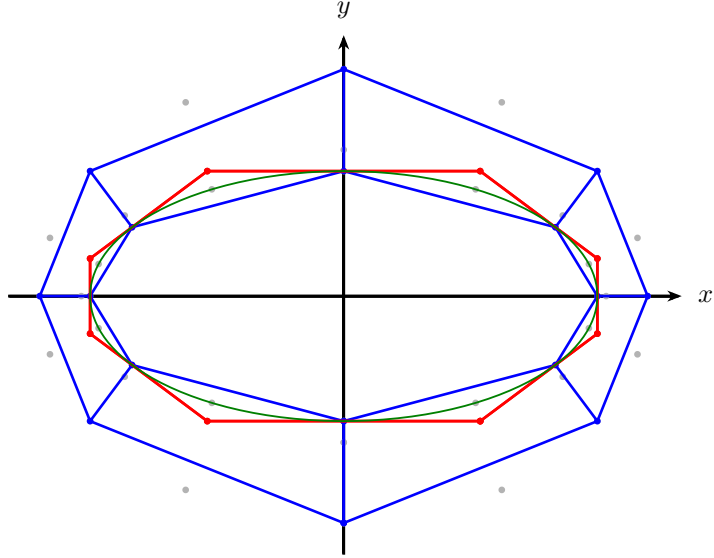
\begin{figure}[H]
  \centering
  \input{ic-symmetric-with-ellipse1_image.pgf}
  \caption{
    Gray points show discrete confocal coordinates \eqref{P10}, \eqref{S13}
    on $\left(\tfrac{1}{2}\Z\right)^2$ with $a=2$, $b=1$, $k=0.9$, $m=2$, $c_1=c_2=0$.
    The two dual sublattices on $\Z^2$ and on $\left( \Z + \tfrac{1}{2} \right)^2$
    are shown in blue and in red, respectively.
    The innermost blue discrete ellipse, corresponding to $n_2=0$, is inscribed in $\hat{\mathcal E}_0$,
    while the neighboring red discrete ellipse, corresponding to $n_2=\frac{1}{2}$, is circumscribed about $\hat{\mathcal E}_0$.
  }
    \label{Fig: discrete ellipses IC detailed}
\end{figure}

Consider the case $c_2=0$. Then a short computation shows that the vertices of the innermost discrete ellipse $\bx(n_1,0)$ lie on the ellipse $\hat{\mathcal E}_0$. The tangent line to $\hat{\mathcal E}_0$  at the point $\bx(n_1,0)$ 
%given by
%$$
%\frac{x x(n_1,0)}{\hat a_0} +\frac{y y(n_1,0)}{\hat b_0} =1
%$$
contains the vertices $\bx(n_1+m,m)$ and $\bx(n_1-m,m)$, $m\in\frac{1}{2}\Z_{\ge 0}$. In particular, this tangent line contains the edge $[\bx(n_1+\frac{1}{2},\frac{1}{2}),\bx(n_1-\frac{1}{2},\frac{1}{2})]$. In other words, the innermost discrete ellipse $n_2=0$ is inscribed in $\hat{\mathcal E}_0$, while the neighboring discrete ellipse $n_2=\frac{1}{2}$ is circumscribed about $\hat{\mathcal E}_0$, with the points of contact being the vertices of the discrete ellipse $n_2=0$ (cf. Figure \ref{Fig: discrete ellipses IC detailed}).
%Furthermore, it turns out that the affine transforms (in the above sense) of the vertices of the circumscribing discrete ellipse also lie on the ellipse \eqref{S15a}.

A further important observation is that the points \eqref{P10} with $f_1,f_2,g_1,g_2$ given by \eqref{S13} upon an affine transformation
\begin{equation}
(x,y)\mapsto(\alpha x,\beta y)\quad {\rm with}\quad 
  \alpha = \dfrac{\jac{cn}(\tfrac{\delta}{2},k)}{\jac{dn}(\tfrac{\delta}{2},k)},\quad \beta =\jac{cn}(\tfrac{\delta}{2},k)
  \label{eq:affine-trafo}
\end{equation}
lie on continuous conics given by the parametrization \eqref{eq: 2d ell par}, i.e., on conics of the original confocal familiy \eqref{eq: 2d confocal}. By the Theorem of Graves-Chasles (see \cite{Akopyan-Bobenko}), all elementary quadrilaterals of the diagonal net upon this affine transformation become circumscribed around circles. This means that the discrete confocal quadrics with the parametrization \eqref{P10}, \eqref{S13}
constitute affine images of ``incircular nets'' (IC-nets) studied in \cite{Akopyan-Bobenko}. An additional computation sketched in \cite{BSST18} shows that, amazingly, the centers of all incircles coincide with the original points of the discrete confocal coordinate system.

\begin{figure}[H]
  \begin{center}
    \scalebox{1.2}{\input{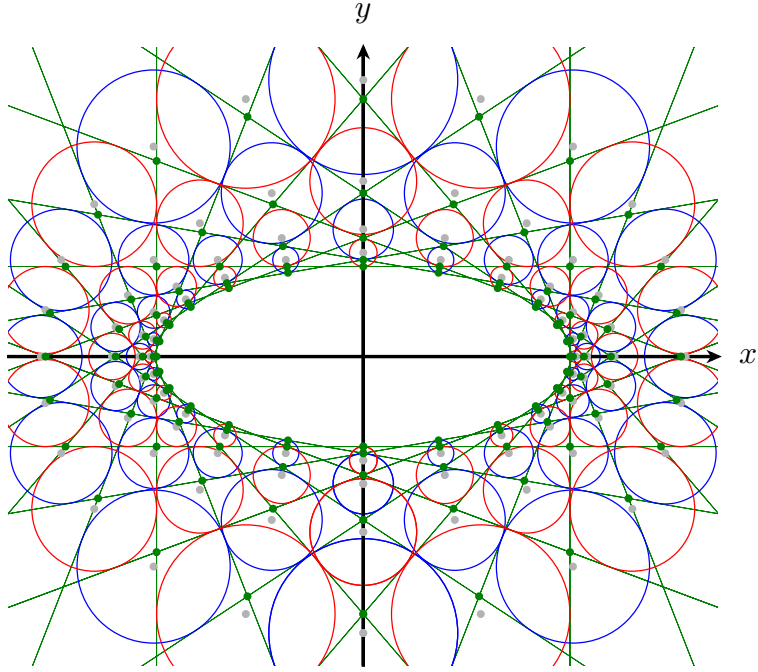}}
  \end{center}
  \caption{
    Grey points show discrete confocal coordinates \eqref{P10}, \eqref{S13}
    on $\left(\tfrac{1}{2}\Z\right)^2$ with $a=2$, $b=1$, $k=0.9$, $m=3$, $c_1 = c_2 = 0$.
    Green points show the net after the affine transformation \eqref{eq:affine-trafo}.
    The straight lines through the transformed net constitute
    two incircular nets. The grey points are the centers of the incircles.
  }
\end{figure}

\subsection{Confocal coordinates outside of a hyperbola, \\ diagonally related to a straight line coordinate system} 
\label{sect hyperbolic IC}

\subsubsection{Continuous case}
Similar properties to those mentioned above has the following coordinate system:

\bela{eq: 2d hyp par}
 \bear{r@{}lr@{}l}
 f_1(s_1) &{} = \alpha_1\jac{sn}(s_1,k),\quad & g_1(s_1) &{} = \beta_1\jac{dn}(s_1,k),\as
f_2(s_2) &{} = \alpha_2\dfrac{1}{\jac{sn}(s_2,k)},\quad & g_2(s_2)  &{} = \beta_2\dfrac{\jac{cn}(s_2,k)}{\jac{sn}(s_2,k)},
 \ear
\ela
where 
$$
 \alpha_1 =  k\textstyle\sqrt{a-b},\quad  \alpha_2=\beta_1=\beta_2 = \sqrt{a-b}.
$$
\begin{proposition}
In the coordinate system \eqref{P3}, \eqref{eq: 2d hyp par}, the points $(x,y)$ with $s_1+s_2=\xi={\rm const}$ lie on straight lines. The same is true for points $(x,y)$ with $s_2-s_1=\eta={\rm const}$. Moreover, all these lines are tangent to the hyperbola (cf. Figure \ref{Fig system outside hyp})
$$
\mathcal H_0: \quad \frac{x^2}{a_0} -\frac{y^2}{c_0} =1,
$$
where
$$
a_0=k^2(a-b), \quad c_0=(1-k^2)(a-b).
$$
This hyperbola belongs to the confocal family \eqref{eq: 2d confocal}.
\end{proposition}
\begin{proof}
Due to the fact that the functions $f_2,g_2$ are odd with respect to $s_2$, it is enough to demonstrate the second statement. We set $s_2=s_1+\eta$ and use addition theorems for elliptic functions to derive:
\begin{eqnarray*}
x(s_1,s_1+\eta) &  = & k\sqrt{a-b}\ \frac{\jac{sn}(s_1)\big(1-k^2\jac{sn}^2(s_1)\jac{sn}^2(\eta)\big)}
                             {\jac{sn}(s_1)\jac{cn}(\eta)\jac{dn}(\eta)+\jac{cn}(s_1)\jac{dn}(s_1)\jac{sn}(\eta)}, \\
y(s_1,s_1+\eta) & = & \sqrt{a-b}\ \frac{\jac{dn}(s_1)\big(\jac{cn}(s_1)\jac{cn}(\eta)-\jac{sn}(s_1)\jac{dn}(s_1)\jac{sn}(\eta)\jac{dn}(\eta)\big)}
                             {\jac{sn}(s_1)\jac{cn}(\eta)\jac{dn}(\eta)+\jac{cn}(s_1)\jac{dn}(s_1)\jac{sn}(\eta)}.
\end{eqnarray*}
For these points, equation $Ax+By=C$ is satisfied with
$$
Ak\sqrt{a-b}\jac{cn}(\eta)=C\jac{dn}(\eta), \quad B\sqrt{a-b}\jac{cn}(\eta)=C\jac{sn}(\eta).
$$
Obviously, for any $\eta$ the coefficients $A,B,C$ satisfy
$$
a_0\left(\frac{A}{C}\right)^2-c_0\left(\frac{B}{C}\right)^2=\frac{\jac{dn}^2(\eta)}{\jac{cn}^2(\eta)}-\frac{(1-k^2)\jac{sn}^2(\eta)}{\jac{cn}^2(\eta)}=1
$$
and the quantities $a_0$, $c_0$ obey $a_0+c_0=a-b$.
\end{proof}

\begin{figure}[H]
  \begin{center}
    \scalebox{1.2}{\input{ic-symmetric-hyperbolic-cont2_image.pgf}}
  \end{center}
  \caption{
    Two-dimensional continuous confocal coordinate system \eqref{P3}, \eqref{eq: 2d hyp par}
    with $a=2$, $b=1$, $k=0.9$.
    Points with $s_1 + s_2 = \text{const}$ as well as points with $s_1 - s_2 = \text{const}$
    lie on straight lines which are tangent to a hyperbola $\mathcal H_0$.
    The parametrization is only defined outside $\mathcal H_0$.
  }
  \label{Fig system outside hyp}
\end{figure}

\subsubsection{Discrete case and ``hyperbolic'' IC-nets}
A solution of \eqref{P12a}, \eqref{P12b} which approximates \eqref{eq: 2d hyp par} in the continuum limit $\delta\to 0$ reads:
\bela{T4}
 \bear{r@{}lr@{}l}
 f_1(n_1) &{} = \hat\alpha_1\jac{sn}(\delta n_1 + c_1,k),\quad & g_1(n_1) &{} = \hat\beta_1\jac{dn}(\delta n_1 + c_1,k),\as
f_2(n_2) &{} = \hat\alpha_2\dfrac{1}{\jac{sn}(\delta n_2 + c_2,k)},\quad & 
g_2(n_2)  &{} = \hat\beta_2\dfrac{\jac{cn}(\delta n_2 + c_2,k)}{\jac{sn}(\delta n_2 + c_2,k)},
 \ear
\ela
where 
$$
 \hat\alpha_2 = \displaystyle \frac{1}{k}\hat\alpha_1=\sqrt{(a-b)\dfrac{\jac{cn}(\frac{\delta}{2},k)}{\jac{dn}(\frac{\delta}{2},k)}},\quad 
 \hat\beta_2 = \hat\beta_1 = \sqrt{\dfrac{a-b}{\jac{dn}(\frac{\delta}{2},k)}}.
$$
\begin{proposition}
The points $(x,y)$ with $n_1+n_2=\xi={\rm const}$ lie on straight lines. The same is true for points $(x,y)$ with $n_2-n_1=\eta={\rm const}$. Moreover, all these lines are tangent to the hyperbola (cf. Figure \ref{Fig hyp IC})
\bela{S15c}
\hat{\mathcal H}_0: \quad \frac{x^2}{\hat a_0} -\frac{y^2}{\hat c_0} =1,
\ela
where
$$
\hat a_0=(a-b)k^2\dfrac{\jac{cn}^2(\frac{\delta}{2},k)}{\jac{dn}^2(\frac{\delta}{2},k)}, \quad 
\hat c_0=(a-b)(1-k^2)\dfrac{1}{\jac{dn}^2(\frac{\delta}{2},k)}.
$$
This hyperbola belongs to the confocal family \eqref{eq: 2d confocal}, since $\hat a_0+\hat c_0=a-b$.
\end{proposition}

\begin{figure}[H]
  \begin{center}
    \scalebox{1.2}{\input{ic-symmetric-hyperbolic2_image.pgf}}
  \end{center}
  \caption{
    Two-dimensional discrete confocal coordinate system \eqref{P10}, \eqref{S13}
    on $\left(\tfrac{1}{2}\Z\right)^2$ with $a=2$, $b=1$, $k=0.9$, $m=3$, $c_1 = c_2 = 0$.
    Points with $n_1 + n_2 = \text{const}$ as well as points with $n_1 - n_2 = \text{const}$
    lie on straight lines which are tangent to a hyperbola.
    The parametrization is only defined outside this hyperbola.
  }
  \label{Fig hyp IC}
\end{figure}

An affine transformation converting the discrete confocal system \eqref{P10} with \eqref{T4} into ``hyperbolic'' IC-nets is characterized by 
\bela{T6}
  (x,y)\mapsto(\alpha x,\beta y) \quad {\rm with}\quad \alpha = \dfrac{\jac{dn}(\tfrac{\delta}{2},k)}{\jac{cn}(\tfrac{\delta}{2},k)},\quad 
  \beta =\jac{dn}(\tfrac{\delta}{2},k).
\ela

\subsection{Confocal coordinates, diagonally related to vertical lines and a hyperbolic pencil of circles}
\label{s: akopyan}

\subsubsection{Continuous case}
Consider a coordinate system \eqref{P3} with
\bela{eq: 2d hyp pencil}
 \bear{r@{}lr@{}l}
 f_1(s_1) &{} = \sqrt{a-b}\ e^{s_1}, \quad & g_1(s_1) &{} = \sqrt{a-b}\ \sqrt{1-e^{2s_1}},\as
 f_2(s_2) &{} = \sqrt{a-b}\ e^{s_2}, \quad & g_2(s_2) &{} = \sqrt{a-b}\ \sqrt{e^{2s_2}-1},
 \ear
\ela
where $s_1<0$ and $s_2>0$.
The image is the first quadrant where the $y$-axis is approached in the limit $s_1 \rightarrow -\infty$.
This parametrization is diagonally related to a hyperbolic pencil of circles
which has the two foci of the confocal conics as limiting points.
The following statement from \cite{Akopyan} shows that the confocal ellipses, confocal hyperbolas, and the pencil of circles 
constitute a \emph{3-web}.

\begin{figure}[H]
  \begin{center}
    \scalebox{1.2}{\input{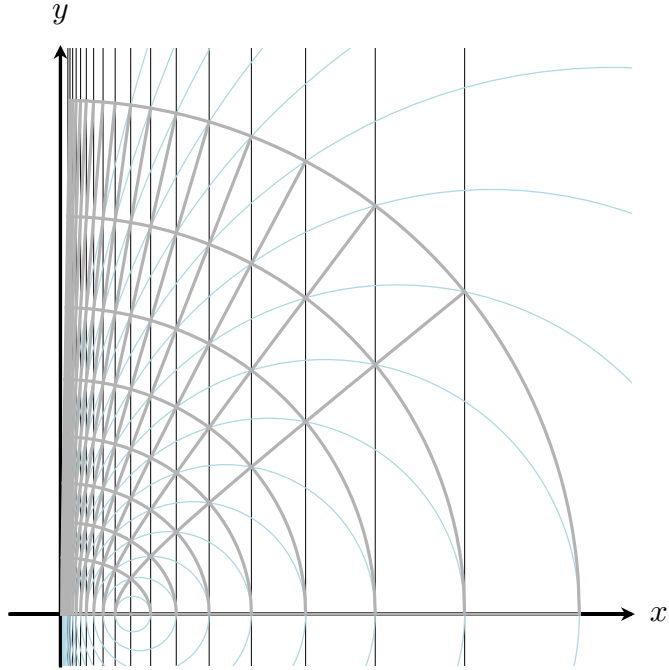}} 
  \end{center}
  \caption{
    Two-dimensional continuous confocal coordinate system \eqref{P3}, \eqref{eq: 2d hyp pencil}
    with $a=2$, $b=1$.
    Points with $s_1 + s_2 = {\rm const}$ lie on vertical lines,
    while points with $s_2 - s_1 = {\rm const}$ lie on circles of a hyperbolic pencil
    which has as limiting points the two foci of the confocal conics.
  }
  \label{fig:hyperbolic-pencil-cont}
\end{figure}

\begin{proposition}
  The points $(x,y)$ with $s_1+s_2={\rm const}$ lie on vertical lines.
  The points $(x,y)$ with $s_2-s_1=\eta={\rm const}$ lie on circles
  with centers $(c(\eta),0) = (\sqrt{a-b}\ \cosh(\eta), 0)$ and radii $r(\eta) = \sqrt{a-b}\ \sinh(\eta)$  (cf.~Figure~\ref{fig:hyperbolic-pencil-cont}).
\end{proposition}

% \begin{proof}
%   The first statement is equivalent to
%   \begin{equation}
%     f_1'f_2 - f_1f_2' = 0.
%   \end{equation}
%   For the second statement one computes that points $(x,y)$ with $s_1-s_2=\eta={\rm const}$ satisfy
%   \begin{equation}
%     x^2 + y^2 - 2c(\eta)y + c(\eta)^2 - r(\eta)^2 = 0.
%   \end{equation}
% \end{proof}

\subsubsection{Discrete case}
A solution of the functional equations \eqref{P12a} and \eqref{P12b} which approximates \eqref{eq: 2d hyp pencil} 
in the continuum limit $\delta\to 0$, is given by
\bela{eq: 2d hyp pencil discr}
 \bear{r@{}lr@{}l}
 f_1(n_1) &{} = \alpha_1~e^{\delta(n_1 + c_1)},\quad & g_1(n_1) &{} = \beta_1~\displaystyle\frac{\Gamma_{e^{-2\delta}}(-n_1-c_1+\frac{3}{4})}{\Gamma_{e^{-2\delta}}(-n_1-c_1+\frac{1}{4})},\as
 f_2(n_2) &{} = \alpha_2~e^{\delta(n_2 + c_2 )},\quad & g_2(n_2) &{} = \beta_2~\displaystyle\frac{\Gamma_{e^{2\delta}}(n_2+c_2+\frac{3}{4})}{\Gamma_{e^{2\delta}}(n_2+c_2+\frac{1}{4})},
 \ear
\ela
where
\begin{equation}\
  \begin{aligned}
    \alpha_1 &= \alpha_2 = \sqrt{a-b},\\
    \beta_1 &= \sqrt{a-b}\ \sqrt{1-e^{-2\delta}}, \quad \beta_2 = \sqrt{a-b}\ \sqrt{e^{2\delta}-1},
  \end{aligned}
\end{equation}
and $\Gamma_q$ denotes the $q$-gamma function, which satisfies
\begin{equation}
  \Gamma_q(z+1) = \frac{1-q^z}{1-q} \Gamma_q(z).
\end{equation}
Boundary conditions $y |_{n_1=0} = 0$, $y |_{n_2=0} = 0$ may be achieved by setting $c_1=\frac{1}{4}$, $c_2=-\frac{1}{4}$.

\begin{figure}[H]
  \begin{center}
    \scalebox{1.2}{\input{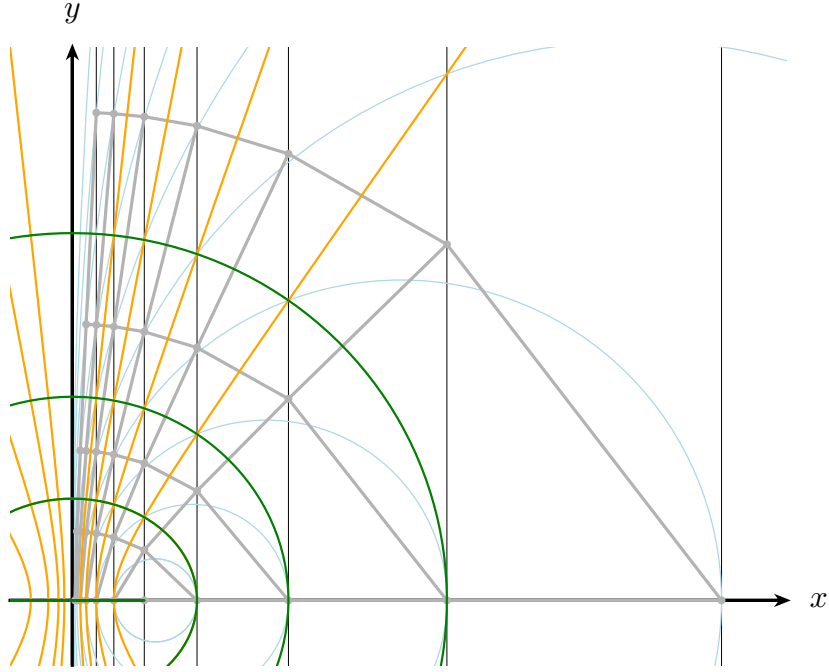}} 
  \end{center}
  \caption{
    Two-dimensional discrete confocal coordinate system \eqref{P10}, \eqref{eq: 2d hyp pencil discr}
    with $a=2$, $b=1$, $c_1 = \frac{1}{4}$, $c_2 = -\frac{1}{4}$.
    Points with $n_1 + n_2 = {\rm const}$ lie on vertical lines.
    In addition to the polarity relation given by the corresponding classical confocal ellipses (green) and hyperbolas (orange)
    two points adjacent to a diagonal $n_2 - n_1 = {\rm const}$ also satisfy a polarity relation
    with respect to circles of a hyperbolic pencil.
  }
  \label{fig:hyperbolic-pencil}
\end{figure}

\begin{proposition}
  The points $(x,y)$ with $n_1 + n_2 = {\rm const}$ lie on vertical lines.
  Pairs of points 
 $$
      (x,y) = \big(x(n_1+\tfrac{1}{2},n_2),y(n_1+\tfrac{1}{2},n_2)\big)
 $$
 and
 $$
 (\tilde x,\tilde y) = \big(x(n_1,n_2+\tfrac{1}{2}),y(n_1,n_2+\tfrac{1}{2})\big),
$$
  which are adjacent to the diagonal $n_2-n_1 = \eta = {\rm const}$, are related by the polarity
  \begin{equation}
    x \tilde x + y \tilde y - c(\eta)\left(x + \tilde x\right) + c(\eta)^2 - r(\eta)^2 = 0
  \end{equation}
  with respect to the circle with the center $(c(\eta), 0)$ and the radius $r(\eta)$, where
  $$ 
   c(\eta)= \sqrt{a-b}\ \cosh(\delta(\eta + c_1 - c_2)), \quad  r(\eta) = \sqrt{a-b}\ \sinh(\delta(\eta + c_1 - c_2))
   $$
 (cf.~Figure~\ref{fig:hyperbolic-pencil}).
 
  The two classical confocal conics corresponding to 
  the parameter values $u(n_1+\tfrac{1}{4})$ and $v(n_2+\tfrac{1}{4})$,
  and the circle with center $(c(\eta), 0)$ and radius $r(\eta)$, $\eta=n_2-n_1$, intersect at a point.
\end{proposition}

This can be checked by a direct computation. There exists an analogous parametrization diagonally related to horizontal lines and an elliptic pencil of circles, both in the continuous case (see \cite{Akopyan}) and in the discrete case.

\subsection{Confocal coordinates, diagonally related to two families of concentric circles}
\label{s: 2 concentric}

\subsubsection{Continuous case}
Consider a coordinate system \eqref{P3} with
\bela{eq: 2d concetric}
 \bear{r@{}lr@{}l}
 f_1(s_1) &{} = \sqrt{a-b}~s_1, \quad & g_1(s_1) &{} = \sqrt{a-b}~\sqrt{1-s_1^2},\as
 f_2(s_2) &{} = \sqrt{a-b}~s_2, \quad & g_2(s_2) &{} = \sqrt{a-b}~\sqrt{s_2^2-1},
 \ear
\ela
where $-1<s_1<1$ and $s_2>1$.
%The image is the upper half-plane.
This parametrization is diagonally related to concentric circles
with centers at the two foci of the confocal conics.

\begin{figure}[H]
  \begin{center}
    \scalebox{0.9}{\input{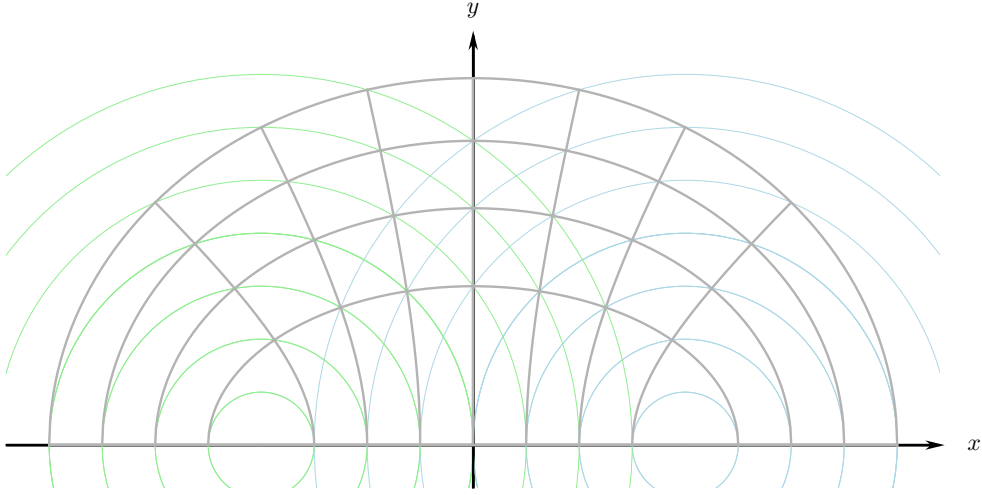}} 
  \end{center}
  \caption{
    Two-dimensional continuous confocal coordinate system \eqref{P3}, \eqref{eq: 2d concetric}
    with $a=2$, $b=1$.
    Points with $s_1 + s_2 = {\rm const}$ and points with $s_2 - s_1 = {\rm const}$ 
    lie on concetric circles with centers at the two foci of the confocal family.
  }
  \label{fig:concetric-cont}
\end{figure}

\begin{proposition}
  The points $(x,y)$ with $s_1+s_2=\xi={\rm const}$ lie on concentric circles
  with the center $\left(-\sqrt{a-b}, 0\right)$ and with the radii $r(\xi) = \sqrt{a-b}\ \xi$.
  The points $(x,y)$ with $s_2-s_1=\eta={\rm const}$ lie on concentric circles
  with the center $\left(\sqrt{a-b}, 0\right)$ and with the radii $r(\eta) = \sqrt{a-b}\ \eta$ (cf. Figure \ref{fig:concetric-cont}).
\end{proposition}

\subsubsection{Discrete case}
A solution of the functional equations \eqref{P12a} and \eqref{P12b} which approximates \eqref{eq: 2d concetric} 
in the continuum limit $\delta\to 0$, is given by
\bela{eq: 2d concetric discr}
 \bear{r@{}lr@{}l}
 f_1(n_1) &{} = \alpha\delta(n_1 + c_1), \quad & g_1(n_1) &{} = \alpha\sqr{n_1+c_1^+}\sqr{-(n_1+c_1^-)+\frac{1}{2}},\as
 f_2(n_2) &{} = \alpha\delta(n_1 + c_1),\quad & g_2(n_2) &{} = \alpha\sqr{n_2+c_2^+}\sqr{n_2+c_2^-},
 \ear
\ela
where $\delta > 0$, $\alpha = \sqrt{a-b}$, $c_1, c_2 \in \R$, and
\begin{equation}
  c_1^\pm = c_1 + \frac{1}{4} \pm \frac{\sqrt{16 + \delta^2}}{4\delta}, \quad c_2^\pm = c_2 + \frac{1}{4} \pm \frac{\sqrt{16 + \delta^2}}{4\delta}.
\end{equation}
If we set $c_1 = c_2 = 0$, and
\begin{equation}
  \delta = \frac{4}{\sqrt{(2l+1)^2 - 1}} ~\Leftrightarrow~ \frac{\sqrt{16 + \delta^2}}{4\delta} = l+\frac{1}{2}
\end{equation}
with some $l \in \N$,
we may let $n_1 \in \left[-\frac{l+1}{2}, \frac{l+1}{2}\right]$, $n_2 \geq \frac{l}{2}$
and achieve boundary conditions
\begin{equation}
  y |_{n_1 = -\frac{l+1}{2}} = y |_{n_1 = \frac{l+1}{2}} = 0, \quad y |_{n_2 = \frac{l}{2}} = 0.
\end{equation}

\begin{figure}[H]
  \begin{center}
    \scalebox{0.9}{\input{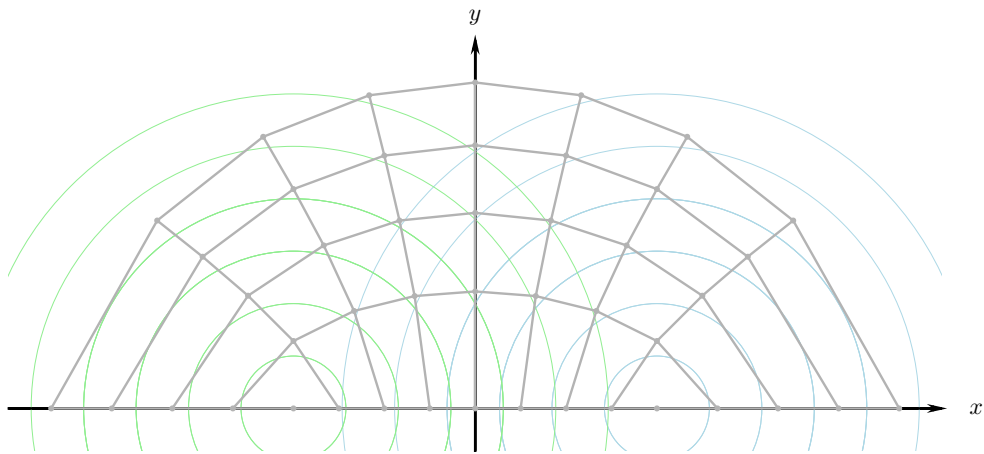}} 
  \end{center}
  \caption{
    Two-dimensional discrete confocal coordinate system \eqref{P10}, \eqref{eq: 2d concetric discr}
    with $a=2$, $b=1$, $c_1 = c_2 = 0$, $2l+1 = 7$.
    Points adjacent to a diagonal $n_2 - n_1 = {\rm const}$ are related by polarity with respect to
    concentric circles with centers in the right focal point.
    Similarly, points adjacent to a diagonal $n_1 + n_2 = {\rm const}$ are related by polarity with respect to
    concentric circles with centers in the left focal point.
  }
  \label{fig:concetric-discr}
\end{figure}

\begin{proposition}
  Pairs of points 
  $$
      (x,y) = \big(x(n_1,n_2-\tfrac{1}{2}),y(n_1,n_2-\tfrac{1}{2})\big)
  $$
and
$$
(\tilde x,\tilde y) = \big(x(n_1+\tfrac{1}{2},n_2),y(n_1+\tfrac{1}{2},n_2)\big),
$$
  which are adjacent to the diagonal $n_1+n_2 = \xi = {\rm const}$, are related by polarity with respect to the circle  with the center $(c, 0) =(\sqrt{a-b}, 0)$ and the radius $r(\xi) = \sqrt{a-b}\,\delta(\xi + c_1 + c_2)$:
  \begin{equation}
    x \tilde x + y \tilde y - c\left(x + \tilde x\right) + c^2 - r(\xi)^2 = 0.
  \end{equation}

Similarly, pairs of points 
$$
      (x,y) = \big(x(n_1+\tfrac{1}{2},n_2),y(n_1+\tfrac{1}{2},n_2)\big)
$$
and
$$
(\tilde x,\tilde y) = \big(x(n_1,n_2+\tfrac{1}{2}),y(n_1,n_2+\tfrac{1}{2})\big),
$$
  which are adjacent to the diagonal $n_2-n_1 = \eta = {\rm const}$, are related by polarity  with respect to the circle
  with the center $(-c, 0)$ and with the radius $\hat{r}(\eta) = \sqrt{a-b}\, \delta(\eta + c_2 - c_1)$:
  \begin{equation}
    x \tilde x + y \tilde y +c\left(x + \tilde x\right) + c^2 - \hat{r}(\eta)^2 = 0
  \end{equation}
(cf. Figure \ref{fig:concetric-discr}).
\end{proposition}

We remark that in this case the corresponding classical ellipses, hyperbolas, and circles participating in the polarity relations are not incident. The proof of all these statements is by direct computation.

%%%%%%%%%%%%%%%%%%%%%%%%%%%%%%%%%%%%%%%%%%%%%%%%%%%%

\section{The case \boldmath $N=3$}\label{section_3d}

\subsection{Classical confocal coordinate systems}

\begin{figure}[H]
  \begin{center}
    \includegraphics[width=0.7\textwidth]{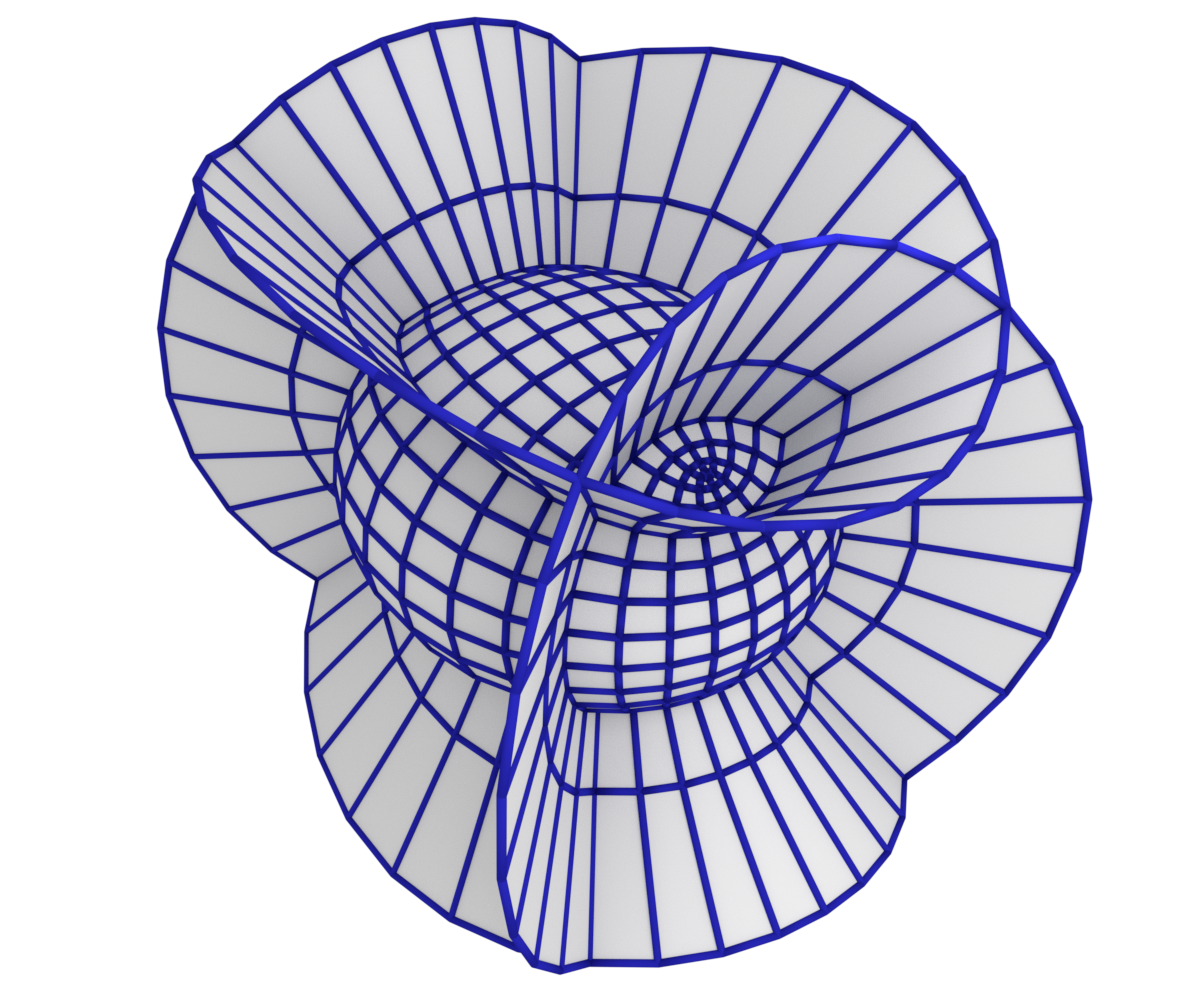}
  \end{center}
  \caption{
    Three-dimensional classical confocal coordinate system \eqref{Q7}
    in terms of Jacobi elliptic functions with $a=8$, $b=4$, $c=0$.
    One quadric of each signature is shown.
  }
  \label{Fig 3d confocal}
\end{figure}

The defining equations  
\bela{Q1}
 \begin{split}
  \frac{x^2}{u+a} + \frac{y^2}{u+b} + \frac{z^2}{u+c} & = 1,\as
  \frac{x^2}{v+a} + \frac{y^2}{v+b} + \frac{z^2}{v+c} &= 1,\as
  \frac{x^2}{w+a} + \frac{y^2}{w+b} + \frac{z^2}{w+c} &= 1
 \end{split}
\ela
of confocal coordinates $\big\{(u,v,w): -a<u<-b<v<-c<w\big\}$ in three dimensions give rise to the expressions
\bela{Q2}
 \begin{split}
  x^2 & = \frac{(u+a)(v+a)(w+a)}{(a-b)(a-c)},\as
  y^2 &= \frac{(u+b)(v+b)(w+b)}{(b-a)(b-c)},\as
  z^2 &= \frac{(u+c)(v+c)(w+c)}{(c-a)(c-b)}.
 \end{split}
\ela
For an arbitrary re-parametrization of the coordinate lines we obtain:
\bela{Q3}
  x = \frac{f_1(s_1)f_2(s_2)f_3(s_3)}{\sqrt{(a-b)(a-c)}},\quad 
  y = \frac{g_1(s_1)g_2(s_2)g_3(s_3)}{\sqrt{(a-b)(b-c)}},\quad 
  z = \frac{h_1(s_1)h_2(s_2)h_3(s_3)}{\sqrt{(a-c)(b-c)}},
\ela
where
\bela{Q4}
\renewcommand{\arraystretch}{1.3}
\left\{\begin{array}{l}
(f_1(s_1))^2 =u + a,\\
( g_1(s_1))^2= -(u+ b),\\
( h_1(s_1))^2 = -(u+ c),
 \end{array}\right. \quad
\left\{\begin{array}{l}
 ( f_2(s_2))^2  =v + a,  \\
  (g_2(s_2))^2 = v+ b,\\
  ( h_2(s_2))^2 = -(v+ c),
   \end{array}\right.  \quad
   \left\{\begin{array}{l}
  (f_3(s_3))^2 = w+a, \\  
  (g_3(s_3))^2  = w+b, \\ 
  (h_3(s_3))^2  = w+c.
   \end{array}\right. 
\ela
Elimination of $u$, $v$ and $w$ leads to functional equations
\begin{eqnarray}\label{Q5}
 & \renewcommand{\arraystretch}{1.4}
 \left\{\begin{array}{l}
  (f_1(s_1))^2 + (g_1(s_1))^2 = a-b, \\ ( f_1(s_1))^2 + (h_1(s_1))^2 = a-c,   
  \end{array}\right.  \quad
\left\{\begin{array}{l}
  (f_2(s_2))^2 - (g_2(s_2))^2 = a-b, \\  (f_2(s_2))^2 + (h_2(s_2))^2= a-c, 
   \end{array}\right. & \nonumber\\
& \renewcommand{\arraystretch}{1.4} 
\left\{\begin{array}{l}   
  (f_3(s_3))^2-(g_3(s_3))^2 =a-b, \\   (f_3(s_3))^2-(h_3(s_3))^2 =a-c. 
   \end{array}\right. & 
\end{eqnarray}
There exists a solution parametrized in terms of Jacobi elliptic functions:
\bela{Q6}
   \bear{r@{}lr@{}lr@{}l}
  f_1(s_1) &{}=\sqrt{a-b}\,\jac{sn}(s_1,k_1),&  g_1(s_1) &{}=\sqrt{a-b}\,\jac{cn}(s_1,k_1),& 
  h_1(s_1) &{}\dis=\sqrt{a-b}\,\frac{\jac{dn}(s_1,k_1)}{k_1},
  \AS
 f_2(s_2) &{}\dis= \sqrt{b-c}\,\frac{\jac{dn}(s_2,k_2)}{k_2}, & g_2(s_2) &{}\dis= \sqrt{b-c}\,\jac{cn}(s_2,k_2), &
 h_2(s_2) &{}\dis= \sqrt{b-c}\,\jac{sn}(s_2,k_2),
  \AS
 f_3(s_3) &{}= \sqrt{a-c}\,\dfrac{1}{\jac{sn}(s_3,k_3)},& g_3(s_3) &{}= \sqrt{a-c}\,\dfrac{\jac{dn}(s_3,k_3)}{\jac{sn}(s_3,k_3)}, &
 h_3(s_3) &{}= \sqrt{a-c}\,\dfrac{\jac{cn}(s_3,k_3)}{\jac{sn}(s_3,k_3)},
 \ear
\ela
where the moduli of the elliptic functions are defined by
\bela{Q6a}
  k_1^2 = \frac{a-b}{a-c},\quad k_2^2 = \frac{b-c}{a-c} = 1-k_1^2,\quad k_3= k_1.
\ela
Hence, we obtain the representation
\bela{Q7}
  \begin{pmatrix} x\\ y\\ z \end{pmatrix}  = \sqrt{a-c} \begin{pmatrix} 
     \jac{sn}(s_1,k_1)\jac{dn}(s_2,k_2)\jac{ns}(s_3,k_3)\\ 
     \jac{cn}(s_1,k_1)\jac{cn}(s_2,k_2)\jac{ds}(s_3,k_3)\\
     \jac{dn}(s_1,k_1)\jac{sn}(s_2,k_2)\jac{cs}(s_3,k_3) \end{pmatrix} 
\ela
of confocal coordinate systems in 3-space. If $\mathsf{K}(k)$ denotes the complete elliptic integral of the first kind, which constitutes the quarter-period of $\jac{sn}(s,k)$, then the parameters may be restricted to $s_1\in [0,4\mathsf{K}(k_1)]$, $s_2\in[0,4\mathsf{K}(k_2)]$ and $s_3\in(0,2\mathsf{K}(k_3))$. Three corresponding coordinate surfaces are depicted in Figure \ref{Fig 3d confocal}.

\subsection{Discrete confocal quadrics}

For any discrete set of confocal quadrics \eqref{Q1}, indexed by $u(n_1+\tfrac{1}{4})$, $v(n_2+\tfrac{1}{4})$, and 
$w(n_3+\tfrac{1}{4})$, we have introduced the discrete confocal quadrics defined by the equations of polarity relating nearest neighbors $\bx(\bn)$ and $\bx(\bn+\tfrac{1}{2}\bsigma)$:

\bela{Q8}
 \begin{split}
  \frac{x(\bn)x(\bn+\tfrac{1}{2}\bsigma)}{u(n_1+\tfrac{1}{4}\sigma_1)+a} + 
  \frac{y(\bn)y(\bn+\tfrac{1}{2}\bsigma)}{u(n_1+\tfrac{1}{4}\sigma_1)+b} + 
  \frac{z(\bn)z(\bn+\tfrac{1}{2}\bsigma)}{u(n_1+\tfrac{1}{4}\sigma_1)+c} & = 1,\as
  \frac{x(\bn)x(\bn+\tfrac{1}{2}\bsigma)}{v(n_2+\tfrac{1}{4}\sigma_2)+a} + 
  \frac{y(\bn)y(\bn+\tfrac{1}{2}\bsigma)}{v(n_2+\tfrac{1}{4}\sigma_2)+b} + 
  \frac{z(\bn)z(\bn+\tfrac{1}{2}\bsigma)}{v(n_2+\tfrac{1}{4}\sigma_2)+c} & = 1,\as
  \frac{x(\bn)x(\bn+\tfrac{1}{2}\bsigma)}{w(n_3+\tfrac{1}{4}\sigma_3)+a} + 
  \frac{y(\bn)y(\bn+\tfrac{1}{2}\bsigma)}{w(n_3+\tfrac{1}{4}\sigma_3)+b} + 
  \frac{z(\bn)z(\bn+\tfrac{1}{2}\bsigma)}{w(n_3+\tfrac{1}{4}\sigma_3)+c} & = 1.
 \end{split}
\ela
This is equivalent to
\bela{Q9}
 \begin{split}
  x(\bn)x(\bn+\tfrac{1}{2}\bsigma) & 
  = \frac{\big(u(n_1+\tfrac{1}{4}\sigma_1)+a\big)\big(v(n_2+\tfrac{1}{4}\sigma_2)+a\big)\big(w(n_3+\tfrac{1}{4}\sigma_3)+a\big)}{(a-b)(a-c)},\as
  y(\bn)y(\bn+\tfrac{1}{2}\bsigma) &
  =  \frac{\big(u(n_1+\tfrac{1}{4}\sigma_1)+b\big)\big(v(n_2+\tfrac{1}{4}\sigma_2)+b\big)\big(w(n_3+\tfrac{1}{4}\sigma_3)+b\big)}{(b-a)(b-c)},\as
  z(\bn)z(\bn+\tfrac{1}{2}\bsigma) &
  =  \frac{\big(u(n_1+\tfrac{1}{4}\sigma_1)+c\big)\big(v(n_2+\tfrac{1}{4}\sigma_2)+c\big)\big(w(n_3+\tfrac{1}{4}\sigma_3)+c\big)}{(c-a)(c-b)}.
 \end{split}
\ela
According to Theorem \ref{th resolution}, these equations can be resolved as follows:
\bela{Q10}
  x(\bn) = \frac{f_1(n_1)f_2(n_2)f_3(n_3)}{\sqrt{(a-b)(a-c)}},\quad 
  y(\bn) = \frac{g_1(n_1)g_2(n_2)g_3(n_3)}{\sqrt{(a-b)(b-c)}},\quad 
  z(\bn) = \frac{h_1(n_1)h_2(n_2)h_3(n_3)}{\sqrt{(a-c)(b-c)}},
\ela
where
\begin{equation*}
\renewcommand{\arraystretch}{1.5}
\left\{\begin{array}{l}
  f_1(n_1)f_1(n_1+\tfrac{1}{2})=u(n_1+\tfrac{1}{4})+a, \\
 g_1(n_1)g_1(n_1+\tfrac{1}{2})=-(u(n_1+\tfrac{1}{4})+b), \\
 h_1(n_1)h_1(n_1+\tfrac{1}{2})=-(u(n_1+\tfrac{1}{4})+c), 
\end{array}\right.  \quad
\left\{\begin{array}{l}
f_2(n_2)f_2(n_2+\tfrac{1}{2})=v(n_2+\tfrac{1}{4})+a,\\
g_2(n_2)g_2(n_2+\tfrac{1}{2})=v(n_2+\tfrac{1}{4})+b,\\
h_2(n_2)h_2(n_2+\tfrac{1}{2})=-(v(n_2+\tfrac{1}{4})+c),
\end{array}\right.
\end{equation*}
\begin{equation}\label{Q11}
\renewcommand{\arraystretch}{1.5}
\left\{\begin{array}{l}
 f_3(n_3)f_3(n_3+\tfrac{1}{2})=w(n_3+\tfrac{1}{4})+a,\\
 g_3(n_3)g_3(n_3+\tfrac{1}{2})=w(n_3+\tfrac{1}{4})+b,\\
 h_3(n_3)h_3(n_3+\tfrac{1}{2})=w(n_3+\tfrac{1}{4})+c.
\end{array}\right.
\end{equation}
 The solution of equations \eqref{Q11} in terms of gamma functions found in \cite{BSST16} and reproduced for general $N$ in Section \ref{sect gamma}, is given (in the first octant) by
\begin{equation*}
\renewcommand{\arraystretch}{1.3}
\left\{ \begin{array}{l}
 f_1(n_1)=\sqr{n_1+\alpha}, \\ g_1(n_1)=\sqr{-n_1-\beta},\\  h_1(n_1)=\sqr{-n_1-\gamma-\frac{1}{2}},
\end{array}\right. \quad
\left\{ \begin{array}{l}
f_2(n_2)=\sqr{n_2+\alpha-\frac{1}{2}},\\ g_2(n_2)=\sqr{n_2+\beta}, \\ h_2(n_2)=\sqr{-n_2-\gamma},
\end{array}\right. 
\end{equation*}
\begin{equation}\label{Q29}
\renewcommand{\arraystretch}{1.3}
\left\{ \begin{array}{l}
 f_3(n_3= \sqr{n_3+\alpha-1}, \\ g_3(n_3)=\sqr{n_3+\beta-\frac{1}{2}}, \\ h_3(n_3)=\sqr{n_3+\gamma},
 \end{array}\right. 
\end{equation}
with $\alpha>\beta>\gamma$ being three integers, and with the identification
\bela{Q30}
  \textstyle a = \alpha + \frac{1}{2},\quad b = \beta+1,\quad c = \gamma + \frac{3}{2}.
\ela

On the other hand, the construction in Theorem \ref{th funct eqs} can be specialized in the case $N=3$ as follows: the nine functions $f_i(n_i)$, $g_i(n_i)$, $h_i(n_i)$ satisfy the functional equations
\begin{equation*}
\renewcommand{\arraystretch}{1.5}
\left\{\begin{array}{l}
f_1(n_1)f_1(n_1+\tfrac{1}{2})+g_1(n_1)g_1(n_1+\tfrac{1}{2})=a-b, \\
f_1(n_1)f_1(n_1+\tfrac{1}{2})+h_1(n_1)h_1(n_1+\tfrac{1}{2})=a-c, 
\end{array}\right.
\end{equation*}
\begin{equation*}
\renewcommand{\arraystretch}{1.5}
\left\{\begin{array}{l}
f_2(n_2)f_2(n_2+\tfrac{1}{2})-g_2(n_2)g_2(n_2+\tfrac{1}{2})=a-b, \\
f_2(n_2)f_2(n_2+\tfrac{1}{2})+h_2(n_2)h_2(n_2+\tfrac{1}{2})=a-c, 
\end{array}\right.
\end{equation*}
\begin{equation}\label{Q12}
\renewcommand{\arraystretch}{1.5}
\left\{\begin{array}{l}
f_3(n_3)f_3(n_3+\tfrac{1}{2})-g_3(n_3)g_3(n_3+\tfrac{1}{2})=a-b, \\
f_3(n_3)f_3(n_3+\tfrac{1}{2})-h_3(n_3)h_3(n_3+\tfrac{1}{2})=a-c.
\end{array}\right.
\end{equation}
A solution of system  (\ref{Q12}) in terms of Jacobi elliptic functions reads: 
\bela{Q15}
   \bear{r@{}lr@{}lr@{}l}
  f_1(n_1) &{}\dis=\alpha_1\jac{sn}(\delta_1 n_1,k_1),&
  g_1(n_1) &{}\dis=\beta_1\jac{cn}(\delta_1 n_1,k_1)& 
   h_1(n_1) &{}\dis=\gamma_1\jac{dn}(\delta_1 n_1,k_1),\\[8mm]
  f_2(n_2) &{}\dis= \alpha_2\jac{dn}(\delta_2n_2,k_2),&
  g_2(n_2) &{}\dis=\beta_2\jac{cn}(\delta_2n_2,k_2),&
   h_2(n_2) &{}\dis= \gamma_2\jac{sn}(\delta_2n_2,k_2),\\[8mm]
  f_3(n_3) &{}= \alpha_3\dfrac{1}{\jac{sn}(\delta_3n_3,k_3)},&
  g_3(n_3) &{}\dis= \beta_3\dfrac{\jac{dn}(\delta_3n_3,k_3)}{\jac{sn}(\delta_3n_3,k_3)},&
 h_3(n_3) &{}\dis= \gamma_3\frac{\jac{cn}(\delta_3n_3,k_3)}{\jac{sn}(\delta_3n_3,k_3)},
 \ear
\ela
where the moduli $k_1$, $k_2$, $k_3$ are defined as solutions of the following transcendental equations:
\bela{Q16}
  k_1^2 = \dfrac{a-b}{a-c}\cdot\dfrac{\jac{dn}^2(\tfrac{\delta_1}{2},k_1)}{\jac{cn}^2(\tfrac{\delta_1}{2},k_1)},\quad
  k_2^2 = \dfrac{b-c}{a-c}\cdot\dfrac{\jac{dn}^2(\tfrac{\delta_2}{2},k_2)}{\jac{cn}^2(\tfrac{\delta_2}{2},k_2)},\quad
  k_3^2 = \dfrac{a-b}{a-c}\cdot\dfrac{\jac{dn}^2(\tfrac{\delta_3}{2},k_3)}{\jac{cn}^2(\tfrac{\delta_3}{2},k_3)},
\ela
and the amplitudes $\alpha_1,\ldots,\gamma_3$ are given by
$$
\alpha_1 = \sqrt{(a-b)\dfrac{\jac{dn}(\frac{\delta_1}{2},k_1)}{\jac{cn}(\frac{\delta_1}{2},k_1)}}, \quad \beta_1=\frac{\alpha_1}{\sqrt{\jac{dn}(\frac{\delta_1}{2},k_1)}}, \quad 
\gamma_1=\frac{\alpha_1}{k_1\sqrt{\jac{cn}(\frac{\delta_1}{2},k_1)}},
$$
$$
\gamma_2 = \sqrt{(b-c)\dfrac{\jac{dn}(\frac{\delta_2}{2},k_2)}{\jac{cn}(\frac{\delta_2}{2},k_2)}}, \quad \alpha_2=\frac{\gamma_2}{k_2\sqrt{\jac{cn}(\frac{\delta_2}{2},k_2)}}, \quad 
\beta_2=\frac{\gamma_2}{\sqrt{\jac{dn}(\frac{\delta_2}{2},k_2)}},
$$
$$
\alpha_3 = \sqrt{(a-c)\dfrac{\jac{cn}(\frac{\delta_3}{2},k_3)}{\jac{dn}(\frac{\delta_3}{2},k_3)}}, \quad \beta_3=\frac{\alpha_3}{\sqrt{\jac{dn}(\frac{\delta_3}{2},k_3)}}, \quad 
\gamma_3=\frac{\alpha_3}{\sqrt{\jac{cn}(\frac{\delta_3}{2},k_3)}}.
$$
%

%
%Hence, we obtain discrete confocal quadrics of the form
%%
%\bela{Q17}
%  \left(\bear{c}x\\ y\\ z\ear\right) = \sqrt{a-c}\left(\bear{c}
%     \kappa^{(x)}\jac{sn}(\delta n_1,k)\jac{dn}(\tilde{\delta}n_2,\tilde{k})\jac{ns}(\hat{\delta}n_3,\hat{k})\\ 
%     \kappa^{(y)}\jac{cn}(\delta n_1,k)\jac{cn}(\tilde{\delta}n_2,\tilde{k})\jac{ds}(\hat{\delta}n_3,\hat{k})\\
%     \kappa^{(z)}\jac{dn}(\delta n_1,k)\jac{sn}(\tilde{\delta}n_2,\tilde{k})\jac{cs}(\hat{\delta}n_3,\hat{k})
%\ear\right),
%\ela
%%
%where the coefficients $\kappa^{(\cdot)}$ are given by
%\note[JT]{$c_i$'s?}
%%
%\bela{Q18}
% \begin{split}
%  \kappa^{(x)} & =\textstyle \sqrt{\jac{dc}(\frac{\delta}{2},k)}\,\sqrt{\jac{nd}(\frac{\tilde{\delta}}{2},\tilde{k})}\,\sqrt{\jac{cd}(\frac{\hat{\delta}}{2},\hat{k})}\\
% \kappa^{(y)} & =\textstyle \sqrt{\jac{nc}(\frac{\delta}{2},k)}\,\sqrt{\jac{nc}(\frac{\tilde{\delta}}{2},\tilde{k})}\,\sqrt{\jac{cn}(\frac{\hat{\delta}}{2},\hat{k})}\,\jac{nd}(\frac{\hat{\delta}}{2},\hat{k})\\
% \kappa^{(z)} & =\textstyle \sqrt{\jac{nd}(\frac{\delta}{2},k)}\,\sqrt{\jac{dc}(\frac{\tilde{\delta}}{2},\tilde{k})}\,\sqrt{\jac{nd}(\frac{\hat{\delta}}{2},\hat{k})}.
% \end{split}
%\ela
%
In order for the discrete confocal quadrics to respect the symmetries of their classical counterparts, we set
\bela{Q19}
  \delta_1 = \frac{\mathsf{K}(k_1)}{m_1},\quad  \delta_2 = \frac{\mathsf{K}(k_2)}{m_2},\quad \delta_3 = \frac{\mathsf{K}(k_3)}{m_3},\qquad m_i\in\N.
\ela
The parameters $n_i$ may then be restricted to $n_1\in[ 0,4m_1]$, $n_2\in[0, 4m_2]$ and $n_3 \in(0,2m_3)$. 

\begin{figure}[H]
  \begin{center}
    \includegraphics[width=0.7\textwidth]{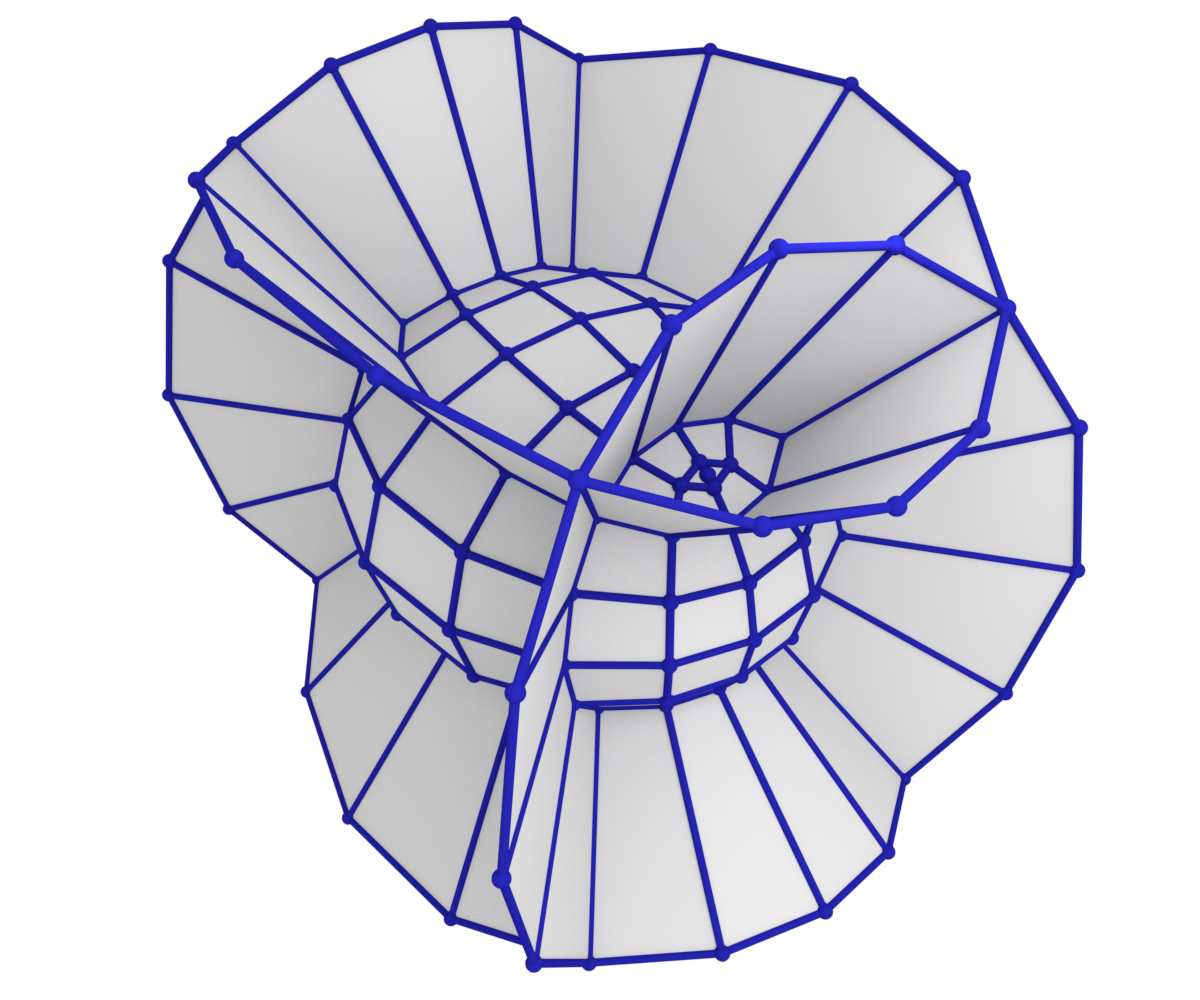}
  \end{center}
  \caption{
    Three-dimensional discrete confocal coordinate system \eqref{Q10}, \eqref{Q15}
    in terms of Jacobi elliptic functions on the sublattice $\Z^2$
    with $a=8$, $b=4$, $c=0$, $m=4$.
    Three discrete quadrics are shown:
    a discrete two-sheeted hyperboloid for $n_1=2$,
    a discrete one-sheeted hyperboloid for $n_2=2$,
    and a discrete ellipsoid for $n_3=2$. A point of valence 2 on a discrete ellipsoid is a discrete analog of an umbilic point.
  }
  \label{fig:confocal3da}
\end{figure}

\begin{figure}[H]
  \begin{center}
    \includegraphics[width=0.7\textwidth]{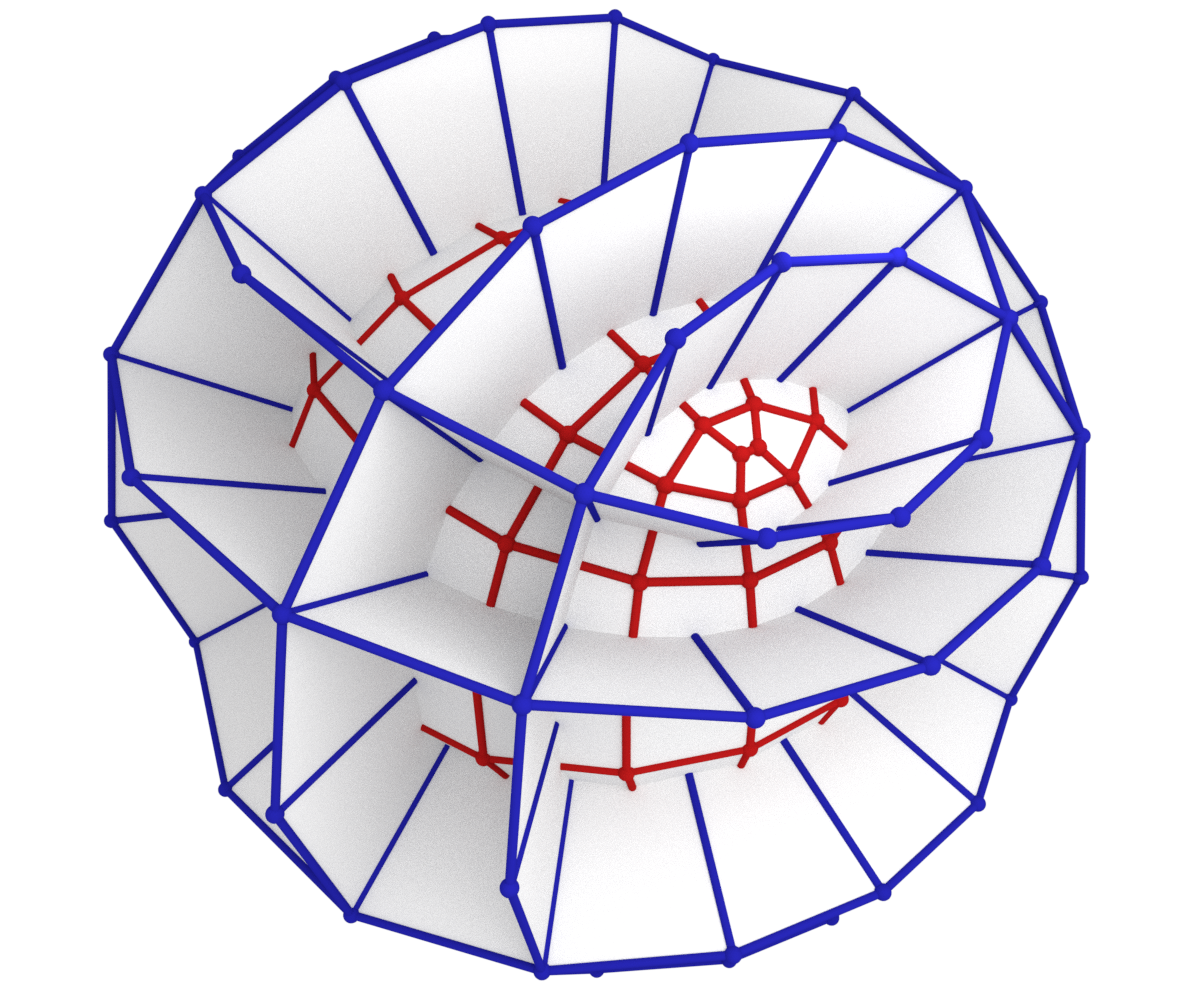}
  \end{center}
  \caption{
    Discrete confocal coordinates with $a=8$, $b=4$, $c=0$, $m=4$.
    Discrete quadrics from the pair of dual orthogonal sublattices
    $\Z^3$ and $\left(\Z+\tfrac{1}{2}\right)^3$ are shown in blue and red respectively:
    two-sheeted hyperboloids $n_1=1,2,~n_2, n_3 \in \Z$,
    one-sheeted hyperboloids $n_2=1,2,~n_1, n_3 \in \Z$,
    ellipse $n_3=1.5,~n_1, n_2 \in \Z+\frac{1}{2}$, }
  \label{fig:confocal3db}
\end{figure}

As in the 2-dimensional case, for arbitrary $m_i$, there exist  special vertices of valence $\neq4$ (cf.\ Figures \ref{fig:confocal3da} and \ref{fig:confocal3db}) which are discrete analogs of the umbilic points on smooth confocal ellipsoids and two-sheeted hyperboloids. In the parametrization \eqref{Q7}, these umbilic points are seen to be
 \begin{gather}\label{Q19a}
  \bx(\epsilon_1\mathsf{K}(k_1),\epsilon_2\mathsf{K}(k_2),s_3),\quad \bx(s_1,\epsilon_2^*\mathsf{K}(k_2),\mathsf{K}(k_3))\as
  \epsilon_1,\epsilon_2\in\{1,3\},\quad \epsilon_2^*\in\{0,2\}
 \end{gather}
respectively. Their discrete analogues are given by the vertices 
 \begin{gather}\label{Q20}
  \bx(\epsilon_1m_1,\epsilon_2m_2,n_3),\quad \bx(n_1,\epsilon_2^*m_2,m_3)\as
  \epsilon_1,\epsilon_2\in\{1,3\},\quad \epsilon_2^*\in\{0,2\}
 \end{gather}
 which have valence 2 (cf.\ Figure \ref{fig:confocal3da}) as may be inferred from the parametrization \eqref{Q15}.

%%%%%%%%%%%%%%%%%%%%%%%%%%%%%%%%%%%%%%%%%%%%%%%%

\begin{appendix}
\section{Euler-Poisson-Darboux equation}

\subsection{Classical Euler-Poisson-Darboux equation}
The discretization of confocal quadrics in \cite{BSST16} was based on
an integrable discretization of the Euler-Poisson-Darboux equation.
We adapt the characterization of confocal coordinates in terms of the Euler-Poisson-Darboux equation
to our present approach by arbitrary re-parametrization of the coordinate lines.

Consider the classical \emph{Euler-Poisson-Darboux system}
\begin{equation*}
  \partial_{u_i}\partial_{u_j}\q = \frac{\gamma}{u_i - u_j} ( \partial_{u_j}\q - \partial_{u_i}\q ), \quad i,j\in\{1,\ldots,N\}
\end{equation*}
with some constant $\gamma \in \R$.
Under re-parametrization $u_i = u_i(s_i)$ this becomes
\begin{equation}
  \label{eq:EPD-system}
  \partial_{s_i}\partial_{s_j} \q = \frac{\gamma}{u_i(s_i)-u_j(s_j)} \left( u_i'(s_i) \partial_{s_j}\q - u_j'(s_j) \partial_{s_i}\q \right).
  \tag{EPD$_\gamma$}
\end{equation}
Confocal coordinates are given by certain factorizable solutions of this equation, and can be characterized as such.

\begin{theorem}
  Let $\q = (x_1, \ldots, x_N)$ be $N$ independent factorizable solutions
  \begin{equation*}
    x_k(s_1, \ldots, s_N) = \prod_{i=1}^N f_i^{k}(s_i)
  \end{equation*}
  of the Euler-Poisson-Darboux system \eqref{eq:EPD-system} with $\gamma = \frac{1}{2}$
  defined on a suitable domain 
   \begin{equation*}
    \mathcal{U} = \{ (s_1, \ldots, s_N) \in \R^N \setbar -a_1 < u_1(s_1) < -a_2 < u_2(s_2)< \cdots < -a_N < u_N(s_N) \}.
  \end{equation*}
  Then the net $\q:\mathcal{U} \rightarrow \R^N$  coincides with confocal coordinates \eqref{eq: x=fff smooth},
  up to independent scaling along the coordinate axes
  $(x_1, \ldots, x_N) \rightarrow (C_1 x_1, \ldots, C_N x_N )$
  with some $C_1, \ldots, C_N~>~0$.
\end{theorem}

\begin{proof}
  A factorizable function
  \begin{equation*}
    x(s_1, \ldots, s_N) = f_1(s_1)\cdots f_N(s_N)
  \end{equation*}
  is a solution of \eqref{eq:EPD-system}, if and only if (compare \cite{BSST16})
  \begin{equation*}
    \frac{f_i'}{f_i} = \frac{u_i'}{2(u_i + a)}, \quad i=1,\ldots,N,
  \end{equation*}
  with some integration constant $a \in \R$.
  The general solution is, up to constant factors, given by
  \begin{equation}
  \renewcommand{\arraystretch}{1.4}
  f_i(s_i) = \left\{\begin{array}{lcl}
     \sqrt{u_i(s_i) + a}, & \; &  \text{if} \quad u_i(s_i) > - a, \\
     \sqrt{-(u_i(s_i) + a)}, & \; & \text{if} \quad u_i(s_i) < - a.
  \end{array}\right.
\end{equation}  
Now $N$ independent separable solutions $x_k(s_1, \ldots, s_N)$, $1\le k\le N$ with constants of integration $a_1 > \ldots > a_N$ are, on the domain $\mathcal U$, given by
\begin{equation}
    \label{eq:separable-solutions}
    x_k(s_1, \ldots, s_N) =  D_k \prod_{i<k} \sqrt{-(u_i(s_i) + a_k)} \prod_{i\ge k} \sqrt{u_i(s_i) + a_k}
\end{equation} 
with some constants $D_k \neq 0$. The choice
\begin{equation}
    \label{eq:confocal-constants}
    D_k^{-1}=\prod_{i<k}\sqrt{a_i-a_k}\prod_{i>k}\sqrt{a_k-a_i}
\end{equation}
is the unique scaling (up to a common factor) for which the parameter curves are pairwise orthogonal (see \cite{BSST16}).
\end{proof}

\subsection{Discrete Euler-Poisson-Darboux equation}
\label{sec:dEPD}

It turns out that discrete confocal coordinate systems may also be characterized in terms of a discrete Euler-Poisson-Darbox equation.

\begin{theorem}
Discrete confocal coordinate systems given by \eqref{eq: x=fff} satisfy the discrete Euler-Poisson-Darboux system with $\gamma = \frac{1}{2}$:
  \begin{equation}    \label{eq:dEPD-system}
    \Delta_i \Delta_j x_k = \frac{1}{u_i - u_j} \left( \halfDelta u_i \Delta_j x_k - \halfDelta u_j \Delta_i x_k \right),
     \tag{dEPD$_{\gamma=\frac{1}{2}}$}
  \end{equation}
  where $x_k=x_k(\bn)$ and $u_i=u_i(n_i+\tfrac{1}{4})$. Here, the difference operator $\Delta$ acts according to
$
     \Delta_i f(\bn) = f(\bn + \be_i) - f(\bn), 
$
%where $\be_i$ is the $i$-th coordinate vector of $\Z^N$, 
and
$$
    \halfDelta u_i=u_i(n_i+\tfrac{3}{4})-u_i(n_i + \tfrac{1}{4}).
$$

Conversely, let $\q = (x_1, \ldots, x_N)$ be $N$ independent factorized solutions
  \begin{equation*}
    x_k(n_1, \ldots, n_N) = \prod_{i=1}^N f_i^{(k)}(n_i)
  \end{equation*}
  of \eqref{eq:dEPD-system} with positive ``discrete squares''
  \begin{equation*}
    f_i^{(k)}(n_i) f_i^{(k)}(n_i + \tfrac{1}{2}) > 0
  \end{equation*}
defined on a suitable domain  
\begin{eqnarray} \label{eq: U discr}
\lefteqn{\mathcal U=\big\{(n_1, \ldots, n_N) \in \left(\tfrac{1}{2}\Z\right)^N:} \nonumber\\
&&    -a_1 < u_1(n_1+\tfrac{1}{4}) < -a_2 < u_2(n_2+\tfrac{1}{4})< \cdots < -a_N < u_N(n_N+\tfrac{1}{4})\big\}.\qquad\qquad
\end{eqnarray}
 Then the net $\q:\mathcal{U} \rightarrow \R^N$  coincides with discrete confocal coordinates
\eqref{eq: x=fff}, \eqref{eq: ff=u+a}  in the first hyperoctant, up to independent scaling along the coordinate axes $(x_1, \ldots, x_N) \rightarrow (C_1 x_1, \ldots, C_N x_N )$ with some constants $C_1, \ldots, C_N>0$.
\end{theorem}

\begin{proof} First, we derive the discrete Euler-Poisson-Darboux equations satisfied by the discrete confocal coordinates given by \eqref{eq: x=fff}.
From equation \eqref{eq: ff=u+a} we obtain
\begin{equation} \label{eq: app 1}
  \frac{f_i^{(k)}(n_i+1)}{f_i^{(k)}(n_i)} = \frac{u_i(n_i+\frac{3}{4}) + a_k}{u_i(n_i+\frac{1}{4}) + a_k},
\end{equation}
which is equivalent to
\begin{equation*}
  \frac{\Delta f_i^{(k)}}{f_i^{(k)}} = \frac{\halfDelta u_i}{u_i + a_k}.
\end{equation*}
%where $f_i^{(k)} = f_i^{(k)}(n_i)$ and $u_i = u_i(n_i + \frac{1}{4})$.
So, for $x_k(\bn)=\prod_{i=1}^N f_i^{(k)}(n_i)$ we obtain:
\begin{equation*}
  \begin{aligned}
    \Delta_i \Delta_j x_k &= \frac{\halfDelta u_i}{u_i + a_k} \frac{\halfDelta u_j}{u_j + a_k} x_k\\
    &= \frac{1}{u_i - u_j} \left(
      \frac{\halfDelta u_i \halfDelta u_j}{u_j + a_k}
      - \frac{\halfDelta u_i \halfDelta u_j}{u_i + a_k}
    \right) x_k\\
    &= \frac{1}{u_i - u_j} \left(
      \halfDelta u_i \Delta_j x_k - \halfDelta u_j \Delta_i x_k
    \right).
  \end{aligned}
\end{equation*}
Conversely, a simple computation shows that a factorizable function
$$
    x(n_1, \ldots, n_N) = f_1(n_1)\cdots f_N(n_N)
$$
  is a solution of \eqref{eq:dEPD-system} if and only if 
$$
\halfDelta u_i\ \frac{f_i}{\Delta f_i}-u_i=\halfDelta u_j\ \frac{f_j}{\Delta f_j}-u_j=a,
$$  
with some constant of integration $a\in\R$.  Equivalently,
\begin{equation*}
\frac{\Delta f_i}{f_i}=\frac{\halfDelta u_i}{u_i+a}\quad\Leftrightarrow\quad   \frac{f_i(n_i+1)}{f_i(n_i)} =\frac{u_i(n_i+\frac{3}{4}) + a}{u_i(n_i+\frac{1}{4}) + a}, \quad i=1,\ldots,N.
\end{equation*}
Here the left-hand sides can be written as
$$
\frac{f_i(n_i+1)}{f_i(n_i)}= \frac{F_i(n_i+\tfrac{3}{4})}{F_i(n_i+\tfrac{1}{4})},
$$
where $F_i(n_i+\tfrac{1}{4}) = f_i(n_i) f_i(n_i + \tfrac{1}{2})$.
Assuming that the discrete squares $F_i(n_i+\tfrac{1}{4})$ are positive,
  the general solution is, up to constant factors, given by
  \begin{equation}
  \renewcommand{\arraystretch}{1.4}
  F_i(n_i+\tfrac{1}{4}) = \left\{\begin{array}{lcl}
     u_i(n_i+\frac{1}{4}) + a, & \; &  \text{if} \quad u_i(n_i+\frac{1}{4}) > - a, \\
     -(u_i(n_i+\frac{1}{4}) + a), & \; & \text{if} \quad u_i(n_i+\frac{1}{4}) < - a.
  \end{array}\right.
\end{equation}  
Now, take $N$ independent factorizable solutions $x_k(n_1, \ldots, n_N)$, $1\le k\le N$,  with the constants of integration $a_1 > \cdots > a_N$.
Define $\mathcal U$ as in \eqref{eq: U discr} with these $a_1,\ldots,a_N$. Then, we find that
  \begin{equation}
    \label{eq:separable-solutions-discrete}
    x_k(\bn)x_k(\bn + \tfrac{1}{2}\bsigma) = D_k^2 \prod_{i<k} -(u_i(n_i+\tfrac{1}{4}\sigma_i) + a_k) \prod_{i\ge k} (u_i(n_i+\tfrac{1}{4}\sigma_i) + a_k) ,
  \end{equation}
where $D_k \neq 0$ are constants. With the choice
  \begin{equation}
    \label{eq:confocal-constants-discrete}
    D_k^{-2}=\prod_{i<k}(a_i-a_k)\prod_{i>k}(a_k-a_i),
  \end{equation}
  expressions \eqref{eq:separable-solutions-discrete} coincide with \eqref{discr confocal gen x thru u}.
\end{proof}

\end{appendix}

%%%%%%%%%%%%%%%%%%%%%%%%%%%%%%%%%%%%%%%%%

\end{document}

%% file: orthogonality.pdf_tex
%% Creator: Inkscape inkscape 0.48.3.1, www.inkscape.org
%% Accompanies image file 'orthogonality.pdf' (pdf, eps, ps)
%%
%% To include the image in your LaTeX document,  write
\begingroup%
  \makeatletter%
  \providecommand\color[2][]{%
    \errmessage{(Inkscape) Color is used for the text in Inkscape, but the package 'color.sty' is not loaded}%
    \renewcommand\color[2][]{}%
  }%
  \providecommand\transparent[1]{%
    \errmessage{(Inkscape) Transparency is used (non-zero) for the text in Inkscape, but the package 'transparent.sty' is not loaded}%
    \renewcommand\transparent[1]{}%
  }%
  \providecommand\rotatebox[2]{#2}%
  \ifx\svgwidth\undefined%
    \setlength{\unitlength}{228.74447013bp}%
    \ifx\svgscale\undefined%
      \relax%
    \else%
      \setlength{\unitlength}{\unitlength * \real{\svgscale}}%
    \fi%
  \else%
    \setlength{\unitlength}{\svgwidth}%
  \fi%
  \global\let\svgwidth\undefined%
  \global\let\svgscale\undefined%
  \makeatother%
  \begin{picture}(1,0.5555234)%
    \put(0,0){\includegraphics[width=\unitlength]{orthogonality.pdf}}%
    \put(0.46,0.58){\color[rgb]{0,0,0}\makebox(0,0)[lb]{\smash{$\q(\bn + \be_i)$}}}%
   % \put(0.08261778,0.3368745){\color[rgb]{0,0,0}\makebox(0,0)[lb]{\smash{$\Delta_i \q(\bn)$}}}%
    \put(-0.245785,0.09743084){\color[rgb]{0,0,0}\makebox(0,0)[lb]{\smash{$\q(\bn+\tfrac{1}{2}\be_i-\tfrac{1}{2}\be_j-\tfrac{1}{2}\be_k)$}}}%
     \put(0.68,0.09743084){\color[rgb]{0,0,0}\makebox(0,0)[lb]{\smash{$\q(\bn+\tfrac{1}{2}\be_i-\tfrac{1}{2}\be_j+\tfrac{1}{2}\be_k)$}}}%
     \put(0.75,0.5){\color[rgb]{0,0,0}\makebox(0,0)[lb]{\smash{$\q(\bn+\tfrac{1}{2}\be_i+\tfrac{1}{2}\be_j+\tfrac{1}{2}\be_k)$}}}%
     \put(-0.05,0.5){\color[rgb]{0,0,0}\makebox(0,0)[lb]{\smash{$\q(\bn+\tfrac{1}{2}\be_i+\tfrac{1}{2}\be_j-\tfrac{1}{2}\be_k)$}}}%
    \put(0.49567522,0.00824815){\color[rgb]{0,0,0}\makebox(0,0)[lb]{\smash{$\q(\bn)$}}}%
  \end{picture}%
\endgroup%

%% file: trigonometric-symmetric-cont1_image.pgf
%% Creator: Matplotlib, PGF backend
%%
%% To include the figure in your LaTeX document, write
%%   \input{<filename>.pgf}
%%
%% Make sure the required packages are loaded in your preamble
%%   \usepackage{pgf}
%%
%% Figures using additional raster images can only be included by \input if
%% they are in the same directory as the main LaTeX file. For loading figures
%% from other directories you can use the `import` package
%% and then include the figures with
%%   \import{<path to file>}{<filename>.pgf}
%%
%% Matplotlib used the following preamble
%%   \usepackage{fontspec}
%%   \setmainfont{DejaVu Serif}
%%   \setsansfont{DejaVu Sans}
%%   \setmonofont{DejaVu Sans Mono}
%%
\begingroup%
\makeatletter%
\begin{pgfpicture}%
\pgfpathrectangle{\pgfpointorigin}{\pgfqpoint{3.316323in}{3.170273in}}%
\pgfusepath{use as bounding box, clip}%
\begin{pgfscope}%
\pgfsetbuttcap%
\pgfsetmiterjoin%
\definecolor{currentfill}{rgb}{1.000000,1.000000,1.000000}%
\pgfsetfillcolor{currentfill}%
\pgfsetlinewidth{0.000000pt}%
\definecolor{currentstroke}{rgb}{1.000000,1.000000,1.000000}%
\pgfsetstrokecolor{currentstroke}%
\pgfsetdash{}{0pt}%
\pgfpathmoveto{\pgfqpoint{0.000000in}{0.000000in}}%
\pgfpathlineto{\pgfqpoint{3.316323in}{0.000000in}}%
\pgfpathlineto{\pgfqpoint{3.316323in}{3.170273in}}%
\pgfpathlineto{\pgfqpoint{0.000000in}{3.170273in}}%
\pgfpathclose%
\pgfusepath{fill}%
\end{pgfscope}%
\begin{pgfscope}%
\pgfsetbuttcap%
\pgfsetmiterjoin%
\definecolor{currentfill}{rgb}{1.000000,1.000000,1.000000}%
\pgfsetfillcolor{currentfill}%
\pgfsetlinewidth{0.000000pt}%
\definecolor{currentstroke}{rgb}{0.000000,0.000000,0.000000}%
\pgfsetstrokecolor{currentstroke}%
\pgfsetstrokeopacity{0.000000}%
\pgfsetdash{}{0pt}%
\pgfpathmoveto{\pgfqpoint{0.135000in}{0.135000in}}%
\pgfpathlineto{\pgfqpoint{3.049509in}{0.135000in}}%
\pgfpathlineto{\pgfqpoint{3.049509in}{2.830000in}}%
\pgfpathlineto{\pgfqpoint{0.135000in}{2.830000in}}%
\pgfpathclose%
\pgfusepath{fill}%
\end{pgfscope}%
\begin{pgfscope}%
\pgfpathrectangle{\pgfqpoint{0.135000in}{0.135000in}}{\pgfqpoint{2.914509in}{2.695000in}} %
\pgfusepath{clip}%
\pgfsetrectcap%
\pgfsetroundjoin%
\pgfsetlinewidth{0.702625pt}%
\definecolor{currentstroke}{rgb}{0.000000,0.000000,0.000000}%
\pgfsetstrokecolor{currentstroke}%
\pgfsetdash{}{0pt}%
\pgfpathmoveto{\pgfqpoint{0.135000in}{1.482500in}}%
\pgfpathlineto{\pgfqpoint{3.049509in}{1.482500in}}%
\pgfusepath{stroke}%
\end{pgfscope}%
\begin{pgfscope}%
\pgfpathrectangle{\pgfqpoint{0.135000in}{0.135000in}}{\pgfqpoint{2.914509in}{2.695000in}} %
\pgfusepath{clip}%
\pgfsetrectcap%
\pgfsetroundjoin%
\pgfsetlinewidth{0.702625pt}%
\definecolor{currentstroke}{rgb}{0.000000,0.000000,0.000000}%
\pgfsetstrokecolor{currentstroke}%
\pgfsetdash{}{0pt}%
\pgfpathmoveto{\pgfqpoint{1.592255in}{0.135000in}}%
\pgfpathlineto{\pgfqpoint{1.592255in}{2.830000in}}%
\pgfusepath{stroke}%
\end{pgfscope}%
\begin{pgfscope}%
\pgfsetbuttcap%
\pgfsetmiterjoin%
\definecolor{currentfill}{rgb}{0.000000,0.000000,0.000000}%
\pgfsetfillcolor{currentfill}%
\pgfsetlinewidth{1.003750pt}%
\definecolor{currentstroke}{rgb}{0.000000,0.000000,0.000000}%
\pgfsetstrokecolor{currentstroke}%
\pgfsetdash{}{0pt}%
\pgfpathmoveto{\pgfqpoint{3.049509in}{1.482500in}}%
\pgfpathlineto{\pgfqpoint{3.020364in}{1.469025in}}%
\pgfpathlineto{\pgfqpoint{3.029108in}{1.482236in}}%
\pgfpathlineto{\pgfqpoint{0.135000in}{1.482236in}}%
\pgfpathlineto{\pgfqpoint{0.135000in}{1.482764in}}%
\pgfpathlineto{\pgfqpoint{3.029108in}{1.482764in}}%
\pgfpathlineto{\pgfqpoint{3.020364in}{1.495975in}}%
\pgfpathclose%
\pgfusepath{stroke,fill}%
\end{pgfscope}%
\begin{pgfscope}%
\pgfsetbuttcap%
\pgfsetmiterjoin%
\definecolor{currentfill}{rgb}{0.000000,0.000000,0.000000}%
\pgfsetfillcolor{currentfill}%
\pgfsetlinewidth{1.003750pt}%
\definecolor{currentstroke}{rgb}{0.000000,0.000000,0.000000}%
\pgfsetstrokecolor{currentstroke}%
\pgfsetdash{}{0pt}%
\pgfpathmoveto{\pgfqpoint{1.592255in}{2.830000in}}%
\pgfpathlineto{\pgfqpoint{1.605730in}{2.800855in}}%
\pgfpathlineto{\pgfqpoint{1.592519in}{2.809598in}}%
\pgfpathlineto{\pgfqpoint{1.592519in}{0.135000in}}%
\pgfpathlineto{\pgfqpoint{1.591991in}{0.135000in}}%
\pgfpathlineto{\pgfqpoint{1.591991in}{2.809598in}}%
\pgfpathlineto{\pgfqpoint{1.578780in}{2.800855in}}%
\pgfpathclose%
\pgfusepath{stroke,fill}%
\end{pgfscope}%
\begin{pgfscope}%
\pgfpathrectangle{\pgfqpoint{0.135000in}{0.135000in}}{\pgfqpoint{2.914509in}{2.695000in}} %
\pgfusepath{clip}%
\pgfsetrectcap%
\pgfsetroundjoin%
\pgfsetlinewidth{0.803000pt}%
\definecolor{currentstroke}{rgb}{0.000000,0.000000,0.000000}%
\pgfsetstrokecolor{currentstroke}%
\pgfsetdash{}{0pt}%
\pgfpathmoveto{\pgfqpoint{1.592255in}{0.135000in}}%
\pgfpathlineto{\pgfqpoint{1.592255in}{2.830000in}}%
\pgfusepath{stroke}%
\end{pgfscope}%
\begin{pgfscope}%
\pgfpathrectangle{\pgfqpoint{0.135000in}{0.135000in}}{\pgfqpoint{2.914509in}{2.695000in}} %
\pgfusepath{clip}%
\pgfsetrectcap%
\pgfsetroundjoin%
\pgfsetlinewidth{0.803000pt}%
\definecolor{currentstroke}{rgb}{0.000000,0.000000,0.000000}%
\pgfsetstrokecolor{currentstroke}%
\pgfsetdash{}{0pt}%
\pgfpathmoveto{\pgfqpoint{0.135000in}{1.482500in}}%
\pgfpathlineto{\pgfqpoint{3.049509in}{1.482500in}}%
\pgfusepath{stroke}%
\end{pgfscope}%
\begin{pgfscope}%
\pgfpathrectangle{\pgfqpoint{0.135000in}{0.135000in}}{\pgfqpoint{2.914509in}{2.695000in}} %
\pgfusepath{clip}%
\pgfsetrectcap%
\pgfsetroundjoin%
\pgfsetlinewidth{0.803000pt}%
\definecolor{currentstroke}{rgb}{0.000000,0.000000,0.000000}%
\pgfsetstrokecolor{currentstroke}%
\pgfsetdash{}{0pt}%
\pgfpathmoveto{\pgfqpoint{1.592255in}{0.135000in}}%
\pgfpathlineto{\pgfqpoint{1.592255in}{2.830000in}}%
\pgfusepath{stroke}%
\end{pgfscope}%
\begin{pgfscope}%
\pgfpathrectangle{\pgfqpoint{0.135000in}{0.135000in}}{\pgfqpoint{2.914509in}{2.695000in}} %
\pgfusepath{clip}%
\pgfsetrectcap%
\pgfsetroundjoin%
\pgfsetlinewidth{0.803000pt}%
\definecolor{currentstroke}{rgb}{0.000000,0.000000,0.000000}%
\pgfsetstrokecolor{currentstroke}%
\pgfsetdash{}{0pt}%
\pgfpathmoveto{\pgfqpoint{0.135000in}{1.482500in}}%
\pgfpathlineto{\pgfqpoint{3.049509in}{1.482500in}}%
\pgfusepath{stroke}%
\end{pgfscope}%
\begin{pgfscope}%
\pgftext[x=1.592255in,y=2.964750in,,base]{\sffamily\fontsize{10.000000}{12.000000}\selectfont \(\displaystyle y\)}%
\end{pgfscope}%
\begin{pgfscope}%
\pgftext[x=3.136945in,y=1.482500in,left,]{\sffamily\fontsize{10.000000}{12.000000}\selectfont \(\displaystyle x\)}%
\end{pgfscope}%
\begin{pgfscope}%
\pgfpathrectangle{\pgfqpoint{0.135000in}{0.135000in}}{\pgfqpoint{2.914509in}{2.695000in}} %
\pgfusepath{clip}%
\pgfsetrectcap%
\pgfsetroundjoin%
\pgfsetlinewidth{1.003750pt}%
\definecolor{currentstroke}{rgb}{0.700000,0.700000,0.700000}%
\pgfsetstrokecolor{currentstroke}%
\pgfsetdash{}{0pt}%
\pgfpathmoveto{\pgfqpoint{2.120227in}{1.482500in}}%
\pgfpathlineto{\pgfqpoint{1.064282in}{1.482500in}}%
\pgfpathlineto{\pgfqpoint{2.120211in}{1.482500in}}%
\pgfpathlineto{\pgfqpoint{2.120211in}{1.482500in}}%
\pgfpathlineto{\pgfqpoint{2.120227in}{1.482500in}}%
\pgfpathlineto{\pgfqpoint{2.120227in}{1.482500in}}%
\pgfusepath{stroke}%
\end{pgfscope}%
\begin{pgfscope}%
\pgfpathrectangle{\pgfqpoint{0.135000in}{0.135000in}}{\pgfqpoint{2.914509in}{2.695000in}} %
\pgfusepath{clip}%
\pgfsetrectcap%
\pgfsetroundjoin%
\pgfsetlinewidth{1.003750pt}%
\definecolor{currentstroke}{rgb}{0.700000,0.700000,0.700000}%
\pgfsetstrokecolor{currentstroke}%
\pgfsetdash{}{0pt}%
\pgfpathmoveto{\pgfqpoint{2.161463in}{1.482500in}}%
\pgfpathlineto{\pgfqpoint{2.160340in}{1.495856in}}%
\pgfpathlineto{\pgfqpoint{2.156975in}{1.509159in}}%
\pgfpathlineto{\pgfqpoint{2.151381in}{1.522357in}}%
\pgfpathlineto{\pgfqpoint{2.143020in}{1.536206in}}%
\pgfpathlineto{\pgfqpoint{2.132206in}{1.549816in}}%
\pgfpathlineto{\pgfqpoint{2.118135in}{1.563898in}}%
\pgfpathlineto{\pgfqpoint{2.100434in}{1.578321in}}%
\pgfpathlineto{\pgfqpoint{2.078749in}{1.592925in}}%
\pgfpathlineto{\pgfqpoint{2.054064in}{1.606848in}}%
\pgfpathlineto{\pgfqpoint{2.025084in}{1.620640in}}%
\pgfpathlineto{\pgfqpoint{1.991572in}{1.634081in}}%
\pgfpathlineto{\pgfqpoint{1.953357in}{1.646923in}}%
\pgfpathlineto{\pgfqpoint{1.912197in}{1.658424in}}%
\pgfpathlineto{\pgfqpoint{1.866473in}{1.668894in}}%
\pgfpathlineto{\pgfqpoint{1.816261in}{1.678041in}}%
\pgfpathlineto{\pgfqpoint{1.763893in}{1.685304in}}%
\pgfpathlineto{\pgfqpoint{1.707683in}{1.690785in}}%
\pgfpathlineto{\pgfqpoint{1.650271in}{1.694097in}}%
\pgfpathlineto{\pgfqpoint{1.592255in}{1.695204in}}%
\pgfpathlineto{\pgfqpoint{1.534238in}{1.694097in}}%
\pgfpathlineto{\pgfqpoint{1.476826in}{1.690785in}}%
\pgfpathlineto{\pgfqpoint{1.420617in}{1.685304in}}%
\pgfpathlineto{\pgfqpoint{1.366195in}{1.677710in}}%
\pgfpathlineto{\pgfqpoint{1.316080in}{1.668490in}}%
\pgfpathlineto{\pgfqpoint{1.268625in}{1.657479in}}%
\pgfpathlineto{\pgfqpoint{1.225994in}{1.645322in}}%
\pgfpathlineto{\pgfqpoint{1.188186in}{1.632313in}}%
\pgfpathlineto{\pgfqpoint{1.153672in}{1.618083in}}%
\pgfpathlineto{\pgfqpoint{1.124001in}{1.603436in}}%
\pgfpathlineto{\pgfqpoint{1.098934in}{1.588611in}}%
\pgfpathlineto{\pgfqpoint{1.078175in}{1.573820in}}%
\pgfpathlineto{\pgfqpoint{1.061389in}{1.559246in}}%
\pgfpathlineto{\pgfqpoint{1.048209in}{1.545044in}}%
\pgfpathlineto{\pgfqpoint{1.037747in}{1.530529in}}%
\pgfpathlineto{\pgfqpoint{1.030408in}{1.516599in}}%
\pgfpathlineto{\pgfqpoint{1.025573in}{1.502517in}}%
\pgfpathlineto{\pgfqpoint{1.023327in}{1.489181in}}%
\pgfpathlineto{\pgfqpoint{1.023327in}{1.475819in}}%
\pgfpathlineto{\pgfqpoint{1.025573in}{1.462483in}}%
\pgfpathlineto{\pgfqpoint{1.030054in}{1.449226in}}%
\pgfpathlineto{\pgfqpoint{1.037247in}{1.435285in}}%
\pgfpathlineto{\pgfqpoint{1.046912in}{1.421555in}}%
\pgfpathlineto{\pgfqpoint{1.059792in}{1.407314in}}%
\pgfpathlineto{\pgfqpoint{1.075332in}{1.393449in}}%
\pgfpathlineto{\pgfqpoint{1.094534in}{1.379298in}}%
\pgfpathlineto{\pgfqpoint{1.117737in}{1.365023in}}%
\pgfpathlineto{\pgfqpoint{1.145246in}{1.350816in}}%
\pgfpathlineto{\pgfqpoint{1.177320in}{1.336894in}}%
\pgfpathlineto{\pgfqpoint{1.212489in}{1.324058in}}%
\pgfpathlineto{\pgfqpoint{1.252281in}{1.311903in}}%
\pgfpathlineto{\pgfqpoint{1.296752in}{1.300705in}}%
\pgfpathlineto{\pgfqpoint{1.343847in}{1.291120in}}%
\pgfpathlineto{\pgfqpoint{1.395242in}{1.282943in}}%
\pgfpathlineto{\pgfqpoint{1.448534in}{1.276687in}}%
\pgfpathlineto{\pgfqpoint{1.505419in}{1.272285in}}%
\pgfpathlineto{\pgfqpoint{1.563209in}{1.270073in}}%
\pgfpathlineto{\pgfqpoint{1.621301in}{1.270073in}}%
\pgfpathlineto{\pgfqpoint{1.679090in}{1.272285in}}%
\pgfpathlineto{\pgfqpoint{1.735975in}{1.276687in}}%
\pgfpathlineto{\pgfqpoint{1.791363in}{1.283233in}}%
\pgfpathlineto{\pgfqpoint{1.842672in}{1.291486in}}%
\pgfpathlineto{\pgfqpoint{1.891569in}{1.301577in}}%
\pgfpathlineto{\pgfqpoint{1.935804in}{1.312906in}}%
\pgfpathlineto{\pgfqpoint{1.976989in}{1.325741in}}%
\pgfpathlineto{\pgfqpoint{2.013259in}{1.339347in}}%
\pgfpathlineto{\pgfqpoint{2.044744in}{1.353456in}}%
\pgfpathlineto{\pgfqpoint{2.071652in}{1.367823in}}%
\pgfpathlineto{\pgfqpoint{2.095301in}{1.382969in}}%
\pgfpathlineto{\pgfqpoint{2.114648in}{1.398025in}}%
\pgfpathlineto{\pgfqpoint{2.130046in}{1.412812in}}%
\pgfpathlineto{\pgfqpoint{2.141874in}{1.427179in}}%
\pgfpathlineto{\pgfqpoint{2.150526in}{1.441003in}}%
\pgfpathlineto{\pgfqpoint{2.156690in}{1.455013in}}%
\pgfpathlineto{\pgfqpoint{2.160195in}{1.468311in}}%
\pgfpathlineto{\pgfqpoint{2.161459in}{1.481665in}}%
\pgfpathlineto{\pgfqpoint{2.161463in}{1.482500in}}%
\pgfpathlineto{\pgfqpoint{2.161463in}{1.482500in}}%
\pgfusepath{stroke}%
\end{pgfscope}%
\begin{pgfscope}%
\pgfpathrectangle{\pgfqpoint{0.135000in}{0.135000in}}{\pgfqpoint{2.914509in}{2.695000in}} %
\pgfusepath{clip}%
\pgfsetrectcap%
\pgfsetroundjoin%
\pgfsetlinewidth{1.003750pt}%
\definecolor{currentstroke}{rgb}{0.700000,0.700000,0.700000}%
\pgfsetstrokecolor{currentstroke}%
\pgfsetdash{}{0pt}%
\pgfpathmoveto{\pgfqpoint{2.291612in}{1.482500in}}%
\pgfpathlineto{\pgfqpoint{2.290555in}{1.507702in}}%
\pgfpathlineto{\pgfqpoint{2.287388in}{1.532828in}}%
\pgfpathlineto{\pgfqpoint{2.282121in}{1.557802in}}%
\pgfpathlineto{\pgfqpoint{2.274769in}{1.582548in}}%
\pgfpathlineto{\pgfqpoint{2.264604in}{1.608724in}}%
\pgfpathlineto{\pgfqpoint{2.252107in}{1.634463in}}%
\pgfpathlineto{\pgfqpoint{2.237320in}{1.659674in}}%
\pgfpathlineto{\pgfqpoint{2.220297in}{1.684271in}}%
\pgfpathlineto{\pgfqpoint{2.201094in}{1.708168in}}%
\pgfpathlineto{\pgfqpoint{2.179780in}{1.731282in}}%
\pgfpathlineto{\pgfqpoint{2.154800in}{1.754984in}}%
\pgfpathlineto{\pgfqpoint{2.127601in}{1.777611in}}%
\pgfpathlineto{\pgfqpoint{2.098288in}{1.799073in}}%
\pgfpathlineto{\pgfqpoint{2.066979in}{1.819285in}}%
\pgfpathlineto{\pgfqpoint{2.033796in}{1.838169in}}%
\pgfpathlineto{\pgfqpoint{1.996632in}{1.856693in}}%
\pgfpathlineto{\pgfqpoint{1.957668in}{1.873550in}}%
\pgfpathlineto{\pgfqpoint{1.917075in}{1.888665in}}%
\pgfpathlineto{\pgfqpoint{1.875035in}{1.901970in}}%
\pgfpathlineto{\pgfqpoint{1.831736in}{1.913407in}}%
\pgfpathlineto{\pgfqpoint{1.787369in}{1.922924in}}%
\pgfpathlineto{\pgfqpoint{1.739450in}{1.930861in}}%
\pgfpathlineto{\pgfqpoint{1.690795in}{1.936559in}}%
\pgfpathlineto{\pgfqpoint{1.641648in}{1.939989in}}%
\pgfpathlineto{\pgfqpoint{1.592255in}{1.941134in}}%
\pgfpathlineto{\pgfqpoint{1.542861in}{1.939989in}}%
\pgfpathlineto{\pgfqpoint{1.493714in}{1.936559in}}%
\pgfpathlineto{\pgfqpoint{1.445060in}{1.930861in}}%
\pgfpathlineto{\pgfqpoint{1.397140in}{1.922924in}}%
\pgfpathlineto{\pgfqpoint{1.350195in}{1.912787in}}%
\pgfpathlineto{\pgfqpoint{1.306964in}{1.901239in}}%
\pgfpathlineto{\pgfqpoint{1.265005in}{1.887825in}}%
\pgfpathlineto{\pgfqpoint{1.224503in}{1.872606in}}%
\pgfpathlineto{\pgfqpoint{1.185639in}{1.855649in}}%
\pgfpathlineto{\pgfqpoint{1.148587in}{1.837029in}}%
\pgfpathlineto{\pgfqpoint{1.113512in}{1.816830in}}%
\pgfpathlineto{\pgfqpoint{1.082445in}{1.796457in}}%
\pgfpathlineto{\pgfqpoint{1.053391in}{1.774844in}}%
\pgfpathlineto{\pgfqpoint{1.026463in}{1.752078in}}%
\pgfpathlineto{\pgfqpoint{1.001768in}{1.728249in}}%
\pgfpathlineto{\pgfqpoint{0.979403in}{1.703449in}}%
\pgfpathlineto{\pgfqpoint{0.960632in}{1.679405in}}%
\pgfpathlineto{\pgfqpoint{0.944051in}{1.654678in}}%
\pgfpathlineto{\pgfqpoint{0.929719in}{1.629355in}}%
\pgfpathlineto{\pgfqpoint{0.917685in}{1.603521in}}%
\pgfpathlineto{\pgfqpoint{0.907991in}{1.577268in}}%
\pgfpathlineto{\pgfqpoint{0.900670in}{1.550686in}}%
\pgfpathlineto{\pgfqpoint{0.896001in}{1.525661in}}%
\pgfpathlineto{\pgfqpoint{0.893437in}{1.500506in}}%
\pgfpathlineto{\pgfqpoint{0.892984in}{1.475296in}}%
\pgfpathlineto{\pgfqpoint{0.894644in}{1.450108in}}%
\pgfpathlineto{\pgfqpoint{0.898412in}{1.425018in}}%
\pgfpathlineto{\pgfqpoint{0.904277in}{1.400101in}}%
\pgfpathlineto{\pgfqpoint{0.912867in}{1.373683in}}%
\pgfpathlineto{\pgfqpoint{0.923814in}{1.347643in}}%
\pgfpathlineto{\pgfqpoint{0.937080in}{1.322070in}}%
\pgfpathlineto{\pgfqpoint{0.952618in}{1.297054in}}%
\pgfpathlineto{\pgfqpoint{0.970375in}{1.272681in}}%
\pgfpathlineto{\pgfqpoint{0.990289in}{1.249036in}}%
\pgfpathlineto{\pgfqpoint{1.012291in}{1.226201in}}%
\pgfpathlineto{\pgfqpoint{1.037975in}{1.202825in}}%
\pgfpathlineto{\pgfqpoint{1.065847in}{1.180553in}}%
\pgfpathlineto{\pgfqpoint{1.095796in}{1.159473in}}%
\pgfpathlineto{\pgfqpoint{1.127705in}{1.139667in}}%
\pgfpathlineto{\pgfqpoint{1.163614in}{1.120108in}}%
\pgfpathlineto{\pgfqpoint{1.201432in}{1.102163in}}%
\pgfpathlineto{\pgfqpoint{1.240992in}{1.085913in}}%
\pgfpathlineto{\pgfqpoint{1.282116in}{1.071429in}}%
\pgfpathlineto{\pgfqpoint{1.324622in}{1.058777in}}%
\pgfpathlineto{\pgfqpoint{1.368321in}{1.048013in}}%
\pgfpathlineto{\pgfqpoint{1.413017in}{1.039184in}}%
\pgfpathlineto{\pgfqpoint{1.461208in}{1.031989in}}%
\pgfpathlineto{\pgfqpoint{1.510054in}{1.027045in}}%
\pgfpathlineto{\pgfqpoint{1.559310in}{1.024375in}}%
\pgfpathlineto{\pgfqpoint{1.608731in}{1.023993in}}%
\pgfpathlineto{\pgfqpoint{1.658070in}{1.025901in}}%
\pgfpathlineto{\pgfqpoint{1.707080in}{1.030090in}}%
\pgfpathlineto{\pgfqpoint{1.755516in}{1.036538in}}%
\pgfpathlineto{\pgfqpoint{1.803137in}{1.045213in}}%
\pgfpathlineto{\pgfqpoint{1.847150in}{1.055413in}}%
\pgfpathlineto{\pgfqpoint{1.890026in}{1.067515in}}%
\pgfpathlineto{\pgfqpoint{1.931577in}{1.081466in}}%
\pgfpathlineto{\pgfqpoint{1.971615in}{1.097204in}}%
\pgfpathlineto{\pgfqpoint{2.009963in}{1.114659in}}%
\pgfpathlineto{\pgfqpoint{2.046451in}{1.133752in}}%
\pgfpathlineto{\pgfqpoint{2.080915in}{1.154399in}}%
\pgfpathlineto{\pgfqpoint{2.111365in}{1.175168in}}%
\pgfpathlineto{\pgfqpoint{2.139767in}{1.197150in}}%
\pgfpathlineto{\pgfqpoint{2.166008in}{1.220258in}}%
\pgfpathlineto{\pgfqpoint{2.189985in}{1.244402in}}%
\pgfpathlineto{\pgfqpoint{2.210322in}{1.267891in}}%
\pgfpathlineto{\pgfqpoint{2.228516in}{1.292125in}}%
\pgfpathlineto{\pgfqpoint{2.244502in}{1.317020in}}%
\pgfpathlineto{\pgfqpoint{2.258226in}{1.342488in}}%
\pgfpathlineto{\pgfqpoint{2.269640in}{1.368442in}}%
\pgfpathlineto{\pgfqpoint{2.278704in}{1.394792in}}%
\pgfpathlineto{\pgfqpoint{2.285387in}{1.421446in}}%
\pgfpathlineto{\pgfqpoint{2.289456in}{1.446516in}}%
\pgfpathlineto{\pgfqpoint{2.291418in}{1.471695in}}%
\pgfpathlineto{\pgfqpoint{2.291612in}{1.482500in}}%
\pgfpathlineto{\pgfqpoint{2.291612in}{1.482500in}}%
\pgfusepath{stroke}%
\end{pgfscope}%
\begin{pgfscope}%
\pgfpathrectangle{\pgfqpoint{0.135000in}{0.135000in}}{\pgfqpoint{2.914509in}{2.695000in}} %
\pgfusepath{clip}%
\pgfsetrectcap%
\pgfsetroundjoin%
\pgfsetlinewidth{1.003750pt}%
\definecolor{currentstroke}{rgb}{0.700000,0.700000,0.700000}%
\pgfsetstrokecolor{currentstroke}%
\pgfsetdash{}{0pt}%
\pgfpathmoveto{\pgfqpoint{2.531003in}{1.482500in}}%
\pgfpathlineto{\pgfqpoint{2.529961in}{1.519064in}}%
\pgfpathlineto{\pgfqpoint{2.526837in}{1.555547in}}%
\pgfpathlineto{\pgfqpoint{2.521638in}{1.591868in}}%
\pgfpathlineto{\pgfqpoint{2.514375in}{1.627946in}}%
\pgfpathlineto{\pgfqpoint{2.505066in}{1.663701in}}%
\pgfpathlineto{\pgfqpoint{2.493729in}{1.699054in}}%
\pgfpathlineto{\pgfqpoint{2.479190in}{1.736808in}}%
\pgfpathlineto{\pgfqpoint{2.462340in}{1.773899in}}%
\pgfpathlineto{\pgfqpoint{2.443223in}{1.810231in}}%
\pgfpathlineto{\pgfqpoint{2.421888in}{1.845709in}}%
\pgfpathlineto{\pgfqpoint{2.398392in}{1.880241in}}%
\pgfpathlineto{\pgfqpoint{2.372796in}{1.913736in}}%
\pgfpathlineto{\pgfqpoint{2.345165in}{1.946108in}}%
\pgfpathlineto{\pgfqpoint{2.315573in}{1.977272in}}%
\pgfpathlineto{\pgfqpoint{2.284096in}{2.007146in}}%
\pgfpathlineto{\pgfqpoint{2.250816in}{2.035654in}}%
\pgfpathlineto{\pgfqpoint{2.215821in}{2.062720in}}%
\pgfpathlineto{\pgfqpoint{2.179200in}{2.088274in}}%
\pgfpathlineto{\pgfqpoint{2.141051in}{2.112250in}}%
\pgfpathlineto{\pgfqpoint{2.101471in}{2.134585in}}%
\pgfpathlineto{\pgfqpoint{2.060564in}{2.155221in}}%
\pgfpathlineto{\pgfqpoint{2.018438in}{2.174104in}}%
\pgfpathlineto{\pgfqpoint{1.975200in}{2.191185in}}%
\pgfpathlineto{\pgfqpoint{1.930965in}{2.206419in}}%
\pgfpathlineto{\pgfqpoint{1.885848in}{2.219767in}}%
\pgfpathlineto{\pgfqpoint{1.836407in}{2.231993in}}%
\pgfpathlineto{\pgfqpoint{1.786229in}{2.241954in}}%
\pgfpathlineto{\pgfqpoint{1.735465in}{2.249619in}}%
\pgfpathlineto{\pgfqpoint{1.684268in}{2.254967in}}%
\pgfpathlineto{\pgfqpoint{1.632793in}{2.257981in}}%
\pgfpathlineto{\pgfqpoint{1.581196in}{2.258651in}}%
\pgfpathlineto{\pgfqpoint{1.529631in}{2.256976in}}%
\pgfpathlineto{\pgfqpoint{1.478257in}{2.252960in}}%
\pgfpathlineto{\pgfqpoint{1.427226in}{2.246617in}}%
\pgfpathlineto{\pgfqpoint{1.376694in}{2.237964in}}%
\pgfpathlineto{\pgfqpoint{1.326814in}{2.227029in}}%
\pgfpathlineto{\pgfqpoint{1.277736in}{2.213843in}}%
\pgfpathlineto{\pgfqpoint{1.229608in}{2.198448in}}%
\pgfpathlineto{\pgfqpoint{1.185896in}{2.182214in}}%
\pgfpathlineto{\pgfqpoint{1.143243in}{2.164157in}}%
\pgfpathlineto{\pgfqpoint{1.101760in}{2.144323in}}%
\pgfpathlineto{\pgfqpoint{1.061555in}{2.122766in}}%
\pgfpathlineto{\pgfqpoint{1.022732in}{2.099540in}}%
\pgfpathlineto{\pgfqpoint{0.985394in}{2.074706in}}%
\pgfpathlineto{\pgfqpoint{0.949637in}{2.048329in}}%
\pgfpathlineto{\pgfqpoint{0.915554in}{2.020478in}}%
\pgfpathlineto{\pgfqpoint{0.883235in}{1.991225in}}%
\pgfpathlineto{\pgfqpoint{0.852763in}{1.960646in}}%
\pgfpathlineto{\pgfqpoint{0.824218in}{1.928822in}}%
\pgfpathlineto{\pgfqpoint{0.797674in}{1.895835in}}%
\pgfpathlineto{\pgfqpoint{0.773200in}{1.861770in}}%
\pgfpathlineto{\pgfqpoint{0.750861in}{1.826718in}}%
\pgfpathlineto{\pgfqpoint{0.730714in}{1.790768in}}%
\pgfpathlineto{\pgfqpoint{0.712812in}{1.754015in}}%
\pgfpathlineto{\pgfqpoint{0.697201in}{1.716555in}}%
\pgfpathlineto{\pgfqpoint{0.683922in}{1.678485in}}%
\pgfpathlineto{\pgfqpoint{0.673765in}{1.642888in}}%
\pgfpathlineto{\pgfqpoint{0.665647in}{1.606935in}}%
\pgfpathlineto{\pgfqpoint{0.659587in}{1.570705in}}%
\pgfpathlineto{\pgfqpoint{0.655597in}{1.534280in}}%
\pgfpathlineto{\pgfqpoint{0.653687in}{1.497740in}}%
\pgfpathlineto{\pgfqpoint{0.653861in}{1.461166in}}%
\pgfpathlineto{\pgfqpoint{0.656118in}{1.424639in}}%
\pgfpathlineto{\pgfqpoint{0.660454in}{1.388241in}}%
\pgfpathlineto{\pgfqpoint{0.666858in}{1.352052in}}%
\pgfpathlineto{\pgfqpoint{0.675317in}{1.316152in}}%
\pgfpathlineto{\pgfqpoint{0.685812in}{1.280622in}}%
\pgfpathlineto{\pgfqpoint{0.698319in}{1.245540in}}%
\pgfpathlineto{\pgfqpoint{0.714108in}{1.208131in}}%
\pgfpathlineto{\pgfqpoint{0.732185in}{1.171437in}}%
\pgfpathlineto{\pgfqpoint{0.752502in}{1.135553in}}%
\pgfpathlineto{\pgfqpoint{0.775008in}{1.100573in}}%
\pgfpathlineto{\pgfqpoint{0.799643in}{1.066589in}}%
\pgfpathlineto{\pgfqpoint{0.826343in}{1.033688in}}%
\pgfpathlineto{\pgfqpoint{0.855040in}{1.001956in}}%
\pgfpathlineto{\pgfqpoint{0.885657in}{0.971477in}}%
\pgfpathlineto{\pgfqpoint{0.918115in}{0.942329in}}%
\pgfpathlineto{\pgfqpoint{0.952329in}{0.914589in}}%
\pgfpathlineto{\pgfqpoint{0.988211in}{0.888328in}}%
\pgfpathlineto{\pgfqpoint{1.025667in}{0.863616in}}%
\pgfpathlineto{\pgfqpoint{1.064600in}{0.840516in}}%
\pgfpathlineto{\pgfqpoint{1.104907in}{0.819089in}}%
\pgfpathlineto{\pgfqpoint{1.146484in}{0.799391in}}%
\pgfpathlineto{\pgfqpoint{1.189222in}{0.781472in}}%
\pgfpathlineto{\pgfqpoint{1.233011in}{0.765380in}}%
\pgfpathlineto{\pgfqpoint{1.277736in}{0.751157in}}%
\pgfpathlineto{\pgfqpoint{1.326814in}{0.737971in}}%
\pgfpathlineto{\pgfqpoint{1.376694in}{0.727036in}}%
\pgfpathlineto{\pgfqpoint{1.427226in}{0.718383in}}%
\pgfpathlineto{\pgfqpoint{1.478257in}{0.712040in}}%
\pgfpathlineto{\pgfqpoint{1.529631in}{0.708024in}}%
\pgfpathlineto{\pgfqpoint{1.581196in}{0.706349in}}%
\pgfpathlineto{\pgfqpoint{1.632793in}{0.707019in}}%
\pgfpathlineto{\pgfqpoint{1.684268in}{0.710033in}}%
\pgfpathlineto{\pgfqpoint{1.735465in}{0.715381in}}%
\pgfpathlineto{\pgfqpoint{1.786229in}{0.723046in}}%
\pgfpathlineto{\pgfqpoint{1.836407in}{0.733007in}}%
\pgfpathlineto{\pgfqpoint{1.885848in}{0.745233in}}%
\pgfpathlineto{\pgfqpoint{1.934401in}{0.759686in}}%
\pgfpathlineto{\pgfqpoint{1.978563in}{0.775064in}}%
\pgfpathlineto{\pgfqpoint{2.021719in}{0.792286in}}%
\pgfpathlineto{\pgfqpoint{2.063756in}{0.811305in}}%
\pgfpathlineto{\pgfqpoint{2.104564in}{0.832074in}}%
\pgfpathlineto{\pgfqpoint{2.144037in}{0.854537in}}%
\pgfpathlineto{\pgfqpoint{2.182073in}{0.878637in}}%
\pgfpathlineto{\pgfqpoint{2.218572in}{0.904310in}}%
\pgfpathlineto{\pgfqpoint{2.253438in}{0.931489in}}%
\pgfpathlineto{\pgfqpoint{2.286582in}{0.960104in}}%
\pgfpathlineto{\pgfqpoint{2.317917in}{0.990081in}}%
\pgfpathlineto{\pgfqpoint{2.347361in}{1.021340in}}%
\pgfpathlineto{\pgfqpoint{2.374838in}{1.053801in}}%
\pgfpathlineto{\pgfqpoint{2.400275in}{1.087380in}}%
\pgfpathlineto{\pgfqpoint{2.423607in}{1.121987in}}%
\pgfpathlineto{\pgfqpoint{2.444773in}{1.157534in}}%
\pgfpathlineto{\pgfqpoint{2.463717in}{1.193928in}}%
\pgfpathlineto{\pgfqpoint{2.480391in}{1.231074in}}%
\pgfpathlineto{\pgfqpoint{2.494751in}{1.268875in}}%
\pgfpathlineto{\pgfqpoint{2.506759in}{1.307232in}}%
\pgfpathlineto{\pgfqpoint{2.515729in}{1.343046in}}%
\pgfpathlineto{\pgfqpoint{2.522648in}{1.379171in}}%
\pgfpathlineto{\pgfqpoint{2.527502in}{1.415524in}}%
\pgfpathlineto{\pgfqpoint{2.530280in}{1.452026in}}%
\pgfpathlineto{\pgfqpoint{2.531003in}{1.482500in}}%
\pgfpathlineto{\pgfqpoint{2.531003in}{1.482500in}}%
\pgfusepath{stroke}%
\end{pgfscope}%
\begin{pgfscope}%
\pgfpathrectangle{\pgfqpoint{0.135000in}{0.135000in}}{\pgfqpoint{2.914509in}{2.695000in}} %
\pgfusepath{clip}%
\pgfsetrectcap%
\pgfsetroundjoin%
\pgfsetlinewidth{1.003750pt}%
\definecolor{currentstroke}{rgb}{0.700000,0.700000,0.700000}%
\pgfsetstrokecolor{currentstroke}%
\pgfsetdash{}{0pt}%
\pgfpathmoveto{\pgfqpoint{2.917032in}{1.482500in}}%
\pgfpathlineto{\pgfqpoint{2.916010in}{1.530202in}}%
\pgfpathlineto{\pgfqpoint{2.912948in}{1.577830in}}%
\pgfpathlineto{\pgfqpoint{2.907849in}{1.625311in}}%
\pgfpathlineto{\pgfqpoint{2.900721in}{1.672571in}}%
\pgfpathlineto{\pgfqpoint{2.891576in}{1.719539in}}%
\pgfpathlineto{\pgfqpoint{2.880428in}{1.766141in}}%
\pgfpathlineto{\pgfqpoint{2.867293in}{1.812306in}}%
\pgfpathlineto{\pgfqpoint{2.852192in}{1.857963in}}%
\pgfpathlineto{\pgfqpoint{2.835149in}{1.903040in}}%
\pgfpathlineto{\pgfqpoint{2.814189in}{1.951873in}}%
\pgfpathlineto{\pgfqpoint{2.790949in}{1.999831in}}%
\pgfpathlineto{\pgfqpoint{2.765472in}{2.046824in}}%
\pgfpathlineto{\pgfqpoint{2.737807in}{2.092764in}}%
\pgfpathlineto{\pgfqpoint{2.708005in}{2.137565in}}%
\pgfpathlineto{\pgfqpoint{2.676121in}{2.181144in}}%
\pgfpathlineto{\pgfqpoint{2.642214in}{2.223420in}}%
\pgfpathlineto{\pgfqpoint{2.606349in}{2.264313in}}%
\pgfpathlineto{\pgfqpoint{2.568592in}{2.303748in}}%
\pgfpathlineto{\pgfqpoint{2.529013in}{2.341651in}}%
\pgfpathlineto{\pgfqpoint{2.487687in}{2.377950in}}%
\pgfpathlineto{\pgfqpoint{2.444690in}{2.412579in}}%
\pgfpathlineto{\pgfqpoint{2.400103in}{2.445473in}}%
\pgfpathlineto{\pgfqpoint{2.354008in}{2.476570in}}%
\pgfpathlineto{\pgfqpoint{2.306493in}{2.505813in}}%
\pgfpathlineto{\pgfqpoint{2.257644in}{2.533146in}}%
\pgfpathlineto{\pgfqpoint{2.207555in}{2.558519in}}%
\pgfpathlineto{\pgfqpoint{2.156317in}{2.581885in}}%
\pgfpathlineto{\pgfqpoint{2.104027in}{2.603200in}}%
\pgfpathlineto{\pgfqpoint{2.050783in}{2.622423in}}%
\pgfpathlineto{\pgfqpoint{1.996682in}{2.639520in}}%
\pgfpathlineto{\pgfqpoint{1.941828in}{2.654459in}}%
\pgfpathlineto{\pgfqpoint{1.886321in}{2.667211in}}%
\pgfpathlineto{\pgfqpoint{1.830265in}{2.677752in}}%
\pgfpathlineto{\pgfqpoint{1.773765in}{2.686064in}}%
\pgfpathlineto{\pgfqpoint{1.716927in}{2.692130in}}%
\pgfpathlineto{\pgfqpoint{1.659856in}{2.695939in}}%
\pgfpathlineto{\pgfqpoint{1.602659in}{2.697485in}}%
\pgfpathlineto{\pgfqpoint{1.545443in}{2.696764in}}%
\pgfpathlineto{\pgfqpoint{1.488314in}{2.693777in}}%
\pgfpathlineto{\pgfqpoint{1.431379in}{2.688530in}}%
\pgfpathlineto{\pgfqpoint{1.374744in}{2.681034in}}%
\pgfpathlineto{\pgfqpoint{1.318515in}{2.671301in}}%
\pgfpathlineto{\pgfqpoint{1.262796in}{2.659350in}}%
\pgfpathlineto{\pgfqpoint{1.207692in}{2.645204in}}%
\pgfpathlineto{\pgfqpoint{1.153306in}{2.628888in}}%
\pgfpathlineto{\pgfqpoint{1.099738in}{2.610434in}}%
\pgfpathlineto{\pgfqpoint{1.047090in}{2.589875in}}%
\pgfpathlineto{\pgfqpoint{0.995459in}{2.567251in}}%
\pgfpathlineto{\pgfqpoint{0.944941in}{2.542602in}}%
\pgfpathlineto{\pgfqpoint{0.895630in}{2.515976in}}%
\pgfpathlineto{\pgfqpoint{0.847620in}{2.487421in}}%
\pgfpathlineto{\pgfqpoint{0.800998in}{2.456992in}}%
\pgfpathlineto{\pgfqpoint{0.755853in}{2.424745in}}%
\pgfpathlineto{\pgfqpoint{0.712268in}{2.390739in}}%
\pgfpathlineto{\pgfqpoint{0.670325in}{2.355039in}}%
\pgfpathlineto{\pgfqpoint{0.630103in}{2.317712in}}%
\pgfpathlineto{\pgfqpoint{0.591675in}{2.278826in}}%
\pgfpathlineto{\pgfqpoint{0.555114in}{2.238454in}}%
\pgfpathlineto{\pgfqpoint{0.520488in}{2.196672in}}%
\pgfpathlineto{\pgfqpoint{0.487861in}{2.153558in}}%
\pgfpathlineto{\pgfqpoint{0.457295in}{2.109191in}}%
\pgfpathlineto{\pgfqpoint{0.428846in}{2.063656in}}%
\pgfpathlineto{\pgfqpoint{0.402568in}{2.017036in}}%
\pgfpathlineto{\pgfqpoint{0.378510in}{1.969419in}}%
\pgfpathlineto{\pgfqpoint{0.356716in}{1.920893in}}%
\pgfpathlineto{\pgfqpoint{0.337227in}{1.871550in}}%
\pgfpathlineto{\pgfqpoint{0.321541in}{1.826060in}}%
\pgfpathlineto{\pgfqpoint{0.307814in}{1.780040in}}%
\pgfpathlineto{\pgfqpoint{0.296068in}{1.733561in}}%
\pgfpathlineto{\pgfqpoint{0.286320in}{1.686696in}}%
\pgfpathlineto{\pgfqpoint{0.278586in}{1.639515in}}%
\pgfpathlineto{\pgfqpoint{0.272878in}{1.592093in}}%
\pgfpathlineto{\pgfqpoint{0.269204in}{1.544501in}}%
\pgfpathlineto{\pgfqpoint{0.267570in}{1.496814in}}%
\pgfpathlineto{\pgfqpoint{0.267978in}{1.449105in}}%
\pgfpathlineto{\pgfqpoint{0.270429in}{1.401447in}}%
\pgfpathlineto{\pgfqpoint{0.274917in}{1.353914in}}%
\pgfpathlineto{\pgfqpoint{0.281437in}{1.306579in}}%
\pgfpathlineto{\pgfqpoint{0.289978in}{1.259516in}}%
\pgfpathlineto{\pgfqpoint{0.300527in}{1.212797in}}%
\pgfpathlineto{\pgfqpoint{0.313068in}{1.166493in}}%
\pgfpathlineto{\pgfqpoint{0.327582in}{1.120677in}}%
\pgfpathlineto{\pgfqpoint{0.344045in}{1.075418in}}%
\pgfpathlineto{\pgfqpoint{0.362433in}{1.030788in}}%
\pgfpathlineto{\pgfqpoint{0.384849in}{0.982501in}}%
\pgfpathlineto{\pgfqpoint{0.409517in}{0.935147in}}%
\pgfpathlineto{\pgfqpoint{0.436392in}{0.888814in}}%
\pgfpathlineto{\pgfqpoint{0.465424in}{0.843589in}}%
\pgfpathlineto{\pgfqpoint{0.496557in}{0.799556in}}%
\pgfpathlineto{\pgfqpoint{0.529735in}{0.756797in}}%
\pgfpathlineto{\pgfqpoint{0.564896in}{0.715392in}}%
\pgfpathlineto{\pgfqpoint{0.601973in}{0.675418in}}%
\pgfpathlineto{\pgfqpoint{0.640897in}{0.636950in}}%
\pgfpathlineto{\pgfqpoint{0.681597in}{0.600060in}}%
\pgfpathlineto{\pgfqpoint{0.723996in}{0.564816in}}%
\pgfpathlineto{\pgfqpoint{0.768014in}{0.531284in}}%
\pgfpathlineto{\pgfqpoint{0.813570in}{0.499526in}}%
\pgfpathlineto{\pgfqpoint{0.860579in}{0.469603in}}%
\pgfpathlineto{\pgfqpoint{0.908953in}{0.441569in}}%
\pgfpathlineto{\pgfqpoint{0.958602in}{0.415477in}}%
\pgfpathlineto{\pgfqpoint{1.009433in}{0.391376in}}%
\pgfpathlineto{\pgfqpoint{1.061352in}{0.369311in}}%
\pgfpathlineto{\pgfqpoint{1.114261in}{0.349323in}}%
\pgfpathlineto{\pgfqpoint{1.168062in}{0.331449in}}%
\pgfpathlineto{\pgfqpoint{1.222654in}{0.315722in}}%
\pgfpathlineto{\pgfqpoint{1.277935in}{0.302172in}}%
\pgfpathlineto{\pgfqpoint{1.333804in}{0.290824in}}%
\pgfpathlineto{\pgfqpoint{1.390154in}{0.281699in}}%
\pgfpathlineto{\pgfqpoint{1.446881in}{0.274815in}}%
\pgfpathlineto{\pgfqpoint{1.503880in}{0.270184in}}%
\pgfpathlineto{\pgfqpoint{1.561043in}{0.267815in}}%
\pgfpathlineto{\pgfqpoint{1.618265in}{0.267712in}}%
\pgfpathlineto{\pgfqpoint{1.675438in}{0.269875in}}%
\pgfpathlineto{\pgfqpoint{1.732456in}{0.274301in}}%
\pgfpathlineto{\pgfqpoint{1.789212in}{0.280981in}}%
\pgfpathlineto{\pgfqpoint{1.845601in}{0.289902in}}%
\pgfpathlineto{\pgfqpoint{1.901518in}{0.301049in}}%
\pgfpathlineto{\pgfqpoint{1.956857in}{0.314400in}}%
\pgfpathlineto{\pgfqpoint{2.011516in}{0.329930in}}%
\pgfpathlineto{\pgfqpoint{2.065393in}{0.347610in}}%
\pgfpathlineto{\pgfqpoint{2.118387in}{0.367408in}}%
\pgfpathlineto{\pgfqpoint{2.170400in}{0.389286in}}%
\pgfpathlineto{\pgfqpoint{2.221333in}{0.413203in}}%
\pgfpathlineto{\pgfqpoint{2.271094in}{0.439116in}}%
\pgfpathlineto{\pgfqpoint{2.319587in}{0.466975in}}%
\pgfpathlineto{\pgfqpoint{2.366724in}{0.496729in}}%
\pgfpathlineto{\pgfqpoint{2.412416in}{0.528322in}}%
\pgfpathlineto{\pgfqpoint{2.456578in}{0.561696in}}%
\pgfpathlineto{\pgfqpoint{2.499127in}{0.596787in}}%
\pgfpathlineto{\pgfqpoint{2.539984in}{0.633531in}}%
\pgfpathlineto{\pgfqpoint{2.579073in}{0.671858in}}%
\pgfpathlineto{\pgfqpoint{2.616321in}{0.711698in}}%
\pgfpathlineto{\pgfqpoint{2.651659in}{0.752976in}}%
\pgfpathlineto{\pgfqpoint{2.685019in}{0.795615in}}%
\pgfpathlineto{\pgfqpoint{2.716341in}{0.839536in}}%
\pgfpathlineto{\pgfqpoint{2.745566in}{0.884656in}}%
\pgfpathlineto{\pgfqpoint{2.772640in}{0.930891in}}%
\pgfpathlineto{\pgfqpoint{2.797511in}{0.978156in}}%
\pgfpathlineto{\pgfqpoint{2.820133in}{1.026362in}}%
\pgfpathlineto{\pgfqpoint{2.840464in}{1.075418in}}%
\pgfpathlineto{\pgfqpoint{2.856928in}{1.120677in}}%
\pgfpathlineto{\pgfqpoint{2.871441in}{1.166493in}}%
\pgfpathlineto{\pgfqpoint{2.883982in}{1.212797in}}%
\pgfpathlineto{\pgfqpoint{2.894531in}{1.259516in}}%
\pgfpathlineto{\pgfqpoint{2.903072in}{1.306579in}}%
\pgfpathlineto{\pgfqpoint{2.909592in}{1.353914in}}%
\pgfpathlineto{\pgfqpoint{2.914081in}{1.401447in}}%
\pgfpathlineto{\pgfqpoint{2.916531in}{1.449105in}}%
\pgfpathlineto{\pgfqpoint{2.917032in}{1.482500in}}%
\pgfpathlineto{\pgfqpoint{2.917032in}{1.482500in}}%
\pgfusepath{stroke}%
\end{pgfscope}%
\begin{pgfscope}%
\pgfpathrectangle{\pgfqpoint{0.135000in}{0.135000in}}{\pgfqpoint{2.914509in}{2.695000in}} %
\pgfusepath{clip}%
\pgfsetrectcap%
\pgfsetroundjoin%
\pgfsetlinewidth{1.003750pt}%
\definecolor{currentstroke}{rgb}{0.700000,0.700000,0.700000}%
\pgfsetstrokecolor{currentstroke}%
\pgfsetdash{}{0pt}%
\pgfpathmoveto{\pgfqpoint{2.120227in}{1.482500in}}%
\pgfpathlineto{\pgfqpoint{2.916065in}{1.482500in}}%
\pgfpathlineto{\pgfqpoint{2.916065in}{1.482500in}}%
\pgfusepath{stroke}%
\end{pgfscope}%
\begin{pgfscope}%
\pgfpathrectangle{\pgfqpoint{0.135000in}{0.135000in}}{\pgfqpoint{2.914509in}{2.695000in}} %
\pgfusepath{clip}%
\pgfsetrectcap%
\pgfsetroundjoin%
\pgfsetlinewidth{1.003750pt}%
\definecolor{currentstroke}{rgb}{0.700000,0.700000,0.700000}%
\pgfsetstrokecolor{currentstroke}%
\pgfsetdash{}{0pt}%
\pgfpathmoveto{\pgfqpoint{2.080038in}{1.482500in}}%
\pgfpathlineto{\pgfqpoint{2.081125in}{1.495998in}}%
\pgfpathlineto{\pgfqpoint{2.084392in}{1.509557in}}%
\pgfpathlineto{\pgfqpoint{2.090243in}{1.524046in}}%
\pgfpathlineto{\pgfqpoint{2.098583in}{1.538742in}}%
\pgfpathlineto{\pgfqpoint{2.110134in}{1.554562in}}%
\pgfpathlineto{\pgfqpoint{2.125410in}{1.571650in}}%
\pgfpathlineto{\pgfqpoint{2.144993in}{1.590185in}}%
\pgfpathlineto{\pgfqpoint{2.170762in}{1.611332in}}%
\pgfpathlineto{\pgfqpoint{2.202691in}{1.634522in}}%
\pgfpathlineto{\pgfqpoint{2.243413in}{1.661177in}}%
\pgfpathlineto{\pgfqpoint{2.292900in}{1.690834in}}%
\pgfpathlineto{\pgfqpoint{2.354851in}{1.725308in}}%
\pgfpathlineto{\pgfqpoint{2.432046in}{1.765659in}}%
\pgfpathlineto{\pgfqpoint{2.531087in}{1.814769in}}%
\pgfpathlineto{\pgfqpoint{2.654266in}{1.873254in}}%
\pgfpathlineto{\pgfqpoint{2.811790in}{1.945482in}}%
\pgfpathlineto{\pgfqpoint{2.815295in}{1.947065in}}%
\pgfpathlineto{\pgfqpoint{2.815295in}{1.947065in}}%
\pgfusepath{stroke}%
\end{pgfscope}%
\begin{pgfscope}%
\pgfpathrectangle{\pgfqpoint{0.135000in}{0.135000in}}{\pgfqpoint{2.914509in}{2.695000in}} %
\pgfusepath{clip}%
\pgfsetrectcap%
\pgfsetroundjoin%
\pgfsetlinewidth{1.003750pt}%
\definecolor{currentstroke}{rgb}{0.700000,0.700000,0.700000}%
\pgfsetstrokecolor{currentstroke}%
\pgfsetdash{}{0pt}%
\pgfpathmoveto{\pgfqpoint{1.965588in}{1.482500in}}%
\pgfpathlineto{\pgfqpoint{1.966740in}{1.511852in}}%
\pgfpathlineto{\pgfqpoint{1.970203in}{1.541384in}}%
\pgfpathlineto{\pgfqpoint{1.975999in}{1.571281in}}%
\pgfpathlineto{\pgfqpoint{1.984163in}{1.601725in}}%
\pgfpathlineto{\pgfqpoint{1.995339in}{1.634486in}}%
\pgfpathlineto{\pgfqpoint{2.009259in}{1.668282in}}%
\pgfpathlineto{\pgfqpoint{2.026887in}{1.705048in}}%
\pgfpathlineto{\pgfqpoint{2.047760in}{1.743475in}}%
\pgfpathlineto{\pgfqpoint{2.073223in}{1.785738in}}%
\pgfpathlineto{\pgfqpoint{2.103982in}{1.832483in}}%
\pgfpathlineto{\pgfqpoint{2.140865in}{1.884490in}}%
\pgfpathlineto{\pgfqpoint{2.186654in}{1.945030in}}%
\pgfpathlineto{\pgfqpoint{2.241212in}{2.013318in}}%
\pgfpathlineto{\pgfqpoint{2.306006in}{2.090829in}}%
\pgfpathlineto{\pgfqpoint{2.385540in}{2.182446in}}%
\pgfpathlineto{\pgfqpoint{2.483171in}{2.291422in}}%
\pgfpathlineto{\pgfqpoint{2.528330in}{2.340905in}}%
\pgfpathlineto{\pgfqpoint{2.528330in}{2.340905in}}%
\pgfusepath{stroke}%
\end{pgfscope}%
\begin{pgfscope}%
\pgfpathrectangle{\pgfqpoint{0.135000in}{0.135000in}}{\pgfqpoint{2.914509in}{2.695000in}} %
\pgfusepath{clip}%
\pgfsetrectcap%
\pgfsetroundjoin%
\pgfsetlinewidth{1.003750pt}%
\definecolor{currentstroke}{rgb}{0.700000,0.700000,0.700000}%
\pgfsetstrokecolor{currentstroke}%
\pgfsetdash{}{0pt}%
\pgfpathmoveto{\pgfqpoint{1.794301in}{1.482500in}}%
\pgfpathlineto{\pgfqpoint{1.795438in}{1.534316in}}%
\pgfpathlineto{\pgfqpoint{1.798861in}{1.586715in}}%
\pgfpathlineto{\pgfqpoint{1.804609in}{1.640287in}}%
\pgfpathlineto{\pgfqpoint{1.812746in}{1.695634in}}%
\pgfpathlineto{\pgfqpoint{1.823807in}{1.755573in}}%
\pgfpathlineto{\pgfqpoint{1.838220in}{1.821146in}}%
\pgfpathlineto{\pgfqpoint{1.856492in}{1.893617in}}%
\pgfpathlineto{\pgfqpoint{1.879234in}{1.974518in}}%
\pgfpathlineto{\pgfqpoint{1.908132in}{2.068691in}}%
\pgfpathlineto{\pgfqpoint{1.943468in}{2.176047in}}%
\pgfpathlineto{\pgfqpoint{1.987731in}{2.303258in}}%
\pgfpathlineto{\pgfqpoint{2.043063in}{2.455418in}}%
\pgfpathlineto{\pgfqpoint{2.098855in}{2.604060in}}%
\pgfpathlineto{\pgfqpoint{2.098855in}{2.604060in}}%
\pgfusepath{stroke}%
\end{pgfscope}%
\begin{pgfscope}%
\pgfpathrectangle{\pgfqpoint{0.135000in}{0.135000in}}{\pgfqpoint{2.914509in}{2.695000in}} %
\pgfusepath{clip}%
\pgfsetrectcap%
\pgfsetroundjoin%
\pgfsetlinewidth{1.003750pt}%
\definecolor{currentstroke}{rgb}{0.700000,0.700000,0.700000}%
\pgfsetstrokecolor{currentstroke}%
\pgfsetdash{}{0pt}%
\pgfpathmoveto{\pgfqpoint{1.592255in}{1.482500in}}%
\pgfpathlineto{\pgfqpoint{1.592255in}{2.696468in}}%
\pgfpathlineto{\pgfqpoint{1.592255in}{2.696468in}}%
\pgfusepath{stroke}%
\end{pgfscope}%
\begin{pgfscope}%
\pgfpathrectangle{\pgfqpoint{0.135000in}{0.135000in}}{\pgfqpoint{2.914509in}{2.695000in}} %
\pgfusepath{clip}%
\pgfsetrectcap%
\pgfsetroundjoin%
\pgfsetlinewidth{1.003750pt}%
\definecolor{currentstroke}{rgb}{0.700000,0.700000,0.700000}%
\pgfsetstrokecolor{currentstroke}%
\pgfsetdash{}{0pt}%
\pgfpathmoveto{\pgfqpoint{1.390208in}{1.482500in}}%
\pgfpathlineto{\pgfqpoint{1.389072in}{1.534316in}}%
\pgfpathlineto{\pgfqpoint{1.385648in}{1.586715in}}%
\pgfpathlineto{\pgfqpoint{1.379901in}{1.640287in}}%
\pgfpathlineto{\pgfqpoint{1.371763in}{1.695634in}}%
\pgfpathlineto{\pgfqpoint{1.360702in}{1.755573in}}%
\pgfpathlineto{\pgfqpoint{1.346290in}{1.821146in}}%
\pgfpathlineto{\pgfqpoint{1.328017in}{1.893617in}}%
\pgfpathlineto{\pgfqpoint{1.305275in}{1.974518in}}%
\pgfpathlineto{\pgfqpoint{1.276377in}{2.068691in}}%
\pgfpathlineto{\pgfqpoint{1.241042in}{2.176047in}}%
\pgfpathlineto{\pgfqpoint{1.196778in}{2.303258in}}%
\pgfpathlineto{\pgfqpoint{1.141446in}{2.455418in}}%
\pgfpathlineto{\pgfqpoint{1.085655in}{2.604060in}}%
\pgfpathlineto{\pgfqpoint{1.085655in}{2.604060in}}%
\pgfusepath{stroke}%
\end{pgfscope}%
\begin{pgfscope}%
\pgfpathrectangle{\pgfqpoint{0.135000in}{0.135000in}}{\pgfqpoint{2.914509in}{2.695000in}} %
\pgfusepath{clip}%
\pgfsetrectcap%
\pgfsetroundjoin%
\pgfsetlinewidth{1.003750pt}%
\definecolor{currentstroke}{rgb}{0.700000,0.700000,0.700000}%
\pgfsetstrokecolor{currentstroke}%
\pgfsetdash{}{0pt}%
\pgfpathmoveto{\pgfqpoint{1.218922in}{1.482500in}}%
\pgfpathlineto{\pgfqpoint{1.217770in}{1.511852in}}%
\pgfpathlineto{\pgfqpoint{1.214307in}{1.541384in}}%
\pgfpathlineto{\pgfqpoint{1.208511in}{1.571281in}}%
\pgfpathlineto{\pgfqpoint{1.200347in}{1.601725in}}%
\pgfpathlineto{\pgfqpoint{1.189170in}{1.634486in}}%
\pgfpathlineto{\pgfqpoint{1.175250in}{1.668282in}}%
\pgfpathlineto{\pgfqpoint{1.157623in}{1.705048in}}%
\pgfpathlineto{\pgfqpoint{1.136749in}{1.743475in}}%
\pgfpathlineto{\pgfqpoint{1.111286in}{1.785738in}}%
\pgfpathlineto{\pgfqpoint{1.080527in}{1.832483in}}%
\pgfpathlineto{\pgfqpoint{1.043644in}{1.884490in}}%
\pgfpathlineto{\pgfqpoint{0.997855in}{1.945030in}}%
\pgfpathlineto{\pgfqpoint{0.943298in}{2.013318in}}%
\pgfpathlineto{\pgfqpoint{0.878503in}{2.090829in}}%
\pgfpathlineto{\pgfqpoint{0.798969in}{2.182446in}}%
\pgfpathlineto{\pgfqpoint{0.701338in}{2.291422in}}%
\pgfpathlineto{\pgfqpoint{0.656180in}{2.340905in}}%
\pgfpathlineto{\pgfqpoint{0.656180in}{2.340905in}}%
\pgfusepath{stroke}%
\end{pgfscope}%
\begin{pgfscope}%
\pgfpathrectangle{\pgfqpoint{0.135000in}{0.135000in}}{\pgfqpoint{2.914509in}{2.695000in}} %
\pgfusepath{clip}%
\pgfsetrectcap%
\pgfsetroundjoin%
\pgfsetlinewidth{1.003750pt}%
\definecolor{currentstroke}{rgb}{0.700000,0.700000,0.700000}%
\pgfsetstrokecolor{currentstroke}%
\pgfsetdash{}{0pt}%
\pgfpathmoveto{\pgfqpoint{1.104472in}{1.482500in}}%
\pgfpathlineto{\pgfqpoint{1.103384in}{1.495998in}}%
\pgfpathlineto{\pgfqpoint{1.100117in}{1.509557in}}%
\pgfpathlineto{\pgfqpoint{1.094266in}{1.524046in}}%
\pgfpathlineto{\pgfqpoint{1.085926in}{1.538742in}}%
\pgfpathlineto{\pgfqpoint{1.074376in}{1.554562in}}%
\pgfpathlineto{\pgfqpoint{1.059100in}{1.571650in}}%
\pgfpathlineto{\pgfqpoint{1.039516in}{1.590185in}}%
\pgfpathlineto{\pgfqpoint{1.013747in}{1.611332in}}%
\pgfpathlineto{\pgfqpoint{0.981819in}{1.634522in}}%
\pgfpathlineto{\pgfqpoint{0.941096in}{1.661177in}}%
\pgfpathlineto{\pgfqpoint{0.891610in}{1.690834in}}%
\pgfpathlineto{\pgfqpoint{0.829659in}{1.725308in}}%
\pgfpathlineto{\pgfqpoint{0.752464in}{1.765659in}}%
\pgfpathlineto{\pgfqpoint{0.653422in}{1.814769in}}%
\pgfpathlineto{\pgfqpoint{0.530243in}{1.873254in}}%
\pgfpathlineto{\pgfqpoint{0.372719in}{1.945482in}}%
\pgfpathlineto{\pgfqpoint{0.369214in}{1.947065in}}%
\pgfpathlineto{\pgfqpoint{0.369214in}{1.947065in}}%
\pgfusepath{stroke}%
\end{pgfscope}%
\begin{pgfscope}%
\pgfpathrectangle{\pgfqpoint{0.135000in}{0.135000in}}{\pgfqpoint{2.914509in}{2.695000in}} %
\pgfusepath{clip}%
\pgfsetrectcap%
\pgfsetroundjoin%
\pgfsetlinewidth{1.003750pt}%
\definecolor{currentstroke}{rgb}{0.700000,0.700000,0.700000}%
\pgfsetstrokecolor{currentstroke}%
\pgfsetdash{}{0pt}%
\pgfpathmoveto{\pgfqpoint{1.064282in}{1.482500in}}%
\pgfpathlineto{\pgfqpoint{0.268445in}{1.482500in}}%
\pgfpathlineto{\pgfqpoint{0.268445in}{1.482500in}}%
\pgfusepath{stroke}%
\end{pgfscope}%
\begin{pgfscope}%
\pgfpathrectangle{\pgfqpoint{0.135000in}{0.135000in}}{\pgfqpoint{2.914509in}{2.695000in}} %
\pgfusepath{clip}%
\pgfsetrectcap%
\pgfsetroundjoin%
\pgfsetlinewidth{1.003750pt}%
\definecolor{currentstroke}{rgb}{0.700000,0.700000,0.700000}%
\pgfsetstrokecolor{currentstroke}%
\pgfsetdash{}{0pt}%
\pgfpathmoveto{\pgfqpoint{1.104472in}{1.482500in}}%
\pgfpathlineto{\pgfqpoint{1.103384in}{1.469002in}}%
\pgfpathlineto{\pgfqpoint{1.100117in}{1.455443in}}%
\pgfpathlineto{\pgfqpoint{1.094266in}{1.440954in}}%
\pgfpathlineto{\pgfqpoint{1.085926in}{1.426258in}}%
\pgfpathlineto{\pgfqpoint{1.074376in}{1.410438in}}%
\pgfpathlineto{\pgfqpoint{1.059100in}{1.393350in}}%
\pgfpathlineto{\pgfqpoint{1.039516in}{1.374815in}}%
\pgfpathlineto{\pgfqpoint{1.013747in}{1.353668in}}%
\pgfpathlineto{\pgfqpoint{0.981819in}{1.330478in}}%
\pgfpathlineto{\pgfqpoint{0.941096in}{1.303823in}}%
\pgfpathlineto{\pgfqpoint{0.891610in}{1.274166in}}%
\pgfpathlineto{\pgfqpoint{0.829659in}{1.239692in}}%
\pgfpathlineto{\pgfqpoint{0.752464in}{1.199341in}}%
\pgfpathlineto{\pgfqpoint{0.653422in}{1.150231in}}%
\pgfpathlineto{\pgfqpoint{0.530243in}{1.091746in}}%
\pgfpathlineto{\pgfqpoint{0.372719in}{1.019518in}}%
\pgfpathlineto{\pgfqpoint{0.369214in}{1.017935in}}%
\pgfpathlineto{\pgfqpoint{0.369214in}{1.017935in}}%
\pgfusepath{stroke}%
\end{pgfscope}%
\begin{pgfscope}%
\pgfpathrectangle{\pgfqpoint{0.135000in}{0.135000in}}{\pgfqpoint{2.914509in}{2.695000in}} %
\pgfusepath{clip}%
\pgfsetrectcap%
\pgfsetroundjoin%
\pgfsetlinewidth{1.003750pt}%
\definecolor{currentstroke}{rgb}{0.700000,0.700000,0.700000}%
\pgfsetstrokecolor{currentstroke}%
\pgfsetdash{}{0pt}%
\pgfpathmoveto{\pgfqpoint{1.218922in}{1.482500in}}%
\pgfpathlineto{\pgfqpoint{1.217770in}{1.453148in}}%
\pgfpathlineto{\pgfqpoint{1.214307in}{1.423616in}}%
\pgfpathlineto{\pgfqpoint{1.208511in}{1.393719in}}%
\pgfpathlineto{\pgfqpoint{1.200347in}{1.363275in}}%
\pgfpathlineto{\pgfqpoint{1.189170in}{1.330514in}}%
\pgfpathlineto{\pgfqpoint{1.175250in}{1.296718in}}%
\pgfpathlineto{\pgfqpoint{1.157623in}{1.259952in}}%
\pgfpathlineto{\pgfqpoint{1.136749in}{1.221525in}}%
\pgfpathlineto{\pgfqpoint{1.111286in}{1.179262in}}%
\pgfpathlineto{\pgfqpoint{1.080527in}{1.132517in}}%
\pgfpathlineto{\pgfqpoint{1.043644in}{1.080510in}}%
\pgfpathlineto{\pgfqpoint{0.997855in}{1.019970in}}%
\pgfpathlineto{\pgfqpoint{0.943298in}{0.951682in}}%
\pgfpathlineto{\pgfqpoint{0.878503in}{0.874171in}}%
\pgfpathlineto{\pgfqpoint{0.798969in}{0.782554in}}%
\pgfpathlineto{\pgfqpoint{0.701338in}{0.673578in}}%
\pgfpathlineto{\pgfqpoint{0.656180in}{0.624095in}}%
\pgfpathlineto{\pgfqpoint{0.656180in}{0.624095in}}%
\pgfusepath{stroke}%
\end{pgfscope}%
\begin{pgfscope}%
\pgfpathrectangle{\pgfqpoint{0.135000in}{0.135000in}}{\pgfqpoint{2.914509in}{2.695000in}} %
\pgfusepath{clip}%
\pgfsetrectcap%
\pgfsetroundjoin%
\pgfsetlinewidth{1.003750pt}%
\definecolor{currentstroke}{rgb}{0.700000,0.700000,0.700000}%
\pgfsetstrokecolor{currentstroke}%
\pgfsetdash{}{0pt}%
\pgfpathmoveto{\pgfqpoint{1.390208in}{1.482500in}}%
\pgfpathlineto{\pgfqpoint{1.389072in}{1.430684in}}%
\pgfpathlineto{\pgfqpoint{1.385648in}{1.378285in}}%
\pgfpathlineto{\pgfqpoint{1.379901in}{1.324713in}}%
\pgfpathlineto{\pgfqpoint{1.371763in}{1.269366in}}%
\pgfpathlineto{\pgfqpoint{1.360702in}{1.209427in}}%
\pgfpathlineto{\pgfqpoint{1.346290in}{1.143854in}}%
\pgfpathlineto{\pgfqpoint{1.328017in}{1.071383in}}%
\pgfpathlineto{\pgfqpoint{1.305275in}{0.990482in}}%
\pgfpathlineto{\pgfqpoint{1.276377in}{0.896309in}}%
\pgfpathlineto{\pgfqpoint{1.241042in}{0.788953in}}%
\pgfpathlineto{\pgfqpoint{1.196778in}{0.661742in}}%
\pgfpathlineto{\pgfqpoint{1.141446in}{0.509582in}}%
\pgfpathlineto{\pgfqpoint{1.085655in}{0.360940in}}%
\pgfpathlineto{\pgfqpoint{1.085655in}{0.360940in}}%
\pgfusepath{stroke}%
\end{pgfscope}%
\begin{pgfscope}%
\pgfpathrectangle{\pgfqpoint{0.135000in}{0.135000in}}{\pgfqpoint{2.914509in}{2.695000in}} %
\pgfusepath{clip}%
\pgfsetrectcap%
\pgfsetroundjoin%
\pgfsetlinewidth{1.003750pt}%
\definecolor{currentstroke}{rgb}{0.700000,0.700000,0.700000}%
\pgfsetstrokecolor{currentstroke}%
\pgfsetdash{}{0pt}%
\pgfpathmoveto{\pgfqpoint{1.592255in}{1.482500in}}%
\pgfpathlineto{\pgfqpoint{1.592255in}{0.268532in}}%
\pgfpathlineto{\pgfqpoint{1.592255in}{0.268532in}}%
\pgfusepath{stroke}%
\end{pgfscope}%
\begin{pgfscope}%
\pgfpathrectangle{\pgfqpoint{0.135000in}{0.135000in}}{\pgfqpoint{2.914509in}{2.695000in}} %
\pgfusepath{clip}%
\pgfsetrectcap%
\pgfsetroundjoin%
\pgfsetlinewidth{1.003750pt}%
\definecolor{currentstroke}{rgb}{0.700000,0.700000,0.700000}%
\pgfsetstrokecolor{currentstroke}%
\pgfsetdash{}{0pt}%
\pgfpathmoveto{\pgfqpoint{1.794301in}{1.482500in}}%
\pgfpathlineto{\pgfqpoint{1.795438in}{1.430684in}}%
\pgfpathlineto{\pgfqpoint{1.798861in}{1.378285in}}%
\pgfpathlineto{\pgfqpoint{1.804609in}{1.324713in}}%
\pgfpathlineto{\pgfqpoint{1.812746in}{1.269366in}}%
\pgfpathlineto{\pgfqpoint{1.823807in}{1.209427in}}%
\pgfpathlineto{\pgfqpoint{1.838220in}{1.143854in}}%
\pgfpathlineto{\pgfqpoint{1.856492in}{1.071383in}}%
\pgfpathlineto{\pgfqpoint{1.879234in}{0.990482in}}%
\pgfpathlineto{\pgfqpoint{1.908132in}{0.896309in}}%
\pgfpathlineto{\pgfqpoint{1.943468in}{0.788953in}}%
\pgfpathlineto{\pgfqpoint{1.987731in}{0.661742in}}%
\pgfpathlineto{\pgfqpoint{2.043063in}{0.509582in}}%
\pgfpathlineto{\pgfqpoint{2.098855in}{0.360940in}}%
\pgfpathlineto{\pgfqpoint{2.098855in}{0.360940in}}%
\pgfusepath{stroke}%
\end{pgfscope}%
\begin{pgfscope}%
\pgfpathrectangle{\pgfqpoint{0.135000in}{0.135000in}}{\pgfqpoint{2.914509in}{2.695000in}} %
\pgfusepath{clip}%
\pgfsetrectcap%
\pgfsetroundjoin%
\pgfsetlinewidth{1.003750pt}%
\definecolor{currentstroke}{rgb}{0.700000,0.700000,0.700000}%
\pgfsetstrokecolor{currentstroke}%
\pgfsetdash{}{0pt}%
\pgfpathmoveto{\pgfqpoint{1.965588in}{1.482500in}}%
\pgfpathlineto{\pgfqpoint{1.966740in}{1.453148in}}%
\pgfpathlineto{\pgfqpoint{1.970203in}{1.423616in}}%
\pgfpathlineto{\pgfqpoint{1.975999in}{1.393719in}}%
\pgfpathlineto{\pgfqpoint{1.984163in}{1.363275in}}%
\pgfpathlineto{\pgfqpoint{1.995339in}{1.330514in}}%
\pgfpathlineto{\pgfqpoint{2.009259in}{1.296718in}}%
\pgfpathlineto{\pgfqpoint{2.026887in}{1.259952in}}%
\pgfpathlineto{\pgfqpoint{2.047760in}{1.221525in}}%
\pgfpathlineto{\pgfqpoint{2.073223in}{1.179262in}}%
\pgfpathlineto{\pgfqpoint{2.103982in}{1.132517in}}%
\pgfpathlineto{\pgfqpoint{2.140865in}{1.080510in}}%
\pgfpathlineto{\pgfqpoint{2.186654in}{1.019970in}}%
\pgfpathlineto{\pgfqpoint{2.241212in}{0.951682in}}%
\pgfpathlineto{\pgfqpoint{2.306006in}{0.874171in}}%
\pgfpathlineto{\pgfqpoint{2.385540in}{0.782554in}}%
\pgfpathlineto{\pgfqpoint{2.483171in}{0.673578in}}%
\pgfpathlineto{\pgfqpoint{2.528330in}{0.624095in}}%
\pgfpathlineto{\pgfqpoint{2.528330in}{0.624095in}}%
\pgfusepath{stroke}%
\end{pgfscope}%
\begin{pgfscope}%
\pgfpathrectangle{\pgfqpoint{0.135000in}{0.135000in}}{\pgfqpoint{2.914509in}{2.695000in}} %
\pgfusepath{clip}%
\pgfsetrectcap%
\pgfsetroundjoin%
\pgfsetlinewidth{1.003750pt}%
\definecolor{currentstroke}{rgb}{0.700000,0.700000,0.700000}%
\pgfsetstrokecolor{currentstroke}%
\pgfsetdash{}{0pt}%
\pgfpathmoveto{\pgfqpoint{2.080038in}{1.482500in}}%
\pgfpathlineto{\pgfqpoint{2.081125in}{1.469002in}}%
\pgfpathlineto{\pgfqpoint{2.084392in}{1.455443in}}%
\pgfpathlineto{\pgfqpoint{2.090243in}{1.440954in}}%
\pgfpathlineto{\pgfqpoint{2.098583in}{1.426258in}}%
\pgfpathlineto{\pgfqpoint{2.110134in}{1.410438in}}%
\pgfpathlineto{\pgfqpoint{2.125410in}{1.393350in}}%
\pgfpathlineto{\pgfqpoint{2.144993in}{1.374815in}}%
\pgfpathlineto{\pgfqpoint{2.170762in}{1.353668in}}%
\pgfpathlineto{\pgfqpoint{2.202691in}{1.330478in}}%
\pgfpathlineto{\pgfqpoint{2.243413in}{1.303823in}}%
\pgfpathlineto{\pgfqpoint{2.292900in}{1.274166in}}%
\pgfpathlineto{\pgfqpoint{2.354851in}{1.239692in}}%
\pgfpathlineto{\pgfqpoint{2.432046in}{1.199341in}}%
\pgfpathlineto{\pgfqpoint{2.531087in}{1.150231in}}%
\pgfpathlineto{\pgfqpoint{2.654266in}{1.091746in}}%
\pgfpathlineto{\pgfqpoint{2.811790in}{1.019518in}}%
\pgfpathlineto{\pgfqpoint{2.815295in}{1.017935in}}%
\pgfpathlineto{\pgfqpoint{2.815295in}{1.017935in}}%
\pgfusepath{stroke}%
\end{pgfscope}%
\begin{pgfscope}%
\pgfpathrectangle{\pgfqpoint{0.135000in}{0.135000in}}{\pgfqpoint{2.914509in}{2.695000in}} %
\pgfusepath{clip}%
\pgfsetrectcap%
\pgfsetroundjoin%
\pgfsetlinewidth{1.003750pt}%
\definecolor{currentstroke}{rgb}{0.700000,0.700000,0.700000}%
\pgfsetstrokecolor{currentstroke}%
\pgfsetdash{}{0pt}%
\pgfpathmoveto{\pgfqpoint{2.120227in}{1.482500in}}%
\pgfpathlineto{\pgfqpoint{2.916065in}{1.482500in}}%
\pgfpathlineto{\pgfqpoint{2.916065in}{1.482500in}}%
\pgfusepath{stroke}%
\end{pgfscope}%
\end{pgfpicture}%
\makeatother%
\endgroup%

%% file: trigonometric-symmetric1_image.pgf
%% Creator: Matplotlib, PGF backend
%%
%% To include the figure in your LaTeX document, write
%%   \input{<filename>.pgf}
%%
%% Make sure the required packages are loaded in your preamble
%%   \usepackage{pgf}
%%
%% Figures using additional raster images can only be included by \input if
%% they are in the same directory as the main LaTeX file. For loading figures
%% from other directories you can use the `import` package
%% and then include the figures with
%%   \import{<path to file>}{<filename>.pgf}
%%
%% Matplotlib used the following preamble
%%   \usepackage{fontspec}
%%   \setmainfont{DejaVu Serif}
%%   \setsansfont{DejaVu Sans}
%%   \setmonofont{DejaVu Sans Mono}
%%
\begingroup%
\makeatletter%
\begin{pgfpicture}%
\pgfpathrectangle{\pgfpointorigin}{\pgfqpoint{3.316323in}{3.170273in}}%
\pgfusepath{use as bounding box, clip}%
\begin{pgfscope}%
\pgfsetbuttcap%
\pgfsetmiterjoin%
\definecolor{currentfill}{rgb}{1.000000,1.000000,1.000000}%
\pgfsetfillcolor{currentfill}%
\pgfsetlinewidth{0.000000pt}%
\definecolor{currentstroke}{rgb}{1.000000,1.000000,1.000000}%
\pgfsetstrokecolor{currentstroke}%
\pgfsetdash{}{0pt}%
\pgfpathmoveto{\pgfqpoint{0.000000in}{0.000000in}}%
\pgfpathlineto{\pgfqpoint{3.316323in}{0.000000in}}%
\pgfpathlineto{\pgfqpoint{3.316323in}{3.170273in}}%
\pgfpathlineto{\pgfqpoint{0.000000in}{3.170273in}}%
\pgfpathclose%
\pgfusepath{fill}%
\end{pgfscope}%
\begin{pgfscope}%
\pgfsetbuttcap%
\pgfsetmiterjoin%
\definecolor{currentfill}{rgb}{1.000000,1.000000,1.000000}%
\pgfsetfillcolor{currentfill}%
\pgfsetlinewidth{0.000000pt}%
\definecolor{currentstroke}{rgb}{0.000000,0.000000,0.000000}%
\pgfsetstrokecolor{currentstroke}%
\pgfsetstrokeopacity{0.000000}%
\pgfsetdash{}{0pt}%
\pgfpathmoveto{\pgfqpoint{0.135000in}{0.135000in}}%
\pgfpathlineto{\pgfqpoint{3.049509in}{0.135000in}}%
\pgfpathlineto{\pgfqpoint{3.049509in}{2.830000in}}%
\pgfpathlineto{\pgfqpoint{0.135000in}{2.830000in}}%
\pgfpathclose%
\pgfusepath{fill}%
\end{pgfscope}%
\begin{pgfscope}%
\pgfpathrectangle{\pgfqpoint{0.135000in}{0.135000in}}{\pgfqpoint{2.914509in}{2.695000in}} %
\pgfusepath{clip}%
\pgfsetrectcap%
\pgfsetroundjoin%
\pgfsetlinewidth{0.702625pt}%
\definecolor{currentstroke}{rgb}{0.000000,0.000000,0.000000}%
\pgfsetstrokecolor{currentstroke}%
\pgfsetdash{}{0pt}%
\pgfpathmoveto{\pgfqpoint{0.135000in}{1.482500in}}%
\pgfpathlineto{\pgfqpoint{3.049509in}{1.482500in}}%
\pgfusepath{stroke}%
\end{pgfscope}%
\begin{pgfscope}%
\pgfpathrectangle{\pgfqpoint{0.135000in}{0.135000in}}{\pgfqpoint{2.914509in}{2.695000in}} %
\pgfusepath{clip}%
\pgfsetrectcap%
\pgfsetroundjoin%
\pgfsetlinewidth{0.702625pt}%
\definecolor{currentstroke}{rgb}{0.000000,0.000000,0.000000}%
\pgfsetstrokecolor{currentstroke}%
\pgfsetdash{}{0pt}%
\pgfpathmoveto{\pgfqpoint{1.592255in}{0.135000in}}%
\pgfpathlineto{\pgfqpoint{1.592255in}{2.830000in}}%
\pgfusepath{stroke}%
\end{pgfscope}%
\begin{pgfscope}%
\pgfsetbuttcap%
\pgfsetmiterjoin%
\definecolor{currentfill}{rgb}{0.000000,0.000000,0.000000}%
\pgfsetfillcolor{currentfill}%
\pgfsetlinewidth{1.003750pt}%
\definecolor{currentstroke}{rgb}{0.000000,0.000000,0.000000}%
\pgfsetstrokecolor{currentstroke}%
\pgfsetdash{}{0pt}%
\pgfpathmoveto{\pgfqpoint{3.049509in}{1.482500in}}%
\pgfpathlineto{\pgfqpoint{3.020364in}{1.469025in}}%
\pgfpathlineto{\pgfqpoint{3.029108in}{1.482237in}}%
\pgfpathlineto{\pgfqpoint{0.135000in}{1.482237in}}%
\pgfpathlineto{\pgfqpoint{0.135000in}{1.482763in}}%
\pgfpathlineto{\pgfqpoint{3.029108in}{1.482763in}}%
\pgfpathlineto{\pgfqpoint{3.020364in}{1.495975in}}%
\pgfpathclose%
\pgfusepath{stroke,fill}%
\end{pgfscope}%
\begin{pgfscope}%
\pgfsetbuttcap%
\pgfsetmiterjoin%
\definecolor{currentfill}{rgb}{0.000000,0.000000,0.000000}%
\pgfsetfillcolor{currentfill}%
\pgfsetlinewidth{1.003750pt}%
\definecolor{currentstroke}{rgb}{0.000000,0.000000,0.000000}%
\pgfsetstrokecolor{currentstroke}%
\pgfsetdash{}{0pt}%
\pgfpathmoveto{\pgfqpoint{1.592255in}{2.830000in}}%
\pgfpathlineto{\pgfqpoint{1.605730in}{2.800855in}}%
\pgfpathlineto{\pgfqpoint{1.592518in}{2.809598in}}%
\pgfpathlineto{\pgfqpoint{1.592518in}{0.135000in}}%
\pgfpathlineto{\pgfqpoint{1.591991in}{0.135000in}}%
\pgfpathlineto{\pgfqpoint{1.591991in}{2.809598in}}%
\pgfpathlineto{\pgfqpoint{1.578780in}{2.800855in}}%
\pgfpathclose%
\pgfusepath{stroke,fill}%
\end{pgfscope}%
\begin{pgfscope}%
\pgfpathrectangle{\pgfqpoint{0.135000in}{0.135000in}}{\pgfqpoint{2.914509in}{2.695000in}} %
\pgfusepath{clip}%
\pgfsetrectcap%
\pgfsetroundjoin%
\pgfsetlinewidth{0.803000pt}%
\definecolor{currentstroke}{rgb}{0.000000,0.000000,0.000000}%
\pgfsetstrokecolor{currentstroke}%
\pgfsetdash{}{0pt}%
\pgfpathmoveto{\pgfqpoint{1.592255in}{0.135000in}}%
\pgfpathlineto{\pgfqpoint{1.592255in}{2.830000in}}%
\pgfusepath{stroke}%
\end{pgfscope}%
\begin{pgfscope}%
\pgfpathrectangle{\pgfqpoint{0.135000in}{0.135000in}}{\pgfqpoint{2.914509in}{2.695000in}} %
\pgfusepath{clip}%
\pgfsetrectcap%
\pgfsetroundjoin%
\pgfsetlinewidth{0.803000pt}%
\definecolor{currentstroke}{rgb}{0.000000,0.000000,0.000000}%
\pgfsetstrokecolor{currentstroke}%
\pgfsetdash{}{0pt}%
\pgfpathmoveto{\pgfqpoint{0.135000in}{1.482500in}}%
\pgfpathlineto{\pgfqpoint{3.049509in}{1.482500in}}%
\pgfusepath{stroke}%
\end{pgfscope}%
\begin{pgfscope}%
\pgfpathrectangle{\pgfqpoint{0.135000in}{0.135000in}}{\pgfqpoint{2.914509in}{2.695000in}} %
\pgfusepath{clip}%
\pgfsetrectcap%
\pgfsetroundjoin%
\pgfsetlinewidth{0.803000pt}%
\definecolor{currentstroke}{rgb}{0.000000,0.000000,0.000000}%
\pgfsetstrokecolor{currentstroke}%
\pgfsetdash{}{0pt}%
\pgfpathmoveto{\pgfqpoint{1.592255in}{0.135000in}}%
\pgfpathlineto{\pgfqpoint{1.592255in}{2.830000in}}%
\pgfusepath{stroke}%
\end{pgfscope}%
\begin{pgfscope}%
\pgfpathrectangle{\pgfqpoint{0.135000in}{0.135000in}}{\pgfqpoint{2.914509in}{2.695000in}} %
\pgfusepath{clip}%
\pgfsetrectcap%
\pgfsetroundjoin%
\pgfsetlinewidth{0.803000pt}%
\definecolor{currentstroke}{rgb}{0.000000,0.000000,0.000000}%
\pgfsetstrokecolor{currentstroke}%
\pgfsetdash{}{0pt}%
\pgfpathmoveto{\pgfqpoint{0.135000in}{1.482500in}}%
\pgfpathlineto{\pgfqpoint{3.049509in}{1.482500in}}%
\pgfusepath{stroke}%
\end{pgfscope}%
\begin{pgfscope}%
\pgftext[x=1.592255in,y=2.964750in,,base]{\sffamily\fontsize{10.000000}{12.000000}\selectfont \(\displaystyle y\)}%
\end{pgfscope}%
\begin{pgfscope}%
\pgftext[x=3.136945in,y=1.482500in,left,]{\sffamily\fontsize{10.000000}{12.000000}\selectfont \(\displaystyle x\)}%
\end{pgfscope}%
\begin{pgfscope}%
\pgfpathrectangle{\pgfqpoint{0.135000in}{0.135000in}}{\pgfqpoint{2.914509in}{2.695000in}} %
\pgfusepath{clip}%
\pgfsetrectcap%
\pgfsetroundjoin%
\pgfsetlinewidth{1.003750pt}%
\definecolor{currentstroke}{rgb}{0.700000,0.700000,0.700000}%
\pgfsetstrokecolor{currentstroke}%
\pgfsetdash{}{0pt}%
\pgfpathmoveto{\pgfqpoint{2.120227in}{1.482500in}}%
\pgfpathlineto{\pgfqpoint{2.080038in}{1.482500in}}%
\pgfpathlineto{\pgfqpoint{1.965588in}{1.482500in}}%
\pgfpathlineto{\pgfqpoint{1.794301in}{1.482500in}}%
\pgfpathlineto{\pgfqpoint{1.592255in}{1.482500in}}%
\pgfpathlineto{\pgfqpoint{1.390208in}{1.482500in}}%
\pgfpathlineto{\pgfqpoint{1.218922in}{1.482500in}}%
\pgfpathlineto{\pgfqpoint{1.104472in}{1.482500in}}%
\pgfpathlineto{\pgfqpoint{1.064282in}{1.482500in}}%
\pgfpathlineto{\pgfqpoint{1.104472in}{1.482500in}}%
\pgfpathlineto{\pgfqpoint{1.218922in}{1.482500in}}%
\pgfpathlineto{\pgfqpoint{1.390208in}{1.482500in}}%
\pgfpathlineto{\pgfqpoint{1.592255in}{1.482500in}}%
\pgfpathlineto{\pgfqpoint{1.794301in}{1.482500in}}%
\pgfpathlineto{\pgfqpoint{1.965588in}{1.482500in}}%
\pgfpathlineto{\pgfqpoint{2.080038in}{1.482500in}}%
\pgfpathlineto{\pgfqpoint{2.120227in}{1.482500in}}%
\pgfusepath{stroke}%
\end{pgfscope}%
\begin{pgfscope}%
\pgfpathrectangle{\pgfqpoint{0.135000in}{0.135000in}}{\pgfqpoint{2.914509in}{2.695000in}} %
\pgfusepath{clip}%
\pgfsetrectcap%
\pgfsetroundjoin%
\pgfsetlinewidth{1.003750pt}%
\definecolor{currentstroke}{rgb}{0.700000,0.700000,0.700000}%
\pgfsetstrokecolor{currentstroke}%
\pgfsetdash{}{0pt}%
\pgfpathmoveto{\pgfqpoint{2.161463in}{1.482500in}}%
\pgfpathlineto{\pgfqpoint{2.118135in}{1.563898in}}%
\pgfpathlineto{\pgfqpoint{1.994746in}{1.632905in}}%
\pgfpathlineto{\pgfqpoint{1.810081in}{1.679013in}}%
\pgfpathlineto{\pgfqpoint{1.592255in}{1.695204in}}%
\pgfpathlineto{\pgfqpoint{1.374428in}{1.679013in}}%
\pgfpathlineto{\pgfqpoint{1.189764in}{1.632905in}}%
\pgfpathlineto{\pgfqpoint{1.066375in}{1.563898in}}%
\pgfpathlineto{\pgfqpoint{1.023046in}{1.482500in}}%
\pgfpathlineto{\pgfqpoint{1.066375in}{1.401102in}}%
\pgfpathlineto{\pgfqpoint{1.189764in}{1.332095in}}%
\pgfpathlineto{\pgfqpoint{1.374428in}{1.285987in}}%
\pgfpathlineto{\pgfqpoint{1.592255in}{1.269796in}}%
\pgfpathlineto{\pgfqpoint{1.810081in}{1.285987in}}%
\pgfpathlineto{\pgfqpoint{1.994746in}{1.332095in}}%
\pgfpathlineto{\pgfqpoint{2.118135in}{1.401102in}}%
\pgfpathlineto{\pgfqpoint{2.161463in}{1.482500in}}%
\pgfusepath{stroke}%
\end{pgfscope}%
\begin{pgfscope}%
\pgfpathrectangle{\pgfqpoint{0.135000in}{0.135000in}}{\pgfqpoint{2.914509in}{2.695000in}} %
\pgfusepath{clip}%
\pgfsetrectcap%
\pgfsetroundjoin%
\pgfsetlinewidth{1.003750pt}%
\definecolor{currentstroke}{rgb}{0.700000,0.700000,0.700000}%
\pgfsetstrokecolor{currentstroke}%
\pgfsetdash{}{0pt}%
\pgfpathmoveto{\pgfqpoint{2.291612in}{1.482500in}}%
\pgfpathlineto{\pgfqpoint{2.238376in}{1.658012in}}%
\pgfpathlineto{\pgfqpoint{2.086775in}{1.806803in}}%
\pgfpathlineto{\pgfqpoint{1.859887in}{1.906223in}}%
\pgfpathlineto{\pgfqpoint{1.592255in}{1.941134in}}%
\pgfpathlineto{\pgfqpoint{1.324622in}{1.906223in}}%
\pgfpathlineto{\pgfqpoint{1.097735in}{1.806803in}}%
\pgfpathlineto{\pgfqpoint{0.946133in}{1.658012in}}%
\pgfpathlineto{\pgfqpoint{0.892898in}{1.482500in}}%
\pgfpathlineto{\pgfqpoint{0.946133in}{1.306988in}}%
\pgfpathlineto{\pgfqpoint{1.097735in}{1.158197in}}%
\pgfpathlineto{\pgfqpoint{1.324622in}{1.058777in}}%
\pgfpathlineto{\pgfqpoint{1.592255in}{1.023866in}}%
\pgfpathlineto{\pgfqpoint{1.859887in}{1.058777in}}%
\pgfpathlineto{\pgfqpoint{2.086775in}{1.158197in}}%
\pgfpathlineto{\pgfqpoint{2.238376in}{1.306988in}}%
\pgfpathlineto{\pgfqpoint{2.291612in}{1.482500in}}%
\pgfusepath{stroke}%
\end{pgfscope}%
\begin{pgfscope}%
\pgfpathrectangle{\pgfqpoint{0.135000in}{0.135000in}}{\pgfqpoint{2.914509in}{2.695000in}} %
\pgfusepath{clip}%
\pgfsetrectcap%
\pgfsetroundjoin%
\pgfsetlinewidth{1.003750pt}%
\definecolor{currentstroke}{rgb}{0.700000,0.700000,0.700000}%
\pgfsetstrokecolor{currentstroke}%
\pgfsetdash{}{0pt}%
\pgfpathmoveto{\pgfqpoint{2.531003in}{1.482500in}}%
\pgfpathlineto{\pgfqpoint{2.459545in}{1.779541in}}%
\pgfpathlineto{\pgfqpoint{2.256050in}{2.031360in}}%
\pgfpathlineto{\pgfqpoint{1.951498in}{2.199620in}}%
\pgfpathlineto{\pgfqpoint{1.592255in}{2.258705in}}%
\pgfpathlineto{\pgfqpoint{1.233011in}{2.199620in}}%
\pgfpathlineto{\pgfqpoint{0.928459in}{2.031360in}}%
\pgfpathlineto{\pgfqpoint{0.724964in}{1.779541in}}%
\pgfpathlineto{\pgfqpoint{0.653506in}{1.482500in}}%
\pgfpathlineto{\pgfqpoint{0.724964in}{1.185459in}}%
\pgfpathlineto{\pgfqpoint{0.928459in}{0.933640in}}%
\pgfpathlineto{\pgfqpoint{1.233011in}{0.765380in}}%
\pgfpathlineto{\pgfqpoint{1.592255in}{0.706295in}}%
\pgfpathlineto{\pgfqpoint{1.951498in}{0.765380in}}%
\pgfpathlineto{\pgfqpoint{2.256050in}{0.933640in}}%
\pgfpathlineto{\pgfqpoint{2.459545in}{1.185459in}}%
\pgfpathlineto{\pgfqpoint{2.531003in}{1.482500in}}%
\pgfusepath{stroke}%
\end{pgfscope}%
\begin{pgfscope}%
\pgfpathrectangle{\pgfqpoint{0.135000in}{0.135000in}}{\pgfqpoint{2.914509in}{2.695000in}} %
\pgfusepath{clip}%
\pgfsetrectcap%
\pgfsetroundjoin%
\pgfsetlinewidth{1.003750pt}%
\definecolor{currentstroke}{rgb}{0.700000,0.700000,0.700000}%
\pgfsetstrokecolor{currentstroke}%
\pgfsetdash{}{0pt}%
\pgfpathmoveto{\pgfqpoint{2.917032in}{1.482500in}}%
\pgfpathlineto{\pgfqpoint{2.816189in}{1.947469in}}%
\pgfpathlineto{\pgfqpoint{2.529013in}{2.341651in}}%
\pgfpathlineto{\pgfqpoint{2.099225in}{2.605034in}}%
\pgfpathlineto{\pgfqpoint{1.592255in}{2.697522in}}%
\pgfpathlineto{\pgfqpoint{1.085284in}{2.605034in}}%
\pgfpathlineto{\pgfqpoint{0.655496in}{2.341651in}}%
\pgfpathlineto{\pgfqpoint{0.368320in}{1.947469in}}%
\pgfpathlineto{\pgfqpoint{0.267478in}{1.482500in}}%
\pgfpathlineto{\pgfqpoint{0.368320in}{1.017531in}}%
\pgfpathlineto{\pgfqpoint{0.655496in}{0.623349in}}%
\pgfpathlineto{\pgfqpoint{1.085284in}{0.359966in}}%
\pgfpathlineto{\pgfqpoint{1.592255in}{0.267478in}}%
\pgfpathlineto{\pgfqpoint{2.099225in}{0.359966in}}%
\pgfpathlineto{\pgfqpoint{2.529013in}{0.623349in}}%
\pgfpathlineto{\pgfqpoint{2.816189in}{1.017531in}}%
\pgfpathlineto{\pgfqpoint{2.917032in}{1.482500in}}%
\pgfusepath{stroke}%
\end{pgfscope}%
\begin{pgfscope}%
\pgfpathrectangle{\pgfqpoint{0.135000in}{0.135000in}}{\pgfqpoint{2.914509in}{2.695000in}} %
\pgfusepath{clip}%
\pgfsetrectcap%
\pgfsetroundjoin%
\pgfsetlinewidth{1.003750pt}%
\definecolor{currentstroke}{rgb}{0.700000,0.700000,0.700000}%
\pgfsetstrokecolor{currentstroke}%
\pgfsetdash{}{0pt}%
\pgfpathmoveto{\pgfqpoint{2.120227in}{1.482500in}}%
\pgfpathlineto{\pgfqpoint{2.161463in}{1.482500in}}%
\pgfpathlineto{\pgfqpoint{2.291612in}{1.482500in}}%
\pgfpathlineto{\pgfqpoint{2.531003in}{1.482500in}}%
\pgfpathlineto{\pgfqpoint{2.917032in}{1.482500in}}%
\pgfusepath{stroke}%
\end{pgfscope}%
\begin{pgfscope}%
\pgfpathrectangle{\pgfqpoint{0.135000in}{0.135000in}}{\pgfqpoint{2.914509in}{2.695000in}} %
\pgfusepath{clip}%
\pgfsetrectcap%
\pgfsetroundjoin%
\pgfsetlinewidth{1.003750pt}%
\definecolor{currentstroke}{rgb}{0.700000,0.700000,0.700000}%
\pgfsetstrokecolor{currentstroke}%
\pgfsetdash{}{0pt}%
\pgfpathmoveto{\pgfqpoint{2.080038in}{1.482500in}}%
\pgfpathlineto{\pgfqpoint{2.118135in}{1.563898in}}%
\pgfpathlineto{\pgfqpoint{2.238376in}{1.658012in}}%
\pgfpathlineto{\pgfqpoint{2.459545in}{1.779541in}}%
\pgfpathlineto{\pgfqpoint{2.816189in}{1.947469in}}%
\pgfusepath{stroke}%
\end{pgfscope}%
\begin{pgfscope}%
\pgfpathrectangle{\pgfqpoint{0.135000in}{0.135000in}}{\pgfqpoint{2.914509in}{2.695000in}} %
\pgfusepath{clip}%
\pgfsetrectcap%
\pgfsetroundjoin%
\pgfsetlinewidth{1.003750pt}%
\definecolor{currentstroke}{rgb}{0.700000,0.700000,0.700000}%
\pgfsetstrokecolor{currentstroke}%
\pgfsetdash{}{0pt}%
\pgfpathmoveto{\pgfqpoint{1.965588in}{1.482500in}}%
\pgfpathlineto{\pgfqpoint{1.994746in}{1.632905in}}%
\pgfpathlineto{\pgfqpoint{2.086775in}{1.806803in}}%
\pgfpathlineto{\pgfqpoint{2.256050in}{2.031360in}}%
\pgfpathlineto{\pgfqpoint{2.529013in}{2.341651in}}%
\pgfusepath{stroke}%
\end{pgfscope}%
\begin{pgfscope}%
\pgfpathrectangle{\pgfqpoint{0.135000in}{0.135000in}}{\pgfqpoint{2.914509in}{2.695000in}} %
\pgfusepath{clip}%
\pgfsetrectcap%
\pgfsetroundjoin%
\pgfsetlinewidth{1.003750pt}%
\definecolor{currentstroke}{rgb}{0.700000,0.700000,0.700000}%
\pgfsetstrokecolor{currentstroke}%
\pgfsetdash{}{0pt}%
\pgfpathmoveto{\pgfqpoint{1.794301in}{1.482500in}}%
\pgfpathlineto{\pgfqpoint{1.810081in}{1.679013in}}%
\pgfpathlineto{\pgfqpoint{1.859887in}{1.906223in}}%
\pgfpathlineto{\pgfqpoint{1.951498in}{2.199620in}}%
\pgfpathlineto{\pgfqpoint{2.099225in}{2.605034in}}%
\pgfusepath{stroke}%
\end{pgfscope}%
\begin{pgfscope}%
\pgfpathrectangle{\pgfqpoint{0.135000in}{0.135000in}}{\pgfqpoint{2.914509in}{2.695000in}} %
\pgfusepath{clip}%
\pgfsetrectcap%
\pgfsetroundjoin%
\pgfsetlinewidth{1.003750pt}%
\definecolor{currentstroke}{rgb}{0.700000,0.700000,0.700000}%
\pgfsetstrokecolor{currentstroke}%
\pgfsetdash{}{0pt}%
\pgfpathmoveto{\pgfqpoint{1.592255in}{1.482500in}}%
\pgfpathlineto{\pgfqpoint{1.592255in}{1.695204in}}%
\pgfpathlineto{\pgfqpoint{1.592255in}{1.941134in}}%
\pgfpathlineto{\pgfqpoint{1.592255in}{2.258705in}}%
\pgfpathlineto{\pgfqpoint{1.592255in}{2.697522in}}%
\pgfusepath{stroke}%
\end{pgfscope}%
\begin{pgfscope}%
\pgfpathrectangle{\pgfqpoint{0.135000in}{0.135000in}}{\pgfqpoint{2.914509in}{2.695000in}} %
\pgfusepath{clip}%
\pgfsetrectcap%
\pgfsetroundjoin%
\pgfsetlinewidth{1.003750pt}%
\definecolor{currentstroke}{rgb}{0.700000,0.700000,0.700000}%
\pgfsetstrokecolor{currentstroke}%
\pgfsetdash{}{0pt}%
\pgfpathmoveto{\pgfqpoint{1.390208in}{1.482500in}}%
\pgfpathlineto{\pgfqpoint{1.374428in}{1.679013in}}%
\pgfpathlineto{\pgfqpoint{1.324622in}{1.906223in}}%
\pgfpathlineto{\pgfqpoint{1.233011in}{2.199620in}}%
\pgfpathlineto{\pgfqpoint{1.085284in}{2.605034in}}%
\pgfusepath{stroke}%
\end{pgfscope}%
\begin{pgfscope}%
\pgfpathrectangle{\pgfqpoint{0.135000in}{0.135000in}}{\pgfqpoint{2.914509in}{2.695000in}} %
\pgfusepath{clip}%
\pgfsetrectcap%
\pgfsetroundjoin%
\pgfsetlinewidth{1.003750pt}%
\definecolor{currentstroke}{rgb}{0.700000,0.700000,0.700000}%
\pgfsetstrokecolor{currentstroke}%
\pgfsetdash{}{0pt}%
\pgfpathmoveto{\pgfqpoint{1.218922in}{1.482500in}}%
\pgfpathlineto{\pgfqpoint{1.189764in}{1.632905in}}%
\pgfpathlineto{\pgfqpoint{1.097735in}{1.806803in}}%
\pgfpathlineto{\pgfqpoint{0.928459in}{2.031360in}}%
\pgfpathlineto{\pgfqpoint{0.655496in}{2.341651in}}%
\pgfusepath{stroke}%
\end{pgfscope}%
\begin{pgfscope}%
\pgfpathrectangle{\pgfqpoint{0.135000in}{0.135000in}}{\pgfqpoint{2.914509in}{2.695000in}} %
\pgfusepath{clip}%
\pgfsetrectcap%
\pgfsetroundjoin%
\pgfsetlinewidth{1.003750pt}%
\definecolor{currentstroke}{rgb}{0.700000,0.700000,0.700000}%
\pgfsetstrokecolor{currentstroke}%
\pgfsetdash{}{0pt}%
\pgfpathmoveto{\pgfqpoint{1.104472in}{1.482500in}}%
\pgfpathlineto{\pgfqpoint{1.066375in}{1.563898in}}%
\pgfpathlineto{\pgfqpoint{0.946133in}{1.658012in}}%
\pgfpathlineto{\pgfqpoint{0.724964in}{1.779541in}}%
\pgfpathlineto{\pgfqpoint{0.368320in}{1.947469in}}%
\pgfusepath{stroke}%
\end{pgfscope}%
\begin{pgfscope}%
\pgfpathrectangle{\pgfqpoint{0.135000in}{0.135000in}}{\pgfqpoint{2.914509in}{2.695000in}} %
\pgfusepath{clip}%
\pgfsetrectcap%
\pgfsetroundjoin%
\pgfsetlinewidth{1.003750pt}%
\definecolor{currentstroke}{rgb}{0.700000,0.700000,0.700000}%
\pgfsetstrokecolor{currentstroke}%
\pgfsetdash{}{0pt}%
\pgfpathmoveto{\pgfqpoint{1.064282in}{1.482500in}}%
\pgfpathlineto{\pgfqpoint{1.023046in}{1.482500in}}%
\pgfpathlineto{\pgfqpoint{0.892898in}{1.482500in}}%
\pgfpathlineto{\pgfqpoint{0.653506in}{1.482500in}}%
\pgfpathlineto{\pgfqpoint{0.267478in}{1.482500in}}%
\pgfusepath{stroke}%
\end{pgfscope}%
\begin{pgfscope}%
\pgfpathrectangle{\pgfqpoint{0.135000in}{0.135000in}}{\pgfqpoint{2.914509in}{2.695000in}} %
\pgfusepath{clip}%
\pgfsetrectcap%
\pgfsetroundjoin%
\pgfsetlinewidth{1.003750pt}%
\definecolor{currentstroke}{rgb}{0.700000,0.700000,0.700000}%
\pgfsetstrokecolor{currentstroke}%
\pgfsetdash{}{0pt}%
\pgfpathmoveto{\pgfqpoint{1.104472in}{1.482500in}}%
\pgfpathlineto{\pgfqpoint{1.066375in}{1.401102in}}%
\pgfpathlineto{\pgfqpoint{0.946133in}{1.306988in}}%
\pgfpathlineto{\pgfqpoint{0.724964in}{1.185459in}}%
\pgfpathlineto{\pgfqpoint{0.368320in}{1.017531in}}%
\pgfusepath{stroke}%
\end{pgfscope}%
\begin{pgfscope}%
\pgfpathrectangle{\pgfqpoint{0.135000in}{0.135000in}}{\pgfqpoint{2.914509in}{2.695000in}} %
\pgfusepath{clip}%
\pgfsetrectcap%
\pgfsetroundjoin%
\pgfsetlinewidth{1.003750pt}%
\definecolor{currentstroke}{rgb}{0.700000,0.700000,0.700000}%
\pgfsetstrokecolor{currentstroke}%
\pgfsetdash{}{0pt}%
\pgfpathmoveto{\pgfqpoint{1.218922in}{1.482500in}}%
\pgfpathlineto{\pgfqpoint{1.189764in}{1.332095in}}%
\pgfpathlineto{\pgfqpoint{1.097735in}{1.158197in}}%
\pgfpathlineto{\pgfqpoint{0.928459in}{0.933640in}}%
\pgfpathlineto{\pgfqpoint{0.655496in}{0.623349in}}%
\pgfusepath{stroke}%
\end{pgfscope}%
\begin{pgfscope}%
\pgfpathrectangle{\pgfqpoint{0.135000in}{0.135000in}}{\pgfqpoint{2.914509in}{2.695000in}} %
\pgfusepath{clip}%
\pgfsetrectcap%
\pgfsetroundjoin%
\pgfsetlinewidth{1.003750pt}%
\definecolor{currentstroke}{rgb}{0.700000,0.700000,0.700000}%
\pgfsetstrokecolor{currentstroke}%
\pgfsetdash{}{0pt}%
\pgfpathmoveto{\pgfqpoint{1.390208in}{1.482500in}}%
\pgfpathlineto{\pgfqpoint{1.374428in}{1.285987in}}%
\pgfpathlineto{\pgfqpoint{1.324622in}{1.058777in}}%
\pgfpathlineto{\pgfqpoint{1.233011in}{0.765380in}}%
\pgfpathlineto{\pgfqpoint{1.085284in}{0.359966in}}%
\pgfusepath{stroke}%
\end{pgfscope}%
\begin{pgfscope}%
\pgfpathrectangle{\pgfqpoint{0.135000in}{0.135000in}}{\pgfqpoint{2.914509in}{2.695000in}} %
\pgfusepath{clip}%
\pgfsetrectcap%
\pgfsetroundjoin%
\pgfsetlinewidth{1.003750pt}%
\definecolor{currentstroke}{rgb}{0.700000,0.700000,0.700000}%
\pgfsetstrokecolor{currentstroke}%
\pgfsetdash{}{0pt}%
\pgfpathmoveto{\pgfqpoint{1.592255in}{1.482500in}}%
\pgfpathlineto{\pgfqpoint{1.592255in}{1.269796in}}%
\pgfpathlineto{\pgfqpoint{1.592255in}{1.023866in}}%
\pgfpathlineto{\pgfqpoint{1.592255in}{0.706295in}}%
\pgfpathlineto{\pgfqpoint{1.592255in}{0.267478in}}%
\pgfusepath{stroke}%
\end{pgfscope}%
\begin{pgfscope}%
\pgfpathrectangle{\pgfqpoint{0.135000in}{0.135000in}}{\pgfqpoint{2.914509in}{2.695000in}} %
\pgfusepath{clip}%
\pgfsetrectcap%
\pgfsetroundjoin%
\pgfsetlinewidth{1.003750pt}%
\definecolor{currentstroke}{rgb}{0.700000,0.700000,0.700000}%
\pgfsetstrokecolor{currentstroke}%
\pgfsetdash{}{0pt}%
\pgfpathmoveto{\pgfqpoint{1.794301in}{1.482500in}}%
\pgfpathlineto{\pgfqpoint{1.810081in}{1.285987in}}%
\pgfpathlineto{\pgfqpoint{1.859887in}{1.058777in}}%
\pgfpathlineto{\pgfqpoint{1.951498in}{0.765380in}}%
\pgfpathlineto{\pgfqpoint{2.099225in}{0.359966in}}%
\pgfusepath{stroke}%
\end{pgfscope}%
\begin{pgfscope}%
\pgfpathrectangle{\pgfqpoint{0.135000in}{0.135000in}}{\pgfqpoint{2.914509in}{2.695000in}} %
\pgfusepath{clip}%
\pgfsetrectcap%
\pgfsetroundjoin%
\pgfsetlinewidth{1.003750pt}%
\definecolor{currentstroke}{rgb}{0.700000,0.700000,0.700000}%
\pgfsetstrokecolor{currentstroke}%
\pgfsetdash{}{0pt}%
\pgfpathmoveto{\pgfqpoint{1.965588in}{1.482500in}}%
\pgfpathlineto{\pgfqpoint{1.994746in}{1.332095in}}%
\pgfpathlineto{\pgfqpoint{2.086775in}{1.158197in}}%
\pgfpathlineto{\pgfqpoint{2.256050in}{0.933640in}}%
\pgfpathlineto{\pgfqpoint{2.529013in}{0.623349in}}%
\pgfusepath{stroke}%
\end{pgfscope}%
\begin{pgfscope}%
\pgfpathrectangle{\pgfqpoint{0.135000in}{0.135000in}}{\pgfqpoint{2.914509in}{2.695000in}} %
\pgfusepath{clip}%
\pgfsetrectcap%
\pgfsetroundjoin%
\pgfsetlinewidth{1.003750pt}%
\definecolor{currentstroke}{rgb}{0.700000,0.700000,0.700000}%
\pgfsetstrokecolor{currentstroke}%
\pgfsetdash{}{0pt}%
\pgfpathmoveto{\pgfqpoint{2.080038in}{1.482500in}}%
\pgfpathlineto{\pgfqpoint{2.118135in}{1.401102in}}%
\pgfpathlineto{\pgfqpoint{2.238376in}{1.306988in}}%
\pgfpathlineto{\pgfqpoint{2.459545in}{1.185459in}}%
\pgfpathlineto{\pgfqpoint{2.816189in}{1.017531in}}%
\pgfusepath{stroke}%
\end{pgfscope}%
\begin{pgfscope}%
\pgfpathrectangle{\pgfqpoint{0.135000in}{0.135000in}}{\pgfqpoint{2.914509in}{2.695000in}} %
\pgfusepath{clip}%
\pgfsetrectcap%
\pgfsetroundjoin%
\pgfsetlinewidth{1.003750pt}%
\definecolor{currentstroke}{rgb}{0.700000,0.700000,0.700000}%
\pgfsetstrokecolor{currentstroke}%
\pgfsetdash{}{0pt}%
\pgfpathmoveto{\pgfqpoint{2.120227in}{1.482500in}}%
\pgfpathlineto{\pgfqpoint{2.161463in}{1.482500in}}%
\pgfpathlineto{\pgfqpoint{2.291612in}{1.482500in}}%
\pgfpathlineto{\pgfqpoint{2.531003in}{1.482500in}}%
\pgfpathlineto{\pgfqpoint{2.917032in}{1.482500in}}%
\pgfusepath{stroke}%
\end{pgfscope}%
\begin{pgfscope}%
\pgfpathrectangle{\pgfqpoint{0.135000in}{0.135000in}}{\pgfqpoint{2.914509in}{2.695000in}} %
\pgfusepath{clip}%
\pgfsetbuttcap%
\pgfsetroundjoin%
\definecolor{currentfill}{rgb}{0.700000,0.700000,0.700000}%
\pgfsetfillcolor{currentfill}%
\pgfsetlinewidth{1.003750pt}%
\definecolor{currentstroke}{rgb}{0.700000,0.700000,0.700000}%
\pgfsetstrokecolor{currentstroke}%
\pgfsetdash{}{0pt}%
\pgfsys@defobject{currentmarker}{\pgfqpoint{-0.010417in}{-0.010417in}}{\pgfqpoint{0.010417in}{0.010417in}}{%
\pgfpathmoveto{\pgfqpoint{0.000000in}{-0.010417in}}%
\pgfpathcurveto{\pgfqpoint{0.002763in}{-0.010417in}}{\pgfqpoint{0.005412in}{-0.009319in}}{\pgfqpoint{0.007366in}{-0.007366in}}%
\pgfpathcurveto{\pgfqpoint{0.009319in}{-0.005412in}}{\pgfqpoint{0.010417in}{-0.002763in}}{\pgfqpoint{0.010417in}{0.000000in}}%
\pgfpathcurveto{\pgfqpoint{0.010417in}{0.002763in}}{\pgfqpoint{0.009319in}{0.005412in}}{\pgfqpoint{0.007366in}{0.007366in}}%
\pgfpathcurveto{\pgfqpoint{0.005412in}{0.009319in}}{\pgfqpoint{0.002763in}{0.010417in}}{\pgfqpoint{0.000000in}{0.010417in}}%
\pgfpathcurveto{\pgfqpoint{-0.002763in}{0.010417in}}{\pgfqpoint{-0.005412in}{0.009319in}}{\pgfqpoint{-0.007366in}{0.007366in}}%
\pgfpathcurveto{\pgfqpoint{-0.009319in}{0.005412in}}{\pgfqpoint{-0.010417in}{0.002763in}}{\pgfqpoint{-0.010417in}{0.000000in}}%
\pgfpathcurveto{\pgfqpoint{-0.010417in}{-0.002763in}}{\pgfqpoint{-0.009319in}{-0.005412in}}{\pgfqpoint{-0.007366in}{-0.007366in}}%
\pgfpathcurveto{\pgfqpoint{-0.005412in}{-0.009319in}}{\pgfqpoint{-0.002763in}{-0.010417in}}{\pgfqpoint{0.000000in}{-0.010417in}}%
\pgfpathclose%
\pgfusepath{stroke,fill}%
}%
\begin{pgfscope}%
\pgfsys@transformshift{2.120227in}{1.482500in}%
\pgfsys@useobject{currentmarker}{}%
\end{pgfscope}%
\begin{pgfscope}%
\pgfsys@transformshift{2.161463in}{1.482500in}%
\pgfsys@useobject{currentmarker}{}%
\end{pgfscope}%
\begin{pgfscope}%
\pgfsys@transformshift{2.291612in}{1.482500in}%
\pgfsys@useobject{currentmarker}{}%
\end{pgfscope}%
\begin{pgfscope}%
\pgfsys@transformshift{2.531003in}{1.482500in}%
\pgfsys@useobject{currentmarker}{}%
\end{pgfscope}%
\begin{pgfscope}%
\pgfsys@transformshift{2.917032in}{1.482500in}%
\pgfsys@useobject{currentmarker}{}%
\end{pgfscope}%
\begin{pgfscope}%
\pgfsys@transformshift{2.080038in}{1.482500in}%
\pgfsys@useobject{currentmarker}{}%
\end{pgfscope}%
\begin{pgfscope}%
\pgfsys@transformshift{2.118135in}{1.563898in}%
\pgfsys@useobject{currentmarker}{}%
\end{pgfscope}%
\begin{pgfscope}%
\pgfsys@transformshift{2.238376in}{1.658012in}%
\pgfsys@useobject{currentmarker}{}%
\end{pgfscope}%
\begin{pgfscope}%
\pgfsys@transformshift{2.459545in}{1.779541in}%
\pgfsys@useobject{currentmarker}{}%
\end{pgfscope}%
\begin{pgfscope}%
\pgfsys@transformshift{2.816189in}{1.947469in}%
\pgfsys@useobject{currentmarker}{}%
\end{pgfscope}%
\begin{pgfscope}%
\pgfsys@transformshift{1.965588in}{1.482500in}%
\pgfsys@useobject{currentmarker}{}%
\end{pgfscope}%
\begin{pgfscope}%
\pgfsys@transformshift{1.994746in}{1.632905in}%
\pgfsys@useobject{currentmarker}{}%
\end{pgfscope}%
\begin{pgfscope}%
\pgfsys@transformshift{2.086775in}{1.806803in}%
\pgfsys@useobject{currentmarker}{}%
\end{pgfscope}%
\begin{pgfscope}%
\pgfsys@transformshift{2.256050in}{2.031360in}%
\pgfsys@useobject{currentmarker}{}%
\end{pgfscope}%
\begin{pgfscope}%
\pgfsys@transformshift{2.529013in}{2.341651in}%
\pgfsys@useobject{currentmarker}{}%
\end{pgfscope}%
\begin{pgfscope}%
\pgfsys@transformshift{1.794301in}{1.482500in}%
\pgfsys@useobject{currentmarker}{}%
\end{pgfscope}%
\begin{pgfscope}%
\pgfsys@transformshift{1.810081in}{1.679013in}%
\pgfsys@useobject{currentmarker}{}%
\end{pgfscope}%
\begin{pgfscope}%
\pgfsys@transformshift{1.859887in}{1.906223in}%
\pgfsys@useobject{currentmarker}{}%
\end{pgfscope}%
\begin{pgfscope}%
\pgfsys@transformshift{1.951498in}{2.199620in}%
\pgfsys@useobject{currentmarker}{}%
\end{pgfscope}%
\begin{pgfscope}%
\pgfsys@transformshift{2.099225in}{2.605034in}%
\pgfsys@useobject{currentmarker}{}%
\end{pgfscope}%
\begin{pgfscope}%
\pgfsys@transformshift{1.592255in}{1.482500in}%
\pgfsys@useobject{currentmarker}{}%
\end{pgfscope}%
\begin{pgfscope}%
\pgfsys@transformshift{1.592255in}{1.695204in}%
\pgfsys@useobject{currentmarker}{}%
\end{pgfscope}%
\begin{pgfscope}%
\pgfsys@transformshift{1.592255in}{1.941134in}%
\pgfsys@useobject{currentmarker}{}%
\end{pgfscope}%
\begin{pgfscope}%
\pgfsys@transformshift{1.592255in}{2.258705in}%
\pgfsys@useobject{currentmarker}{}%
\end{pgfscope}%
\begin{pgfscope}%
\pgfsys@transformshift{1.592255in}{2.697522in}%
\pgfsys@useobject{currentmarker}{}%
\end{pgfscope}%
\begin{pgfscope}%
\pgfsys@transformshift{1.390208in}{1.482500in}%
\pgfsys@useobject{currentmarker}{}%
\end{pgfscope}%
\begin{pgfscope}%
\pgfsys@transformshift{1.374428in}{1.679013in}%
\pgfsys@useobject{currentmarker}{}%
\end{pgfscope}%
\begin{pgfscope}%
\pgfsys@transformshift{1.324622in}{1.906223in}%
\pgfsys@useobject{currentmarker}{}%
\end{pgfscope}%
\begin{pgfscope}%
\pgfsys@transformshift{1.233011in}{2.199620in}%
\pgfsys@useobject{currentmarker}{}%
\end{pgfscope}%
\begin{pgfscope}%
\pgfsys@transformshift{1.085284in}{2.605034in}%
\pgfsys@useobject{currentmarker}{}%
\end{pgfscope}%
\begin{pgfscope}%
\pgfsys@transformshift{1.218922in}{1.482500in}%
\pgfsys@useobject{currentmarker}{}%
\end{pgfscope}%
\begin{pgfscope}%
\pgfsys@transformshift{1.189764in}{1.632905in}%
\pgfsys@useobject{currentmarker}{}%
\end{pgfscope}%
\begin{pgfscope}%
\pgfsys@transformshift{1.097735in}{1.806803in}%
\pgfsys@useobject{currentmarker}{}%
\end{pgfscope}%
\begin{pgfscope}%
\pgfsys@transformshift{0.928459in}{2.031360in}%
\pgfsys@useobject{currentmarker}{}%
\end{pgfscope}%
\begin{pgfscope}%
\pgfsys@transformshift{0.655496in}{2.341651in}%
\pgfsys@useobject{currentmarker}{}%
\end{pgfscope}%
\begin{pgfscope}%
\pgfsys@transformshift{1.104472in}{1.482500in}%
\pgfsys@useobject{currentmarker}{}%
\end{pgfscope}%
\begin{pgfscope}%
\pgfsys@transformshift{1.066375in}{1.563898in}%
\pgfsys@useobject{currentmarker}{}%
\end{pgfscope}%
\begin{pgfscope}%
\pgfsys@transformshift{0.946133in}{1.658012in}%
\pgfsys@useobject{currentmarker}{}%
\end{pgfscope}%
\begin{pgfscope}%
\pgfsys@transformshift{0.724964in}{1.779541in}%
\pgfsys@useobject{currentmarker}{}%
\end{pgfscope}%
\begin{pgfscope}%
\pgfsys@transformshift{0.368320in}{1.947469in}%
\pgfsys@useobject{currentmarker}{}%
\end{pgfscope}%
\begin{pgfscope}%
\pgfsys@transformshift{1.064282in}{1.482500in}%
\pgfsys@useobject{currentmarker}{}%
\end{pgfscope}%
\begin{pgfscope}%
\pgfsys@transformshift{1.023046in}{1.482500in}%
\pgfsys@useobject{currentmarker}{}%
\end{pgfscope}%
\begin{pgfscope}%
\pgfsys@transformshift{0.892898in}{1.482500in}%
\pgfsys@useobject{currentmarker}{}%
\end{pgfscope}%
\begin{pgfscope}%
\pgfsys@transformshift{0.653506in}{1.482500in}%
\pgfsys@useobject{currentmarker}{}%
\end{pgfscope}%
\begin{pgfscope}%
\pgfsys@transformshift{0.267478in}{1.482500in}%
\pgfsys@useobject{currentmarker}{}%
\end{pgfscope}%
\begin{pgfscope}%
\pgfsys@transformshift{1.104472in}{1.482500in}%
\pgfsys@useobject{currentmarker}{}%
\end{pgfscope}%
\begin{pgfscope}%
\pgfsys@transformshift{1.066375in}{1.401102in}%
\pgfsys@useobject{currentmarker}{}%
\end{pgfscope}%
\begin{pgfscope}%
\pgfsys@transformshift{0.946133in}{1.306988in}%
\pgfsys@useobject{currentmarker}{}%
\end{pgfscope}%
\begin{pgfscope}%
\pgfsys@transformshift{0.724964in}{1.185459in}%
\pgfsys@useobject{currentmarker}{}%
\end{pgfscope}%
\begin{pgfscope}%
\pgfsys@transformshift{0.368320in}{1.017531in}%
\pgfsys@useobject{currentmarker}{}%
\end{pgfscope}%
\begin{pgfscope}%
\pgfsys@transformshift{1.218922in}{1.482500in}%
\pgfsys@useobject{currentmarker}{}%
\end{pgfscope}%
\begin{pgfscope}%
\pgfsys@transformshift{1.189764in}{1.332095in}%
\pgfsys@useobject{currentmarker}{}%
\end{pgfscope}%
\begin{pgfscope}%
\pgfsys@transformshift{1.097735in}{1.158197in}%
\pgfsys@useobject{currentmarker}{}%
\end{pgfscope}%
\begin{pgfscope}%
\pgfsys@transformshift{0.928459in}{0.933640in}%
\pgfsys@useobject{currentmarker}{}%
\end{pgfscope}%
\begin{pgfscope}%
\pgfsys@transformshift{0.655496in}{0.623349in}%
\pgfsys@useobject{currentmarker}{}%
\end{pgfscope}%
\begin{pgfscope}%
\pgfsys@transformshift{1.390208in}{1.482500in}%
\pgfsys@useobject{currentmarker}{}%
\end{pgfscope}%
\begin{pgfscope}%
\pgfsys@transformshift{1.374428in}{1.285987in}%
\pgfsys@useobject{currentmarker}{}%
\end{pgfscope}%
\begin{pgfscope}%
\pgfsys@transformshift{1.324622in}{1.058777in}%
\pgfsys@useobject{currentmarker}{}%
\end{pgfscope}%
\begin{pgfscope}%
\pgfsys@transformshift{1.233011in}{0.765380in}%
\pgfsys@useobject{currentmarker}{}%
\end{pgfscope}%
\begin{pgfscope}%
\pgfsys@transformshift{1.085284in}{0.359966in}%
\pgfsys@useobject{currentmarker}{}%
\end{pgfscope}%
\begin{pgfscope}%
\pgfsys@transformshift{1.592255in}{1.482500in}%
\pgfsys@useobject{currentmarker}{}%
\end{pgfscope}%
\begin{pgfscope}%
\pgfsys@transformshift{1.592255in}{1.269796in}%
\pgfsys@useobject{currentmarker}{}%
\end{pgfscope}%
\begin{pgfscope}%
\pgfsys@transformshift{1.592255in}{1.023866in}%
\pgfsys@useobject{currentmarker}{}%
\end{pgfscope}%
\begin{pgfscope}%
\pgfsys@transformshift{1.592255in}{0.706295in}%
\pgfsys@useobject{currentmarker}{}%
\end{pgfscope}%
\begin{pgfscope}%
\pgfsys@transformshift{1.592255in}{0.267478in}%
\pgfsys@useobject{currentmarker}{}%
\end{pgfscope}%
\begin{pgfscope}%
\pgfsys@transformshift{1.794301in}{1.482500in}%
\pgfsys@useobject{currentmarker}{}%
\end{pgfscope}%
\begin{pgfscope}%
\pgfsys@transformshift{1.810081in}{1.285987in}%
\pgfsys@useobject{currentmarker}{}%
\end{pgfscope}%
\begin{pgfscope}%
\pgfsys@transformshift{1.859887in}{1.058777in}%
\pgfsys@useobject{currentmarker}{}%
\end{pgfscope}%
\begin{pgfscope}%
\pgfsys@transformshift{1.951498in}{0.765380in}%
\pgfsys@useobject{currentmarker}{}%
\end{pgfscope}%
\begin{pgfscope}%
\pgfsys@transformshift{2.099225in}{0.359966in}%
\pgfsys@useobject{currentmarker}{}%
\end{pgfscope}%
\begin{pgfscope}%
\pgfsys@transformshift{1.965588in}{1.482500in}%
\pgfsys@useobject{currentmarker}{}%
\end{pgfscope}%
\begin{pgfscope}%
\pgfsys@transformshift{1.994746in}{1.332095in}%
\pgfsys@useobject{currentmarker}{}%
\end{pgfscope}%
\begin{pgfscope}%
\pgfsys@transformshift{2.086775in}{1.158197in}%
\pgfsys@useobject{currentmarker}{}%
\end{pgfscope}%
\begin{pgfscope}%
\pgfsys@transformshift{2.256050in}{0.933640in}%
\pgfsys@useobject{currentmarker}{}%
\end{pgfscope}%
\begin{pgfscope}%
\pgfsys@transformshift{2.529013in}{0.623349in}%
\pgfsys@useobject{currentmarker}{}%
\end{pgfscope}%
\begin{pgfscope}%
\pgfsys@transformshift{2.080038in}{1.482500in}%
\pgfsys@useobject{currentmarker}{}%
\end{pgfscope}%
\begin{pgfscope}%
\pgfsys@transformshift{2.118135in}{1.401102in}%
\pgfsys@useobject{currentmarker}{}%
\end{pgfscope}%
\begin{pgfscope}%
\pgfsys@transformshift{2.238376in}{1.306988in}%
\pgfsys@useobject{currentmarker}{}%
\end{pgfscope}%
\begin{pgfscope}%
\pgfsys@transformshift{2.459545in}{1.185459in}%
\pgfsys@useobject{currentmarker}{}%
\end{pgfscope}%
\begin{pgfscope}%
\pgfsys@transformshift{2.816189in}{1.017531in}%
\pgfsys@useobject{currentmarker}{}%
\end{pgfscope}%
\begin{pgfscope}%
\pgfsys@transformshift{2.120227in}{1.482500in}%
\pgfsys@useobject{currentmarker}{}%
\end{pgfscope}%
\begin{pgfscope}%
\pgfsys@transformshift{2.161463in}{1.482500in}%
\pgfsys@useobject{currentmarker}{}%
\end{pgfscope}%
\begin{pgfscope}%
\pgfsys@transformshift{2.291612in}{1.482500in}%
\pgfsys@useobject{currentmarker}{}%
\end{pgfscope}%
\begin{pgfscope}%
\pgfsys@transformshift{2.531003in}{1.482500in}%
\pgfsys@useobject{currentmarker}{}%
\end{pgfscope}%
\begin{pgfscope}%
\pgfsys@transformshift{2.917032in}{1.482500in}%
\pgfsys@useobject{currentmarker}{}%
\end{pgfscope}%
\end{pgfscope}%
\end{pgfpicture}%
\makeatother%
\endgroup%

%% file: trigonometric-symmetric5_image.pgf
%% Creator: Matplotlib, PGF backend
%%
%% To include the figure in your LaTeX document, write
%%   \input{<filename>.pgf}
%%
%% Make sure the required packages are loaded in your preamble
%%   \usepackage{pgf}
%%
%% Figures using additional raster images can only be included by \input if
%% they are in the same directory as the main LaTeX file. For loading figures
%% from other directories you can use the `import` package
%% and then include the figures with
%%   \import{<path to file>}{<filename>.pgf}
%%
%% Matplotlib used the following preamble
%%   \usepackage{fontspec}
%%   \setmainfont{DejaVu Serif}
%%   \setsansfont{DejaVu Sans}
%%   \setmonofont{DejaVu Sans Mono}
%%
\begingroup%
\makeatletter%
\begin{pgfpicture}%
\pgfpathrectangle{\pgfpointorigin}{\pgfqpoint{3.316323in}{3.170273in}}%
\pgfusepath{use as bounding box, clip}%
\begin{pgfscope}%
\pgfsetbuttcap%
\pgfsetmiterjoin%
\definecolor{currentfill}{rgb}{1.000000,1.000000,1.000000}%
\pgfsetfillcolor{currentfill}%
\pgfsetlinewidth{0.000000pt}%
\definecolor{currentstroke}{rgb}{1.000000,1.000000,1.000000}%
\pgfsetstrokecolor{currentstroke}%
\pgfsetdash{}{0pt}%
\pgfpathmoveto{\pgfqpoint{0.000000in}{0.000000in}}%
\pgfpathlineto{\pgfqpoint{3.316323in}{0.000000in}}%
\pgfpathlineto{\pgfqpoint{3.316323in}{3.170273in}}%
\pgfpathlineto{\pgfqpoint{0.000000in}{3.170273in}}%
\pgfpathclose%
\pgfusepath{fill}%
\end{pgfscope}%
\begin{pgfscope}%
\pgfsetbuttcap%
\pgfsetmiterjoin%
\definecolor{currentfill}{rgb}{1.000000,1.000000,1.000000}%
\pgfsetfillcolor{currentfill}%
\pgfsetlinewidth{0.000000pt}%
\definecolor{currentstroke}{rgb}{0.000000,0.000000,0.000000}%
\pgfsetstrokecolor{currentstroke}%
\pgfsetstrokeopacity{0.000000}%
\pgfsetdash{}{0pt}%
\pgfpathmoveto{\pgfqpoint{0.135000in}{0.135000in}}%
\pgfpathlineto{\pgfqpoint{3.049509in}{0.135000in}}%
\pgfpathlineto{\pgfqpoint{3.049509in}{2.830000in}}%
\pgfpathlineto{\pgfqpoint{0.135000in}{2.830000in}}%
\pgfpathclose%
\pgfusepath{fill}%
\end{pgfscope}%
\begin{pgfscope}%
\pgfpathrectangle{\pgfqpoint{0.135000in}{0.135000in}}{\pgfqpoint{2.914509in}{2.695000in}} %
\pgfusepath{clip}%
\pgfsetrectcap%
\pgfsetroundjoin%
\pgfsetlinewidth{0.702625pt}%
\definecolor{currentstroke}{rgb}{0.000000,0.000000,0.000000}%
\pgfsetstrokecolor{currentstroke}%
\pgfsetdash{}{0pt}%
\pgfpathmoveto{\pgfqpoint{0.135000in}{1.482500in}}%
\pgfpathlineto{\pgfqpoint{3.049509in}{1.482500in}}%
\pgfusepath{stroke}%
\end{pgfscope}%
\begin{pgfscope}%
\pgfpathrectangle{\pgfqpoint{0.135000in}{0.135000in}}{\pgfqpoint{2.914509in}{2.695000in}} %
\pgfusepath{clip}%
\pgfsetrectcap%
\pgfsetroundjoin%
\pgfsetlinewidth{0.702625pt}%
\definecolor{currentstroke}{rgb}{0.000000,0.000000,0.000000}%
\pgfsetstrokecolor{currentstroke}%
\pgfsetdash{}{0pt}%
\pgfpathmoveto{\pgfqpoint{1.592255in}{0.135000in}}%
\pgfpathlineto{\pgfqpoint{1.592255in}{2.830000in}}%
\pgfusepath{stroke}%
\end{pgfscope}%
\begin{pgfscope}%
\pgfsetbuttcap%
\pgfsetmiterjoin%
\definecolor{currentfill}{rgb}{0.000000,0.000000,0.000000}%
\pgfsetfillcolor{currentfill}%
\pgfsetlinewidth{1.003750pt}%
\definecolor{currentstroke}{rgb}{0.000000,0.000000,0.000000}%
\pgfsetstrokecolor{currentstroke}%
\pgfsetdash{}{0pt}%
\pgfpathmoveto{\pgfqpoint{3.049509in}{1.482500in}}%
\pgfpathlineto{\pgfqpoint{3.020364in}{1.469025in}}%
\pgfpathlineto{\pgfqpoint{3.029108in}{1.482237in}}%
\pgfpathlineto{\pgfqpoint{0.135000in}{1.482237in}}%
\pgfpathlineto{\pgfqpoint{0.135000in}{1.482763in}}%
\pgfpathlineto{\pgfqpoint{3.029108in}{1.482763in}}%
\pgfpathlineto{\pgfqpoint{3.020364in}{1.495975in}}%
\pgfpathclose%
\pgfusepath{stroke,fill}%
\end{pgfscope}%
\begin{pgfscope}%
\pgfsetbuttcap%
\pgfsetmiterjoin%
\definecolor{currentfill}{rgb}{0.000000,0.000000,0.000000}%
\pgfsetfillcolor{currentfill}%
\pgfsetlinewidth{1.003750pt}%
\definecolor{currentstroke}{rgb}{0.000000,0.000000,0.000000}%
\pgfsetstrokecolor{currentstroke}%
\pgfsetdash{}{0pt}%
\pgfpathmoveto{\pgfqpoint{1.592255in}{2.830000in}}%
\pgfpathlineto{\pgfqpoint{1.605730in}{2.800855in}}%
\pgfpathlineto{\pgfqpoint{1.592518in}{2.809598in}}%
\pgfpathlineto{\pgfqpoint{1.592518in}{0.135000in}}%
\pgfpathlineto{\pgfqpoint{1.591991in}{0.135000in}}%
\pgfpathlineto{\pgfqpoint{1.591991in}{2.809598in}}%
\pgfpathlineto{\pgfqpoint{1.578780in}{2.800855in}}%
\pgfpathclose%
\pgfusepath{stroke,fill}%
\end{pgfscope}%
\begin{pgfscope}%
\pgfpathrectangle{\pgfqpoint{0.135000in}{0.135000in}}{\pgfqpoint{2.914509in}{2.695000in}} %
\pgfusepath{clip}%
\pgfsetrectcap%
\pgfsetroundjoin%
\pgfsetlinewidth{0.803000pt}%
\definecolor{currentstroke}{rgb}{0.000000,0.000000,0.000000}%
\pgfsetstrokecolor{currentstroke}%
\pgfsetdash{}{0pt}%
\pgfpathmoveto{\pgfqpoint{1.592255in}{0.135000in}}%
\pgfpathlineto{\pgfqpoint{1.592255in}{2.830000in}}%
\pgfusepath{stroke}%
\end{pgfscope}%
\begin{pgfscope}%
\pgfpathrectangle{\pgfqpoint{0.135000in}{0.135000in}}{\pgfqpoint{2.914509in}{2.695000in}} %
\pgfusepath{clip}%
\pgfsetrectcap%
\pgfsetroundjoin%
\pgfsetlinewidth{0.803000pt}%
\definecolor{currentstroke}{rgb}{0.000000,0.000000,0.000000}%
\pgfsetstrokecolor{currentstroke}%
\pgfsetdash{}{0pt}%
\pgfpathmoveto{\pgfqpoint{0.135000in}{1.482500in}}%
\pgfpathlineto{\pgfqpoint{3.049509in}{1.482500in}}%
\pgfusepath{stroke}%
\end{pgfscope}%
\begin{pgfscope}%
\pgfpathrectangle{\pgfqpoint{0.135000in}{0.135000in}}{\pgfqpoint{2.914509in}{2.695000in}} %
\pgfusepath{clip}%
\pgfsetrectcap%
\pgfsetroundjoin%
\pgfsetlinewidth{0.803000pt}%
\definecolor{currentstroke}{rgb}{0.000000,0.000000,0.000000}%
\pgfsetstrokecolor{currentstroke}%
\pgfsetdash{}{0pt}%
\pgfpathmoveto{\pgfqpoint{1.592255in}{0.135000in}}%
\pgfpathlineto{\pgfqpoint{1.592255in}{2.830000in}}%
\pgfusepath{stroke}%
\end{pgfscope}%
\begin{pgfscope}%
\pgfpathrectangle{\pgfqpoint{0.135000in}{0.135000in}}{\pgfqpoint{2.914509in}{2.695000in}} %
\pgfusepath{clip}%
\pgfsetrectcap%
\pgfsetroundjoin%
\pgfsetlinewidth{0.803000pt}%
\definecolor{currentstroke}{rgb}{0.000000,0.000000,0.000000}%
\pgfsetstrokecolor{currentstroke}%
\pgfsetdash{}{0pt}%
\pgfpathmoveto{\pgfqpoint{0.135000in}{1.482500in}}%
\pgfpathlineto{\pgfqpoint{3.049509in}{1.482500in}}%
\pgfusepath{stroke}%
\end{pgfscope}%
\begin{pgfscope}%
\pgftext[x=1.592255in,y=2.964750in,,base]{\sffamily\fontsize{10.000000}{12.000000}\selectfont \(\displaystyle y\)}%
\end{pgfscope}%
\begin{pgfscope}%
\pgftext[x=3.136945in,y=1.482500in,left,]{\sffamily\fontsize{10.000000}{12.000000}\selectfont \(\displaystyle x\)}%
\end{pgfscope}%
\begin{pgfscope}%
\pgfpathrectangle{\pgfqpoint{0.135000in}{0.135000in}}{\pgfqpoint{2.914509in}{2.695000in}} %
\pgfusepath{clip}%
\pgfsetbuttcap%
\pgfsetroundjoin%
\definecolor{currentfill}{rgb}{0.700000,0.700000,0.700000}%
\pgfsetfillcolor{currentfill}%
\pgfsetlinewidth{1.003750pt}%
\definecolor{currentstroke}{rgb}{0.700000,0.700000,0.700000}%
\pgfsetstrokecolor{currentstroke}%
\pgfsetdash{}{0pt}%
\pgfsys@defobject{currentmarker}{\pgfqpoint{-0.010417in}{-0.010417in}}{\pgfqpoint{0.010417in}{0.010417in}}{%
\pgfpathmoveto{\pgfqpoint{0.000000in}{-0.010417in}}%
\pgfpathcurveto{\pgfqpoint{0.002763in}{-0.010417in}}{\pgfqpoint{0.005412in}{-0.009319in}}{\pgfqpoint{0.007366in}{-0.007366in}}%
\pgfpathcurveto{\pgfqpoint{0.009319in}{-0.005412in}}{\pgfqpoint{0.010417in}{-0.002763in}}{\pgfqpoint{0.010417in}{0.000000in}}%
\pgfpathcurveto{\pgfqpoint{0.010417in}{0.002763in}}{\pgfqpoint{0.009319in}{0.005412in}}{\pgfqpoint{0.007366in}{0.007366in}}%
\pgfpathcurveto{\pgfqpoint{0.005412in}{0.009319in}}{\pgfqpoint{0.002763in}{0.010417in}}{\pgfqpoint{0.000000in}{0.010417in}}%
\pgfpathcurveto{\pgfqpoint{-0.002763in}{0.010417in}}{\pgfqpoint{-0.005412in}{0.009319in}}{\pgfqpoint{-0.007366in}{0.007366in}}%
\pgfpathcurveto{\pgfqpoint{-0.009319in}{0.005412in}}{\pgfqpoint{-0.010417in}{0.002763in}}{\pgfqpoint{-0.010417in}{0.000000in}}%
\pgfpathcurveto{\pgfqpoint{-0.010417in}{-0.002763in}}{\pgfqpoint{-0.009319in}{-0.005412in}}{\pgfqpoint{-0.007366in}{-0.007366in}}%
\pgfpathcurveto{\pgfqpoint{-0.005412in}{-0.009319in}}{\pgfqpoint{-0.002763in}{-0.010417in}}{\pgfqpoint{0.000000in}{-0.010417in}}%
\pgfpathclose%
\pgfusepath{stroke,fill}%
}%
\begin{pgfscope}%
\pgfsys@transformshift{2.120227in}{1.482500in}%
\pgfsys@useobject{currentmarker}{}%
\end{pgfscope}%
\begin{pgfscope}%
\pgfsys@transformshift{2.161463in}{1.482500in}%
\pgfsys@useobject{currentmarker}{}%
\end{pgfscope}%
\begin{pgfscope}%
\pgfsys@transformshift{2.291612in}{1.482500in}%
\pgfsys@useobject{currentmarker}{}%
\end{pgfscope}%
\begin{pgfscope}%
\pgfsys@transformshift{2.531003in}{1.482500in}%
\pgfsys@useobject{currentmarker}{}%
\end{pgfscope}%
\begin{pgfscope}%
\pgfsys@transformshift{2.917032in}{1.482500in}%
\pgfsys@useobject{currentmarker}{}%
\end{pgfscope}%
\begin{pgfscope}%
\pgfsys@transformshift{2.080038in}{1.482500in}%
\pgfsys@useobject{currentmarker}{}%
\end{pgfscope}%
\begin{pgfscope}%
\pgfsys@transformshift{2.118135in}{1.563898in}%
\pgfsys@useobject{currentmarker}{}%
\end{pgfscope}%
\begin{pgfscope}%
\pgfsys@transformshift{2.238376in}{1.658012in}%
\pgfsys@useobject{currentmarker}{}%
\end{pgfscope}%
\begin{pgfscope}%
\pgfsys@transformshift{2.459545in}{1.779541in}%
\pgfsys@useobject{currentmarker}{}%
\end{pgfscope}%
\begin{pgfscope}%
\pgfsys@transformshift{2.816189in}{1.947469in}%
\pgfsys@useobject{currentmarker}{}%
\end{pgfscope}%
\begin{pgfscope}%
\pgfsys@transformshift{1.965588in}{1.482500in}%
\pgfsys@useobject{currentmarker}{}%
\end{pgfscope}%
\begin{pgfscope}%
\pgfsys@transformshift{1.994746in}{1.632905in}%
\pgfsys@useobject{currentmarker}{}%
\end{pgfscope}%
\begin{pgfscope}%
\pgfsys@transformshift{2.086775in}{1.806803in}%
\pgfsys@useobject{currentmarker}{}%
\end{pgfscope}%
\begin{pgfscope}%
\pgfsys@transformshift{2.256050in}{2.031360in}%
\pgfsys@useobject{currentmarker}{}%
\end{pgfscope}%
\begin{pgfscope}%
\pgfsys@transformshift{2.529013in}{2.341651in}%
\pgfsys@useobject{currentmarker}{}%
\end{pgfscope}%
\begin{pgfscope}%
\pgfsys@transformshift{1.794301in}{1.482500in}%
\pgfsys@useobject{currentmarker}{}%
\end{pgfscope}%
\begin{pgfscope}%
\pgfsys@transformshift{1.810081in}{1.679013in}%
\pgfsys@useobject{currentmarker}{}%
\end{pgfscope}%
\begin{pgfscope}%
\pgfsys@transformshift{1.859887in}{1.906223in}%
\pgfsys@useobject{currentmarker}{}%
\end{pgfscope}%
\begin{pgfscope}%
\pgfsys@transformshift{1.951498in}{2.199620in}%
\pgfsys@useobject{currentmarker}{}%
\end{pgfscope}%
\begin{pgfscope}%
\pgfsys@transformshift{2.099225in}{2.605034in}%
\pgfsys@useobject{currentmarker}{}%
\end{pgfscope}%
\begin{pgfscope}%
\pgfsys@transformshift{1.592255in}{1.482500in}%
\pgfsys@useobject{currentmarker}{}%
\end{pgfscope}%
\begin{pgfscope}%
\pgfsys@transformshift{1.592255in}{1.695204in}%
\pgfsys@useobject{currentmarker}{}%
\end{pgfscope}%
\begin{pgfscope}%
\pgfsys@transformshift{1.592255in}{1.941134in}%
\pgfsys@useobject{currentmarker}{}%
\end{pgfscope}%
\begin{pgfscope}%
\pgfsys@transformshift{1.592255in}{2.258705in}%
\pgfsys@useobject{currentmarker}{}%
\end{pgfscope}%
\begin{pgfscope}%
\pgfsys@transformshift{1.592255in}{2.697522in}%
\pgfsys@useobject{currentmarker}{}%
\end{pgfscope}%
\begin{pgfscope}%
\pgfsys@transformshift{1.390208in}{1.482500in}%
\pgfsys@useobject{currentmarker}{}%
\end{pgfscope}%
\begin{pgfscope}%
\pgfsys@transformshift{1.374428in}{1.679013in}%
\pgfsys@useobject{currentmarker}{}%
\end{pgfscope}%
\begin{pgfscope}%
\pgfsys@transformshift{1.324622in}{1.906223in}%
\pgfsys@useobject{currentmarker}{}%
\end{pgfscope}%
\begin{pgfscope}%
\pgfsys@transformshift{1.233011in}{2.199620in}%
\pgfsys@useobject{currentmarker}{}%
\end{pgfscope}%
\begin{pgfscope}%
\pgfsys@transformshift{1.085284in}{2.605034in}%
\pgfsys@useobject{currentmarker}{}%
\end{pgfscope}%
\begin{pgfscope}%
\pgfsys@transformshift{1.218922in}{1.482500in}%
\pgfsys@useobject{currentmarker}{}%
\end{pgfscope}%
\begin{pgfscope}%
\pgfsys@transformshift{1.189764in}{1.632905in}%
\pgfsys@useobject{currentmarker}{}%
\end{pgfscope}%
\begin{pgfscope}%
\pgfsys@transformshift{1.097735in}{1.806803in}%
\pgfsys@useobject{currentmarker}{}%
\end{pgfscope}%
\begin{pgfscope}%
\pgfsys@transformshift{0.928459in}{2.031360in}%
\pgfsys@useobject{currentmarker}{}%
\end{pgfscope}%
\begin{pgfscope}%
\pgfsys@transformshift{0.655496in}{2.341651in}%
\pgfsys@useobject{currentmarker}{}%
\end{pgfscope}%
\begin{pgfscope}%
\pgfsys@transformshift{1.104472in}{1.482500in}%
\pgfsys@useobject{currentmarker}{}%
\end{pgfscope}%
\begin{pgfscope}%
\pgfsys@transformshift{1.066375in}{1.563898in}%
\pgfsys@useobject{currentmarker}{}%
\end{pgfscope}%
\begin{pgfscope}%
\pgfsys@transformshift{0.946133in}{1.658012in}%
\pgfsys@useobject{currentmarker}{}%
\end{pgfscope}%
\begin{pgfscope}%
\pgfsys@transformshift{0.724964in}{1.779541in}%
\pgfsys@useobject{currentmarker}{}%
\end{pgfscope}%
\begin{pgfscope}%
\pgfsys@transformshift{0.368320in}{1.947469in}%
\pgfsys@useobject{currentmarker}{}%
\end{pgfscope}%
\begin{pgfscope}%
\pgfsys@transformshift{1.064282in}{1.482500in}%
\pgfsys@useobject{currentmarker}{}%
\end{pgfscope}%
\begin{pgfscope}%
\pgfsys@transformshift{1.023046in}{1.482500in}%
\pgfsys@useobject{currentmarker}{}%
\end{pgfscope}%
\begin{pgfscope}%
\pgfsys@transformshift{0.892898in}{1.482500in}%
\pgfsys@useobject{currentmarker}{}%
\end{pgfscope}%
\begin{pgfscope}%
\pgfsys@transformshift{0.653506in}{1.482500in}%
\pgfsys@useobject{currentmarker}{}%
\end{pgfscope}%
\begin{pgfscope}%
\pgfsys@transformshift{0.267478in}{1.482500in}%
\pgfsys@useobject{currentmarker}{}%
\end{pgfscope}%
\begin{pgfscope}%
\pgfsys@transformshift{1.104472in}{1.482500in}%
\pgfsys@useobject{currentmarker}{}%
\end{pgfscope}%
\begin{pgfscope}%
\pgfsys@transformshift{1.066375in}{1.401102in}%
\pgfsys@useobject{currentmarker}{}%
\end{pgfscope}%
\begin{pgfscope}%
\pgfsys@transformshift{0.946133in}{1.306988in}%
\pgfsys@useobject{currentmarker}{}%
\end{pgfscope}%
\begin{pgfscope}%
\pgfsys@transformshift{0.724964in}{1.185459in}%
\pgfsys@useobject{currentmarker}{}%
\end{pgfscope}%
\begin{pgfscope}%
\pgfsys@transformshift{0.368320in}{1.017531in}%
\pgfsys@useobject{currentmarker}{}%
\end{pgfscope}%
\begin{pgfscope}%
\pgfsys@transformshift{1.218922in}{1.482500in}%
\pgfsys@useobject{currentmarker}{}%
\end{pgfscope}%
\begin{pgfscope}%
\pgfsys@transformshift{1.189764in}{1.332095in}%
\pgfsys@useobject{currentmarker}{}%
\end{pgfscope}%
\begin{pgfscope}%
\pgfsys@transformshift{1.097735in}{1.158197in}%
\pgfsys@useobject{currentmarker}{}%
\end{pgfscope}%
\begin{pgfscope}%
\pgfsys@transformshift{0.928459in}{0.933640in}%
\pgfsys@useobject{currentmarker}{}%
\end{pgfscope}%
\begin{pgfscope}%
\pgfsys@transformshift{0.655496in}{0.623349in}%
\pgfsys@useobject{currentmarker}{}%
\end{pgfscope}%
\begin{pgfscope}%
\pgfsys@transformshift{1.390208in}{1.482500in}%
\pgfsys@useobject{currentmarker}{}%
\end{pgfscope}%
\begin{pgfscope}%
\pgfsys@transformshift{1.374428in}{1.285987in}%
\pgfsys@useobject{currentmarker}{}%
\end{pgfscope}%
\begin{pgfscope}%
\pgfsys@transformshift{1.324622in}{1.058777in}%
\pgfsys@useobject{currentmarker}{}%
\end{pgfscope}%
\begin{pgfscope}%
\pgfsys@transformshift{1.233011in}{0.765380in}%
\pgfsys@useobject{currentmarker}{}%
\end{pgfscope}%
\begin{pgfscope}%
\pgfsys@transformshift{1.085284in}{0.359966in}%
\pgfsys@useobject{currentmarker}{}%
\end{pgfscope}%
\begin{pgfscope}%
\pgfsys@transformshift{1.592255in}{1.482500in}%
\pgfsys@useobject{currentmarker}{}%
\end{pgfscope}%
\begin{pgfscope}%
\pgfsys@transformshift{1.592255in}{1.269796in}%
\pgfsys@useobject{currentmarker}{}%
\end{pgfscope}%
\begin{pgfscope}%
\pgfsys@transformshift{1.592255in}{1.023866in}%
\pgfsys@useobject{currentmarker}{}%
\end{pgfscope}%
\begin{pgfscope}%
\pgfsys@transformshift{1.592255in}{0.706295in}%
\pgfsys@useobject{currentmarker}{}%
\end{pgfscope}%
\begin{pgfscope}%
\pgfsys@transformshift{1.592255in}{0.267478in}%
\pgfsys@useobject{currentmarker}{}%
\end{pgfscope}%
\begin{pgfscope}%
\pgfsys@transformshift{1.794301in}{1.482500in}%
\pgfsys@useobject{currentmarker}{}%
\end{pgfscope}%
\begin{pgfscope}%
\pgfsys@transformshift{1.810081in}{1.285987in}%
\pgfsys@useobject{currentmarker}{}%
\end{pgfscope}%
\begin{pgfscope}%
\pgfsys@transformshift{1.859887in}{1.058777in}%
\pgfsys@useobject{currentmarker}{}%
\end{pgfscope}%
\begin{pgfscope}%
\pgfsys@transformshift{1.951498in}{0.765380in}%
\pgfsys@useobject{currentmarker}{}%
\end{pgfscope}%
\begin{pgfscope}%
\pgfsys@transformshift{2.099225in}{0.359966in}%
\pgfsys@useobject{currentmarker}{}%
\end{pgfscope}%
\begin{pgfscope}%
\pgfsys@transformshift{1.965588in}{1.482500in}%
\pgfsys@useobject{currentmarker}{}%
\end{pgfscope}%
\begin{pgfscope}%
\pgfsys@transformshift{1.994746in}{1.332095in}%
\pgfsys@useobject{currentmarker}{}%
\end{pgfscope}%
\begin{pgfscope}%
\pgfsys@transformshift{2.086775in}{1.158197in}%
\pgfsys@useobject{currentmarker}{}%
\end{pgfscope}%
\begin{pgfscope}%
\pgfsys@transformshift{2.256050in}{0.933640in}%
\pgfsys@useobject{currentmarker}{}%
\end{pgfscope}%
\begin{pgfscope}%
\pgfsys@transformshift{2.529013in}{0.623349in}%
\pgfsys@useobject{currentmarker}{}%
\end{pgfscope}%
\begin{pgfscope}%
\pgfsys@transformshift{2.080038in}{1.482500in}%
\pgfsys@useobject{currentmarker}{}%
\end{pgfscope}%
\begin{pgfscope}%
\pgfsys@transformshift{2.118135in}{1.401102in}%
\pgfsys@useobject{currentmarker}{}%
\end{pgfscope}%
\begin{pgfscope}%
\pgfsys@transformshift{2.238376in}{1.306988in}%
\pgfsys@useobject{currentmarker}{}%
\end{pgfscope}%
\begin{pgfscope}%
\pgfsys@transformshift{2.459545in}{1.185459in}%
\pgfsys@useobject{currentmarker}{}%
\end{pgfscope}%
\begin{pgfscope}%
\pgfsys@transformshift{2.816189in}{1.017531in}%
\pgfsys@useobject{currentmarker}{}%
\end{pgfscope}%
\begin{pgfscope}%
\pgfsys@transformshift{2.120227in}{1.482500in}%
\pgfsys@useobject{currentmarker}{}%
\end{pgfscope}%
\begin{pgfscope}%
\pgfsys@transformshift{2.161463in}{1.482500in}%
\pgfsys@useobject{currentmarker}{}%
\end{pgfscope}%
\begin{pgfscope}%
\pgfsys@transformshift{2.291612in}{1.482500in}%
\pgfsys@useobject{currentmarker}{}%
\end{pgfscope}%
\begin{pgfscope}%
\pgfsys@transformshift{2.531003in}{1.482500in}%
\pgfsys@useobject{currentmarker}{}%
\end{pgfscope}%
\begin{pgfscope}%
\pgfsys@transformshift{2.917032in}{1.482500in}%
\pgfsys@useobject{currentmarker}{}%
\end{pgfscope}%
\end{pgfscope}%
\begin{pgfscope}%
\pgfpathrectangle{\pgfqpoint{0.135000in}{0.135000in}}{\pgfqpoint{2.914509in}{2.695000in}} %
\pgfusepath{clip}%
\pgfsetrectcap%
\pgfsetroundjoin%
\pgfsetlinewidth{1.003750pt}%
\definecolor{currentstroke}{rgb}{0.000000,0.000000,1.000000}%
\pgfsetstrokecolor{currentstroke}%
\pgfsetdash{}{0pt}%
\pgfpathmoveto{\pgfqpoint{2.120227in}{1.482500in}}%
\pgfpathlineto{\pgfqpoint{1.965588in}{1.482500in}}%
\pgfpathlineto{\pgfqpoint{1.592255in}{1.482500in}}%
\pgfpathlineto{\pgfqpoint{1.218922in}{1.482500in}}%
\pgfpathlineto{\pgfqpoint{1.064282in}{1.482500in}}%
\pgfpathlineto{\pgfqpoint{1.218922in}{1.482500in}}%
\pgfpathlineto{\pgfqpoint{1.592255in}{1.482500in}}%
\pgfpathlineto{\pgfqpoint{1.965588in}{1.482500in}}%
\pgfpathlineto{\pgfqpoint{2.120227in}{1.482500in}}%
\pgfusepath{stroke}%
\end{pgfscope}%
\begin{pgfscope}%
\pgfpathrectangle{\pgfqpoint{0.135000in}{0.135000in}}{\pgfqpoint{2.914509in}{2.695000in}} %
\pgfusepath{clip}%
\pgfsetrectcap%
\pgfsetroundjoin%
\pgfsetlinewidth{1.003750pt}%
\definecolor{currentstroke}{rgb}{0.000000,0.000000,1.000000}%
\pgfsetstrokecolor{currentstroke}%
\pgfsetdash{}{0pt}%
\pgfpathmoveto{\pgfqpoint{2.291612in}{1.482500in}}%
\pgfpathlineto{\pgfqpoint{2.086775in}{1.806803in}}%
\pgfpathlineto{\pgfqpoint{1.592255in}{1.941134in}}%
\pgfpathlineto{\pgfqpoint{1.097735in}{1.806803in}}%
\pgfpathlineto{\pgfqpoint{0.892898in}{1.482500in}}%
\pgfpathlineto{\pgfqpoint{1.097735in}{1.158197in}}%
\pgfpathlineto{\pgfqpoint{1.592255in}{1.023866in}}%
\pgfpathlineto{\pgfqpoint{2.086775in}{1.158197in}}%
\pgfpathlineto{\pgfqpoint{2.291612in}{1.482500in}}%
\pgfusepath{stroke}%
\end{pgfscope}%
\begin{pgfscope}%
\pgfpathrectangle{\pgfqpoint{0.135000in}{0.135000in}}{\pgfqpoint{2.914509in}{2.695000in}} %
\pgfusepath{clip}%
\pgfsetrectcap%
\pgfsetroundjoin%
\pgfsetlinewidth{1.003750pt}%
\definecolor{currentstroke}{rgb}{0.000000,0.000000,1.000000}%
\pgfsetstrokecolor{currentstroke}%
\pgfsetdash{}{0pt}%
\pgfpathmoveto{\pgfqpoint{2.917032in}{1.482500in}}%
\pgfpathlineto{\pgfqpoint{2.529013in}{2.341651in}}%
\pgfpathlineto{\pgfqpoint{1.592255in}{2.697522in}}%
\pgfpathlineto{\pgfqpoint{0.655496in}{2.341651in}}%
\pgfpathlineto{\pgfqpoint{0.267478in}{1.482500in}}%
\pgfpathlineto{\pgfqpoint{0.655496in}{0.623349in}}%
\pgfpathlineto{\pgfqpoint{1.592255in}{0.267478in}}%
\pgfpathlineto{\pgfqpoint{2.529013in}{0.623349in}}%
\pgfpathlineto{\pgfqpoint{2.917032in}{1.482500in}}%
\pgfusepath{stroke}%
\end{pgfscope}%
\begin{pgfscope}%
\pgfpathrectangle{\pgfqpoint{0.135000in}{0.135000in}}{\pgfqpoint{2.914509in}{2.695000in}} %
\pgfusepath{clip}%
\pgfsetrectcap%
\pgfsetroundjoin%
\pgfsetlinewidth{1.003750pt}%
\definecolor{currentstroke}{rgb}{0.000000,0.000000,1.000000}%
\pgfsetstrokecolor{currentstroke}%
\pgfsetdash{}{0pt}%
\pgfpathmoveto{\pgfqpoint{2.120227in}{1.482500in}}%
\pgfpathlineto{\pgfqpoint{2.291612in}{1.482500in}}%
\pgfpathlineto{\pgfqpoint{2.917032in}{1.482500in}}%
\pgfusepath{stroke}%
\end{pgfscope}%
\begin{pgfscope}%
\pgfpathrectangle{\pgfqpoint{0.135000in}{0.135000in}}{\pgfqpoint{2.914509in}{2.695000in}} %
\pgfusepath{clip}%
\pgfsetrectcap%
\pgfsetroundjoin%
\pgfsetlinewidth{1.003750pt}%
\definecolor{currentstroke}{rgb}{0.000000,0.000000,1.000000}%
\pgfsetstrokecolor{currentstroke}%
\pgfsetdash{}{0pt}%
\pgfpathmoveto{\pgfqpoint{1.965588in}{1.482500in}}%
\pgfpathlineto{\pgfqpoint{2.086775in}{1.806803in}}%
\pgfpathlineto{\pgfqpoint{2.529013in}{2.341651in}}%
\pgfusepath{stroke}%
\end{pgfscope}%
\begin{pgfscope}%
\pgfpathrectangle{\pgfqpoint{0.135000in}{0.135000in}}{\pgfqpoint{2.914509in}{2.695000in}} %
\pgfusepath{clip}%
\pgfsetrectcap%
\pgfsetroundjoin%
\pgfsetlinewidth{1.003750pt}%
\definecolor{currentstroke}{rgb}{0.000000,0.000000,1.000000}%
\pgfsetstrokecolor{currentstroke}%
\pgfsetdash{}{0pt}%
\pgfpathmoveto{\pgfqpoint{1.592255in}{1.482500in}}%
\pgfpathlineto{\pgfqpoint{1.592255in}{1.941134in}}%
\pgfpathlineto{\pgfqpoint{1.592255in}{2.697522in}}%
\pgfusepath{stroke}%
\end{pgfscope}%
\begin{pgfscope}%
\pgfpathrectangle{\pgfqpoint{0.135000in}{0.135000in}}{\pgfqpoint{2.914509in}{2.695000in}} %
\pgfusepath{clip}%
\pgfsetrectcap%
\pgfsetroundjoin%
\pgfsetlinewidth{1.003750pt}%
\definecolor{currentstroke}{rgb}{0.000000,0.000000,1.000000}%
\pgfsetstrokecolor{currentstroke}%
\pgfsetdash{}{0pt}%
\pgfpathmoveto{\pgfqpoint{1.218922in}{1.482500in}}%
\pgfpathlineto{\pgfqpoint{1.097735in}{1.806803in}}%
\pgfpathlineto{\pgfqpoint{0.655496in}{2.341651in}}%
\pgfusepath{stroke}%
\end{pgfscope}%
\begin{pgfscope}%
\pgfpathrectangle{\pgfqpoint{0.135000in}{0.135000in}}{\pgfqpoint{2.914509in}{2.695000in}} %
\pgfusepath{clip}%
\pgfsetrectcap%
\pgfsetroundjoin%
\pgfsetlinewidth{1.003750pt}%
\definecolor{currentstroke}{rgb}{0.000000,0.000000,1.000000}%
\pgfsetstrokecolor{currentstroke}%
\pgfsetdash{}{0pt}%
\pgfpathmoveto{\pgfqpoint{1.064282in}{1.482500in}}%
\pgfpathlineto{\pgfqpoint{0.892898in}{1.482500in}}%
\pgfpathlineto{\pgfqpoint{0.267478in}{1.482500in}}%
\pgfusepath{stroke}%
\end{pgfscope}%
\begin{pgfscope}%
\pgfpathrectangle{\pgfqpoint{0.135000in}{0.135000in}}{\pgfqpoint{2.914509in}{2.695000in}} %
\pgfusepath{clip}%
\pgfsetrectcap%
\pgfsetroundjoin%
\pgfsetlinewidth{1.003750pt}%
\definecolor{currentstroke}{rgb}{0.000000,0.000000,1.000000}%
\pgfsetstrokecolor{currentstroke}%
\pgfsetdash{}{0pt}%
\pgfpathmoveto{\pgfqpoint{1.218922in}{1.482500in}}%
\pgfpathlineto{\pgfqpoint{1.097735in}{1.158197in}}%
\pgfpathlineto{\pgfqpoint{0.655496in}{0.623349in}}%
\pgfusepath{stroke}%
\end{pgfscope}%
\begin{pgfscope}%
\pgfpathrectangle{\pgfqpoint{0.135000in}{0.135000in}}{\pgfqpoint{2.914509in}{2.695000in}} %
\pgfusepath{clip}%
\pgfsetrectcap%
\pgfsetroundjoin%
\pgfsetlinewidth{1.003750pt}%
\definecolor{currentstroke}{rgb}{0.000000,0.000000,1.000000}%
\pgfsetstrokecolor{currentstroke}%
\pgfsetdash{}{0pt}%
\pgfpathmoveto{\pgfqpoint{1.592255in}{1.482500in}}%
\pgfpathlineto{\pgfqpoint{1.592255in}{1.023866in}}%
\pgfpathlineto{\pgfqpoint{1.592255in}{0.267478in}}%
\pgfusepath{stroke}%
\end{pgfscope}%
\begin{pgfscope}%
\pgfpathrectangle{\pgfqpoint{0.135000in}{0.135000in}}{\pgfqpoint{2.914509in}{2.695000in}} %
\pgfusepath{clip}%
\pgfsetrectcap%
\pgfsetroundjoin%
\pgfsetlinewidth{1.003750pt}%
\definecolor{currentstroke}{rgb}{0.000000,0.000000,1.000000}%
\pgfsetstrokecolor{currentstroke}%
\pgfsetdash{}{0pt}%
\pgfpathmoveto{\pgfqpoint{1.965588in}{1.482500in}}%
\pgfpathlineto{\pgfqpoint{2.086775in}{1.158197in}}%
\pgfpathlineto{\pgfqpoint{2.529013in}{0.623349in}}%
\pgfusepath{stroke}%
\end{pgfscope}%
\begin{pgfscope}%
\pgfpathrectangle{\pgfqpoint{0.135000in}{0.135000in}}{\pgfqpoint{2.914509in}{2.695000in}} %
\pgfusepath{clip}%
\pgfsetrectcap%
\pgfsetroundjoin%
\pgfsetlinewidth{1.003750pt}%
\definecolor{currentstroke}{rgb}{0.000000,0.000000,1.000000}%
\pgfsetstrokecolor{currentstroke}%
\pgfsetdash{}{0pt}%
\pgfpathmoveto{\pgfqpoint{2.120227in}{1.482500in}}%
\pgfpathlineto{\pgfqpoint{2.291612in}{1.482500in}}%
\pgfpathlineto{\pgfqpoint{2.917032in}{1.482500in}}%
\pgfusepath{stroke}%
\end{pgfscope}%
\begin{pgfscope}%
\pgfpathrectangle{\pgfqpoint{0.135000in}{0.135000in}}{\pgfqpoint{2.914509in}{2.695000in}} %
\pgfusepath{clip}%
\pgfsetbuttcap%
\pgfsetroundjoin%
\definecolor{currentfill}{rgb}{0.000000,0.000000,1.000000}%
\pgfsetfillcolor{currentfill}%
\pgfsetlinewidth{1.003750pt}%
\definecolor{currentstroke}{rgb}{0.000000,0.000000,1.000000}%
\pgfsetstrokecolor{currentstroke}%
\pgfsetdash{}{0pt}%
\pgfsys@defobject{currentmarker}{\pgfqpoint{-0.010417in}{-0.010417in}}{\pgfqpoint{0.010417in}{0.010417in}}{%
\pgfpathmoveto{\pgfqpoint{0.000000in}{-0.010417in}}%
\pgfpathcurveto{\pgfqpoint{0.002763in}{-0.010417in}}{\pgfqpoint{0.005412in}{-0.009319in}}{\pgfqpoint{0.007366in}{-0.007366in}}%
\pgfpathcurveto{\pgfqpoint{0.009319in}{-0.005412in}}{\pgfqpoint{0.010417in}{-0.002763in}}{\pgfqpoint{0.010417in}{0.000000in}}%
\pgfpathcurveto{\pgfqpoint{0.010417in}{0.002763in}}{\pgfqpoint{0.009319in}{0.005412in}}{\pgfqpoint{0.007366in}{0.007366in}}%
\pgfpathcurveto{\pgfqpoint{0.005412in}{0.009319in}}{\pgfqpoint{0.002763in}{0.010417in}}{\pgfqpoint{0.000000in}{0.010417in}}%
\pgfpathcurveto{\pgfqpoint{-0.002763in}{0.010417in}}{\pgfqpoint{-0.005412in}{0.009319in}}{\pgfqpoint{-0.007366in}{0.007366in}}%
\pgfpathcurveto{\pgfqpoint{-0.009319in}{0.005412in}}{\pgfqpoint{-0.010417in}{0.002763in}}{\pgfqpoint{-0.010417in}{0.000000in}}%
\pgfpathcurveto{\pgfqpoint{-0.010417in}{-0.002763in}}{\pgfqpoint{-0.009319in}{-0.005412in}}{\pgfqpoint{-0.007366in}{-0.007366in}}%
\pgfpathcurveto{\pgfqpoint{-0.005412in}{-0.009319in}}{\pgfqpoint{-0.002763in}{-0.010417in}}{\pgfqpoint{0.000000in}{-0.010417in}}%
\pgfpathclose%
\pgfusepath{stroke,fill}%
}%
\begin{pgfscope}%
\pgfsys@transformshift{2.120227in}{1.482500in}%
\pgfsys@useobject{currentmarker}{}%
\end{pgfscope}%
\begin{pgfscope}%
\pgfsys@transformshift{2.291612in}{1.482500in}%
\pgfsys@useobject{currentmarker}{}%
\end{pgfscope}%
\begin{pgfscope}%
\pgfsys@transformshift{2.917032in}{1.482500in}%
\pgfsys@useobject{currentmarker}{}%
\end{pgfscope}%
\begin{pgfscope}%
\pgfsys@transformshift{1.965588in}{1.482500in}%
\pgfsys@useobject{currentmarker}{}%
\end{pgfscope}%
\begin{pgfscope}%
\pgfsys@transformshift{2.086775in}{1.806803in}%
\pgfsys@useobject{currentmarker}{}%
\end{pgfscope}%
\begin{pgfscope}%
\pgfsys@transformshift{2.529013in}{2.341651in}%
\pgfsys@useobject{currentmarker}{}%
\end{pgfscope}%
\begin{pgfscope}%
\pgfsys@transformshift{1.592255in}{1.482500in}%
\pgfsys@useobject{currentmarker}{}%
\end{pgfscope}%
\begin{pgfscope}%
\pgfsys@transformshift{1.592255in}{1.941134in}%
\pgfsys@useobject{currentmarker}{}%
\end{pgfscope}%
\begin{pgfscope}%
\pgfsys@transformshift{1.592255in}{2.697522in}%
\pgfsys@useobject{currentmarker}{}%
\end{pgfscope}%
\begin{pgfscope}%
\pgfsys@transformshift{1.218922in}{1.482500in}%
\pgfsys@useobject{currentmarker}{}%
\end{pgfscope}%
\begin{pgfscope}%
\pgfsys@transformshift{1.097735in}{1.806803in}%
\pgfsys@useobject{currentmarker}{}%
\end{pgfscope}%
\begin{pgfscope}%
\pgfsys@transformshift{0.655496in}{2.341651in}%
\pgfsys@useobject{currentmarker}{}%
\end{pgfscope}%
\begin{pgfscope}%
\pgfsys@transformshift{1.064282in}{1.482500in}%
\pgfsys@useobject{currentmarker}{}%
\end{pgfscope}%
\begin{pgfscope}%
\pgfsys@transformshift{0.892898in}{1.482500in}%
\pgfsys@useobject{currentmarker}{}%
\end{pgfscope}%
\begin{pgfscope}%
\pgfsys@transformshift{0.267478in}{1.482500in}%
\pgfsys@useobject{currentmarker}{}%
\end{pgfscope}%
\begin{pgfscope}%
\pgfsys@transformshift{1.218922in}{1.482500in}%
\pgfsys@useobject{currentmarker}{}%
\end{pgfscope}%
\begin{pgfscope}%
\pgfsys@transformshift{1.097735in}{1.158197in}%
\pgfsys@useobject{currentmarker}{}%
\end{pgfscope}%
\begin{pgfscope}%
\pgfsys@transformshift{0.655496in}{0.623349in}%
\pgfsys@useobject{currentmarker}{}%
\end{pgfscope}%
\begin{pgfscope}%
\pgfsys@transformshift{1.592255in}{1.482500in}%
\pgfsys@useobject{currentmarker}{}%
\end{pgfscope}%
\begin{pgfscope}%
\pgfsys@transformshift{1.592255in}{1.023866in}%
\pgfsys@useobject{currentmarker}{}%
\end{pgfscope}%
\begin{pgfscope}%
\pgfsys@transformshift{1.592255in}{0.267478in}%
\pgfsys@useobject{currentmarker}{}%
\end{pgfscope}%
\begin{pgfscope}%
\pgfsys@transformshift{1.965588in}{1.482500in}%
\pgfsys@useobject{currentmarker}{}%
\end{pgfscope}%
\begin{pgfscope}%
\pgfsys@transformshift{2.086775in}{1.158197in}%
\pgfsys@useobject{currentmarker}{}%
\end{pgfscope}%
\begin{pgfscope}%
\pgfsys@transformshift{2.529013in}{0.623349in}%
\pgfsys@useobject{currentmarker}{}%
\end{pgfscope}%
\begin{pgfscope}%
\pgfsys@transformshift{2.120227in}{1.482500in}%
\pgfsys@useobject{currentmarker}{}%
\end{pgfscope}%
\begin{pgfscope}%
\pgfsys@transformshift{2.291612in}{1.482500in}%
\pgfsys@useobject{currentmarker}{}%
\end{pgfscope}%
\begin{pgfscope}%
\pgfsys@transformshift{2.917032in}{1.482500in}%
\pgfsys@useobject{currentmarker}{}%
\end{pgfscope}%
\end{pgfscope}%
\begin{pgfscope}%
\pgfpathrectangle{\pgfqpoint{0.135000in}{0.135000in}}{\pgfqpoint{2.914509in}{2.695000in}} %
\pgfusepath{clip}%
\pgfsetrectcap%
\pgfsetroundjoin%
\pgfsetlinewidth{1.003750pt}%
\definecolor{currentstroke}{rgb}{1.000000,0.000000,0.000000}%
\pgfsetstrokecolor{currentstroke}%
\pgfsetdash{}{0pt}%
\pgfpathmoveto{\pgfqpoint{2.118135in}{1.563898in}}%
\pgfpathlineto{\pgfqpoint{2.118135in}{1.401102in}}%
\pgfpathlineto{\pgfqpoint{1.810081in}{1.285987in}}%
\pgfpathlineto{\pgfqpoint{1.374428in}{1.285987in}}%
\pgfpathlineto{\pgfqpoint{1.066375in}{1.401102in}}%
\pgfpathlineto{\pgfqpoint{1.066375in}{1.563898in}}%
\pgfpathlineto{\pgfqpoint{1.374428in}{1.679013in}}%
\pgfpathlineto{\pgfqpoint{1.810081in}{1.679013in}}%
\pgfpathlineto{\pgfqpoint{2.118135in}{1.563898in}}%
\pgfusepath{stroke}%
\end{pgfscope}%
\begin{pgfscope}%
\pgfpathrectangle{\pgfqpoint{0.135000in}{0.135000in}}{\pgfqpoint{2.914509in}{2.695000in}} %
\pgfusepath{clip}%
\pgfsetrectcap%
\pgfsetroundjoin%
\pgfsetlinewidth{1.003750pt}%
\definecolor{currentstroke}{rgb}{1.000000,0.000000,0.000000}%
\pgfsetstrokecolor{currentstroke}%
\pgfsetdash{}{0pt}%
\pgfpathmoveto{\pgfqpoint{2.118135in}{1.401102in}}%
\pgfpathlineto{\pgfqpoint{2.118135in}{1.563898in}}%
\pgfpathlineto{\pgfqpoint{1.810081in}{1.679013in}}%
\pgfpathlineto{\pgfqpoint{1.374428in}{1.679013in}}%
\pgfpathlineto{\pgfqpoint{1.066375in}{1.563898in}}%
\pgfpathlineto{\pgfqpoint{1.066375in}{1.401102in}}%
\pgfpathlineto{\pgfqpoint{1.374428in}{1.285987in}}%
\pgfpathlineto{\pgfqpoint{1.810081in}{1.285987in}}%
\pgfpathlineto{\pgfqpoint{2.118135in}{1.401102in}}%
\pgfusepath{stroke}%
\end{pgfscope}%
\begin{pgfscope}%
\pgfpathrectangle{\pgfqpoint{0.135000in}{0.135000in}}{\pgfqpoint{2.914509in}{2.695000in}} %
\pgfusepath{clip}%
\pgfsetrectcap%
\pgfsetroundjoin%
\pgfsetlinewidth{1.003750pt}%
\definecolor{currentstroke}{rgb}{1.000000,0.000000,0.000000}%
\pgfsetstrokecolor{currentstroke}%
\pgfsetdash{}{0pt}%
\pgfpathmoveto{\pgfqpoint{2.459545in}{1.185459in}}%
\pgfpathlineto{\pgfqpoint{2.459545in}{1.779541in}}%
\pgfpathlineto{\pgfqpoint{1.951498in}{2.199620in}}%
\pgfpathlineto{\pgfqpoint{1.233011in}{2.199620in}}%
\pgfpathlineto{\pgfqpoint{0.724964in}{1.779541in}}%
\pgfpathlineto{\pgfqpoint{0.724964in}{1.185459in}}%
\pgfpathlineto{\pgfqpoint{1.233011in}{0.765380in}}%
\pgfpathlineto{\pgfqpoint{1.951498in}{0.765380in}}%
\pgfpathlineto{\pgfqpoint{2.459545in}{1.185459in}}%
\pgfusepath{stroke}%
\end{pgfscope}%
\begin{pgfscope}%
\pgfpathrectangle{\pgfqpoint{0.135000in}{0.135000in}}{\pgfqpoint{2.914509in}{2.695000in}} %
\pgfusepath{clip}%
\pgfsetrectcap%
\pgfsetroundjoin%
\pgfsetlinewidth{1.003750pt}%
\definecolor{currentstroke}{rgb}{1.000000,0.000000,0.000000}%
\pgfsetstrokecolor{currentstroke}%
\pgfsetdash{}{0pt}%
\pgfpathmoveto{\pgfqpoint{2.118135in}{1.563898in}}%
\pgfpathlineto{\pgfqpoint{2.118135in}{1.401102in}}%
\pgfpathlineto{\pgfqpoint{2.459545in}{1.185459in}}%
\pgfusepath{stroke}%
\end{pgfscope}%
\begin{pgfscope}%
\pgfpathrectangle{\pgfqpoint{0.135000in}{0.135000in}}{\pgfqpoint{2.914509in}{2.695000in}} %
\pgfusepath{clip}%
\pgfsetrectcap%
\pgfsetroundjoin%
\pgfsetlinewidth{1.003750pt}%
\definecolor{currentstroke}{rgb}{1.000000,0.000000,0.000000}%
\pgfsetstrokecolor{currentstroke}%
\pgfsetdash{}{0pt}%
\pgfpathmoveto{\pgfqpoint{2.118135in}{1.401102in}}%
\pgfpathlineto{\pgfqpoint{2.118135in}{1.563898in}}%
\pgfpathlineto{\pgfqpoint{2.459545in}{1.779541in}}%
\pgfusepath{stroke}%
\end{pgfscope}%
\begin{pgfscope}%
\pgfpathrectangle{\pgfqpoint{0.135000in}{0.135000in}}{\pgfqpoint{2.914509in}{2.695000in}} %
\pgfusepath{clip}%
\pgfsetrectcap%
\pgfsetroundjoin%
\pgfsetlinewidth{1.003750pt}%
\definecolor{currentstroke}{rgb}{1.000000,0.000000,0.000000}%
\pgfsetstrokecolor{currentstroke}%
\pgfsetdash{}{0pt}%
\pgfpathmoveto{\pgfqpoint{1.810081in}{1.285987in}}%
\pgfpathlineto{\pgfqpoint{1.810081in}{1.679013in}}%
\pgfpathlineto{\pgfqpoint{1.951498in}{2.199620in}}%
\pgfusepath{stroke}%
\end{pgfscope}%
\begin{pgfscope}%
\pgfpathrectangle{\pgfqpoint{0.135000in}{0.135000in}}{\pgfqpoint{2.914509in}{2.695000in}} %
\pgfusepath{clip}%
\pgfsetrectcap%
\pgfsetroundjoin%
\pgfsetlinewidth{1.003750pt}%
\definecolor{currentstroke}{rgb}{1.000000,0.000000,0.000000}%
\pgfsetstrokecolor{currentstroke}%
\pgfsetdash{}{0pt}%
\pgfpathmoveto{\pgfqpoint{1.374428in}{1.285987in}}%
\pgfpathlineto{\pgfqpoint{1.374428in}{1.679013in}}%
\pgfpathlineto{\pgfqpoint{1.233011in}{2.199620in}}%
\pgfusepath{stroke}%
\end{pgfscope}%
\begin{pgfscope}%
\pgfpathrectangle{\pgfqpoint{0.135000in}{0.135000in}}{\pgfqpoint{2.914509in}{2.695000in}} %
\pgfusepath{clip}%
\pgfsetrectcap%
\pgfsetroundjoin%
\pgfsetlinewidth{1.003750pt}%
\definecolor{currentstroke}{rgb}{1.000000,0.000000,0.000000}%
\pgfsetstrokecolor{currentstroke}%
\pgfsetdash{}{0pt}%
\pgfpathmoveto{\pgfqpoint{1.066375in}{1.401102in}}%
\pgfpathlineto{\pgfqpoint{1.066375in}{1.563898in}}%
\pgfpathlineto{\pgfqpoint{0.724964in}{1.779541in}}%
\pgfusepath{stroke}%
\end{pgfscope}%
\begin{pgfscope}%
\pgfpathrectangle{\pgfqpoint{0.135000in}{0.135000in}}{\pgfqpoint{2.914509in}{2.695000in}} %
\pgfusepath{clip}%
\pgfsetrectcap%
\pgfsetroundjoin%
\pgfsetlinewidth{1.003750pt}%
\definecolor{currentstroke}{rgb}{1.000000,0.000000,0.000000}%
\pgfsetstrokecolor{currentstroke}%
\pgfsetdash{}{0pt}%
\pgfpathmoveto{\pgfqpoint{1.066375in}{1.563898in}}%
\pgfpathlineto{\pgfqpoint{1.066375in}{1.401102in}}%
\pgfpathlineto{\pgfqpoint{0.724964in}{1.185459in}}%
\pgfusepath{stroke}%
\end{pgfscope}%
\begin{pgfscope}%
\pgfpathrectangle{\pgfqpoint{0.135000in}{0.135000in}}{\pgfqpoint{2.914509in}{2.695000in}} %
\pgfusepath{clip}%
\pgfsetrectcap%
\pgfsetroundjoin%
\pgfsetlinewidth{1.003750pt}%
\definecolor{currentstroke}{rgb}{1.000000,0.000000,0.000000}%
\pgfsetstrokecolor{currentstroke}%
\pgfsetdash{}{0pt}%
\pgfpathmoveto{\pgfqpoint{1.374428in}{1.679013in}}%
\pgfpathlineto{\pgfqpoint{1.374428in}{1.285987in}}%
\pgfpathlineto{\pgfqpoint{1.233011in}{0.765380in}}%
\pgfusepath{stroke}%
\end{pgfscope}%
\begin{pgfscope}%
\pgfpathrectangle{\pgfqpoint{0.135000in}{0.135000in}}{\pgfqpoint{2.914509in}{2.695000in}} %
\pgfusepath{clip}%
\pgfsetrectcap%
\pgfsetroundjoin%
\pgfsetlinewidth{1.003750pt}%
\definecolor{currentstroke}{rgb}{1.000000,0.000000,0.000000}%
\pgfsetstrokecolor{currentstroke}%
\pgfsetdash{}{0pt}%
\pgfpathmoveto{\pgfqpoint{1.810081in}{1.679013in}}%
\pgfpathlineto{\pgfqpoint{1.810081in}{1.285987in}}%
\pgfpathlineto{\pgfqpoint{1.951498in}{0.765380in}}%
\pgfusepath{stroke}%
\end{pgfscope}%
\begin{pgfscope}%
\pgfpathrectangle{\pgfqpoint{0.135000in}{0.135000in}}{\pgfqpoint{2.914509in}{2.695000in}} %
\pgfusepath{clip}%
\pgfsetrectcap%
\pgfsetroundjoin%
\pgfsetlinewidth{1.003750pt}%
\definecolor{currentstroke}{rgb}{1.000000,0.000000,0.000000}%
\pgfsetstrokecolor{currentstroke}%
\pgfsetdash{}{0pt}%
\pgfpathmoveto{\pgfqpoint{2.118135in}{1.563898in}}%
\pgfpathlineto{\pgfqpoint{2.118135in}{1.401102in}}%
\pgfpathlineto{\pgfqpoint{2.459545in}{1.185459in}}%
\pgfusepath{stroke}%
\end{pgfscope}%
\begin{pgfscope}%
\pgfpathrectangle{\pgfqpoint{0.135000in}{0.135000in}}{\pgfqpoint{2.914509in}{2.695000in}} %
\pgfusepath{clip}%
\pgfsetbuttcap%
\pgfsetroundjoin%
\definecolor{currentfill}{rgb}{1.000000,0.000000,0.000000}%
\pgfsetfillcolor{currentfill}%
\pgfsetlinewidth{1.003750pt}%
\definecolor{currentstroke}{rgb}{1.000000,0.000000,0.000000}%
\pgfsetstrokecolor{currentstroke}%
\pgfsetdash{}{0pt}%
\pgfsys@defobject{currentmarker}{\pgfqpoint{-0.010417in}{-0.010417in}}{\pgfqpoint{0.010417in}{0.010417in}}{%
\pgfpathmoveto{\pgfqpoint{0.000000in}{-0.010417in}}%
\pgfpathcurveto{\pgfqpoint{0.002763in}{-0.010417in}}{\pgfqpoint{0.005412in}{-0.009319in}}{\pgfqpoint{0.007366in}{-0.007366in}}%
\pgfpathcurveto{\pgfqpoint{0.009319in}{-0.005412in}}{\pgfqpoint{0.010417in}{-0.002763in}}{\pgfqpoint{0.010417in}{0.000000in}}%
\pgfpathcurveto{\pgfqpoint{0.010417in}{0.002763in}}{\pgfqpoint{0.009319in}{0.005412in}}{\pgfqpoint{0.007366in}{0.007366in}}%
\pgfpathcurveto{\pgfqpoint{0.005412in}{0.009319in}}{\pgfqpoint{0.002763in}{0.010417in}}{\pgfqpoint{0.000000in}{0.010417in}}%
\pgfpathcurveto{\pgfqpoint{-0.002763in}{0.010417in}}{\pgfqpoint{-0.005412in}{0.009319in}}{\pgfqpoint{-0.007366in}{0.007366in}}%
\pgfpathcurveto{\pgfqpoint{-0.009319in}{0.005412in}}{\pgfqpoint{-0.010417in}{0.002763in}}{\pgfqpoint{-0.010417in}{0.000000in}}%
\pgfpathcurveto{\pgfqpoint{-0.010417in}{-0.002763in}}{\pgfqpoint{-0.009319in}{-0.005412in}}{\pgfqpoint{-0.007366in}{-0.007366in}}%
\pgfpathcurveto{\pgfqpoint{-0.005412in}{-0.009319in}}{\pgfqpoint{-0.002763in}{-0.010417in}}{\pgfqpoint{0.000000in}{-0.010417in}}%
\pgfpathclose%
\pgfusepath{stroke,fill}%
}%
\begin{pgfscope}%
\pgfsys@transformshift{2.118135in}{1.563898in}%
\pgfsys@useobject{currentmarker}{}%
\end{pgfscope}%
\begin{pgfscope}%
\pgfsys@transformshift{2.118135in}{1.401102in}%
\pgfsys@useobject{currentmarker}{}%
\end{pgfscope}%
\begin{pgfscope}%
\pgfsys@transformshift{2.459545in}{1.185459in}%
\pgfsys@useobject{currentmarker}{}%
\end{pgfscope}%
\begin{pgfscope}%
\pgfsys@transformshift{2.118135in}{1.401102in}%
\pgfsys@useobject{currentmarker}{}%
\end{pgfscope}%
\begin{pgfscope}%
\pgfsys@transformshift{2.118135in}{1.563898in}%
\pgfsys@useobject{currentmarker}{}%
\end{pgfscope}%
\begin{pgfscope}%
\pgfsys@transformshift{2.459545in}{1.779541in}%
\pgfsys@useobject{currentmarker}{}%
\end{pgfscope}%
\begin{pgfscope}%
\pgfsys@transformshift{1.810081in}{1.285987in}%
\pgfsys@useobject{currentmarker}{}%
\end{pgfscope}%
\begin{pgfscope}%
\pgfsys@transformshift{1.810081in}{1.679013in}%
\pgfsys@useobject{currentmarker}{}%
\end{pgfscope}%
\begin{pgfscope}%
\pgfsys@transformshift{1.951498in}{2.199620in}%
\pgfsys@useobject{currentmarker}{}%
\end{pgfscope}%
\begin{pgfscope}%
\pgfsys@transformshift{1.374428in}{1.285987in}%
\pgfsys@useobject{currentmarker}{}%
\end{pgfscope}%
\begin{pgfscope}%
\pgfsys@transformshift{1.374428in}{1.679013in}%
\pgfsys@useobject{currentmarker}{}%
\end{pgfscope}%
\begin{pgfscope}%
\pgfsys@transformshift{1.233011in}{2.199620in}%
\pgfsys@useobject{currentmarker}{}%
\end{pgfscope}%
\begin{pgfscope}%
\pgfsys@transformshift{1.066375in}{1.401102in}%
\pgfsys@useobject{currentmarker}{}%
\end{pgfscope}%
\begin{pgfscope}%
\pgfsys@transformshift{1.066375in}{1.563898in}%
\pgfsys@useobject{currentmarker}{}%
\end{pgfscope}%
\begin{pgfscope}%
\pgfsys@transformshift{0.724964in}{1.779541in}%
\pgfsys@useobject{currentmarker}{}%
\end{pgfscope}%
\begin{pgfscope}%
\pgfsys@transformshift{1.066375in}{1.563898in}%
\pgfsys@useobject{currentmarker}{}%
\end{pgfscope}%
\begin{pgfscope}%
\pgfsys@transformshift{1.066375in}{1.401102in}%
\pgfsys@useobject{currentmarker}{}%
\end{pgfscope}%
\begin{pgfscope}%
\pgfsys@transformshift{0.724964in}{1.185459in}%
\pgfsys@useobject{currentmarker}{}%
\end{pgfscope}%
\begin{pgfscope}%
\pgfsys@transformshift{1.374428in}{1.679013in}%
\pgfsys@useobject{currentmarker}{}%
\end{pgfscope}%
\begin{pgfscope}%
\pgfsys@transformshift{1.374428in}{1.285987in}%
\pgfsys@useobject{currentmarker}{}%
\end{pgfscope}%
\begin{pgfscope}%
\pgfsys@transformshift{1.233011in}{0.765380in}%
\pgfsys@useobject{currentmarker}{}%
\end{pgfscope}%
\begin{pgfscope}%
\pgfsys@transformshift{1.810081in}{1.679013in}%
\pgfsys@useobject{currentmarker}{}%
\end{pgfscope}%
\begin{pgfscope}%
\pgfsys@transformshift{1.810081in}{1.285987in}%
\pgfsys@useobject{currentmarker}{}%
\end{pgfscope}%
\begin{pgfscope}%
\pgfsys@transformshift{1.951498in}{0.765380in}%
\pgfsys@useobject{currentmarker}{}%
\end{pgfscope}%
\begin{pgfscope}%
\pgfsys@transformshift{2.118135in}{1.563898in}%
\pgfsys@useobject{currentmarker}{}%
\end{pgfscope}%
\begin{pgfscope}%
\pgfsys@transformshift{2.118135in}{1.401102in}%
\pgfsys@useobject{currentmarker}{}%
\end{pgfscope}%
\begin{pgfscope}%
\pgfsys@transformshift{2.459545in}{1.185459in}%
\pgfsys@useobject{currentmarker}{}%
\end{pgfscope}%
\end{pgfscope}%
\end{pgfpicture}%
\makeatother%
\endgroup%

%% file: trigonometric-symmetric7_image.pgf
%% Creator: Matplotlib, PGF backend
%%
%% To include the figure in your LaTeX document, write
%%   \input{<filename>.pgf}
%%
%% Make sure the required packages are loaded in your preamble
%%   \usepackage{pgf}
%%
%% Figures using additional raster images can only be included by \input if
%% they are in the same directory as the main LaTeX file. For loading figures
%% from other directories you can use the `import` package
%% and then include the figures with
%%   \import{<path to file>}{<filename>.pgf}
%%
%% Matplotlib used the following preamble
%%   \usepackage{fontspec}
%%   \setmainfont{DejaVu Serif}
%%   \setsansfont{DejaVu Sans}
%%   \setmonofont{DejaVu Sans Mono}
%%
\begingroup%
\makeatletter%
\begin{pgfpicture}%
\pgfpathrectangle{\pgfpointorigin}{\pgfqpoint{3.316323in}{3.170273in}}%
\pgfusepath{use as bounding box, clip}%
\begin{pgfscope}%
\pgfsetbuttcap%
\pgfsetmiterjoin%
\definecolor{currentfill}{rgb}{1.000000,1.000000,1.000000}%
\pgfsetfillcolor{currentfill}%
\pgfsetlinewidth{0.000000pt}%
\definecolor{currentstroke}{rgb}{1.000000,1.000000,1.000000}%
\pgfsetstrokecolor{currentstroke}%
\pgfsetdash{}{0pt}%
\pgfpathmoveto{\pgfqpoint{0.000000in}{0.000000in}}%
\pgfpathlineto{\pgfqpoint{3.316323in}{0.000000in}}%
\pgfpathlineto{\pgfqpoint{3.316323in}{3.170273in}}%
\pgfpathlineto{\pgfqpoint{0.000000in}{3.170273in}}%
\pgfpathclose%
\pgfusepath{fill}%
\end{pgfscope}%
\begin{pgfscope}%
\pgfsetbuttcap%
\pgfsetmiterjoin%
\definecolor{currentfill}{rgb}{1.000000,1.000000,1.000000}%
\pgfsetfillcolor{currentfill}%
\pgfsetlinewidth{0.000000pt}%
\definecolor{currentstroke}{rgb}{0.000000,0.000000,0.000000}%
\pgfsetstrokecolor{currentstroke}%
\pgfsetstrokeopacity{0.000000}%
\pgfsetdash{}{0pt}%
\pgfpathmoveto{\pgfqpoint{0.135000in}{0.135000in}}%
\pgfpathlineto{\pgfqpoint{3.049509in}{0.135000in}}%
\pgfpathlineto{\pgfqpoint{3.049509in}{2.830000in}}%
\pgfpathlineto{\pgfqpoint{0.135000in}{2.830000in}}%
\pgfpathclose%
\pgfusepath{fill}%
\end{pgfscope}%
\begin{pgfscope}%
\pgfpathrectangle{\pgfqpoint{0.135000in}{0.135000in}}{\pgfqpoint{2.914509in}{2.695000in}} %
\pgfusepath{clip}%
\pgfsetrectcap%
\pgfsetroundjoin%
\pgfsetlinewidth{0.702625pt}%
\definecolor{currentstroke}{rgb}{0.000000,0.000000,0.000000}%
\pgfsetstrokecolor{currentstroke}%
\pgfsetdash{}{0pt}%
\pgfpathmoveto{\pgfqpoint{0.135000in}{1.482500in}}%
\pgfpathlineto{\pgfqpoint{3.049509in}{1.482500in}}%
\pgfusepath{stroke}%
\end{pgfscope}%
\begin{pgfscope}%
\pgfpathrectangle{\pgfqpoint{0.135000in}{0.135000in}}{\pgfqpoint{2.914509in}{2.695000in}} %
\pgfusepath{clip}%
\pgfsetrectcap%
\pgfsetroundjoin%
\pgfsetlinewidth{0.702625pt}%
\definecolor{currentstroke}{rgb}{0.000000,0.000000,0.000000}%
\pgfsetstrokecolor{currentstroke}%
\pgfsetdash{}{0pt}%
\pgfpathmoveto{\pgfqpoint{1.592255in}{0.135000in}}%
\pgfpathlineto{\pgfqpoint{1.592255in}{2.830000in}}%
\pgfusepath{stroke}%
\end{pgfscope}%
\begin{pgfscope}%
\pgfsetbuttcap%
\pgfsetmiterjoin%
\definecolor{currentfill}{rgb}{0.000000,0.000000,0.000000}%
\pgfsetfillcolor{currentfill}%
\pgfsetlinewidth{1.003750pt}%
\definecolor{currentstroke}{rgb}{0.000000,0.000000,0.000000}%
\pgfsetstrokecolor{currentstroke}%
\pgfsetdash{}{0pt}%
\pgfpathmoveto{\pgfqpoint{3.049509in}{1.482500in}}%
\pgfpathlineto{\pgfqpoint{3.020364in}{1.469025in}}%
\pgfpathlineto{\pgfqpoint{3.029108in}{1.482237in}}%
\pgfpathlineto{\pgfqpoint{0.135000in}{1.482237in}}%
\pgfpathlineto{\pgfqpoint{0.135000in}{1.482763in}}%
\pgfpathlineto{\pgfqpoint{3.029108in}{1.482763in}}%
\pgfpathlineto{\pgfqpoint{3.020364in}{1.495975in}}%
\pgfpathclose%
\pgfusepath{stroke,fill}%
\end{pgfscope}%
\begin{pgfscope}%
\pgfsetbuttcap%
\pgfsetmiterjoin%
\definecolor{currentfill}{rgb}{0.000000,0.000000,0.000000}%
\pgfsetfillcolor{currentfill}%
\pgfsetlinewidth{1.003750pt}%
\definecolor{currentstroke}{rgb}{0.000000,0.000000,0.000000}%
\pgfsetstrokecolor{currentstroke}%
\pgfsetdash{}{0pt}%
\pgfpathmoveto{\pgfqpoint{1.592255in}{2.830000in}}%
\pgfpathlineto{\pgfqpoint{1.605730in}{2.800855in}}%
\pgfpathlineto{\pgfqpoint{1.592518in}{2.809598in}}%
\pgfpathlineto{\pgfqpoint{1.592518in}{0.135000in}}%
\pgfpathlineto{\pgfqpoint{1.591991in}{0.135000in}}%
\pgfpathlineto{\pgfqpoint{1.591991in}{2.809598in}}%
\pgfpathlineto{\pgfqpoint{1.578780in}{2.800855in}}%
\pgfpathclose%
\pgfusepath{stroke,fill}%
\end{pgfscope}%
\begin{pgfscope}%
\pgfpathrectangle{\pgfqpoint{0.135000in}{0.135000in}}{\pgfqpoint{2.914509in}{2.695000in}} %
\pgfusepath{clip}%
\pgfsetrectcap%
\pgfsetroundjoin%
\pgfsetlinewidth{0.803000pt}%
\definecolor{currentstroke}{rgb}{0.000000,0.000000,0.000000}%
\pgfsetstrokecolor{currentstroke}%
\pgfsetdash{}{0pt}%
\pgfpathmoveto{\pgfqpoint{1.592255in}{0.135000in}}%
\pgfpathlineto{\pgfqpoint{1.592255in}{2.830000in}}%
\pgfusepath{stroke}%
\end{pgfscope}%
\begin{pgfscope}%
\pgfpathrectangle{\pgfqpoint{0.135000in}{0.135000in}}{\pgfqpoint{2.914509in}{2.695000in}} %
\pgfusepath{clip}%
\pgfsetrectcap%
\pgfsetroundjoin%
\pgfsetlinewidth{0.803000pt}%
\definecolor{currentstroke}{rgb}{0.000000,0.000000,0.000000}%
\pgfsetstrokecolor{currentstroke}%
\pgfsetdash{}{0pt}%
\pgfpathmoveto{\pgfqpoint{0.135000in}{1.482500in}}%
\pgfpathlineto{\pgfqpoint{3.049509in}{1.482500in}}%
\pgfusepath{stroke}%
\end{pgfscope}%
\begin{pgfscope}%
\pgfpathrectangle{\pgfqpoint{0.135000in}{0.135000in}}{\pgfqpoint{2.914509in}{2.695000in}} %
\pgfusepath{clip}%
\pgfsetrectcap%
\pgfsetroundjoin%
\pgfsetlinewidth{0.803000pt}%
\definecolor{currentstroke}{rgb}{0.000000,0.000000,0.000000}%
\pgfsetstrokecolor{currentstroke}%
\pgfsetdash{}{0pt}%
\pgfpathmoveto{\pgfqpoint{1.592255in}{0.135000in}}%
\pgfpathlineto{\pgfqpoint{1.592255in}{2.830000in}}%
\pgfusepath{stroke}%
\end{pgfscope}%
\begin{pgfscope}%
\pgfpathrectangle{\pgfqpoint{0.135000in}{0.135000in}}{\pgfqpoint{2.914509in}{2.695000in}} %
\pgfusepath{clip}%
\pgfsetrectcap%
\pgfsetroundjoin%
\pgfsetlinewidth{0.803000pt}%
\definecolor{currentstroke}{rgb}{0.000000,0.000000,0.000000}%
\pgfsetstrokecolor{currentstroke}%
\pgfsetdash{}{0pt}%
\pgfpathmoveto{\pgfqpoint{0.135000in}{1.482500in}}%
\pgfpathlineto{\pgfqpoint{3.049509in}{1.482500in}}%
\pgfusepath{stroke}%
\end{pgfscope}%
\begin{pgfscope}%
\pgftext[x=1.592255in,y=2.964750in,,base]{\sffamily\fontsize{10.000000}{12.000000}\selectfont \(\displaystyle y\)}%
\end{pgfscope}%
\begin{pgfscope}%
\pgftext[x=3.136945in,y=1.482500in,left,]{\sffamily\fontsize{10.000000}{12.000000}\selectfont \(\displaystyle x\)}%
\end{pgfscope}%
\begin{pgfscope}%
\pgfpathrectangle{\pgfqpoint{0.135000in}{0.135000in}}{\pgfqpoint{2.914509in}{2.695000in}} %
\pgfusepath{clip}%
\pgfsetbuttcap%
\pgfsetroundjoin%
\definecolor{currentfill}{rgb}{0.700000,0.700000,0.700000}%
\pgfsetfillcolor{currentfill}%
\pgfsetlinewidth{1.003750pt}%
\definecolor{currentstroke}{rgb}{0.700000,0.700000,0.700000}%
\pgfsetstrokecolor{currentstroke}%
\pgfsetdash{}{0pt}%
\pgfsys@defobject{currentmarker}{\pgfqpoint{-0.010417in}{-0.010417in}}{\pgfqpoint{0.010417in}{0.010417in}}{%
\pgfpathmoveto{\pgfqpoint{0.000000in}{-0.010417in}}%
\pgfpathcurveto{\pgfqpoint{0.002763in}{-0.010417in}}{\pgfqpoint{0.005412in}{-0.009319in}}{\pgfqpoint{0.007366in}{-0.007366in}}%
\pgfpathcurveto{\pgfqpoint{0.009319in}{-0.005412in}}{\pgfqpoint{0.010417in}{-0.002763in}}{\pgfqpoint{0.010417in}{0.000000in}}%
\pgfpathcurveto{\pgfqpoint{0.010417in}{0.002763in}}{\pgfqpoint{0.009319in}{0.005412in}}{\pgfqpoint{0.007366in}{0.007366in}}%
\pgfpathcurveto{\pgfqpoint{0.005412in}{0.009319in}}{\pgfqpoint{0.002763in}{0.010417in}}{\pgfqpoint{0.000000in}{0.010417in}}%
\pgfpathcurveto{\pgfqpoint{-0.002763in}{0.010417in}}{\pgfqpoint{-0.005412in}{0.009319in}}{\pgfqpoint{-0.007366in}{0.007366in}}%
\pgfpathcurveto{\pgfqpoint{-0.009319in}{0.005412in}}{\pgfqpoint{-0.010417in}{0.002763in}}{\pgfqpoint{-0.010417in}{0.000000in}}%
\pgfpathcurveto{\pgfqpoint{-0.010417in}{-0.002763in}}{\pgfqpoint{-0.009319in}{-0.005412in}}{\pgfqpoint{-0.007366in}{-0.007366in}}%
\pgfpathcurveto{\pgfqpoint{-0.005412in}{-0.009319in}}{\pgfqpoint{-0.002763in}{-0.010417in}}{\pgfqpoint{0.000000in}{-0.010417in}}%
\pgfpathclose%
\pgfusepath{stroke,fill}%
}%
\begin{pgfscope}%
\pgfsys@transformshift{2.120227in}{1.482500in}%
\pgfsys@useobject{currentmarker}{}%
\end{pgfscope}%
\begin{pgfscope}%
\pgfsys@transformshift{2.161463in}{1.482500in}%
\pgfsys@useobject{currentmarker}{}%
\end{pgfscope}%
\begin{pgfscope}%
\pgfsys@transformshift{2.291612in}{1.482500in}%
\pgfsys@useobject{currentmarker}{}%
\end{pgfscope}%
\begin{pgfscope}%
\pgfsys@transformshift{2.531003in}{1.482500in}%
\pgfsys@useobject{currentmarker}{}%
\end{pgfscope}%
\begin{pgfscope}%
\pgfsys@transformshift{2.917032in}{1.482500in}%
\pgfsys@useobject{currentmarker}{}%
\end{pgfscope}%
\begin{pgfscope}%
\pgfsys@transformshift{2.080038in}{1.482500in}%
\pgfsys@useobject{currentmarker}{}%
\end{pgfscope}%
\begin{pgfscope}%
\pgfsys@transformshift{2.118135in}{1.563898in}%
\pgfsys@useobject{currentmarker}{}%
\end{pgfscope}%
\begin{pgfscope}%
\pgfsys@transformshift{2.238376in}{1.658012in}%
\pgfsys@useobject{currentmarker}{}%
\end{pgfscope}%
\begin{pgfscope}%
\pgfsys@transformshift{2.459545in}{1.779541in}%
\pgfsys@useobject{currentmarker}{}%
\end{pgfscope}%
\begin{pgfscope}%
\pgfsys@transformshift{2.816189in}{1.947469in}%
\pgfsys@useobject{currentmarker}{}%
\end{pgfscope}%
\begin{pgfscope}%
\pgfsys@transformshift{1.965588in}{1.482500in}%
\pgfsys@useobject{currentmarker}{}%
\end{pgfscope}%
\begin{pgfscope}%
\pgfsys@transformshift{1.994746in}{1.632905in}%
\pgfsys@useobject{currentmarker}{}%
\end{pgfscope}%
\begin{pgfscope}%
\pgfsys@transformshift{2.086775in}{1.806803in}%
\pgfsys@useobject{currentmarker}{}%
\end{pgfscope}%
\begin{pgfscope}%
\pgfsys@transformshift{2.256050in}{2.031360in}%
\pgfsys@useobject{currentmarker}{}%
\end{pgfscope}%
\begin{pgfscope}%
\pgfsys@transformshift{2.529013in}{2.341651in}%
\pgfsys@useobject{currentmarker}{}%
\end{pgfscope}%
\begin{pgfscope}%
\pgfsys@transformshift{1.794301in}{1.482500in}%
\pgfsys@useobject{currentmarker}{}%
\end{pgfscope}%
\begin{pgfscope}%
\pgfsys@transformshift{1.810081in}{1.679013in}%
\pgfsys@useobject{currentmarker}{}%
\end{pgfscope}%
\begin{pgfscope}%
\pgfsys@transformshift{1.859887in}{1.906223in}%
\pgfsys@useobject{currentmarker}{}%
\end{pgfscope}%
\begin{pgfscope}%
\pgfsys@transformshift{1.951498in}{2.199620in}%
\pgfsys@useobject{currentmarker}{}%
\end{pgfscope}%
\begin{pgfscope}%
\pgfsys@transformshift{2.099225in}{2.605034in}%
\pgfsys@useobject{currentmarker}{}%
\end{pgfscope}%
\begin{pgfscope}%
\pgfsys@transformshift{1.592255in}{1.482500in}%
\pgfsys@useobject{currentmarker}{}%
\end{pgfscope}%
\begin{pgfscope}%
\pgfsys@transformshift{1.592255in}{1.695204in}%
\pgfsys@useobject{currentmarker}{}%
\end{pgfscope}%
\begin{pgfscope}%
\pgfsys@transformshift{1.592255in}{1.941134in}%
\pgfsys@useobject{currentmarker}{}%
\end{pgfscope}%
\begin{pgfscope}%
\pgfsys@transformshift{1.592255in}{2.258705in}%
\pgfsys@useobject{currentmarker}{}%
\end{pgfscope}%
\begin{pgfscope}%
\pgfsys@transformshift{1.592255in}{2.697522in}%
\pgfsys@useobject{currentmarker}{}%
\end{pgfscope}%
\begin{pgfscope}%
\pgfsys@transformshift{1.390208in}{1.482500in}%
\pgfsys@useobject{currentmarker}{}%
\end{pgfscope}%
\begin{pgfscope}%
\pgfsys@transformshift{1.374428in}{1.679013in}%
\pgfsys@useobject{currentmarker}{}%
\end{pgfscope}%
\begin{pgfscope}%
\pgfsys@transformshift{1.324622in}{1.906223in}%
\pgfsys@useobject{currentmarker}{}%
\end{pgfscope}%
\begin{pgfscope}%
\pgfsys@transformshift{1.233011in}{2.199620in}%
\pgfsys@useobject{currentmarker}{}%
\end{pgfscope}%
\begin{pgfscope}%
\pgfsys@transformshift{1.085284in}{2.605034in}%
\pgfsys@useobject{currentmarker}{}%
\end{pgfscope}%
\begin{pgfscope}%
\pgfsys@transformshift{1.218922in}{1.482500in}%
\pgfsys@useobject{currentmarker}{}%
\end{pgfscope}%
\begin{pgfscope}%
\pgfsys@transformshift{1.189764in}{1.632905in}%
\pgfsys@useobject{currentmarker}{}%
\end{pgfscope}%
\begin{pgfscope}%
\pgfsys@transformshift{1.097735in}{1.806803in}%
\pgfsys@useobject{currentmarker}{}%
\end{pgfscope}%
\begin{pgfscope}%
\pgfsys@transformshift{0.928459in}{2.031360in}%
\pgfsys@useobject{currentmarker}{}%
\end{pgfscope}%
\begin{pgfscope}%
\pgfsys@transformshift{0.655496in}{2.341651in}%
\pgfsys@useobject{currentmarker}{}%
\end{pgfscope}%
\begin{pgfscope}%
\pgfsys@transformshift{1.104472in}{1.482500in}%
\pgfsys@useobject{currentmarker}{}%
\end{pgfscope}%
\begin{pgfscope}%
\pgfsys@transformshift{1.066375in}{1.563898in}%
\pgfsys@useobject{currentmarker}{}%
\end{pgfscope}%
\begin{pgfscope}%
\pgfsys@transformshift{0.946133in}{1.658012in}%
\pgfsys@useobject{currentmarker}{}%
\end{pgfscope}%
\begin{pgfscope}%
\pgfsys@transformshift{0.724964in}{1.779541in}%
\pgfsys@useobject{currentmarker}{}%
\end{pgfscope}%
\begin{pgfscope}%
\pgfsys@transformshift{0.368320in}{1.947469in}%
\pgfsys@useobject{currentmarker}{}%
\end{pgfscope}%
\begin{pgfscope}%
\pgfsys@transformshift{1.064282in}{1.482500in}%
\pgfsys@useobject{currentmarker}{}%
\end{pgfscope}%
\begin{pgfscope}%
\pgfsys@transformshift{1.023046in}{1.482500in}%
\pgfsys@useobject{currentmarker}{}%
\end{pgfscope}%
\begin{pgfscope}%
\pgfsys@transformshift{0.892898in}{1.482500in}%
\pgfsys@useobject{currentmarker}{}%
\end{pgfscope}%
\begin{pgfscope}%
\pgfsys@transformshift{0.653506in}{1.482500in}%
\pgfsys@useobject{currentmarker}{}%
\end{pgfscope}%
\begin{pgfscope}%
\pgfsys@transformshift{0.267478in}{1.482500in}%
\pgfsys@useobject{currentmarker}{}%
\end{pgfscope}%
\begin{pgfscope}%
\pgfsys@transformshift{1.104472in}{1.482500in}%
\pgfsys@useobject{currentmarker}{}%
\end{pgfscope}%
\begin{pgfscope}%
\pgfsys@transformshift{1.066375in}{1.401102in}%
\pgfsys@useobject{currentmarker}{}%
\end{pgfscope}%
\begin{pgfscope}%
\pgfsys@transformshift{0.946133in}{1.306988in}%
\pgfsys@useobject{currentmarker}{}%
\end{pgfscope}%
\begin{pgfscope}%
\pgfsys@transformshift{0.724964in}{1.185459in}%
\pgfsys@useobject{currentmarker}{}%
\end{pgfscope}%
\begin{pgfscope}%
\pgfsys@transformshift{0.368320in}{1.017531in}%
\pgfsys@useobject{currentmarker}{}%
\end{pgfscope}%
\begin{pgfscope}%
\pgfsys@transformshift{1.218922in}{1.482500in}%
\pgfsys@useobject{currentmarker}{}%
\end{pgfscope}%
\begin{pgfscope}%
\pgfsys@transformshift{1.189764in}{1.332095in}%
\pgfsys@useobject{currentmarker}{}%
\end{pgfscope}%
\begin{pgfscope}%
\pgfsys@transformshift{1.097735in}{1.158197in}%
\pgfsys@useobject{currentmarker}{}%
\end{pgfscope}%
\begin{pgfscope}%
\pgfsys@transformshift{0.928459in}{0.933640in}%
\pgfsys@useobject{currentmarker}{}%
\end{pgfscope}%
\begin{pgfscope}%
\pgfsys@transformshift{0.655496in}{0.623349in}%
\pgfsys@useobject{currentmarker}{}%
\end{pgfscope}%
\begin{pgfscope}%
\pgfsys@transformshift{1.390208in}{1.482500in}%
\pgfsys@useobject{currentmarker}{}%
\end{pgfscope}%
\begin{pgfscope}%
\pgfsys@transformshift{1.374428in}{1.285987in}%
\pgfsys@useobject{currentmarker}{}%
\end{pgfscope}%
\begin{pgfscope}%
\pgfsys@transformshift{1.324622in}{1.058777in}%
\pgfsys@useobject{currentmarker}{}%
\end{pgfscope}%
\begin{pgfscope}%
\pgfsys@transformshift{1.233011in}{0.765380in}%
\pgfsys@useobject{currentmarker}{}%
\end{pgfscope}%
\begin{pgfscope}%
\pgfsys@transformshift{1.085284in}{0.359966in}%
\pgfsys@useobject{currentmarker}{}%
\end{pgfscope}%
\begin{pgfscope}%
\pgfsys@transformshift{1.592255in}{1.482500in}%
\pgfsys@useobject{currentmarker}{}%
\end{pgfscope}%
\begin{pgfscope}%
\pgfsys@transformshift{1.592255in}{1.269796in}%
\pgfsys@useobject{currentmarker}{}%
\end{pgfscope}%
\begin{pgfscope}%
\pgfsys@transformshift{1.592255in}{1.023866in}%
\pgfsys@useobject{currentmarker}{}%
\end{pgfscope}%
\begin{pgfscope}%
\pgfsys@transformshift{1.592255in}{0.706295in}%
\pgfsys@useobject{currentmarker}{}%
\end{pgfscope}%
\begin{pgfscope}%
\pgfsys@transformshift{1.592255in}{0.267478in}%
\pgfsys@useobject{currentmarker}{}%
\end{pgfscope}%
\begin{pgfscope}%
\pgfsys@transformshift{1.794301in}{1.482500in}%
\pgfsys@useobject{currentmarker}{}%
\end{pgfscope}%
\begin{pgfscope}%
\pgfsys@transformshift{1.810081in}{1.285987in}%
\pgfsys@useobject{currentmarker}{}%
\end{pgfscope}%
\begin{pgfscope}%
\pgfsys@transformshift{1.859887in}{1.058777in}%
\pgfsys@useobject{currentmarker}{}%
\end{pgfscope}%
\begin{pgfscope}%
\pgfsys@transformshift{1.951498in}{0.765380in}%
\pgfsys@useobject{currentmarker}{}%
\end{pgfscope}%
\begin{pgfscope}%
\pgfsys@transformshift{2.099225in}{0.359966in}%
\pgfsys@useobject{currentmarker}{}%
\end{pgfscope}%
\begin{pgfscope}%
\pgfsys@transformshift{1.965588in}{1.482500in}%
\pgfsys@useobject{currentmarker}{}%
\end{pgfscope}%
\begin{pgfscope}%
\pgfsys@transformshift{1.994746in}{1.332095in}%
\pgfsys@useobject{currentmarker}{}%
\end{pgfscope}%
\begin{pgfscope}%
\pgfsys@transformshift{2.086775in}{1.158197in}%
\pgfsys@useobject{currentmarker}{}%
\end{pgfscope}%
\begin{pgfscope}%
\pgfsys@transformshift{2.256050in}{0.933640in}%
\pgfsys@useobject{currentmarker}{}%
\end{pgfscope}%
\begin{pgfscope}%
\pgfsys@transformshift{2.529013in}{0.623349in}%
\pgfsys@useobject{currentmarker}{}%
\end{pgfscope}%
\begin{pgfscope}%
\pgfsys@transformshift{2.080038in}{1.482500in}%
\pgfsys@useobject{currentmarker}{}%
\end{pgfscope}%
\begin{pgfscope}%
\pgfsys@transformshift{2.118135in}{1.401102in}%
\pgfsys@useobject{currentmarker}{}%
\end{pgfscope}%
\begin{pgfscope}%
\pgfsys@transformshift{2.238376in}{1.306988in}%
\pgfsys@useobject{currentmarker}{}%
\end{pgfscope}%
\begin{pgfscope}%
\pgfsys@transformshift{2.459545in}{1.185459in}%
\pgfsys@useobject{currentmarker}{}%
\end{pgfscope}%
\begin{pgfscope}%
\pgfsys@transformshift{2.816189in}{1.017531in}%
\pgfsys@useobject{currentmarker}{}%
\end{pgfscope}%
\begin{pgfscope}%
\pgfsys@transformshift{2.120227in}{1.482500in}%
\pgfsys@useobject{currentmarker}{}%
\end{pgfscope}%
\begin{pgfscope}%
\pgfsys@transformshift{2.161463in}{1.482500in}%
\pgfsys@useobject{currentmarker}{}%
\end{pgfscope}%
\begin{pgfscope}%
\pgfsys@transformshift{2.291612in}{1.482500in}%
\pgfsys@useobject{currentmarker}{}%
\end{pgfscope}%
\begin{pgfscope}%
\pgfsys@transformshift{2.531003in}{1.482500in}%
\pgfsys@useobject{currentmarker}{}%
\end{pgfscope}%
\begin{pgfscope}%
\pgfsys@transformshift{2.917032in}{1.482500in}%
\pgfsys@useobject{currentmarker}{}%
\end{pgfscope}%
\end{pgfscope}%
\begin{pgfscope}%
\pgfpathrectangle{\pgfqpoint{0.135000in}{0.135000in}}{\pgfqpoint{2.914509in}{2.695000in}} %
\pgfusepath{clip}%
\pgfsetrectcap%
\pgfsetroundjoin%
\pgfsetlinewidth{1.003750pt}%
\definecolor{currentstroke}{rgb}{0.392157,0.584314,0.929412}%
\pgfsetstrokecolor{currentstroke}%
\pgfsetdash{}{0pt}%
\pgfpathmoveto{\pgfqpoint{2.080038in}{1.482500in}}%
\pgfpathlineto{\pgfqpoint{2.080038in}{1.482500in}}%
\pgfpathlineto{\pgfqpoint{1.794301in}{1.482500in}}%
\pgfpathlineto{\pgfqpoint{1.390208in}{1.482500in}}%
\pgfpathlineto{\pgfqpoint{1.104472in}{1.482500in}}%
\pgfpathlineto{\pgfqpoint{1.104472in}{1.482500in}}%
\pgfpathlineto{\pgfqpoint{1.390208in}{1.482500in}}%
\pgfpathlineto{\pgfqpoint{1.794301in}{1.482500in}}%
\pgfpathlineto{\pgfqpoint{2.080038in}{1.482500in}}%
\pgfusepath{stroke}%
\end{pgfscope}%
\begin{pgfscope}%
\pgfpathrectangle{\pgfqpoint{0.135000in}{0.135000in}}{\pgfqpoint{2.914509in}{2.695000in}} %
\pgfusepath{clip}%
\pgfsetrectcap%
\pgfsetroundjoin%
\pgfsetlinewidth{1.003750pt}%
\definecolor{currentstroke}{rgb}{0.392157,0.584314,0.929412}%
\pgfsetstrokecolor{currentstroke}%
\pgfsetdash{}{0pt}%
\pgfpathmoveto{\pgfqpoint{2.238376in}{1.306988in}}%
\pgfpathlineto{\pgfqpoint{2.238376in}{1.658012in}}%
\pgfpathlineto{\pgfqpoint{1.859887in}{1.906223in}}%
\pgfpathlineto{\pgfqpoint{1.324622in}{1.906223in}}%
\pgfpathlineto{\pgfqpoint{0.946133in}{1.658012in}}%
\pgfpathlineto{\pgfqpoint{0.946133in}{1.306988in}}%
\pgfpathlineto{\pgfqpoint{1.324622in}{1.058777in}}%
\pgfpathlineto{\pgfqpoint{1.859887in}{1.058777in}}%
\pgfpathlineto{\pgfqpoint{2.238376in}{1.306988in}}%
\pgfusepath{stroke}%
\end{pgfscope}%
\begin{pgfscope}%
\pgfpathrectangle{\pgfqpoint{0.135000in}{0.135000in}}{\pgfqpoint{2.914509in}{2.695000in}} %
\pgfusepath{clip}%
\pgfsetrectcap%
\pgfsetroundjoin%
\pgfsetlinewidth{1.003750pt}%
\definecolor{currentstroke}{rgb}{0.392157,0.584314,0.929412}%
\pgfsetstrokecolor{currentstroke}%
\pgfsetdash{}{0pt}%
\pgfpathmoveto{\pgfqpoint{2.816189in}{1.017531in}}%
\pgfpathlineto{\pgfqpoint{2.816189in}{1.947469in}}%
\pgfpathlineto{\pgfqpoint{2.099225in}{2.605034in}}%
\pgfpathlineto{\pgfqpoint{1.085284in}{2.605034in}}%
\pgfpathlineto{\pgfqpoint{0.368320in}{1.947469in}}%
\pgfpathlineto{\pgfqpoint{0.368320in}{1.017531in}}%
\pgfpathlineto{\pgfqpoint{1.085284in}{0.359966in}}%
\pgfpathlineto{\pgfqpoint{2.099225in}{0.359966in}}%
\pgfpathlineto{\pgfqpoint{2.816189in}{1.017531in}}%
\pgfusepath{stroke}%
\end{pgfscope}%
\begin{pgfscope}%
\pgfpathrectangle{\pgfqpoint{0.135000in}{0.135000in}}{\pgfqpoint{2.914509in}{2.695000in}} %
\pgfusepath{clip}%
\pgfsetrectcap%
\pgfsetroundjoin%
\pgfsetlinewidth{1.003750pt}%
\definecolor{currentstroke}{rgb}{0.392157,0.584314,0.929412}%
\pgfsetstrokecolor{currentstroke}%
\pgfsetdash{}{0pt}%
\pgfpathmoveto{\pgfqpoint{2.080038in}{1.482500in}}%
\pgfpathlineto{\pgfqpoint{2.238376in}{1.306988in}}%
\pgfpathlineto{\pgfqpoint{2.816189in}{1.017531in}}%
\pgfusepath{stroke}%
\end{pgfscope}%
\begin{pgfscope}%
\pgfpathrectangle{\pgfqpoint{0.135000in}{0.135000in}}{\pgfqpoint{2.914509in}{2.695000in}} %
\pgfusepath{clip}%
\pgfsetrectcap%
\pgfsetroundjoin%
\pgfsetlinewidth{1.003750pt}%
\definecolor{currentstroke}{rgb}{0.392157,0.584314,0.929412}%
\pgfsetstrokecolor{currentstroke}%
\pgfsetdash{}{0pt}%
\pgfpathmoveto{\pgfqpoint{2.080038in}{1.482500in}}%
\pgfpathlineto{\pgfqpoint{2.238376in}{1.658012in}}%
\pgfpathlineto{\pgfqpoint{2.816189in}{1.947469in}}%
\pgfusepath{stroke}%
\end{pgfscope}%
\begin{pgfscope}%
\pgfpathrectangle{\pgfqpoint{0.135000in}{0.135000in}}{\pgfqpoint{2.914509in}{2.695000in}} %
\pgfusepath{clip}%
\pgfsetrectcap%
\pgfsetroundjoin%
\pgfsetlinewidth{1.003750pt}%
\definecolor{currentstroke}{rgb}{0.392157,0.584314,0.929412}%
\pgfsetstrokecolor{currentstroke}%
\pgfsetdash{}{0pt}%
\pgfpathmoveto{\pgfqpoint{1.794301in}{1.482500in}}%
\pgfpathlineto{\pgfqpoint{1.859887in}{1.906223in}}%
\pgfpathlineto{\pgfqpoint{2.099225in}{2.605034in}}%
\pgfusepath{stroke}%
\end{pgfscope}%
\begin{pgfscope}%
\pgfpathrectangle{\pgfqpoint{0.135000in}{0.135000in}}{\pgfqpoint{2.914509in}{2.695000in}} %
\pgfusepath{clip}%
\pgfsetrectcap%
\pgfsetroundjoin%
\pgfsetlinewidth{1.003750pt}%
\definecolor{currentstroke}{rgb}{0.392157,0.584314,0.929412}%
\pgfsetstrokecolor{currentstroke}%
\pgfsetdash{}{0pt}%
\pgfpathmoveto{\pgfqpoint{1.390208in}{1.482500in}}%
\pgfpathlineto{\pgfqpoint{1.324622in}{1.906223in}}%
\pgfpathlineto{\pgfqpoint{1.085284in}{2.605034in}}%
\pgfusepath{stroke}%
\end{pgfscope}%
\begin{pgfscope}%
\pgfpathrectangle{\pgfqpoint{0.135000in}{0.135000in}}{\pgfqpoint{2.914509in}{2.695000in}} %
\pgfusepath{clip}%
\pgfsetrectcap%
\pgfsetroundjoin%
\pgfsetlinewidth{1.003750pt}%
\definecolor{currentstroke}{rgb}{0.392157,0.584314,0.929412}%
\pgfsetstrokecolor{currentstroke}%
\pgfsetdash{}{0pt}%
\pgfpathmoveto{\pgfqpoint{1.104472in}{1.482500in}}%
\pgfpathlineto{\pgfqpoint{0.946133in}{1.658012in}}%
\pgfpathlineto{\pgfqpoint{0.368320in}{1.947469in}}%
\pgfusepath{stroke}%
\end{pgfscope}%
\begin{pgfscope}%
\pgfpathrectangle{\pgfqpoint{0.135000in}{0.135000in}}{\pgfqpoint{2.914509in}{2.695000in}} %
\pgfusepath{clip}%
\pgfsetrectcap%
\pgfsetroundjoin%
\pgfsetlinewidth{1.003750pt}%
\definecolor{currentstroke}{rgb}{0.392157,0.584314,0.929412}%
\pgfsetstrokecolor{currentstroke}%
\pgfsetdash{}{0pt}%
\pgfpathmoveto{\pgfqpoint{1.104472in}{1.482500in}}%
\pgfpathlineto{\pgfqpoint{0.946133in}{1.306988in}}%
\pgfpathlineto{\pgfqpoint{0.368320in}{1.017531in}}%
\pgfusepath{stroke}%
\end{pgfscope}%
\begin{pgfscope}%
\pgfpathrectangle{\pgfqpoint{0.135000in}{0.135000in}}{\pgfqpoint{2.914509in}{2.695000in}} %
\pgfusepath{clip}%
\pgfsetrectcap%
\pgfsetroundjoin%
\pgfsetlinewidth{1.003750pt}%
\definecolor{currentstroke}{rgb}{0.392157,0.584314,0.929412}%
\pgfsetstrokecolor{currentstroke}%
\pgfsetdash{}{0pt}%
\pgfpathmoveto{\pgfqpoint{1.390208in}{1.482500in}}%
\pgfpathlineto{\pgfqpoint{1.324622in}{1.058777in}}%
\pgfpathlineto{\pgfqpoint{1.085284in}{0.359966in}}%
\pgfusepath{stroke}%
\end{pgfscope}%
\begin{pgfscope}%
\pgfpathrectangle{\pgfqpoint{0.135000in}{0.135000in}}{\pgfqpoint{2.914509in}{2.695000in}} %
\pgfusepath{clip}%
\pgfsetrectcap%
\pgfsetroundjoin%
\pgfsetlinewidth{1.003750pt}%
\definecolor{currentstroke}{rgb}{0.392157,0.584314,0.929412}%
\pgfsetstrokecolor{currentstroke}%
\pgfsetdash{}{0pt}%
\pgfpathmoveto{\pgfqpoint{1.794301in}{1.482500in}}%
\pgfpathlineto{\pgfqpoint{1.859887in}{1.058777in}}%
\pgfpathlineto{\pgfqpoint{2.099225in}{0.359966in}}%
\pgfusepath{stroke}%
\end{pgfscope}%
\begin{pgfscope}%
\pgfpathrectangle{\pgfqpoint{0.135000in}{0.135000in}}{\pgfqpoint{2.914509in}{2.695000in}} %
\pgfusepath{clip}%
\pgfsetrectcap%
\pgfsetroundjoin%
\pgfsetlinewidth{1.003750pt}%
\definecolor{currentstroke}{rgb}{0.392157,0.584314,0.929412}%
\pgfsetstrokecolor{currentstroke}%
\pgfsetdash{}{0pt}%
\pgfpathmoveto{\pgfqpoint{2.080038in}{1.482500in}}%
\pgfpathlineto{\pgfqpoint{2.238376in}{1.306988in}}%
\pgfpathlineto{\pgfqpoint{2.816189in}{1.017531in}}%
\pgfusepath{stroke}%
\end{pgfscope}%
\begin{pgfscope}%
\pgfpathrectangle{\pgfqpoint{0.135000in}{0.135000in}}{\pgfqpoint{2.914509in}{2.695000in}} %
\pgfusepath{clip}%
\pgfsetbuttcap%
\pgfsetroundjoin%
\definecolor{currentfill}{rgb}{0.392157,0.584314,0.929412}%
\pgfsetfillcolor{currentfill}%
\pgfsetlinewidth{1.003750pt}%
\definecolor{currentstroke}{rgb}{0.392157,0.584314,0.929412}%
\pgfsetstrokecolor{currentstroke}%
\pgfsetdash{}{0pt}%
\pgfsys@defobject{currentmarker}{\pgfqpoint{-0.010417in}{-0.010417in}}{\pgfqpoint{0.010417in}{0.010417in}}{%
\pgfpathmoveto{\pgfqpoint{0.000000in}{-0.010417in}}%
\pgfpathcurveto{\pgfqpoint{0.002763in}{-0.010417in}}{\pgfqpoint{0.005412in}{-0.009319in}}{\pgfqpoint{0.007366in}{-0.007366in}}%
\pgfpathcurveto{\pgfqpoint{0.009319in}{-0.005412in}}{\pgfqpoint{0.010417in}{-0.002763in}}{\pgfqpoint{0.010417in}{0.000000in}}%
\pgfpathcurveto{\pgfqpoint{0.010417in}{0.002763in}}{\pgfqpoint{0.009319in}{0.005412in}}{\pgfqpoint{0.007366in}{0.007366in}}%
\pgfpathcurveto{\pgfqpoint{0.005412in}{0.009319in}}{\pgfqpoint{0.002763in}{0.010417in}}{\pgfqpoint{0.000000in}{0.010417in}}%
\pgfpathcurveto{\pgfqpoint{-0.002763in}{0.010417in}}{\pgfqpoint{-0.005412in}{0.009319in}}{\pgfqpoint{-0.007366in}{0.007366in}}%
\pgfpathcurveto{\pgfqpoint{-0.009319in}{0.005412in}}{\pgfqpoint{-0.010417in}{0.002763in}}{\pgfqpoint{-0.010417in}{0.000000in}}%
\pgfpathcurveto{\pgfqpoint{-0.010417in}{-0.002763in}}{\pgfqpoint{-0.009319in}{-0.005412in}}{\pgfqpoint{-0.007366in}{-0.007366in}}%
\pgfpathcurveto{\pgfqpoint{-0.005412in}{-0.009319in}}{\pgfqpoint{-0.002763in}{-0.010417in}}{\pgfqpoint{0.000000in}{-0.010417in}}%
\pgfpathclose%
\pgfusepath{stroke,fill}%
}%
\begin{pgfscope}%
\pgfsys@transformshift{2.080038in}{1.482500in}%
\pgfsys@useobject{currentmarker}{}%
\end{pgfscope}%
\begin{pgfscope}%
\pgfsys@transformshift{2.238376in}{1.306988in}%
\pgfsys@useobject{currentmarker}{}%
\end{pgfscope}%
\begin{pgfscope}%
\pgfsys@transformshift{2.816189in}{1.017531in}%
\pgfsys@useobject{currentmarker}{}%
\end{pgfscope}%
\begin{pgfscope}%
\pgfsys@transformshift{2.080038in}{1.482500in}%
\pgfsys@useobject{currentmarker}{}%
\end{pgfscope}%
\begin{pgfscope}%
\pgfsys@transformshift{2.238376in}{1.658012in}%
\pgfsys@useobject{currentmarker}{}%
\end{pgfscope}%
\begin{pgfscope}%
\pgfsys@transformshift{2.816189in}{1.947469in}%
\pgfsys@useobject{currentmarker}{}%
\end{pgfscope}%
\begin{pgfscope}%
\pgfsys@transformshift{1.794301in}{1.482500in}%
\pgfsys@useobject{currentmarker}{}%
\end{pgfscope}%
\begin{pgfscope}%
\pgfsys@transformshift{1.859887in}{1.906223in}%
\pgfsys@useobject{currentmarker}{}%
\end{pgfscope}%
\begin{pgfscope}%
\pgfsys@transformshift{2.099225in}{2.605034in}%
\pgfsys@useobject{currentmarker}{}%
\end{pgfscope}%
\begin{pgfscope}%
\pgfsys@transformshift{1.390208in}{1.482500in}%
\pgfsys@useobject{currentmarker}{}%
\end{pgfscope}%
\begin{pgfscope}%
\pgfsys@transformshift{1.324622in}{1.906223in}%
\pgfsys@useobject{currentmarker}{}%
\end{pgfscope}%
\begin{pgfscope}%
\pgfsys@transformshift{1.085284in}{2.605034in}%
\pgfsys@useobject{currentmarker}{}%
\end{pgfscope}%
\begin{pgfscope}%
\pgfsys@transformshift{1.104472in}{1.482500in}%
\pgfsys@useobject{currentmarker}{}%
\end{pgfscope}%
\begin{pgfscope}%
\pgfsys@transformshift{0.946133in}{1.658012in}%
\pgfsys@useobject{currentmarker}{}%
\end{pgfscope}%
\begin{pgfscope}%
\pgfsys@transformshift{0.368320in}{1.947469in}%
\pgfsys@useobject{currentmarker}{}%
\end{pgfscope}%
\begin{pgfscope}%
\pgfsys@transformshift{1.104472in}{1.482500in}%
\pgfsys@useobject{currentmarker}{}%
\end{pgfscope}%
\begin{pgfscope}%
\pgfsys@transformshift{0.946133in}{1.306988in}%
\pgfsys@useobject{currentmarker}{}%
\end{pgfscope}%
\begin{pgfscope}%
\pgfsys@transformshift{0.368320in}{1.017531in}%
\pgfsys@useobject{currentmarker}{}%
\end{pgfscope}%
\begin{pgfscope}%
\pgfsys@transformshift{1.390208in}{1.482500in}%
\pgfsys@useobject{currentmarker}{}%
\end{pgfscope}%
\begin{pgfscope}%
\pgfsys@transformshift{1.324622in}{1.058777in}%
\pgfsys@useobject{currentmarker}{}%
\end{pgfscope}%
\begin{pgfscope}%
\pgfsys@transformshift{1.085284in}{0.359966in}%
\pgfsys@useobject{currentmarker}{}%
\end{pgfscope}%
\begin{pgfscope}%
\pgfsys@transformshift{1.794301in}{1.482500in}%
\pgfsys@useobject{currentmarker}{}%
\end{pgfscope}%
\begin{pgfscope}%
\pgfsys@transformshift{1.859887in}{1.058777in}%
\pgfsys@useobject{currentmarker}{}%
\end{pgfscope}%
\begin{pgfscope}%
\pgfsys@transformshift{2.099225in}{0.359966in}%
\pgfsys@useobject{currentmarker}{}%
\end{pgfscope}%
\begin{pgfscope}%
\pgfsys@transformshift{2.080038in}{1.482500in}%
\pgfsys@useobject{currentmarker}{}%
\end{pgfscope}%
\begin{pgfscope}%
\pgfsys@transformshift{2.238376in}{1.306988in}%
\pgfsys@useobject{currentmarker}{}%
\end{pgfscope}%
\begin{pgfscope}%
\pgfsys@transformshift{2.816189in}{1.017531in}%
\pgfsys@useobject{currentmarker}{}%
\end{pgfscope}%
\end{pgfscope}%
\begin{pgfscope}%
\pgfpathrectangle{\pgfqpoint{0.135000in}{0.135000in}}{\pgfqpoint{2.914509in}{2.695000in}} %
\pgfusepath{clip}%
\pgfsetrectcap%
\pgfsetroundjoin%
\pgfsetlinewidth{1.003750pt}%
\definecolor{currentstroke}{rgb}{0.854902,0.439216,0.839216}%
\pgfsetstrokecolor{currentstroke}%
\pgfsetdash{}{0pt}%
\pgfpathmoveto{\pgfqpoint{2.161463in}{1.482500in}}%
\pgfpathlineto{\pgfqpoint{1.994746in}{1.332095in}}%
\pgfpathlineto{\pgfqpoint{1.592255in}{1.269796in}}%
\pgfpathlineto{\pgfqpoint{1.189764in}{1.332095in}}%
\pgfpathlineto{\pgfqpoint{1.023046in}{1.482500in}}%
\pgfpathlineto{\pgfqpoint{1.189764in}{1.632905in}}%
\pgfpathlineto{\pgfqpoint{1.592255in}{1.695204in}}%
\pgfpathlineto{\pgfqpoint{1.994746in}{1.632905in}}%
\pgfpathlineto{\pgfqpoint{2.161463in}{1.482500in}}%
\pgfusepath{stroke}%
\end{pgfscope}%
\begin{pgfscope}%
\pgfpathrectangle{\pgfqpoint{0.135000in}{0.135000in}}{\pgfqpoint{2.914509in}{2.695000in}} %
\pgfusepath{clip}%
\pgfsetrectcap%
\pgfsetroundjoin%
\pgfsetlinewidth{1.003750pt}%
\definecolor{currentstroke}{rgb}{0.854902,0.439216,0.839216}%
\pgfsetstrokecolor{currentstroke}%
\pgfsetdash{}{0pt}%
\pgfpathmoveto{\pgfqpoint{2.161463in}{1.482500in}}%
\pgfpathlineto{\pgfqpoint{1.994746in}{1.632905in}}%
\pgfpathlineto{\pgfqpoint{1.592255in}{1.695204in}}%
\pgfpathlineto{\pgfqpoint{1.189764in}{1.632905in}}%
\pgfpathlineto{\pgfqpoint{1.023046in}{1.482500in}}%
\pgfpathlineto{\pgfqpoint{1.189764in}{1.332095in}}%
\pgfpathlineto{\pgfqpoint{1.592255in}{1.269796in}}%
\pgfpathlineto{\pgfqpoint{1.994746in}{1.332095in}}%
\pgfpathlineto{\pgfqpoint{2.161463in}{1.482500in}}%
\pgfusepath{stroke}%
\end{pgfscope}%
\begin{pgfscope}%
\pgfpathrectangle{\pgfqpoint{0.135000in}{0.135000in}}{\pgfqpoint{2.914509in}{2.695000in}} %
\pgfusepath{clip}%
\pgfsetrectcap%
\pgfsetroundjoin%
\pgfsetlinewidth{1.003750pt}%
\definecolor{currentstroke}{rgb}{0.854902,0.439216,0.839216}%
\pgfsetstrokecolor{currentstroke}%
\pgfsetdash{}{0pt}%
\pgfpathmoveto{\pgfqpoint{2.531003in}{1.482500in}}%
\pgfpathlineto{\pgfqpoint{2.256050in}{2.031360in}}%
\pgfpathlineto{\pgfqpoint{1.592255in}{2.258705in}}%
\pgfpathlineto{\pgfqpoint{0.928459in}{2.031360in}}%
\pgfpathlineto{\pgfqpoint{0.653506in}{1.482500in}}%
\pgfpathlineto{\pgfqpoint{0.928459in}{0.933640in}}%
\pgfpathlineto{\pgfqpoint{1.592255in}{0.706295in}}%
\pgfpathlineto{\pgfqpoint{2.256050in}{0.933640in}}%
\pgfpathlineto{\pgfqpoint{2.531003in}{1.482500in}}%
\pgfusepath{stroke}%
\end{pgfscope}%
\begin{pgfscope}%
\pgfpathrectangle{\pgfqpoint{0.135000in}{0.135000in}}{\pgfqpoint{2.914509in}{2.695000in}} %
\pgfusepath{clip}%
\pgfsetrectcap%
\pgfsetroundjoin%
\pgfsetlinewidth{1.003750pt}%
\definecolor{currentstroke}{rgb}{0.854902,0.439216,0.839216}%
\pgfsetstrokecolor{currentstroke}%
\pgfsetdash{}{0pt}%
\pgfpathmoveto{\pgfqpoint{2.161463in}{1.482500in}}%
\pgfpathlineto{\pgfqpoint{2.161463in}{1.482500in}}%
\pgfpathlineto{\pgfqpoint{2.531003in}{1.482500in}}%
\pgfusepath{stroke}%
\end{pgfscope}%
\begin{pgfscope}%
\pgfpathrectangle{\pgfqpoint{0.135000in}{0.135000in}}{\pgfqpoint{2.914509in}{2.695000in}} %
\pgfusepath{clip}%
\pgfsetrectcap%
\pgfsetroundjoin%
\pgfsetlinewidth{1.003750pt}%
\definecolor{currentstroke}{rgb}{0.854902,0.439216,0.839216}%
\pgfsetstrokecolor{currentstroke}%
\pgfsetdash{}{0pt}%
\pgfpathmoveto{\pgfqpoint{1.994746in}{1.332095in}}%
\pgfpathlineto{\pgfqpoint{1.994746in}{1.632905in}}%
\pgfpathlineto{\pgfqpoint{2.256050in}{2.031360in}}%
\pgfusepath{stroke}%
\end{pgfscope}%
\begin{pgfscope}%
\pgfpathrectangle{\pgfqpoint{0.135000in}{0.135000in}}{\pgfqpoint{2.914509in}{2.695000in}} %
\pgfusepath{clip}%
\pgfsetrectcap%
\pgfsetroundjoin%
\pgfsetlinewidth{1.003750pt}%
\definecolor{currentstroke}{rgb}{0.854902,0.439216,0.839216}%
\pgfsetstrokecolor{currentstroke}%
\pgfsetdash{}{0pt}%
\pgfpathmoveto{\pgfqpoint{1.592255in}{1.269796in}}%
\pgfpathlineto{\pgfqpoint{1.592255in}{1.695204in}}%
\pgfpathlineto{\pgfqpoint{1.592255in}{2.258705in}}%
\pgfusepath{stroke}%
\end{pgfscope}%
\begin{pgfscope}%
\pgfpathrectangle{\pgfqpoint{0.135000in}{0.135000in}}{\pgfqpoint{2.914509in}{2.695000in}} %
\pgfusepath{clip}%
\pgfsetrectcap%
\pgfsetroundjoin%
\pgfsetlinewidth{1.003750pt}%
\definecolor{currentstroke}{rgb}{0.854902,0.439216,0.839216}%
\pgfsetstrokecolor{currentstroke}%
\pgfsetdash{}{0pt}%
\pgfpathmoveto{\pgfqpoint{1.189764in}{1.332095in}}%
\pgfpathlineto{\pgfqpoint{1.189764in}{1.632905in}}%
\pgfpathlineto{\pgfqpoint{0.928459in}{2.031360in}}%
\pgfusepath{stroke}%
\end{pgfscope}%
\begin{pgfscope}%
\pgfpathrectangle{\pgfqpoint{0.135000in}{0.135000in}}{\pgfqpoint{2.914509in}{2.695000in}} %
\pgfusepath{clip}%
\pgfsetrectcap%
\pgfsetroundjoin%
\pgfsetlinewidth{1.003750pt}%
\definecolor{currentstroke}{rgb}{0.854902,0.439216,0.839216}%
\pgfsetstrokecolor{currentstroke}%
\pgfsetdash{}{0pt}%
\pgfpathmoveto{\pgfqpoint{1.023046in}{1.482500in}}%
\pgfpathlineto{\pgfqpoint{1.023046in}{1.482500in}}%
\pgfpathlineto{\pgfqpoint{0.653506in}{1.482500in}}%
\pgfusepath{stroke}%
\end{pgfscope}%
\begin{pgfscope}%
\pgfpathrectangle{\pgfqpoint{0.135000in}{0.135000in}}{\pgfqpoint{2.914509in}{2.695000in}} %
\pgfusepath{clip}%
\pgfsetrectcap%
\pgfsetroundjoin%
\pgfsetlinewidth{1.003750pt}%
\definecolor{currentstroke}{rgb}{0.854902,0.439216,0.839216}%
\pgfsetstrokecolor{currentstroke}%
\pgfsetdash{}{0pt}%
\pgfpathmoveto{\pgfqpoint{1.189764in}{1.632905in}}%
\pgfpathlineto{\pgfqpoint{1.189764in}{1.332095in}}%
\pgfpathlineto{\pgfqpoint{0.928459in}{0.933640in}}%
\pgfusepath{stroke}%
\end{pgfscope}%
\begin{pgfscope}%
\pgfpathrectangle{\pgfqpoint{0.135000in}{0.135000in}}{\pgfqpoint{2.914509in}{2.695000in}} %
\pgfusepath{clip}%
\pgfsetrectcap%
\pgfsetroundjoin%
\pgfsetlinewidth{1.003750pt}%
\definecolor{currentstroke}{rgb}{0.854902,0.439216,0.839216}%
\pgfsetstrokecolor{currentstroke}%
\pgfsetdash{}{0pt}%
\pgfpathmoveto{\pgfqpoint{1.592255in}{1.695204in}}%
\pgfpathlineto{\pgfqpoint{1.592255in}{1.269796in}}%
\pgfpathlineto{\pgfqpoint{1.592255in}{0.706295in}}%
\pgfusepath{stroke}%
\end{pgfscope}%
\begin{pgfscope}%
\pgfpathrectangle{\pgfqpoint{0.135000in}{0.135000in}}{\pgfqpoint{2.914509in}{2.695000in}} %
\pgfusepath{clip}%
\pgfsetrectcap%
\pgfsetroundjoin%
\pgfsetlinewidth{1.003750pt}%
\definecolor{currentstroke}{rgb}{0.854902,0.439216,0.839216}%
\pgfsetstrokecolor{currentstroke}%
\pgfsetdash{}{0pt}%
\pgfpathmoveto{\pgfqpoint{1.994746in}{1.632905in}}%
\pgfpathlineto{\pgfqpoint{1.994746in}{1.332095in}}%
\pgfpathlineto{\pgfqpoint{2.256050in}{0.933640in}}%
\pgfusepath{stroke}%
\end{pgfscope}%
\begin{pgfscope}%
\pgfpathrectangle{\pgfqpoint{0.135000in}{0.135000in}}{\pgfqpoint{2.914509in}{2.695000in}} %
\pgfusepath{clip}%
\pgfsetrectcap%
\pgfsetroundjoin%
\pgfsetlinewidth{1.003750pt}%
\definecolor{currentstroke}{rgb}{0.854902,0.439216,0.839216}%
\pgfsetstrokecolor{currentstroke}%
\pgfsetdash{}{0pt}%
\pgfpathmoveto{\pgfqpoint{2.161463in}{1.482500in}}%
\pgfpathlineto{\pgfqpoint{2.161463in}{1.482500in}}%
\pgfpathlineto{\pgfqpoint{2.531003in}{1.482500in}}%
\pgfusepath{stroke}%
\end{pgfscope}%
\begin{pgfscope}%
\pgfpathrectangle{\pgfqpoint{0.135000in}{0.135000in}}{\pgfqpoint{2.914509in}{2.695000in}} %
\pgfusepath{clip}%
\pgfsetbuttcap%
\pgfsetroundjoin%
\definecolor{currentfill}{rgb}{0.854902,0.439216,0.839216}%
\pgfsetfillcolor{currentfill}%
\pgfsetlinewidth{1.003750pt}%
\definecolor{currentstroke}{rgb}{0.854902,0.439216,0.839216}%
\pgfsetstrokecolor{currentstroke}%
\pgfsetdash{}{0pt}%
\pgfsys@defobject{currentmarker}{\pgfqpoint{-0.010417in}{-0.010417in}}{\pgfqpoint{0.010417in}{0.010417in}}{%
\pgfpathmoveto{\pgfqpoint{0.000000in}{-0.010417in}}%
\pgfpathcurveto{\pgfqpoint{0.002763in}{-0.010417in}}{\pgfqpoint{0.005412in}{-0.009319in}}{\pgfqpoint{0.007366in}{-0.007366in}}%
\pgfpathcurveto{\pgfqpoint{0.009319in}{-0.005412in}}{\pgfqpoint{0.010417in}{-0.002763in}}{\pgfqpoint{0.010417in}{0.000000in}}%
\pgfpathcurveto{\pgfqpoint{0.010417in}{0.002763in}}{\pgfqpoint{0.009319in}{0.005412in}}{\pgfqpoint{0.007366in}{0.007366in}}%
\pgfpathcurveto{\pgfqpoint{0.005412in}{0.009319in}}{\pgfqpoint{0.002763in}{0.010417in}}{\pgfqpoint{0.000000in}{0.010417in}}%
\pgfpathcurveto{\pgfqpoint{-0.002763in}{0.010417in}}{\pgfqpoint{-0.005412in}{0.009319in}}{\pgfqpoint{-0.007366in}{0.007366in}}%
\pgfpathcurveto{\pgfqpoint{-0.009319in}{0.005412in}}{\pgfqpoint{-0.010417in}{0.002763in}}{\pgfqpoint{-0.010417in}{0.000000in}}%
\pgfpathcurveto{\pgfqpoint{-0.010417in}{-0.002763in}}{\pgfqpoint{-0.009319in}{-0.005412in}}{\pgfqpoint{-0.007366in}{-0.007366in}}%
\pgfpathcurveto{\pgfqpoint{-0.005412in}{-0.009319in}}{\pgfqpoint{-0.002763in}{-0.010417in}}{\pgfqpoint{0.000000in}{-0.010417in}}%
\pgfpathclose%
\pgfusepath{stroke,fill}%
}%
\begin{pgfscope}%
\pgfsys@transformshift{2.161463in}{1.482500in}%
\pgfsys@useobject{currentmarker}{}%
\end{pgfscope}%
\begin{pgfscope}%
\pgfsys@transformshift{2.161463in}{1.482500in}%
\pgfsys@useobject{currentmarker}{}%
\end{pgfscope}%
\begin{pgfscope}%
\pgfsys@transformshift{2.531003in}{1.482500in}%
\pgfsys@useobject{currentmarker}{}%
\end{pgfscope}%
\begin{pgfscope}%
\pgfsys@transformshift{1.994746in}{1.332095in}%
\pgfsys@useobject{currentmarker}{}%
\end{pgfscope}%
\begin{pgfscope}%
\pgfsys@transformshift{1.994746in}{1.632905in}%
\pgfsys@useobject{currentmarker}{}%
\end{pgfscope}%
\begin{pgfscope}%
\pgfsys@transformshift{2.256050in}{2.031360in}%
\pgfsys@useobject{currentmarker}{}%
\end{pgfscope}%
\begin{pgfscope}%
\pgfsys@transformshift{1.592255in}{1.269796in}%
\pgfsys@useobject{currentmarker}{}%
\end{pgfscope}%
\begin{pgfscope}%
\pgfsys@transformshift{1.592255in}{1.695204in}%
\pgfsys@useobject{currentmarker}{}%
\end{pgfscope}%
\begin{pgfscope}%
\pgfsys@transformshift{1.592255in}{2.258705in}%
\pgfsys@useobject{currentmarker}{}%
\end{pgfscope}%
\begin{pgfscope}%
\pgfsys@transformshift{1.189764in}{1.332095in}%
\pgfsys@useobject{currentmarker}{}%
\end{pgfscope}%
\begin{pgfscope}%
\pgfsys@transformshift{1.189764in}{1.632905in}%
\pgfsys@useobject{currentmarker}{}%
\end{pgfscope}%
\begin{pgfscope}%
\pgfsys@transformshift{0.928459in}{2.031360in}%
\pgfsys@useobject{currentmarker}{}%
\end{pgfscope}%
\begin{pgfscope}%
\pgfsys@transformshift{1.023046in}{1.482500in}%
\pgfsys@useobject{currentmarker}{}%
\end{pgfscope}%
\begin{pgfscope}%
\pgfsys@transformshift{1.023046in}{1.482500in}%
\pgfsys@useobject{currentmarker}{}%
\end{pgfscope}%
\begin{pgfscope}%
\pgfsys@transformshift{0.653506in}{1.482500in}%
\pgfsys@useobject{currentmarker}{}%
\end{pgfscope}%
\begin{pgfscope}%
\pgfsys@transformshift{1.189764in}{1.632905in}%
\pgfsys@useobject{currentmarker}{}%
\end{pgfscope}%
\begin{pgfscope}%
\pgfsys@transformshift{1.189764in}{1.332095in}%
\pgfsys@useobject{currentmarker}{}%
\end{pgfscope}%
\begin{pgfscope}%
\pgfsys@transformshift{0.928459in}{0.933640in}%
\pgfsys@useobject{currentmarker}{}%
\end{pgfscope}%
\begin{pgfscope}%
\pgfsys@transformshift{1.592255in}{1.695204in}%
\pgfsys@useobject{currentmarker}{}%
\end{pgfscope}%
\begin{pgfscope}%
\pgfsys@transformshift{1.592255in}{1.269796in}%
\pgfsys@useobject{currentmarker}{}%
\end{pgfscope}%
\begin{pgfscope}%
\pgfsys@transformshift{1.592255in}{0.706295in}%
\pgfsys@useobject{currentmarker}{}%
\end{pgfscope}%
\begin{pgfscope}%
\pgfsys@transformshift{1.994746in}{1.632905in}%
\pgfsys@useobject{currentmarker}{}%
\end{pgfscope}%
\begin{pgfscope}%
\pgfsys@transformshift{1.994746in}{1.332095in}%
\pgfsys@useobject{currentmarker}{}%
\end{pgfscope}%
\begin{pgfscope}%
\pgfsys@transformshift{2.256050in}{0.933640in}%
\pgfsys@useobject{currentmarker}{}%
\end{pgfscope}%
\begin{pgfscope}%
\pgfsys@transformshift{2.161463in}{1.482500in}%
\pgfsys@useobject{currentmarker}{}%
\end{pgfscope}%
\begin{pgfscope}%
\pgfsys@transformshift{2.161463in}{1.482500in}%
\pgfsys@useobject{currentmarker}{}%
\end{pgfscope}%
\begin{pgfscope}%
\pgfsys@transformshift{2.531003in}{1.482500in}%
\pgfsys@useobject{currentmarker}{}%
\end{pgfscope}%
\end{pgfscope}%
\end{pgfpicture}%
\makeatother%
\endgroup%

%% file: ic-symmetric-with-ellipse1_image.pgf
%% Creator: Matplotlib, PGF backend
%%
%% To include the figure in your LaTeX document, write
%%   \input{<filename>.pgf}
%%
%% Make sure the required packages are loaded in your preamble
%%   \usepackage{pgf}
%%
%% Figures using additional raster images can only be included by \input if
%% they are in the same directory as the main LaTeX file. For loading figures
%% from other directories you can use the `import` package
%% and then include the figures with
%%   \import{<path to file>}{<filename>.pgf}
%%
%% Matplotlib used the following preamble
%%   \usepackage{fontspec}
%%   \setmainfont{DejaVu Serif}
%%   \setsansfont{DejaVu Sans}
%%   \setmonofont{DejaVu Sans Mono}
%%
\begingroup%
\makeatletter%
\begin{pgfpicture}%
\pgfpathrectangle{\pgfpointorigin}{\pgfqpoint{3.918914in}{3.170273in}}%
\pgfusepath{use as bounding box, clip}%
\begin{pgfscope}%
\pgfsetbuttcap%
\pgfsetmiterjoin%
\definecolor{currentfill}{rgb}{1.000000,1.000000,1.000000}%
\pgfsetfillcolor{currentfill}%
\pgfsetlinewidth{0.000000pt}%
\definecolor{currentstroke}{rgb}{1.000000,1.000000,1.000000}%
\pgfsetstrokecolor{currentstroke}%
\pgfsetdash{}{0pt}%
\pgfpathmoveto{\pgfqpoint{0.000000in}{0.000000in}}%
\pgfpathlineto{\pgfqpoint{3.918914in}{0.000000in}}%
\pgfpathlineto{\pgfqpoint{3.918914in}{3.170273in}}%
\pgfpathlineto{\pgfqpoint{0.000000in}{3.170273in}}%
\pgfpathclose%
\pgfusepath{fill}%
\end{pgfscope}%
\begin{pgfscope}%
\pgfsetbuttcap%
\pgfsetmiterjoin%
\definecolor{currentfill}{rgb}{1.000000,1.000000,1.000000}%
\pgfsetfillcolor{currentfill}%
\pgfsetlinewidth{0.000000pt}%
\definecolor{currentstroke}{rgb}{0.000000,0.000000,0.000000}%
\pgfsetstrokecolor{currentstroke}%
\pgfsetstrokeopacity{0.000000}%
\pgfsetdash{}{0pt}%
\pgfpathmoveto{\pgfqpoint{0.135000in}{0.135000in}}%
\pgfpathlineto{\pgfqpoint{3.634549in}{0.135000in}}%
\pgfpathlineto{\pgfqpoint{3.634549in}{2.830000in}}%
\pgfpathlineto{\pgfqpoint{0.135000in}{2.830000in}}%
\pgfpathclose%
\pgfusepath{fill}%
\end{pgfscope}%
\begin{pgfscope}%
\pgfpathrectangle{\pgfqpoint{0.135000in}{0.135000in}}{\pgfqpoint{3.499549in}{2.695000in}} %
\pgfusepath{clip}%
\pgfsetrectcap%
\pgfsetroundjoin%
\pgfsetlinewidth{0.702625pt}%
\definecolor{currentstroke}{rgb}{0.000000,0.000000,0.000000}%
\pgfsetstrokecolor{currentstroke}%
\pgfsetdash{}{0pt}%
\pgfpathmoveto{\pgfqpoint{0.135000in}{1.482500in}}%
\pgfpathlineto{\pgfqpoint{3.634549in}{1.482500in}}%
\pgfusepath{stroke}%
\end{pgfscope}%
\begin{pgfscope}%
\pgfpathrectangle{\pgfqpoint{0.135000in}{0.135000in}}{\pgfqpoint{3.499549in}{2.695000in}} %
\pgfusepath{clip}%
\pgfsetrectcap%
\pgfsetroundjoin%
\pgfsetlinewidth{0.702625pt}%
\definecolor{currentstroke}{rgb}{0.000000,0.000000,0.000000}%
\pgfsetstrokecolor{currentstroke}%
\pgfsetdash{}{0pt}%
\pgfpathmoveto{\pgfqpoint{1.884775in}{0.135000in}}%
\pgfpathlineto{\pgfqpoint{1.884775in}{2.830000in}}%
\pgfusepath{stroke}%
\end{pgfscope}%
\begin{pgfscope}%
\pgfsetbuttcap%
\pgfsetmiterjoin%
\definecolor{currentfill}{rgb}{0.000000,0.000000,0.000000}%
\pgfsetfillcolor{currentfill}%
\pgfsetlinewidth{1.003750pt}%
\definecolor{currentstroke}{rgb}{0.000000,0.000000,0.000000}%
\pgfsetstrokecolor{currentstroke}%
\pgfsetdash{}{0pt}%
\pgfpathmoveto{\pgfqpoint{3.634549in}{1.482500in}}%
\pgfpathlineto{\pgfqpoint{3.599554in}{1.469025in}}%
\pgfpathlineto{\pgfqpoint{3.610052in}{1.481923in}}%
\pgfpathlineto{\pgfqpoint{0.135000in}{1.481923in}}%
\pgfpathlineto{\pgfqpoint{0.135000in}{1.483077in}}%
\pgfpathlineto{\pgfqpoint{3.610052in}{1.483077in}}%
\pgfpathlineto{\pgfqpoint{3.599554in}{1.495975in}}%
\pgfpathclose%
\pgfusepath{stroke,fill}%
\end{pgfscope}%
\begin{pgfscope}%
\pgfsetbuttcap%
\pgfsetmiterjoin%
\definecolor{currentfill}{rgb}{0.000000,0.000000,0.000000}%
\pgfsetfillcolor{currentfill}%
\pgfsetlinewidth{1.003750pt}%
\definecolor{currentstroke}{rgb}{0.000000,0.000000,0.000000}%
\pgfsetstrokecolor{currentstroke}%
\pgfsetdash{}{0pt}%
\pgfpathmoveto{\pgfqpoint{1.884775in}{2.830000in}}%
\pgfpathlineto{\pgfqpoint{1.898250in}{2.795005in}}%
\pgfpathlineto{\pgfqpoint{1.885352in}{2.805503in}}%
\pgfpathlineto{\pgfqpoint{1.885352in}{0.135000in}}%
\pgfpathlineto{\pgfqpoint{1.884197in}{0.135000in}}%
\pgfpathlineto{\pgfqpoint{1.884197in}{2.805503in}}%
\pgfpathlineto{\pgfqpoint{1.871300in}{2.795005in}}%
\pgfpathclose%
\pgfusepath{stroke,fill}%
\end{pgfscope}%
\begin{pgfscope}%
\pgfpathrectangle{\pgfqpoint{0.135000in}{0.135000in}}{\pgfqpoint{3.499549in}{2.695000in}} %
\pgfusepath{clip}%
\pgfsetrectcap%
\pgfsetroundjoin%
\pgfsetlinewidth{0.803000pt}%
\definecolor{currentstroke}{rgb}{0.000000,0.000000,0.000000}%
\pgfsetstrokecolor{currentstroke}%
\pgfsetdash{}{0pt}%
\pgfpathmoveto{\pgfqpoint{0.135000in}{1.482500in}}%
\pgfpathlineto{\pgfqpoint{3.634549in}{1.482500in}}%
\pgfusepath{stroke}%
\end{pgfscope}%
\begin{pgfscope}%
\pgfpathrectangle{\pgfqpoint{0.135000in}{0.135000in}}{\pgfqpoint{3.499549in}{2.695000in}} %
\pgfusepath{clip}%
\pgfsetrectcap%
\pgfsetroundjoin%
\pgfsetlinewidth{0.803000pt}%
\definecolor{currentstroke}{rgb}{0.000000,0.000000,0.000000}%
\pgfsetstrokecolor{currentstroke}%
\pgfsetdash{}{0pt}%
\pgfpathmoveto{\pgfqpoint{1.884775in}{0.135000in}}%
\pgfpathlineto{\pgfqpoint{1.884775in}{2.830000in}}%
\pgfusepath{stroke}%
\end{pgfscope}%
\begin{pgfscope}%
\pgfpathrectangle{\pgfqpoint{0.135000in}{0.135000in}}{\pgfqpoint{3.499549in}{2.695000in}} %
\pgfusepath{clip}%
\pgfsetrectcap%
\pgfsetroundjoin%
\pgfsetlinewidth{0.803000pt}%
\definecolor{currentstroke}{rgb}{0.000000,0.000000,0.000000}%
\pgfsetstrokecolor{currentstroke}%
\pgfsetdash{}{0pt}%
\pgfpathmoveto{\pgfqpoint{0.135000in}{1.482500in}}%
\pgfpathlineto{\pgfqpoint{3.634549in}{1.482500in}}%
\pgfusepath{stroke}%
\end{pgfscope}%
\begin{pgfscope}%
\pgfpathrectangle{\pgfqpoint{0.135000in}{0.135000in}}{\pgfqpoint{3.499549in}{2.695000in}} %
\pgfusepath{clip}%
\pgfsetrectcap%
\pgfsetroundjoin%
\pgfsetlinewidth{0.803000pt}%
\definecolor{currentstroke}{rgb}{0.000000,0.000000,0.000000}%
\pgfsetstrokecolor{currentstroke}%
\pgfsetdash{}{0pt}%
\pgfpathmoveto{\pgfqpoint{1.884775in}{0.135000in}}%
\pgfpathlineto{\pgfqpoint{1.884775in}{2.830000in}}%
\pgfusepath{stroke}%
\end{pgfscope}%
\begin{pgfscope}%
\pgftext[x=1.884775in,y=2.964750in,,base]{\sffamily\fontsize{10.000000}{12.000000}\selectfont \(\displaystyle y\)}%
\end{pgfscope}%
\begin{pgfscope}%
\pgftext[x=3.739536in,y=1.482500in,left,]{\sffamily\fontsize{10.000000}{12.000000}\selectfont \(\displaystyle x\)}%
\end{pgfscope}%
\begin{pgfscope}%
\pgfpathrectangle{\pgfqpoint{0.135000in}{0.135000in}}{\pgfqpoint{3.499549in}{2.695000in}} %
\pgfusepath{clip}%
\pgfsetbuttcap%
\pgfsetroundjoin%
\definecolor{currentfill}{rgb}{0.700000,0.700000,0.700000}%
\pgfsetfillcolor{currentfill}%
\pgfsetlinewidth{1.003750pt}%
\definecolor{currentstroke}{rgb}{0.700000,0.700000,0.700000}%
\pgfsetstrokecolor{currentstroke}%
\pgfsetdash{}{0pt}%
\pgfsys@defobject{currentmarker}{\pgfqpoint{-0.010417in}{-0.010417in}}{\pgfqpoint{0.010417in}{0.010417in}}{%
\pgfpathmoveto{\pgfqpoint{0.000000in}{-0.010417in}}%
\pgfpathcurveto{\pgfqpoint{0.002763in}{-0.010417in}}{\pgfqpoint{0.005412in}{-0.009319in}}{\pgfqpoint{0.007366in}{-0.007366in}}%
\pgfpathcurveto{\pgfqpoint{0.009319in}{-0.005412in}}{\pgfqpoint{0.010417in}{-0.002763in}}{\pgfqpoint{0.010417in}{0.000000in}}%
\pgfpathcurveto{\pgfqpoint{0.010417in}{0.002763in}}{\pgfqpoint{0.009319in}{0.005412in}}{\pgfqpoint{0.007366in}{0.007366in}}%
\pgfpathcurveto{\pgfqpoint{0.005412in}{0.009319in}}{\pgfqpoint{0.002763in}{0.010417in}}{\pgfqpoint{0.000000in}{0.010417in}}%
\pgfpathcurveto{\pgfqpoint{-0.002763in}{0.010417in}}{\pgfqpoint{-0.005412in}{0.009319in}}{\pgfqpoint{-0.007366in}{0.007366in}}%
\pgfpathcurveto{\pgfqpoint{-0.009319in}{0.005412in}}{\pgfqpoint{-0.010417in}{0.002763in}}{\pgfqpoint{-0.010417in}{0.000000in}}%
\pgfpathcurveto{\pgfqpoint{-0.010417in}{-0.002763in}}{\pgfqpoint{-0.009319in}{-0.005412in}}{\pgfqpoint{-0.007366in}{-0.007366in}}%
\pgfpathcurveto{\pgfqpoint{-0.005412in}{-0.009319in}}{\pgfqpoint{-0.002763in}{-0.010417in}}{\pgfqpoint{0.000000in}{-0.010417in}}%
\pgfpathclose%
\pgfusepath{stroke,fill}%
}%
\begin{pgfscope}%
\pgfsys@transformshift{1.884775in}{0.827711in}%
\pgfsys@useobject{currentmarker}{}%
\end{pgfscope}%
\begin{pgfscope}%
\pgfsys@transformshift{1.884775in}{0.715884in}%
\pgfsys@useobject{currentmarker}{}%
\end{pgfscope}%
\begin{pgfscope}%
\pgfsys@transformshift{1.884775in}{0.294070in}%
\pgfsys@useobject{currentmarker}{}%
\end{pgfscope}%
\begin{pgfscope}%
\pgfsys@transformshift{1.194402in}{0.923226in}%
\pgfsys@useobject{currentmarker}{}%
\end{pgfscope}%
\begin{pgfscope}%
\pgfsys@transformshift{1.170500in}{0.827711in}%
\pgfsys@useobject{currentmarker}{}%
\end{pgfscope}%
\begin{pgfscope}%
\pgfsys@transformshift{1.057510in}{0.467428in}%
\pgfsys@useobject{currentmarker}{}%
\end{pgfscope}%
\begin{pgfscope}%
\pgfsys@transformshift{0.776957in}{1.121731in}%
\pgfsys@useobject{currentmarker}{}%
\end{pgfscope}%
\begin{pgfscope}%
\pgfsys@transformshift{0.738604in}{1.060118in}%
\pgfsys@useobject{currentmarker}{}%
\end{pgfscope}%
\begin{pgfscope}%
\pgfsys@transformshift{0.557292in}{0.827711in}%
\pgfsys@useobject{currentmarker}{}%
\end{pgfscope}%
\begin{pgfscope}%
\pgfsys@transformshift{0.601712in}{1.314531in}%
\pgfsys@useobject{currentmarker}{}%
\end{pgfscope}%
\begin{pgfscope}%
\pgfsys@transformshift{0.557292in}{1.285845in}%
\pgfsys@useobject{currentmarker}{}%
\end{pgfscope}%
\begin{pgfscope}%
\pgfsys@transformshift{0.347299in}{1.177640in}%
\pgfsys@useobject{currentmarker}{}%
\end{pgfscope}%
\begin{pgfscope}%
\pgfsys@transformshift{0.557292in}{1.482500in}%
\pgfsys@useobject{currentmarker}{}%
\end{pgfscope}%
\begin{pgfscope}%
\pgfsys@transformshift{0.511334in}{1.482500in}%
\pgfsys@useobject{currentmarker}{}%
\end{pgfscope}%
\begin{pgfscope}%
\pgfsys@transformshift{0.294070in}{1.482500in}%
\pgfsys@useobject{currentmarker}{}%
\end{pgfscope}%
\begin{pgfscope}%
\pgfsys@transformshift{0.601712in}{1.650469in}%
\pgfsys@useobject{currentmarker}{}%
\end{pgfscope}%
\begin{pgfscope}%
\pgfsys@transformshift{0.557292in}{1.679155in}%
\pgfsys@useobject{currentmarker}{}%
\end{pgfscope}%
\begin{pgfscope}%
\pgfsys@transformshift{0.347299in}{1.787360in}%
\pgfsys@useobject{currentmarker}{}%
\end{pgfscope}%
\begin{pgfscope}%
\pgfsys@transformshift{0.776957in}{1.843269in}%
\pgfsys@useobject{currentmarker}{}%
\end{pgfscope}%
\begin{pgfscope}%
\pgfsys@transformshift{0.738604in}{1.904882in}%
\pgfsys@useobject{currentmarker}{}%
\end{pgfscope}%
\begin{pgfscope}%
\pgfsys@transformshift{0.557292in}{2.137289in}%
\pgfsys@useobject{currentmarker}{}%
\end{pgfscope}%
\begin{pgfscope}%
\pgfsys@transformshift{1.194402in}{2.041774in}%
\pgfsys@useobject{currentmarker}{}%
\end{pgfscope}%
\begin{pgfscope}%
\pgfsys@transformshift{1.170500in}{2.137289in}%
\pgfsys@useobject{currentmarker}{}%
\end{pgfscope}%
\begin{pgfscope}%
\pgfsys@transformshift{1.057510in}{2.497572in}%
\pgfsys@useobject{currentmarker}{}%
\end{pgfscope}%
\begin{pgfscope}%
\pgfsys@transformshift{1.884775in}{2.137289in}%
\pgfsys@useobject{currentmarker}{}%
\end{pgfscope}%
\begin{pgfscope}%
\pgfsys@transformshift{1.884775in}{2.249116in}%
\pgfsys@useobject{currentmarker}{}%
\end{pgfscope}%
\begin{pgfscope}%
\pgfsys@transformshift{1.884775in}{2.670930in}%
\pgfsys@useobject{currentmarker}{}%
\end{pgfscope}%
\begin{pgfscope}%
\pgfsys@transformshift{2.575148in}{2.041774in}%
\pgfsys@useobject{currentmarker}{}%
\end{pgfscope}%
\begin{pgfscope}%
\pgfsys@transformshift{2.599049in}{2.137289in}%
\pgfsys@useobject{currentmarker}{}%
\end{pgfscope}%
\begin{pgfscope}%
\pgfsys@transformshift{2.712039in}{2.497572in}%
\pgfsys@useobject{currentmarker}{}%
\end{pgfscope}%
\begin{pgfscope}%
\pgfsys@transformshift{2.992592in}{1.843269in}%
\pgfsys@useobject{currentmarker}{}%
\end{pgfscope}%
\begin{pgfscope}%
\pgfsys@transformshift{3.030945in}{1.904882in}%
\pgfsys@useobject{currentmarker}{}%
\end{pgfscope}%
\begin{pgfscope}%
\pgfsys@transformshift{3.212257in}{2.137289in}%
\pgfsys@useobject{currentmarker}{}%
\end{pgfscope}%
\begin{pgfscope}%
\pgfsys@transformshift{3.167837in}{1.650469in}%
\pgfsys@useobject{currentmarker}{}%
\end{pgfscope}%
\begin{pgfscope}%
\pgfsys@transformshift{3.212257in}{1.679155in}%
\pgfsys@useobject{currentmarker}{}%
\end{pgfscope}%
\begin{pgfscope}%
\pgfsys@transformshift{3.422250in}{1.787360in}%
\pgfsys@useobject{currentmarker}{}%
\end{pgfscope}%
\begin{pgfscope}%
\pgfsys@transformshift{3.212257in}{1.482500in}%
\pgfsys@useobject{currentmarker}{}%
\end{pgfscope}%
\begin{pgfscope}%
\pgfsys@transformshift{3.258215in}{1.482500in}%
\pgfsys@useobject{currentmarker}{}%
\end{pgfscope}%
\begin{pgfscope}%
\pgfsys@transformshift{3.475479in}{1.482500in}%
\pgfsys@useobject{currentmarker}{}%
\end{pgfscope}%
\begin{pgfscope}%
\pgfsys@transformshift{3.167837in}{1.314531in}%
\pgfsys@useobject{currentmarker}{}%
\end{pgfscope}%
\begin{pgfscope}%
\pgfsys@transformshift{3.212257in}{1.285845in}%
\pgfsys@useobject{currentmarker}{}%
\end{pgfscope}%
\begin{pgfscope}%
\pgfsys@transformshift{3.422250in}{1.177640in}%
\pgfsys@useobject{currentmarker}{}%
\end{pgfscope}%
\begin{pgfscope}%
\pgfsys@transformshift{2.992592in}{1.121731in}%
\pgfsys@useobject{currentmarker}{}%
\end{pgfscope}%
\begin{pgfscope}%
\pgfsys@transformshift{3.030945in}{1.060118in}%
\pgfsys@useobject{currentmarker}{}%
\end{pgfscope}%
\begin{pgfscope}%
\pgfsys@transformshift{3.212257in}{0.827711in}%
\pgfsys@useobject{currentmarker}{}%
\end{pgfscope}%
\begin{pgfscope}%
\pgfsys@transformshift{2.575148in}{0.923226in}%
\pgfsys@useobject{currentmarker}{}%
\end{pgfscope}%
\begin{pgfscope}%
\pgfsys@transformshift{2.599049in}{0.827711in}%
\pgfsys@useobject{currentmarker}{}%
\end{pgfscope}%
\begin{pgfscope}%
\pgfsys@transformshift{2.712039in}{0.467428in}%
\pgfsys@useobject{currentmarker}{}%
\end{pgfscope}%
\begin{pgfscope}%
\pgfsys@transformshift{1.884775in}{0.827711in}%
\pgfsys@useobject{currentmarker}{}%
\end{pgfscope}%
\begin{pgfscope}%
\pgfsys@transformshift{1.884775in}{0.715884in}%
\pgfsys@useobject{currentmarker}{}%
\end{pgfscope}%
\begin{pgfscope}%
\pgfsys@transformshift{1.884775in}{0.294070in}%
\pgfsys@useobject{currentmarker}{}%
\end{pgfscope}%
\end{pgfscope}%
\begin{pgfscope}%
\pgfpathrectangle{\pgfqpoint{0.135000in}{0.135000in}}{\pgfqpoint{3.499549in}{2.695000in}} %
\pgfusepath{clip}%
\pgfsetrectcap%
\pgfsetroundjoin%
\pgfsetlinewidth{1.003750pt}%
\definecolor{currentstroke}{rgb}{0.000000,0.000000,1.000000}%
\pgfsetstrokecolor{currentstroke}%
\pgfsetdash{}{0pt}%
\pgfpathmoveto{\pgfqpoint{1.884775in}{0.827711in}}%
\pgfpathlineto{\pgfqpoint{0.776957in}{1.121731in}}%
\pgfpathlineto{\pgfqpoint{0.557292in}{1.482500in}}%
\pgfpathlineto{\pgfqpoint{0.776957in}{1.843269in}}%
\pgfpathlineto{\pgfqpoint{1.884775in}{2.137289in}}%
\pgfpathlineto{\pgfqpoint{2.992592in}{1.843269in}}%
\pgfpathlineto{\pgfqpoint{3.212257in}{1.482500in}}%
\pgfpathlineto{\pgfqpoint{2.992592in}{1.121731in}}%
\pgfpathlineto{\pgfqpoint{1.884775in}{0.827711in}}%
\pgfusepath{stroke}%
\end{pgfscope}%
\begin{pgfscope}%
\pgfpathrectangle{\pgfqpoint{0.135000in}{0.135000in}}{\pgfqpoint{3.499549in}{2.695000in}} %
\pgfusepath{clip}%
\pgfsetrectcap%
\pgfsetroundjoin%
\pgfsetlinewidth{1.003750pt}%
\definecolor{currentstroke}{rgb}{0.000000,0.000000,1.000000}%
\pgfsetstrokecolor{currentstroke}%
\pgfsetdash{}{0pt}%
\pgfpathmoveto{\pgfqpoint{1.884775in}{0.294070in}}%
\pgfpathlineto{\pgfqpoint{0.557292in}{0.827711in}}%
\pgfpathlineto{\pgfqpoint{0.294070in}{1.482500in}}%
\pgfpathlineto{\pgfqpoint{0.557292in}{2.137289in}}%
\pgfpathlineto{\pgfqpoint{1.884775in}{2.670930in}}%
\pgfpathlineto{\pgfqpoint{3.212257in}{2.137289in}}%
\pgfpathlineto{\pgfqpoint{3.475479in}{1.482500in}}%
\pgfpathlineto{\pgfqpoint{3.212257in}{0.827711in}}%
\pgfpathlineto{\pgfqpoint{1.884775in}{0.294070in}}%
\pgfusepath{stroke}%
\end{pgfscope}%
\begin{pgfscope}%
\pgfpathrectangle{\pgfqpoint{0.135000in}{0.135000in}}{\pgfqpoint{3.499549in}{2.695000in}} %
\pgfusepath{clip}%
\pgfsetrectcap%
\pgfsetroundjoin%
\pgfsetlinewidth{1.003750pt}%
\definecolor{currentstroke}{rgb}{0.000000,0.000000,1.000000}%
\pgfsetstrokecolor{currentstroke}%
\pgfsetdash{}{0pt}%
\pgfpathmoveto{\pgfqpoint{1.884775in}{0.827711in}}%
\pgfpathlineto{\pgfqpoint{1.884775in}{0.294070in}}%
\pgfusepath{stroke}%
\end{pgfscope}%
\begin{pgfscope}%
\pgfpathrectangle{\pgfqpoint{0.135000in}{0.135000in}}{\pgfqpoint{3.499549in}{2.695000in}} %
\pgfusepath{clip}%
\pgfsetrectcap%
\pgfsetroundjoin%
\pgfsetlinewidth{1.003750pt}%
\definecolor{currentstroke}{rgb}{0.000000,0.000000,1.000000}%
\pgfsetstrokecolor{currentstroke}%
\pgfsetdash{}{0pt}%
\pgfpathmoveto{\pgfqpoint{0.776957in}{1.121731in}}%
\pgfpathlineto{\pgfqpoint{0.557292in}{0.827711in}}%
\pgfusepath{stroke}%
\end{pgfscope}%
\begin{pgfscope}%
\pgfpathrectangle{\pgfqpoint{0.135000in}{0.135000in}}{\pgfqpoint{3.499549in}{2.695000in}} %
\pgfusepath{clip}%
\pgfsetrectcap%
\pgfsetroundjoin%
\pgfsetlinewidth{1.003750pt}%
\definecolor{currentstroke}{rgb}{0.000000,0.000000,1.000000}%
\pgfsetstrokecolor{currentstroke}%
\pgfsetdash{}{0pt}%
\pgfpathmoveto{\pgfqpoint{0.557292in}{1.482500in}}%
\pgfpathlineto{\pgfqpoint{0.294070in}{1.482500in}}%
\pgfusepath{stroke}%
\end{pgfscope}%
\begin{pgfscope}%
\pgfpathrectangle{\pgfqpoint{0.135000in}{0.135000in}}{\pgfqpoint{3.499549in}{2.695000in}} %
\pgfusepath{clip}%
\pgfsetrectcap%
\pgfsetroundjoin%
\pgfsetlinewidth{1.003750pt}%
\definecolor{currentstroke}{rgb}{0.000000,0.000000,1.000000}%
\pgfsetstrokecolor{currentstroke}%
\pgfsetdash{}{0pt}%
\pgfpathmoveto{\pgfqpoint{0.776957in}{1.843269in}}%
\pgfpathlineto{\pgfqpoint{0.557292in}{2.137289in}}%
\pgfusepath{stroke}%
\end{pgfscope}%
\begin{pgfscope}%
\pgfpathrectangle{\pgfqpoint{0.135000in}{0.135000in}}{\pgfqpoint{3.499549in}{2.695000in}} %
\pgfusepath{clip}%
\pgfsetrectcap%
\pgfsetroundjoin%
\pgfsetlinewidth{1.003750pt}%
\definecolor{currentstroke}{rgb}{0.000000,0.000000,1.000000}%
\pgfsetstrokecolor{currentstroke}%
\pgfsetdash{}{0pt}%
\pgfpathmoveto{\pgfqpoint{1.884775in}{2.137289in}}%
\pgfpathlineto{\pgfqpoint{1.884775in}{2.670930in}}%
\pgfusepath{stroke}%
\end{pgfscope}%
\begin{pgfscope}%
\pgfpathrectangle{\pgfqpoint{0.135000in}{0.135000in}}{\pgfqpoint{3.499549in}{2.695000in}} %
\pgfusepath{clip}%
\pgfsetrectcap%
\pgfsetroundjoin%
\pgfsetlinewidth{1.003750pt}%
\definecolor{currentstroke}{rgb}{0.000000,0.000000,1.000000}%
\pgfsetstrokecolor{currentstroke}%
\pgfsetdash{}{0pt}%
\pgfpathmoveto{\pgfqpoint{2.992592in}{1.843269in}}%
\pgfpathlineto{\pgfqpoint{3.212257in}{2.137289in}}%
\pgfusepath{stroke}%
\end{pgfscope}%
\begin{pgfscope}%
\pgfpathrectangle{\pgfqpoint{0.135000in}{0.135000in}}{\pgfqpoint{3.499549in}{2.695000in}} %
\pgfusepath{clip}%
\pgfsetrectcap%
\pgfsetroundjoin%
\pgfsetlinewidth{1.003750pt}%
\definecolor{currentstroke}{rgb}{0.000000,0.000000,1.000000}%
\pgfsetstrokecolor{currentstroke}%
\pgfsetdash{}{0pt}%
\pgfpathmoveto{\pgfqpoint{3.212257in}{1.482500in}}%
\pgfpathlineto{\pgfqpoint{3.475479in}{1.482500in}}%
\pgfusepath{stroke}%
\end{pgfscope}%
\begin{pgfscope}%
\pgfpathrectangle{\pgfqpoint{0.135000in}{0.135000in}}{\pgfqpoint{3.499549in}{2.695000in}} %
\pgfusepath{clip}%
\pgfsetrectcap%
\pgfsetroundjoin%
\pgfsetlinewidth{1.003750pt}%
\definecolor{currentstroke}{rgb}{0.000000,0.000000,1.000000}%
\pgfsetstrokecolor{currentstroke}%
\pgfsetdash{}{0pt}%
\pgfpathmoveto{\pgfqpoint{2.992592in}{1.121731in}}%
\pgfpathlineto{\pgfqpoint{3.212257in}{0.827711in}}%
\pgfusepath{stroke}%
\end{pgfscope}%
\begin{pgfscope}%
\pgfpathrectangle{\pgfqpoint{0.135000in}{0.135000in}}{\pgfqpoint{3.499549in}{2.695000in}} %
\pgfusepath{clip}%
\pgfsetrectcap%
\pgfsetroundjoin%
\pgfsetlinewidth{1.003750pt}%
\definecolor{currentstroke}{rgb}{0.000000,0.000000,1.000000}%
\pgfsetstrokecolor{currentstroke}%
\pgfsetdash{}{0pt}%
\pgfpathmoveto{\pgfqpoint{1.884775in}{0.827711in}}%
\pgfpathlineto{\pgfqpoint{1.884775in}{0.294070in}}%
\pgfusepath{stroke}%
\end{pgfscope}%
\begin{pgfscope}%
\pgfpathrectangle{\pgfqpoint{0.135000in}{0.135000in}}{\pgfqpoint{3.499549in}{2.695000in}} %
\pgfusepath{clip}%
\pgfsetbuttcap%
\pgfsetroundjoin%
\definecolor{currentfill}{rgb}{0.000000,0.000000,1.000000}%
\pgfsetfillcolor{currentfill}%
\pgfsetlinewidth{1.003750pt}%
\definecolor{currentstroke}{rgb}{0.000000,0.000000,1.000000}%
\pgfsetstrokecolor{currentstroke}%
\pgfsetdash{}{0pt}%
\pgfsys@defobject{currentmarker}{\pgfqpoint{-0.010417in}{-0.010417in}}{\pgfqpoint{0.010417in}{0.010417in}}{%
\pgfpathmoveto{\pgfqpoint{0.000000in}{-0.010417in}}%
\pgfpathcurveto{\pgfqpoint{0.002763in}{-0.010417in}}{\pgfqpoint{0.005412in}{-0.009319in}}{\pgfqpoint{0.007366in}{-0.007366in}}%
\pgfpathcurveto{\pgfqpoint{0.009319in}{-0.005412in}}{\pgfqpoint{0.010417in}{-0.002763in}}{\pgfqpoint{0.010417in}{0.000000in}}%
\pgfpathcurveto{\pgfqpoint{0.010417in}{0.002763in}}{\pgfqpoint{0.009319in}{0.005412in}}{\pgfqpoint{0.007366in}{0.007366in}}%
\pgfpathcurveto{\pgfqpoint{0.005412in}{0.009319in}}{\pgfqpoint{0.002763in}{0.010417in}}{\pgfqpoint{0.000000in}{0.010417in}}%
\pgfpathcurveto{\pgfqpoint{-0.002763in}{0.010417in}}{\pgfqpoint{-0.005412in}{0.009319in}}{\pgfqpoint{-0.007366in}{0.007366in}}%
\pgfpathcurveto{\pgfqpoint{-0.009319in}{0.005412in}}{\pgfqpoint{-0.010417in}{0.002763in}}{\pgfqpoint{-0.010417in}{0.000000in}}%
\pgfpathcurveto{\pgfqpoint{-0.010417in}{-0.002763in}}{\pgfqpoint{-0.009319in}{-0.005412in}}{\pgfqpoint{-0.007366in}{-0.007366in}}%
\pgfpathcurveto{\pgfqpoint{-0.005412in}{-0.009319in}}{\pgfqpoint{-0.002763in}{-0.010417in}}{\pgfqpoint{0.000000in}{-0.010417in}}%
\pgfpathclose%
\pgfusepath{stroke,fill}%
}%
\begin{pgfscope}%
\pgfsys@transformshift{1.884775in}{0.827711in}%
\pgfsys@useobject{currentmarker}{}%
\end{pgfscope}%
\begin{pgfscope}%
\pgfsys@transformshift{1.884775in}{0.294070in}%
\pgfsys@useobject{currentmarker}{}%
\end{pgfscope}%
\begin{pgfscope}%
\pgfsys@transformshift{0.776957in}{1.121731in}%
\pgfsys@useobject{currentmarker}{}%
\end{pgfscope}%
\begin{pgfscope}%
\pgfsys@transformshift{0.557292in}{0.827711in}%
\pgfsys@useobject{currentmarker}{}%
\end{pgfscope}%
\begin{pgfscope}%
\pgfsys@transformshift{0.557292in}{1.482500in}%
\pgfsys@useobject{currentmarker}{}%
\end{pgfscope}%
\begin{pgfscope}%
\pgfsys@transformshift{0.294070in}{1.482500in}%
\pgfsys@useobject{currentmarker}{}%
\end{pgfscope}%
\begin{pgfscope}%
\pgfsys@transformshift{0.776957in}{1.843269in}%
\pgfsys@useobject{currentmarker}{}%
\end{pgfscope}%
\begin{pgfscope}%
\pgfsys@transformshift{0.557292in}{2.137289in}%
\pgfsys@useobject{currentmarker}{}%
\end{pgfscope}%
\begin{pgfscope}%
\pgfsys@transformshift{1.884775in}{2.137289in}%
\pgfsys@useobject{currentmarker}{}%
\end{pgfscope}%
\begin{pgfscope}%
\pgfsys@transformshift{1.884775in}{2.670930in}%
\pgfsys@useobject{currentmarker}{}%
\end{pgfscope}%
\begin{pgfscope}%
\pgfsys@transformshift{2.992592in}{1.843269in}%
\pgfsys@useobject{currentmarker}{}%
\end{pgfscope}%
\begin{pgfscope}%
\pgfsys@transformshift{3.212257in}{2.137289in}%
\pgfsys@useobject{currentmarker}{}%
\end{pgfscope}%
\begin{pgfscope}%
\pgfsys@transformshift{3.212257in}{1.482500in}%
\pgfsys@useobject{currentmarker}{}%
\end{pgfscope}%
\begin{pgfscope}%
\pgfsys@transformshift{3.475479in}{1.482500in}%
\pgfsys@useobject{currentmarker}{}%
\end{pgfscope}%
\begin{pgfscope}%
\pgfsys@transformshift{2.992592in}{1.121731in}%
\pgfsys@useobject{currentmarker}{}%
\end{pgfscope}%
\begin{pgfscope}%
\pgfsys@transformshift{3.212257in}{0.827711in}%
\pgfsys@useobject{currentmarker}{}%
\end{pgfscope}%
\begin{pgfscope}%
\pgfsys@transformshift{1.884775in}{0.827711in}%
\pgfsys@useobject{currentmarker}{}%
\end{pgfscope}%
\begin{pgfscope}%
\pgfsys@transformshift{1.884775in}{0.294070in}%
\pgfsys@useobject{currentmarker}{}%
\end{pgfscope}%
\end{pgfscope}%
\begin{pgfscope}%
\pgfpathrectangle{\pgfqpoint{0.135000in}{0.135000in}}{\pgfqpoint{3.499549in}{2.695000in}} %
\pgfusepath{clip}%
\pgfsetrectcap%
\pgfsetroundjoin%
\pgfsetlinewidth{1.003750pt}%
\definecolor{currentstroke}{rgb}{1.000000,0.000000,0.000000}%
\pgfsetstrokecolor{currentstroke}%
\pgfsetdash{}{0pt}%
\pgfpathmoveto{\pgfqpoint{2.599049in}{0.827711in}}%
\pgfpathlineto{\pgfqpoint{1.170500in}{0.827711in}}%
\pgfpathlineto{\pgfqpoint{0.557292in}{1.285845in}}%
\pgfpathlineto{\pgfqpoint{0.557292in}{1.679155in}}%
\pgfpathlineto{\pgfqpoint{1.170500in}{2.137289in}}%
\pgfpathlineto{\pgfqpoint{2.599049in}{2.137289in}}%
\pgfpathlineto{\pgfqpoint{3.212257in}{1.679155in}}%
\pgfpathlineto{\pgfqpoint{3.212257in}{1.285845in}}%
\pgfpathlineto{\pgfqpoint{2.599049in}{0.827711in}}%
\pgfusepath{stroke}%
\end{pgfscope}%
\begin{pgfscope}%
\pgfpathrectangle{\pgfqpoint{0.135000in}{0.135000in}}{\pgfqpoint{3.499549in}{2.695000in}} %
\pgfusepath{clip}%
\pgfsetrectcap%
\pgfsetroundjoin%
\pgfsetlinewidth{1.003750pt}%
\definecolor{currentstroke}{rgb}{1.000000,0.000000,0.000000}%
\pgfsetstrokecolor{currentstroke}%
\pgfsetdash{}{0pt}%
\pgfpathmoveto{\pgfqpoint{2.599049in}{0.827711in}}%
\pgfpathlineto{\pgfqpoint{1.170500in}{0.827711in}}%
\pgfpathlineto{\pgfqpoint{0.557292in}{1.285845in}}%
\pgfpathlineto{\pgfqpoint{0.557292in}{1.679155in}}%
\pgfpathlineto{\pgfqpoint{1.170500in}{2.137289in}}%
\pgfpathlineto{\pgfqpoint{2.599049in}{2.137289in}}%
\pgfpathlineto{\pgfqpoint{3.212257in}{1.679155in}}%
\pgfpathlineto{\pgfqpoint{3.212257in}{1.285845in}}%
\pgfpathlineto{\pgfqpoint{2.599049in}{0.827711in}}%
\pgfusepath{stroke}%
\end{pgfscope}%
\begin{pgfscope}%
\pgfpathrectangle{\pgfqpoint{0.135000in}{0.135000in}}{\pgfqpoint{3.499549in}{2.695000in}} %
\pgfusepath{clip}%
\pgfsetrectcap%
\pgfsetroundjoin%
\pgfsetlinewidth{1.003750pt}%
\definecolor{currentstroke}{rgb}{1.000000,0.000000,0.000000}%
\pgfsetstrokecolor{currentstroke}%
\pgfsetdash{}{0pt}%
\pgfpathmoveto{\pgfqpoint{2.599049in}{0.827711in}}%
\pgfpathlineto{\pgfqpoint{2.599049in}{0.827711in}}%
\pgfusepath{stroke}%
\end{pgfscope}%
\begin{pgfscope}%
\pgfpathrectangle{\pgfqpoint{0.135000in}{0.135000in}}{\pgfqpoint{3.499549in}{2.695000in}} %
\pgfusepath{clip}%
\pgfsetrectcap%
\pgfsetroundjoin%
\pgfsetlinewidth{1.003750pt}%
\definecolor{currentstroke}{rgb}{1.000000,0.000000,0.000000}%
\pgfsetstrokecolor{currentstroke}%
\pgfsetdash{}{0pt}%
\pgfpathmoveto{\pgfqpoint{1.170500in}{0.827711in}}%
\pgfpathlineto{\pgfqpoint{1.170500in}{0.827711in}}%
\pgfusepath{stroke}%
\end{pgfscope}%
\begin{pgfscope}%
\pgfpathrectangle{\pgfqpoint{0.135000in}{0.135000in}}{\pgfqpoint{3.499549in}{2.695000in}} %
\pgfusepath{clip}%
\pgfsetrectcap%
\pgfsetroundjoin%
\pgfsetlinewidth{1.003750pt}%
\definecolor{currentstroke}{rgb}{1.000000,0.000000,0.000000}%
\pgfsetstrokecolor{currentstroke}%
\pgfsetdash{}{0pt}%
\pgfpathmoveto{\pgfqpoint{0.557292in}{1.285845in}}%
\pgfpathlineto{\pgfqpoint{0.557292in}{1.285845in}}%
\pgfusepath{stroke}%
\end{pgfscope}%
\begin{pgfscope}%
\pgfpathrectangle{\pgfqpoint{0.135000in}{0.135000in}}{\pgfqpoint{3.499549in}{2.695000in}} %
\pgfusepath{clip}%
\pgfsetrectcap%
\pgfsetroundjoin%
\pgfsetlinewidth{1.003750pt}%
\definecolor{currentstroke}{rgb}{1.000000,0.000000,0.000000}%
\pgfsetstrokecolor{currentstroke}%
\pgfsetdash{}{0pt}%
\pgfpathmoveto{\pgfqpoint{0.557292in}{1.679155in}}%
\pgfpathlineto{\pgfqpoint{0.557292in}{1.679155in}}%
\pgfusepath{stroke}%
\end{pgfscope}%
\begin{pgfscope}%
\pgfpathrectangle{\pgfqpoint{0.135000in}{0.135000in}}{\pgfqpoint{3.499549in}{2.695000in}} %
\pgfusepath{clip}%
\pgfsetrectcap%
\pgfsetroundjoin%
\pgfsetlinewidth{1.003750pt}%
\definecolor{currentstroke}{rgb}{1.000000,0.000000,0.000000}%
\pgfsetstrokecolor{currentstroke}%
\pgfsetdash{}{0pt}%
\pgfpathmoveto{\pgfqpoint{1.170500in}{2.137289in}}%
\pgfpathlineto{\pgfqpoint{1.170500in}{2.137289in}}%
\pgfusepath{stroke}%
\end{pgfscope}%
\begin{pgfscope}%
\pgfpathrectangle{\pgfqpoint{0.135000in}{0.135000in}}{\pgfqpoint{3.499549in}{2.695000in}} %
\pgfusepath{clip}%
\pgfsetrectcap%
\pgfsetroundjoin%
\pgfsetlinewidth{1.003750pt}%
\definecolor{currentstroke}{rgb}{1.000000,0.000000,0.000000}%
\pgfsetstrokecolor{currentstroke}%
\pgfsetdash{}{0pt}%
\pgfpathmoveto{\pgfqpoint{2.599049in}{2.137289in}}%
\pgfpathlineto{\pgfqpoint{2.599049in}{2.137289in}}%
\pgfusepath{stroke}%
\end{pgfscope}%
\begin{pgfscope}%
\pgfpathrectangle{\pgfqpoint{0.135000in}{0.135000in}}{\pgfqpoint{3.499549in}{2.695000in}} %
\pgfusepath{clip}%
\pgfsetrectcap%
\pgfsetroundjoin%
\pgfsetlinewidth{1.003750pt}%
\definecolor{currentstroke}{rgb}{1.000000,0.000000,0.000000}%
\pgfsetstrokecolor{currentstroke}%
\pgfsetdash{}{0pt}%
\pgfpathmoveto{\pgfqpoint{3.212257in}{1.679155in}}%
\pgfpathlineto{\pgfqpoint{3.212257in}{1.679155in}}%
\pgfusepath{stroke}%
\end{pgfscope}%
\begin{pgfscope}%
\pgfpathrectangle{\pgfqpoint{0.135000in}{0.135000in}}{\pgfqpoint{3.499549in}{2.695000in}} %
\pgfusepath{clip}%
\pgfsetrectcap%
\pgfsetroundjoin%
\pgfsetlinewidth{1.003750pt}%
\definecolor{currentstroke}{rgb}{1.000000,0.000000,0.000000}%
\pgfsetstrokecolor{currentstroke}%
\pgfsetdash{}{0pt}%
\pgfpathmoveto{\pgfqpoint{3.212257in}{1.285845in}}%
\pgfpathlineto{\pgfqpoint{3.212257in}{1.285845in}}%
\pgfusepath{stroke}%
\end{pgfscope}%
\begin{pgfscope}%
\pgfpathrectangle{\pgfqpoint{0.135000in}{0.135000in}}{\pgfqpoint{3.499549in}{2.695000in}} %
\pgfusepath{clip}%
\pgfsetrectcap%
\pgfsetroundjoin%
\pgfsetlinewidth{1.003750pt}%
\definecolor{currentstroke}{rgb}{1.000000,0.000000,0.000000}%
\pgfsetstrokecolor{currentstroke}%
\pgfsetdash{}{0pt}%
\pgfpathmoveto{\pgfqpoint{2.599049in}{0.827711in}}%
\pgfpathlineto{\pgfqpoint{2.599049in}{0.827711in}}%
\pgfusepath{stroke}%
\end{pgfscope}%
\begin{pgfscope}%
\pgfpathrectangle{\pgfqpoint{0.135000in}{0.135000in}}{\pgfqpoint{3.499549in}{2.695000in}} %
\pgfusepath{clip}%
\pgfsetbuttcap%
\pgfsetroundjoin%
\definecolor{currentfill}{rgb}{1.000000,0.000000,0.000000}%
\pgfsetfillcolor{currentfill}%
\pgfsetlinewidth{1.003750pt}%
\definecolor{currentstroke}{rgb}{1.000000,0.000000,0.000000}%
\pgfsetstrokecolor{currentstroke}%
\pgfsetdash{}{0pt}%
\pgfsys@defobject{currentmarker}{\pgfqpoint{-0.010417in}{-0.010417in}}{\pgfqpoint{0.010417in}{0.010417in}}{%
\pgfpathmoveto{\pgfqpoint{0.000000in}{-0.010417in}}%
\pgfpathcurveto{\pgfqpoint{0.002763in}{-0.010417in}}{\pgfqpoint{0.005412in}{-0.009319in}}{\pgfqpoint{0.007366in}{-0.007366in}}%
\pgfpathcurveto{\pgfqpoint{0.009319in}{-0.005412in}}{\pgfqpoint{0.010417in}{-0.002763in}}{\pgfqpoint{0.010417in}{0.000000in}}%
\pgfpathcurveto{\pgfqpoint{0.010417in}{0.002763in}}{\pgfqpoint{0.009319in}{0.005412in}}{\pgfqpoint{0.007366in}{0.007366in}}%
\pgfpathcurveto{\pgfqpoint{0.005412in}{0.009319in}}{\pgfqpoint{0.002763in}{0.010417in}}{\pgfqpoint{0.000000in}{0.010417in}}%
\pgfpathcurveto{\pgfqpoint{-0.002763in}{0.010417in}}{\pgfqpoint{-0.005412in}{0.009319in}}{\pgfqpoint{-0.007366in}{0.007366in}}%
\pgfpathcurveto{\pgfqpoint{-0.009319in}{0.005412in}}{\pgfqpoint{-0.010417in}{0.002763in}}{\pgfqpoint{-0.010417in}{0.000000in}}%
\pgfpathcurveto{\pgfqpoint{-0.010417in}{-0.002763in}}{\pgfqpoint{-0.009319in}{-0.005412in}}{\pgfqpoint{-0.007366in}{-0.007366in}}%
\pgfpathcurveto{\pgfqpoint{-0.005412in}{-0.009319in}}{\pgfqpoint{-0.002763in}{-0.010417in}}{\pgfqpoint{0.000000in}{-0.010417in}}%
\pgfpathclose%
\pgfusepath{stroke,fill}%
}%
\begin{pgfscope}%
\pgfsys@transformshift{2.599049in}{0.827711in}%
\pgfsys@useobject{currentmarker}{}%
\end{pgfscope}%
\begin{pgfscope}%
\pgfsys@transformshift{2.599049in}{0.827711in}%
\pgfsys@useobject{currentmarker}{}%
\end{pgfscope}%
\begin{pgfscope}%
\pgfsys@transformshift{1.170500in}{0.827711in}%
\pgfsys@useobject{currentmarker}{}%
\end{pgfscope}%
\begin{pgfscope}%
\pgfsys@transformshift{1.170500in}{0.827711in}%
\pgfsys@useobject{currentmarker}{}%
\end{pgfscope}%
\begin{pgfscope}%
\pgfsys@transformshift{0.557292in}{1.285845in}%
\pgfsys@useobject{currentmarker}{}%
\end{pgfscope}%
\begin{pgfscope}%
\pgfsys@transformshift{0.557292in}{1.285845in}%
\pgfsys@useobject{currentmarker}{}%
\end{pgfscope}%
\begin{pgfscope}%
\pgfsys@transformshift{0.557292in}{1.679155in}%
\pgfsys@useobject{currentmarker}{}%
\end{pgfscope}%
\begin{pgfscope}%
\pgfsys@transformshift{0.557292in}{1.679155in}%
\pgfsys@useobject{currentmarker}{}%
\end{pgfscope}%
\begin{pgfscope}%
\pgfsys@transformshift{1.170500in}{2.137289in}%
\pgfsys@useobject{currentmarker}{}%
\end{pgfscope}%
\begin{pgfscope}%
\pgfsys@transformshift{1.170500in}{2.137289in}%
\pgfsys@useobject{currentmarker}{}%
\end{pgfscope}%
\begin{pgfscope}%
\pgfsys@transformshift{2.599049in}{2.137289in}%
\pgfsys@useobject{currentmarker}{}%
\end{pgfscope}%
\begin{pgfscope}%
\pgfsys@transformshift{2.599049in}{2.137289in}%
\pgfsys@useobject{currentmarker}{}%
\end{pgfscope}%
\begin{pgfscope}%
\pgfsys@transformshift{3.212257in}{1.679155in}%
\pgfsys@useobject{currentmarker}{}%
\end{pgfscope}%
\begin{pgfscope}%
\pgfsys@transformshift{3.212257in}{1.679155in}%
\pgfsys@useobject{currentmarker}{}%
\end{pgfscope}%
\begin{pgfscope}%
\pgfsys@transformshift{3.212257in}{1.285845in}%
\pgfsys@useobject{currentmarker}{}%
\end{pgfscope}%
\begin{pgfscope}%
\pgfsys@transformshift{3.212257in}{1.285845in}%
\pgfsys@useobject{currentmarker}{}%
\end{pgfscope}%
\begin{pgfscope}%
\pgfsys@transformshift{2.599049in}{0.827711in}%
\pgfsys@useobject{currentmarker}{}%
\end{pgfscope}%
\begin{pgfscope}%
\pgfsys@transformshift{2.599049in}{0.827711in}%
\pgfsys@useobject{currentmarker}{}%
\end{pgfscope}%
\end{pgfscope}%
\begin{pgfscope}%
\pgfpathrectangle{\pgfqpoint{0.135000in}{0.135000in}}{\pgfqpoint{3.499549in}{2.695000in}} %
\pgfusepath{clip}%
\pgfsetrectcap%
\pgfsetroundjoin%
\pgfsetlinewidth{0.702625pt}%
\definecolor{currentstroke}{rgb}{0.000000,0.501961,0.000000}%
\pgfsetstrokecolor{currentstroke}%
\pgfsetdash{}{0pt}%
\pgfpathmoveto{\pgfqpoint{1.884775in}{2.137289in}}%
\pgfpathlineto{\pgfqpoint{1.972834in}{2.135847in}}%
\pgfpathlineto{\pgfqpoint{2.050585in}{2.132161in}}%
\pgfpathlineto{\pgfqpoint{2.126665in}{2.126327in}}%
\pgfpathlineto{\pgfqpoint{2.202278in}{2.118284in}}%
\pgfpathlineto{\pgfqpoint{2.275895in}{2.108223in}}%
\pgfpathlineto{\pgfqpoint{2.347707in}{2.096183in}}%
\pgfpathlineto{\pgfqpoint{2.417782in}{2.082189in}}%
\pgfpathlineto{\pgfqpoint{2.485555in}{2.066394in}}%
\pgfpathlineto{\pgfqpoint{2.550670in}{2.048949in}}%
\pgfpathlineto{\pgfqpoint{2.612872in}{2.030011in}}%
\pgfpathlineto{\pgfqpoint{2.672404in}{2.009581in}}%
\pgfpathlineto{\pgfqpoint{2.728996in}{1.987817in}}%
\pgfpathlineto{\pgfqpoint{2.782418in}{1.964895in}}%
\pgfpathlineto{\pgfqpoint{2.832462in}{1.941015in}}%
\pgfpathlineto{\pgfqpoint{2.879285in}{1.916217in}}%
\pgfpathlineto{\pgfqpoint{2.922643in}{1.890754in}}%
\pgfpathlineto{\pgfqpoint{2.962627in}{1.864716in}}%
\pgfpathlineto{\pgfqpoint{2.999283in}{1.838219in}}%
\pgfpathlineto{\pgfqpoint{3.032625in}{1.811413in}}%
\pgfpathlineto{\pgfqpoint{3.062932in}{1.784216in}}%
\pgfpathlineto{\pgfqpoint{3.090142in}{1.756823in}}%
\pgfpathlineto{\pgfqpoint{3.114174in}{1.729516in}}%
\pgfpathlineto{\pgfqpoint{3.135492in}{1.701939in}}%
\pgfpathlineto{\pgfqpoint{3.153956in}{1.674420in}}%
\pgfpathlineto{\pgfqpoint{3.169688in}{1.646989in}}%
\pgfpathlineto{\pgfqpoint{3.182789in}{1.619700in}}%
\pgfpathlineto{\pgfqpoint{3.193337in}{1.592658in}}%
\pgfpathlineto{\pgfqpoint{3.201393in}{1.566102in}}%
\pgfpathlineto{\pgfqpoint{3.207269in}{1.539209in}}%
\pgfpathlineto{\pgfqpoint{3.210730in}{1.513904in}}%
\pgfpathlineto{\pgfqpoint{3.212058in}{1.493842in}}%
\pgfpathlineto{\pgfqpoint{3.211394in}{1.458890in}}%
\pgfpathlineto{\pgfqpoint{3.208202in}{1.431357in}}%
\pgfpathlineto{\pgfqpoint{3.202864in}{1.404744in}}%
\pgfpathlineto{\pgfqpoint{3.195086in}{1.377524in}}%
\pgfpathlineto{\pgfqpoint{3.185095in}{1.350720in}}%
\pgfpathlineto{\pgfqpoint{3.172566in}{1.323579in}}%
\pgfpathlineto{\pgfqpoint{3.157700in}{1.296711in}}%
\pgfpathlineto{\pgfqpoint{3.140133in}{1.269607in}}%
\pgfpathlineto{\pgfqpoint{3.120037in}{1.242707in}}%
\pgfpathlineto{\pgfqpoint{3.096995in}{1.215625in}}%
\pgfpathlineto{\pgfqpoint{3.071131in}{1.188705in}}%
\pgfpathlineto{\pgfqpoint{3.042257in}{1.161905in}}%
\pgfpathlineto{\pgfqpoint{3.010455in}{1.135445in}}%
\pgfpathlineto{\pgfqpoint{2.975467in}{1.109252in}}%
\pgfpathlineto{\pgfqpoint{2.937312in}{1.083488in}}%
\pgfpathlineto{\pgfqpoint{2.895980in}{1.058279in}}%
\pgfpathlineto{\pgfqpoint{2.851427in}{1.033720in}}%
\pgfpathlineto{\pgfqpoint{2.803569in}{1.009893in}}%
\pgfpathlineto{\pgfqpoint{2.752265in}{0.986864in}}%
\pgfpathlineto{\pgfqpoint{2.698162in}{0.965025in}}%
\pgfpathlineto{\pgfqpoint{2.640662in}{0.944230in}}%
\pgfpathlineto{\pgfqpoint{2.580441in}{0.924824in}}%
\pgfpathlineto{\pgfqpoint{2.517279in}{0.906815in}}%
\pgfpathlineto{\pgfqpoint{2.451437in}{0.890366in}}%
\pgfpathlineto{\pgfqpoint{2.383269in}{0.875632in}}%
\pgfpathlineto{\pgfqpoint{2.312484in}{0.862629in}}%
\pgfpathlineto{\pgfqpoint{2.240498in}{0.851658in}}%
\pgfpathlineto{\pgfqpoint{2.166390in}{0.842615in}}%
\pgfpathlineto{\pgfqpoint{2.091293in}{0.835683in}}%
\pgfpathlineto{\pgfqpoint{2.016195in}{0.830928in}}%
\pgfpathlineto{\pgfqpoint{1.941098in}{0.828301in}}%
\pgfpathlineto{\pgfqpoint{1.866000in}{0.827777in}}%
\pgfpathlineto{\pgfqpoint{1.783671in}{0.829613in}}%
\pgfpathlineto{\pgfqpoint{1.705678in}{0.833698in}}%
\pgfpathlineto{\pgfqpoint{1.629415in}{0.839940in}}%
\pgfpathlineto{\pgfqpoint{1.554751in}{0.848269in}}%
\pgfpathlineto{\pgfqpoint{1.480798in}{0.858767in}}%
\pgfpathlineto{\pgfqpoint{1.409445in}{0.871127in}}%
\pgfpathlineto{\pgfqpoint{1.339671in}{0.885462in}}%
\pgfpathlineto{\pgfqpoint{1.272661in}{0.901477in}}%
\pgfpathlineto{\pgfqpoint{1.208115in}{0.919163in}}%
\pgfpathlineto{\pgfqpoint{1.146342in}{0.938367in}}%
\pgfpathlineto{\pgfqpoint{1.087140in}{0.959093in}}%
\pgfpathlineto{\pgfqpoint{1.031211in}{0.981016in}}%
\pgfpathlineto{\pgfqpoint{0.978340in}{1.004122in}}%
\pgfpathlineto{\pgfqpoint{0.928755in}{1.028212in}}%
\pgfpathlineto{\pgfqpoint{0.882673in}{1.053055in}}%
\pgfpathlineto{\pgfqpoint{0.839967in}{1.078575in}}%
\pgfpathlineto{\pgfqpoint{0.800564in}{1.104684in}}%
\pgfpathlineto{\pgfqpoint{0.764430in}{1.131269in}}%
\pgfpathlineto{\pgfqpoint{0.731563in}{1.158181in}}%
\pgfpathlineto{\pgfqpoint{0.701989in}{1.185223in}}%
\pgfpathlineto{\pgfqpoint{0.675466in}{1.212431in}}%
\pgfpathlineto{\pgfqpoint{0.651797in}{1.239862in}}%
\pgfpathlineto{\pgfqpoint{0.631102in}{1.267204in}}%
\pgfpathlineto{\pgfqpoint{0.613235in}{1.294417in}}%
\pgfpathlineto{\pgfqpoint{0.597805in}{1.321969in}}%
\pgfpathlineto{\pgfqpoint{0.584996in}{1.349425in}}%
\pgfpathlineto{\pgfqpoint{0.574732in}{1.376710in}}%
\pgfpathlineto{\pgfqpoint{0.566952in}{1.403649in}}%
\pgfpathlineto{\pgfqpoint{0.561614in}{1.429707in}}%
\pgfpathlineto{\pgfqpoint{0.558421in}{1.455501in}}%
\pgfpathlineto{\pgfqpoint{0.557358in}{1.475952in}}%
\pgfpathlineto{\pgfqpoint{0.558554in}{1.511043in}}%
\pgfpathlineto{\pgfqpoint{0.562280in}{1.539209in}}%
\pgfpathlineto{\pgfqpoint{0.568156in}{1.566102in}}%
\pgfpathlineto{\pgfqpoint{0.576482in}{1.593434in}}%
\pgfpathlineto{\pgfqpoint{0.587032in}{1.620324in}}%
\pgfpathlineto{\pgfqpoint{0.600135in}{1.647510in}}%
\pgfpathlineto{\pgfqpoint{0.615593in}{1.674420in}}%
\pgfpathlineto{\pgfqpoint{0.633775in}{1.701547in}}%
\pgfpathlineto{\pgfqpoint{0.654515in}{1.728472in}}%
\pgfpathlineto{\pgfqpoint{0.678238in}{1.755569in}}%
\pgfpathlineto{\pgfqpoint{0.704823in}{1.782506in}}%
\pgfpathlineto{\pgfqpoint{0.734164in}{1.809058in}}%
\pgfpathlineto{\pgfqpoint{0.766792in}{1.835558in}}%
\pgfpathlineto{\pgfqpoint{0.802354in}{1.861562in}}%
\pgfpathlineto{\pgfqpoint{0.841149in}{1.887167in}}%
\pgfpathlineto{\pgfqpoint{0.883201in}{1.912245in}}%
\pgfpathlineto{\pgfqpoint{0.928571in}{1.936693in}}%
\pgfpathlineto{\pgfqpoint{0.977368in}{1.960430in}}%
\pgfpathlineto{\pgfqpoint{1.029355in}{1.983214in}}%
\pgfpathlineto{\pgfqpoint{1.084272in}{2.004841in}}%
\pgfpathlineto{\pgfqpoint{1.142296in}{2.025291in}}%
\pgfpathlineto{\pgfqpoint{1.203184in}{2.044389in}}%
\pgfpathlineto{\pgfqpoint{1.267215in}{2.062119in}}%
\pgfpathlineto{\pgfqpoint{1.333562in}{2.078172in}}%
\pgfpathlineto{\pgfqpoint{1.402452in}{2.092539in}}%
\pgfpathlineto{\pgfqpoint{1.473449in}{2.105063in}}%
\pgfpathlineto{\pgfqpoint{1.545794in}{2.115581in}}%
\pgfpathlineto{\pgfqpoint{1.620596in}{2.124192in}}%
\pgfpathlineto{\pgfqpoint{1.695163in}{2.130575in}}%
\pgfpathlineto{\pgfqpoint{1.772128in}{2.134927in}}%
\pgfpathlineto{\pgfqpoint{1.847226in}{2.137027in}}%
\pgfpathlineto{\pgfqpoint{1.884775in}{2.137289in}}%
\pgfpathlineto{\pgfqpoint{1.884775in}{2.137289in}}%
\pgfusepath{stroke}%
\end{pgfscope}%
\end{pgfpicture}%
\makeatother%
\endgroup%